%% file: sigma15-060.tex
\pdfoutput=1
\RequirePackage{ifpdf}
\ifpdf 
\documentclass[pdftex]{sigma}
\else
\documentclass{sigma}
\fi

\numberwithin{equation}{section}
\numberwithin{theorem}{section}

\newtheorem{lem}[theorem]{Lemma}
\newtheorem{prop}[theorem]{Proposition}
\newtheorem{cor}[theorem]{Corollary}
{\theoremstyle{definition}
\newtheorem{rem}[theorem]{Remark}
\newtheorem{exam}[theorem]{Example}
\newtheorem{defin}[theorem]{Definition}
\newtheorem{notation}[theorem]{Notation}
}

\usepackage{amsbsy}
\usepackage{subfigure}

\usepackage{mdwlist} 
\usepackage{enumerate} 
\usepackage{multicol} 

\newcommand\notch{^{(p)}}
\newcommand \tilblacktri[1]{\widetilde{\triangle}_{#1}}
\newcommand \blacktri[1]{\triangle_{#1}}
\newcommand \tiltaui[1] {\widetilde{\tau}_{i_{#1}}}
\newcommand \taui[1] {\tau_{i_{#1}}}

\newcommand \al {\alpha} \newcommand \be {\beta}
\newcommand\ga{\gamma}  \newcommand\Ga{\Gamma}
 \newcommand\si{\sigma}  \newcommand\Si{\Sigma}
\newcommand \om {w} 
\newcommand\Cn{\mathcal{C}_n}
\newcommand\Surf{S}
\newcommand\Cnn[1]{\mathcal{C}_{#1}}
\newcommand\Z{\mathbb{Z}}
\newcommand\B{\mathcal B}
\newcommand\AAA {\mathcal A}

\newcommand\tgk{\tau_{[\ga_k]}}
\newcommand\tbk{\tau_{[\be_k]}}

\newcommand \llrrl {\ell, \ell, r, r, \ell} \newcommand \lrrll {\ell, r, r, \ell, \ell}
\newcommand\rrBak{\underline{r,r}}
\newcommand\rrNon{\mathbf{r,r}}
\newcommand\ttBak{\underline{\tau, \tau}}
\newcommand\ttNon{\boldsymbol{\tau, \tau}}
\newcommand\ttQuasi{\overline{\tau, \tau}}
\newcommand\dotRR{\dot{\textbf{r}}, \textbf{r}}
\newcommand\RddotR{\textbf{r}, \ddot{\textbf{r}}}
\newcommand\GT {G_{T^o,\ga}}
\newcommand\overGT {\overline{G}_{T^o,\ga}}
\newcommand\puncture{\text{\tiny $P$}}
\newcommand\Quad[1]{{\rm Quad}{(#1)}}

\newcommand \tilp {\widetilde{p}} \newcommand \tilP {\widetilde{P}}
  
\newcommand\tilt{\tilde{t}} \newcommand\tils{\tilde{s}}  
\newcommand\tilS{\widetilde{\Sga}}  \newcommand\Sga{S_\ga}
\newcommand\tilA{\widetilde{A}} \newcommand\tilD{\widetilde{D}}
\newcommand\tilga{\widetilde{\ga}}
\newcommand\tilsi{\widetilde{\si}}
\newcommand\tilw{\widetilde{\om}}

\newcommand\tilbe{\widetilde{\be}}

 \newcommand\tilom{\widetilde{\om}}
\newcommand\overpi{\overline{\pi}}

\newcommand\tilTo{\widetilde{T^o}}
\def\a{\aleph}
\def\b{\text{\tiny $\beth$}}
\newcommand\ccr{r^{cc}}
\newcommand\bijection{\longleftrightarrow}

\newcommand\ts[1]{\tau_{[\si_{#1}]}}
\newcommand\tg[1]{\tau_{[\ga_{#1}]}}

\newcommand\tiltg[1]{\widetilde{\tau}_{[\ga_{#1}]}}
\newcommand \pkpkone[1] {[p_k,p_{k+1}]_{#1}}

\newcommand\overom{\bar{\om}}
\newcommand\disk[1]{\textbf{Disk}_{#1}}
\newcommand\tildisk[1]{\textbf{Disk}_{\widetilde{#1}}}
\newcommand\tildiskplus[2]{\textbf{Disk}_{\widetilde{#1}}^{#2}}
\newcommand\outsidedisk[1]{\textbf{C}_{#1}}

\newcommand\PMOne{(b_4, 	 	b_2, 	 	b_3, 		 		r, 	\ell, 	b_3)}

\newcommand\PMTwo{(b_4, 	1, 	\ell, 	r, 	\ell, 	b_3)}

\newcommand\PMThree{(b_1,	2,2,\ell, r, 2)}

\newcommand\PMFour{(b_1, 2, b_3, r,\ell, b_3)}

\newcommand\PMFive{(b_4, b_2, 2, \ell,r,2)}

\newcommand\PMSix{(b_4, 	b_2,2,	\ell,	\ddot{r},	\ell, 	b_3)}

\newcommand\PMSeven{(b_4, b_2, 2, \dot{r}, \ell, b_3)}

\newcommand\PMEight{(b_1, 2,2, \ell,	\ddot{r}, b_3)}

\newcommand\PMNine{(b_1, 2, 2, \dot{r}, \ell, b_3)}

\input{figures}

\begin{document}

\newcommand{\arXivNumber}{1409.3610}

\allowdisplaybreaks

\renewcommand{\thefootnote}{$\star$}

\renewcommand{\PaperNumber}{060}

\FirstPageHeading

\ShortArticleName{$T$-Path Formula and Atomic Bases for Cluster Algebras of Type $D$}

\ArticleName{$\boldsymbol{T}$-Path Formula and Atomic Bases\\ for Cluster Algebras of Type~$\boldsymbol{D}$\footnote{This paper is a~contribution to the Special Issue on New Directions in Lie Theory.
The full collection is available at \href{http://www.emis.de/journals/SIGMA/LieTheory2014.html}{http://www.emis.de/journals/SIGMA/LieTheory2014.html}}}

\Author{Emily GUNAWAN and Gregg MUSIKER}
\AuthorNameForHeading{E.~Gunawan and G.~Musiker}
\Address{School of Mathematics,  University of Minnesota, Minneapolis, MN 55455, USA}
\Email{\href{mailto:egunawan@umn.edu}{egunawan@umn.edu}, \href{mailto:musiker@umn.edu}{musiker@umn.edu}}
\URLaddress{\url{http://umn.edu/home/egunawan}, \url{http://math.umn.edu/~musiker}}

\ArticleDates{Received September 24, 2014, in f\/inal form July 09, 2015; Published online July 28, 2015}

\Abstract{We extend a $T$-path expansion formula for arcs on an unpunctured surface
to the case of arcs on a once-punctured polygon
and use this formula to give a combinatorial proof
that cluster monomials form the atomic basis of a cluster algebra of type $D$.}


\Classification{13F60; 05E15; 16S99}

\renewcommand{\thefootnote}{\arabic{footnote}}
\setcounter{footnote}{0}

\section{Introduction}

Cluster algebras, introduced by Fomin and Zelevinsky~\cite{FZ02} in~2000,
 are commutative algebras equipped with a distinguished set of generators, the \emph{cluster variables}. The cluster variables are grouped into sets of constant cardinality $n$, the \emph{clusters}, and the integer~$n$ is called the \emph{rank} of the cluster algebra. Starting with an initial \emph{seed} $(\mathbf{x},B)$, that is, an initial cluster~$\mathbf{x}$ together with a~skew-symmetrizable integer~$n\times n$ \emph{exchange matrix} $B=(b_{ij})$, the set of cluster variables is obtained by repeated application of so called \emph{mutations}.
To be more precise, let $\{x_1,x_2,\dots,x_n\}$ be indeterminates over~$\mathbb{Z}$ and
let ${\mathcal{F}}=\mathbb{Q}(x_1,x_2,\dots,x_n)$.
For every $k=1,2,\dots,n$, the mutation~$\mu_k(\mathbf{x})$ of the cluster $\mathbf{x}=\{ x_1,x_2,\dots,x_n\}$
is a new cluster $\mu(\mathbf{x}) = ( \mathbf{x} \backslash \{ x_k \} ) \cup \{ x_k'\}$
obtained from~$\mathbf{x}$ by replacing the cluster variable~$x_k$ with the new cluster variable
\begin{gather*}
x_k' = \frac{1}{x_k} \left( \prod_{b_{ik>0}} x_i^{b_{ik}} + \prod_{b_{ik<0}} x_i^{-b_{ik}}\right)
\end{gather*}
in ${\mathcal{F}}$. Mutations also change the attached matrix $B$, see~\cite{FZ02}.

The set of all cluster variables is the union of all clusters obtained from an initial cluster $\mathbf{x}$ by repeated mutations. Note that this set may be inf\/inite.

It is clear from the construction that every cluster variable is a rational function in the initial cluster variables $x_1,x_2,\dots,x_n$.
In \cite{FZ02} it is shown that every cluster variable $u$ is actually a~Laurent polynomial in the~$x_i$,
that is, $u$ can be written as a reduced fraction
\begin{gather}\label{eq:expansion_of_u_in_x}
u=\frac{f(x_1,x_2,\dots,x_n)}{\prod\limits_{i=1}^n x_i^{d_i}},
\end{gather}
where $f\in \mathbb{Z}[x_1,x_2,\dots,x_n]$ and $d_i\geq 0$.
The right hand side of equation (\ref{eq:expansion_of_u_in_x}) is called the \emph{cluster expansion of} $u$ in $\mathbf{x}$.

The \emph{coefficient-free cluster algebra} $\AAA(\mathbf{x},B)$ is the subring of ${\mathcal{F}}=\mathbb{Q}(x_1,x_2,\dots,x_n)$
generated by the cluster variables. When the set of cluster variables is f\/inite, we say that $\AAA(\mathbf{x},B)$ is of \emph{finite type}.

We are interested in cluster algebras arising from \emph{bordered surfaces with marked points}  \cite{FG06,FG09,FST08,GSV05}. In particular, we study cluster algebras of type $D_n$ (type $D$ for short), which are of f\/inite type (as classif\/ied in \cite{FZ03}) and correspond to once-punctured $n$-gons, as also described in detail in \cite{Sch08}. Other, related work on type $D$ cluster algebra combinatorial models include \cite{BM09, CP14, FZ03}.

Our f\/irst result is a Laurent polynomial expansion formula for cluster variables arising from a~once-punctured polygon
in terms of certain paths (called $T^o$-paths) on an ideal triangula\-tion~$T^o$ of the surface.
This is an extension of the $T$-path formula (which we call the $T^o$-path) given for any unpunctured surface by Schif\/f\/ler and Thomas \cite{Sch10, ST09}. Our proof takes advantage of two facts proven by the second author, Schif\/f\/ler, and Williams: (1)~a~Laurent polynomial expansion formula for cluster variables arising from any surface in terms of perfect matchings of a \emph{snake graph}~\cite{MSW11}, and (2)~a~bijection between these perfect matchings and $T^o$-paths arising from any unpunctured surface \cite{MS10}. An application of the $T^o$-path formula for type~$D$ is discussed in the next paragraph.

Our second result is a specif\/ic case of a result of \cite{Cer11,CL12}, proven by representation theoretic methods,
that the basis consisting of all \emph{cluster monomials} is in fact the \emph{atomic basis} for any skew-symmetric cluster algebra of f\/inite type (see Section \ref{subsection:atomic_bases_D}). We give a combinatorial proof of this fact for \emph{coefficient-free} cluster algebras of type $D$. Our proof relies heavily on the $T^o$-path formula for type $D$ and
is inspired by Dupont and Thomas' work in \cite{DT13} on atomic bases for cluster algebras of type $A$ and $\tilA$.

In Section~\ref{sec:background}, we provide background material on ideal triangulations and tagged triangulations, focusing on the case of once-punctured polygons. Section~\ref{sec:tpaths} presents our f\/irst result (Theorem~\ref{thm:tpath_expansion_formula}), an extension of the $T^o$-path formula of~\cite{Sch10, ST09} to once-punctured polygons. We give the proof of this $T^o$-path formula in Section~\ref{sec:proof_of_tpath_expansion_formula}.
Finally, in Section~\ref{sec:proof_atomic_basis}, we give our second and main result (Theorem~\ref{thm:main_thm}), which is a combinatorial proof, using the $T^o$-path formula, that the cluster monomials form the atomic basis for a coef\/f\/icient-free type $D$ cluster algebra.

\section{Background: cluster algebras arising\\ from once-punctured disks}  \label{sec:background}
For the reader's convenience, we begin by reviewing terminology arising in the theory of cluster algebras from marked surfaces from \cite[Sections~2 and~7]{FST08}.  We restrict our attention to the case of a once-punctured polygon, which often simplif\/ies the notation.  Let $\Cn$ denote a once-punctured $n$-gon, i.e., a disk with a marked point (called the puncture) in the interior and $n$ marked points on the boundary.

\begin{figure}[t!]
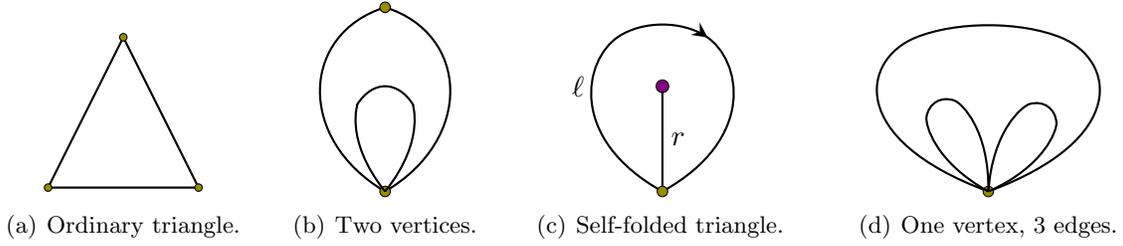

\centering
\mbox{\subfigure[Ordinary triangle.]
{
\TikzTriangle{0.5}
\label{fig:TikzTriangle}
}\quad
\subfigure[Two vertices.]
{
\TikzTwoVertexTriangle{0.7}
\label{fig:TikzTwoVertexTriangle}
}\quad
\subfigure[Self-folded triangle.]
{\quad
\TikzSelfFoldedTriangle{0.7}\quad
\label{fig:TikzSelfFoldedTriangle}
}\quad
\subfigure[One vertex, 3 edges.]
{
\TikzOneVertexTriangle{0.7}
\label{fig:TikzOneVertexTriangle}
}
}
\caption{Possible types of ideal triangles.}
\label{fig:possible_types_of_ideal_triangles}
\end{figure}

\begin{defin}[ordinary arcs]\label{arcs}
A \emph{boundary edge} of $\Cn$ is a segment of the boundary between two consecutive boundary marked points.
An \emph{ordinary arc} $\ga$ of $\Cn$ is a curve (considered up to isotopy) in $\Cn$ such that
the endpoints of $\ga$ are marked points,
$\ga$ does not cross itself except possibly at its endpoints,
$\ga$ does not cross the boundary of $\Cn$ except possibly at its endpoints,
and $\ga$ is not contractible to a marked point
or homotopic to a boundary edge.

A \emph{radius} is an arc between a boundary marked point and the puncture.
Following \cite{DT13}, a~\emph{peripheral} arc is an arc with both endpoints on the boundary.
An \emph{$\ell$-loop} is a loop cutting out a~monogon with a sole puncture inside it
(i.e., as illustrated by loop $\ell$ in Fig.~\ref{fig:TikzSelfFoldedTriangle}).
An $\ell$-loop is considered a peripheral, ordinary arc.
\end{defin}

\begin{defin}[compatibility of ordinary arcs, ideal triangulations]
Two distinct ordinary arcs are said to be \emph{compatible} if they do not intersect
except possibly at endpoints. Also, each arc is compatible with
itself. A maximal (by inclusion) collection of distinct, pairwise compatible ordinary arcs is
called an \emph{ideal triangulation}.
The ordinary arcs of an ideal triangulation cut the surface into \emph{ideal triangles} (see Fig.~\ref{fig:possible_types_of_ideal_triangles}).
We call an $\ell$-loop and the radius it encloses a \emph{self-folded triangle}
(Fig.~\ref{fig:TikzSelfFoldedTriangle}).
\end{defin}

\begin{figure}[t!]
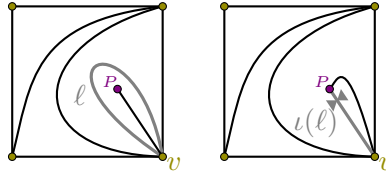
\centering
\TikzTtwoTagged{0.5}
\caption{An ideal triangulation ${T}^o$ and the corresponding tagged triangulation $T=\iota({T}^o)$.
The $\ell$-loop~$\ell$ and the corresponding tagged radius $\iota(\ell)$ (tagged notched at the puncture) are both in gray.}
\label{fig:T2_tagged}
 \end{figure}

\begin{rem}
Three possible types of ideal triangles can appear in an ideal triangulation of~$\Cn$:
an ordinary triangle with~3 distinct vertices and~3 distinct sides (Fig.~\ref{fig:TikzTriangle}),
an ideal triangle with~2 distinct vertices and 3 distinct sides (Fig.~\ref{fig:TikzTwoVertexTriangle}),
and f\/inally a self-folded triangle (Fig.~\ref{fig:TikzSelfFoldedTriangle}).
The ideal triangulation of Fig.~\ref{fig:T2_tagged} (left) contains all 3 types of $\Cn$-ideal triangles.
Note that the one-vertex ideal triangle (Fig.~\ref{fig:TikzOneVertexTriangle}) cannot appear.
\end{rem}

\begin{defin}[tagged arcs] \label{def:tagged}
A \emph{tagged arc} of $\Cn$ is obtained by marking (``tagging'') each endpoint of an ordinary arc (that is not an $\ell$-loop) $\be$ either \emph{plain} or \emph{notched} such that the endpoints of $\be$ on the boundary must be tagged plain. A notching is indicated by a bow tie (see Fig.~\ref{fig:T2_tagged}). Note that a tagged arc never cuts out a once-punctured monogon, i.e., an $\ell$-loop is \emph{not} a tagged arc (even for other punctured surfaces).
\end{defin}

\begin{rem}\label{remark:every_tagged_arc_belongs_to_1of3_classes}
Every tagged arc $\be$ of $\Cn$ belongs to one of the following three classes:
\begin{itemize}\itemsep=0pt
\item $\be$ is a radius tagged plain at both endpoints (which we call a plain radius).
\item $\be$ is a radius tagged notched at the puncture and plain at the boundary (which we call a~notched radius).
\item $\be$ is a peripheral arc connecting distinct endpoints tagged plain at both endpoints.
\end{itemize}
\end{rem}

\begin{defin}[compatibility of tagged arcs, tagged triangulations, and multi-tagged triangulations of $\Cn$]\label{def:compatibility_of_tagged_arcs}
The following is a complete list of possible compatible pairs $\{ \al, \be \}$ of tagged arcs of~$\Cn$:
\begin{itemize}\itemsep=0pt
\item $\al$ and $\be$ are two peripheral arcs (tagged plain at all endpoints)
 that do not intersect in the interior of $\Cn$.
\item $\al$ and $\be$ are a peripheral arc and a radius (tagged plain at boundary endpoints)
that do not intersect in the interior of $\Cn$.
\item $\al$ and $\be$ are two radii both adjacent to the same boundary marked point (tagged plain) but $\al$ is tagged plain at the puncture and $\be$ is tagged notched at the puncture.
\item $\al$ and $\be$ are two radii with distinct boundary endpoints (tagged plain at boundary endpoints) and tagged the same way at the puncture.
\item $\al$ and $\be$ are equal.
\end{itemize}
A maximal (by inclusion) collection of distinct, pairwise compatible tagged arcs is called a \emph{tagged triangulation}.
A collection of pairwise compatible tagged arcs (considered with multiplicity) is called a \emph{multi-tagged triangulation}.
A multi-tagged triangulation $\Si$ is \emph{compatible} with $T$ if $\si \in T$ for every tagged arc $\si \in \Si$.
\end{defin}

\begin{defin}[representing ordinary arcs as tagged arcs]
\label{def:representing_ordinary_arcs_as_tagged_arcs}
Any ordinary arc $\be$ can be represented by a tagged arc $\iota (\be)$ as follows.
Suppose $\be$ is an $\ell$-loop
(based at marked point $v$)
which encloses a radius $r$, where
$r$ is the unique (ordinary) arc connecting~$v$ and the puncture~$\puncture$.
Then~$\iota(\be)$ is obtained by tagging $r$ plain at $v$ and notched at~$\puncture$.
Otherwise, $\iota(\be)$ is simply~$\be$ with both endpoints tagged plain.
Fig.~\ref{fig:T2_tagged} shows an ideal triangulation~${T}^o$ of a once-punctured quadrilateral
and its corresponding tagged triangulation~$\iota({T}^o)$.
\end{defin}

 As a convention, we will usually denote a tagged triangulation by~$T$ and an ideal triangulation by~$T^o$.  Unless otherwise stated, $T = \iota(T^o)$.

\begin{theorem}[\protect{\cite[Theorem 7.11, Example 6.7]{FST08}}]
\label{thm:fst_thm}
A cluster algebra is associated to $\Cn$ as follows.
Choose a tagged triangulation $T_{\rm init}=\{\tau_1,\dots,\tau_n\}$ of $\Cn$.
Let $\AAA$ be the cluster algebra given by the initial seed $(\mathbf{x}_{T_{\rm init}},B_{T_{\rm init}})$ where $\mathbf{x}_{T_{\rm init}}=\{x_{\tau_1},\dots,x_{\tau_n}\}$ is the cluster seed associated to $T_{\rm init}$ and $B_{T_{\rm init}}$ is the exchange matrix corresponding to $T_{\rm init}$ $($see Definition~{\rm 4.1} in {\rm~\cite{FST08})}.
Then the tagged triangulations of $\Cn$ are in bijection with the $($unlabeled$)$ seeds $(\mathbf{x}_T,B_T)$ of~$\AAA$,
and the tagged arcs~$\ga$ of~$\Cn$ are in bijection with the cluster variables $($so we can denote each cluster variable by $x_\ga$ or $x(\ga)$, where $\ga$ is a tagged arc$)$.
Therefore, a multi-tagged triangula\-tion~$\Ga$ corresponds to a cluster monomial $($denoted $x_\Ga$, see Definition~{\rm \ref{def:cluster_monomial})}.
\end{theorem}

If $r$ is a plain radius and $\ell$ is the $\ell$-loop enclosing $r$, denote $x_\ell := x_r x_{r\notch}$.
If $\be$ is a boundary edge, set $x_\be :=1$.
When we say the \emph{$T$-expansion} of a cluster variable $x(\ga)$, we mean the cluster expansion of $x(\ga)$ in the variables of the seed $\mathbf{x}_T$, i.e., a~Laurent polynomial in the variables of~$\mathbf{x}_T$, see equation~(\ref{eq:expansion_of_u_in_x}).

\begin{figure}[t!]
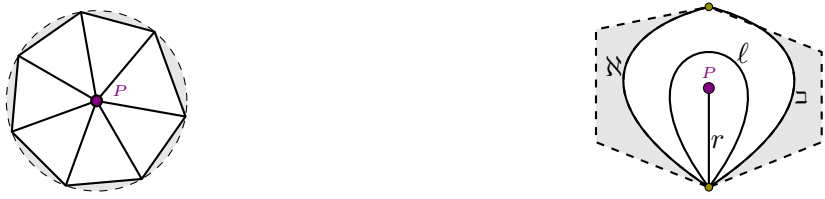

\centering
\mbox{
\subfigure[Wheel-like triangulation around the puncture.]
{\hspace{14mm}
\WheelLocalTriangulation{0.6}\hspace{14mm}
\label{fig:wheel}
}\qquad
\subfigure[A self-folded triangle around the puncture where $\a$, $\b$ are peripheral arcs or boundary edges.]
{\hspace{10mm}
\SelfFoldedLocalTriangulation{0.6}\hspace{10mm}
\label{fig:SelfFoldedLocalTriangulation}
}}
\caption{Local triangulations around the puncture.
The shaded gray area consists only of peripheral arcs and boundary edges, each sector forming a polygon.}
\label{fig:ideal_triangulations}
 \end{figure}

\begin{rem}\label{remark:only_need_to_prove_two_cases}
Due to the following proposition, it is enough to work with only two types of tagged triangulations~$T$: one where $T$ has all plain-tagged radii (so that $T^o$ has a local wheel-like triangulation as in Fig.~\ref{fig:wheel}),
and one where $T$ has two parallel radii, one tagged plain and the other tagged notched at the puncture (so that $T^o$ has a self-folded triangle as in Fig.~\ref{fig:SelfFoldedLocalTriangulation}).
\end{rem}

\begin{prop}[\protect{\cite[Proposition 3.15]{MSW11}}]\label{prop:Prop3.15MSW}
Suppose $T=\{ \tau_1,\dots,\tau_n \}$ is a tagged triangulation of~$\Cn$. Let~$\ga\notch$ denote the arc obtained from~$\ga$ by changing the notching at the puncture~$\puncture$. Let~$T\notch$ denote the tagged triangulation that is obtained from~$T$ by replacing each $\tau \in T$ by~$\tau\notch$. Let~$x(\ga)$ be the $T$-expansion of the cluster variable corresponding to~$\ga$.
Then
\begin{gather*}
x\big({\ga\notch}\big)=x(\ga) |_{x({\tau_i}) \leftarrow x({{\tau_i}\notch})}
\end{gather*}
is the $T\notch$-expansion of the cluster variable corresponding to $\ga\notch$.
\end{prop}

\begin{figure}[t!]
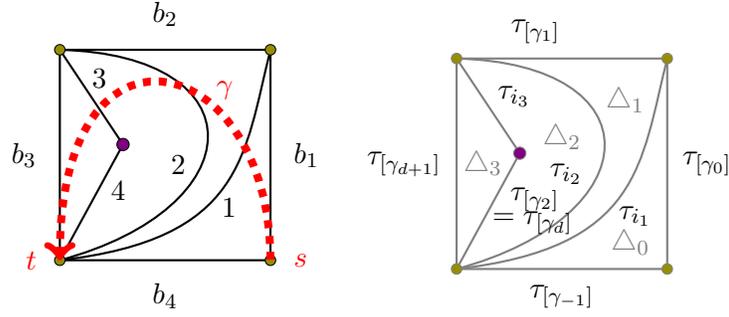

\centering
\TikzToneWithGamma{0.7}
\TikzToneBlackTriangles{0.7}
\caption{A triangulation $T^o$ and an arc $\ga$ of a once-punctured square.
The f\/irst, second, third, and fourth ideal triangles crossed by~$\ga$ are denoted by $\blacktri{0}$, $\blacktri{1}$, $\blacktri{2}$, and~$\blacktri{3}$.}
\label{fig:T1_triangulation}
\end{figure}

\section[$(T^o,\ga)$-path expansion formula]{$\boldsymbol{(T^o,\ga)}$-path expansion formula} \label{sec:tpaths}
We extend Schif\/f\/ler and Thomas' work \cite{Sch10, ST09} to once-punctured disks $\Cn$. Following \cite[Section~3]{Sch10}, we will use the following setup throughout the rest of the paper.
\begin{itemize}\itemsep=0pt
\item
Let $\Surf$ be an unpunctured surface or $\Surf=\Cn$. Let $T^o$ be an ideal triangulation of $\Surf$ and let $\ga\notin T^o$ be an ordinary arc of $\Surf$. Recall that an $\ell$-loop is considered an ordinary arc.

\item Choose an orientation on $\ga$, and let $s^\ga$ and $t^\ga$ be the starting point and the f\/inishing point of $\ga$. Denote by
\begin{gather*}
s^\ga=p_0^\ga, p_1^\ga, p_2^\ga, \ldots, p_{d+1}^\ga = t^\ga
\end{gather*}
the points of intersection of $\ga$ and $T^o$ in order.  Since $\ga$ is considered up to homotopy, we pick a representative so that $d$ is minimal. Let $i_1, i_2, \ldots, i_d$ be such that $\taui{k}^\ga$ is the arc of~$T^o$ containing~$p_k^\ga$. See Fig.~\ref{fig:T1_triangulation}, where $\taui{1}^\ga=1$, $\taui{2}^\ga=2$, and $\taui{3}^\ga=3$, and see Fig.~\ref{fig:T2_triangulated_polygon_original}, where~$\taui{k}^\ga$ ($k=1,\ldots,5$) are labeled $1$, $2$, $\ell$, $r$ and $\ell$.

\item For $k=0, 1,\ldots, d$, let $\ga_k$ denote the segment of $\ga$ from the point $p_k$ to the point $p_{k+1}$, and let $\blacktri{k}^\ga$ denote the (unique) ideal triangle of $T^o$ that $\ga_k$ crosses. When it is clear from the context which arc we mean, we simply write $p_k^\ga$ as $p_k$, $\taui{k}^\ga$ as $\taui{k}$, $\blacktri{k}^\ga$ as~$\blacktri{k}$.

\item
The side/s of $\blacktri{k}$ that is not labeled $\taui{k}$ or $\taui{k+1}$ is labeled as in Figs.~\ref{fig:ways_to_cross_if_k_0} and~\ref{fig:ways_to_cross_if_k_not_0}. In particular, for $k=1,\ldots,d-1$, def\/ine arc $\tgk$ to be
\begin{gather*}
\tgk =
\begin{cases}
\textnormal{the 3rd arc in } \blacktri{k} & \textnormal{if $ \blacktri{k}$ is not self-folded,}\\
\textnormal{the radius in } \blacktri{k}  &\textnormal{if $\blacktri{k}$ is self-folded.}
\end{cases}
\end{gather*}

Def\/ine $\tg{0}$, $\tg{-1}$, $\tg{d}$, and $\tg{d+1}$
as follows:
\begin{itemize}\itemsep=0pt
\item If the ideal triangle $\blacktri{0}$
has three distinct edges, then
 $\blacktri{0}$ is formed by arc $\taui{1}$ and two distinct arcs/boundary edges $\tg{0}$, $\tg{-1}$ (not equal to $\taui{1}$) such that $\tg{0}$, $\tg{-1}$, $\taui{1}$ are arranged in clockwise order around~$\blacktri{0}$.
Similarly, if the ideal triangle $\blacktri{d}$ has three distinct edges, then
it is formed by the arcs~$\taui{d}$
and two distinct arcs/boundary edges~$\tg{d}$,~$\tg{d+1}$ (not equal to~$\taui{d}$) such that~$\tg{d}$,~$\tg{d+1}$,~$\taui{d}$ are arranged in clockwise order around~$\blacktri{d}$.

\item If $\blacktri{0}$ (respectively, $\blacktri{d}$) is self-folded,  then $\tg{0}=\tg{-1}$ (respectively, $\tg{d}=\tg{d+1}$) is the radius.
\end{itemize}
In Fig.~\ref{fig:T2_triangulated_polygon_original}, $\tg{0}=b_1$, $\tg{-1}=b_4$, $\tg{d}=2$, $\tg{d+1}=b_3$.
\end{itemize}

\begin{defin}[Quasi-arc]\label{def:quasi_arc}
If $\tau$ is an ordinary radius of $\Cn$ between a marked point $v$ on the boundary and the puncture $\puncture$,
let an associated \emph{quasi-arc} $\tau'$ be a curve (not passing through $\puncture$) which satisf\/ies the following:
\begin{enumerate}\itemsep=0pt
\item $\tau'$ is between $v$ and a (non-marked) point $\puncture'$ in the vicinity of $\puncture$. (Note that another quasi-arc associated to $\tau$ may use a dif\/ferent point $\puncture''$.)
\item $\tau'$ agrees with the arc $\tau$ outside of a radius-$\epsilon$ disk $D_\epsilon$ around $\puncture$, where $\epsilon$ is chosen small enough so that the intersection of $\tau$ with any other arc is outside of $D_\epsilon$.
\end{enumerate}
If $\tau$ is a peripheral arc, we let the associated quasi-arc be $\tau$ itself.
We label a quasi-arc with the label of the arc that it is associated to.
\end{defin}

By abuse of notation, whenever we say we are going along an arc or a side of an ideal triangle~$\blacktri{k}$ (as part of a~$T^o$-path), we mean traversing an associated quasi-arc.  Note that except for the case when $\tau$ is a radius, no abuse of notation is actually needed.  The following is an extension of the complete $(T^o,\ga)$-path def\/inition as stated in {\cite[Def\/inition~2]{Sch10}}, {\cite[Section~4.1]{MS10}}.
\begin{defin}[complete $(T^o,\ga)$-path]\label{def:Tpath}
A \emph{path} $\om=(\om_1, \om_2, \ldots, \om_{\text{length}(\om)})$ on $T^o$ is a concatenation of \emph{steps},  i.e., oriented quasi-arcs and boundary edges of the ideal triangulation~$T^o$ of~$\Surf$, such that the starting point of a~step~$\om_i$ is the f\/inishing point of the previous step~$\om_{i-1}$. We say that $\om=(\om_1,\dots,\om_{2d+1})$ is a \emph{complete $(T^o,\ga)$-path} if the following axioms hold:
\begin{itemize}\itemsep=0pt
\item [(T1)]
Each even step $\om_{2k}$ ($k=1,\ldots,d$) is a quasi-arc associated to arc $\taui{k}$.
Recall that $\taui{1},\ldots,\taui{d}$ is the sequence of arcs crossed by~$\ga$ in order.
\item[(T2)]
For $k=0,\ldots,d$, each $\om_{2k+1}$ traverses a side of the ideal triangle~$\blacktri{k}$.
In addition, $\om_1$~traverses the edge~$\tg{0}$ or~$\tg{-1}$ (which is adjacent to~$s$)
and $\om_{2d+1}$ traverses the edge~$\tg{d}$ or~$\tg{d+1}$ (which is adjacent to~$t$).
\begin{enumerate}[i)]\itemsep=0pt
\item
Moreover, for $k=1,2, \ldots, d-1$, let $\pkpkone{\om}$
denote the segment of $\om$ starting at the point $p_k$
following $\om_{2k}$, continuing along $\om_{2k+1}$ and $\om_{2(k+1)}$ until the point $p_{k+1}$.
Then the segment $\ga_k$ is homotopic to $\pkpkone{\om}$.
If $S=\Cn$, then we mean homotopy in the disk minus the puncture.

\item
The segment $\ga_0$ is homotopic to the segment $[s,p_1]_{\om}$
of the path starting at
the point $s=p_0$ following $\om_1$ and $\om_{2}$ until the point $p_1$;

\item
The segment $\ga_d$ is homotopic to the segment $[p_d, t]_{\om}$
of the path starting at
the point $p_d$ following $\om_{2d}$ and $\om_{2d+1}$ until the point $p_{d+1}=t$.
\end{enumerate}
\item[(T3)]
The step $\om_{2k+1}$ starts and f\/inishes in the interior of $\blacktri{k}$ or at a boundary marked point. This means that, if $\om_{2k+1}$ goes along a quasi-arc $\tau'$ associated to a radius, $\tau'$ must be chosen so that its endpoint $\puncture'$ near the puncture is located in the interior of $\blacktri{k}$.
\end{itemize}
\end{defin}
See Examples \ref{example:T1_Tpath_computation} and~\ref{example:nine_paths}.

It is clear that this def\/inition agrees with the complete $T^o$-paths of \cite{MS10, Sch10} for unpunctured surfaces.
For short, we will refer to a complete $(T^o,\ga)$-path as simply a $(T^o,\ga)$-path (or a~$T^o$-path) for the rest of this paper.

\begin{rem}
Per (T2), a $(T^o,\ga)$-path is homotopic to $\ga$. It is possible to have $\om_j = \om_{j+1}$ and, if $j$ is odd, to have $\om_j = \om_{j+1} = \om_{j+2}$. However, since for each $k$ we have $\om_{2k}=\taui{k}$ and $\om_{2(k+1)}=\taui{k+1}$ by (T1) but $\taui{k} \neq \taui{k+1}$, no more than three consecutive steps can coincide.
 \end{rem}

\begin{figure}[t!]
\centering

\mbox{\subfigure[Ordinary triangle.]
{
\TikzTriangleCrossesOneEdge{0.4}
\label{fig:TikzTriangleCrossesOneEdge}
}\quad
\subfigure[Two vertices where arc $1^\ga$ or $d^\ga$ is a loop.]
{
\TikzTwoVertexTriangleCrossesOneEdgeLoop{0.8}
\label{fig:TikzTwoVertexTriangleCrossesOneEdgeLoop}
}\quad
\subfigure[Two vertices where arc $1^\ga$ or $d^\ga$ is not a loop.]
{
\TikzTwoVertexTriangleCrossesOneEdgeNotLoop{0.8}
\label{fig:TikzTwoVertexTriangleCrossesOneEdgeNotLoop}
}\quad
\subfigure[Self-folded triangle.]
{
\TikzSelfFoldedTriangleCrossesNoose{0.8}
\label{fig:TikzSelfFoldedTriangleCrossesNoose}
}}
\caption{Ways for $\ga$ to cross the ideal triangle $\blacktri{0}$ or $\blacktri{d}$.}
\label{fig:ways_to_cross_if_k_0}

\mbox{\subfigure[Ordinary triangle.]
{
\TikzTriangleCrossesTwoEdges{0.4}
\label{fig:TikzTriangleCrossesTwoEdges}
}\quad
\subfigure[Two vertices where arc $\tgk$ is a loop.]
{
\TikzTwoVertexTriangleCrossesTwoEdgesNoLoop{0.8}
\label{fig:TikzTwoVertexTriangleCrossesTwoEdgesNoLoop}
}\quad
\subfigure[Two vertices where arc $k+1$ is a loop.]
{
\TikzTwoVertexTriangleCrossesTwoEdgesKOneisLoop{0.8}
\label{fig:TikzTwoVertexTriangleCrossesTwoEdgesKOneisLoop}
}\quad
\subfigure[$\ga$ crossing a self-folded triangle's radius in the counterclockwise direction.]
{\hspace{4mm}
\TikzSelfFoldedTriangleCrossesRadiusNoose{0.8}\hspace{4mm}
\label{fig:TikzSelfFoldedTriangleCrossesRadiusNoose}
}}
\caption{Ways for $\ga$ to cross an ideal triangle $\blacktri{k}$ for $k=1,\dots,k-1$.}
\label{fig:ways_to_cross_if_k_not_0}

\end{figure}

\begin{defin}[backtrack cycle, non-backtrack cycle, quasi-backtrack]
\label{def:backtrack}
Let $\ga\notin T^o$ and
let $\om=(\om_1,\dots,\om_{2d+1})$ be a $(T^o,\ga)$-path. Let $\tau$ be an arc of $T^o$, and let $(\om_j,\om_{j+1})$ be a pair of consecutive quasi-arcs going along $\tau$.
We say that $(\om_j,\om_{j+1})$ is a~\emph{cycle} if the starting point of~$\om_j$ coincides with the ending point of~$\om_{j+1}$.
\begin{enumerate}[i)]\itemsep=0pt
\item A cycle $(\om_j,\om_{j+1})$ is called a \emph{backtrack}, denoted by $(\ttBak)$, if it is contractible.
\item A cycle $(\om_j,\om_{j+1})$ is called a \emph{non-backtrack}, denoted by $(\ttNon)$, otherwise.
\end{enumerate}

In the case that $\tau$ is a radius between the puncture $\puncture$ and a marked point $v$ on the boundary,
we say that $(\om_j,\om_{j+1})$, denoted by $(\ttQuasi)$, is a \emph{quasi-backtrack} if it is a concatenation of two quasi-arcs $\puncture' \leadsto v \leadsto \puncture''$ where $\puncture'$ and $\puncture''$ are distinct points in the vicinity of $\puncture$.
\end{defin}

\begin{rem}\label{rem:pair}
Assume $(\om_j,\om_{j+1})$ is a pair of steps going along $\tau$, as in Def\/inition \ref{def:backtrack}.
\begin{enumerate}[1)]\itemsep=0pt
\item \label{rem:pair:peripheral} Suppose $\tau$ is peripheral.
	\begin{enumerate}[i)]\itemsep=0pt
	\item Then $(\om_j,\om_{j+1})$ is a backtrack if and only if $\om_j$ and $\om_{j+1}$ are opposite orientations of~$\tau$. In particular, if $\tau$ has two distinct endpoints, $(\om_j,\om_{j+1})$ must be a backtrack cycle.
	\item If $\tau$ is an $\ell$-loop, the concatenation of two steps going along the same orientation of~$\tau$ would form a non-contractible cycle.
However, we show in Proposition~\ref{prop:pair}(\ref{prop:pair:peripheral}) that it is impossible for $(\om_j,\om_{j+1})$ to go along the same orientation of an $\ell$-loop twice as part of a $(T^o,\ga)$-path.
	\end{enumerate}
\item Suppose $\tau$ is a radius between the puncture $\puncture$ and a marked point $v$ on the boundary.
\begin{enumerate}[a)]\itemsep=0pt
\item
\label{rem:pair:v_p_v}
Suppose $\om_j$ begins at $v$, so that $(\om_j,\om_{j+1})$ is a concatenation of two associated quasi-arcs
$v\leadsto \puncture'\leadsto v$.
	\begin{enumerate}[i)]\itemsep=0pt
	\item \label{rem:pair:v_p_v:backtrack}
	Then $(\om_j,\om_{j+1})$ is a backtrack cycle if it is contractible to $v$.
	See Fig.~\ref{fig:backtrack_nonbacktrack_backtrack} for the case where $T^o$ has a self-folded triangle and $\tau$ is its radius and Fig.~\ref{fig:TikzBacktracks} for the case where $\tau$ is not the only radius of $T^o$.
	\item \label{rem:pair:v_p_v:non}
	If $(\om_j,\om_{j+1})$ is not contractible, then due to (T2) it must be homotopic to a loop which goes around the puncture once.
	We say that $(\om_j,\om_{j+1})$ is a \emph{counterclockwise non-backtrack} if it goes counterclockwise (Fig.~\ref{fig:backtrack_nonbacktrack_counterclock}), and a \emph{clockwise non-backtrack} if it goes clockwise (Fig.~\ref{fig:backtrack_nonbacktrack_clock}).
	We show in Proposition \ref{prop:pair}(\ref{prop:pair:v_p_v}) that in this case, $\tau$ must be the only radius of $T^o$, i.e., $T^o$ contains a self-folded triangle.
	\end{enumerate}

\item \label{rem:pair:p_v_p}
Suppose $\om_j$ ends at $v$, so that $(\om_j,\om_{j+1})$ is a concatenation of two associated quasi-arcs
$\puncture'\leadsto v \leadsto  \puncture''$,
where $\puncture'$ and $\puncture''$ are in the vicinity of the puncture.
We show in Proposi\-tion~\ref{prop:pair}(\ref{prop:pair:p_v_p}) that~$\puncture'$ and~$\puncture''$ must be distinct, so $(\om_j,\om_{j+1})$ is a~quasi-backtrak.
See Fig.~\ref{fig:TikzQuasiBacktracks}.
\end{enumerate}
\end{enumerate}
\end{rem}

\begin{figure}[t!]
\centering

\subfigure[Backtrack cycle.]
{
\TikzBacktrackNonbacktrackBacktrack{0.7}
\label{fig:backtrack_nonbacktrack_backtrack}
}
\quad
\subfigure[Counterclockwise non-backtrack cycle.]
{\hspace{2mm}
\TikzBacktrackNonbacktrackCounterclock{0.7}\hspace{2mm}
\label{fig:backtrack_nonbacktrack_counterclock}
}\quad
\subfigure[Clockwise \mbox{non-backtrack} \mbox{cycle.}]
{
\TikzBacktrackNonbacktrackClock{0.7}
\label{fig:backtrack_nonbacktrack_clock}
}
\vspace{1mm}

\caption{The three possibilities of a pair $(\om_j,\om_{j+1})$ of consecutive steps along a radius $\tau$ if $T^o$ has a~self-folded triangle with $r$ as the radius.}
\label{fig:backtrack_nonbacktrack}

\subfigure[Backtrack cycle $(\om_{2k+1},\om_{2k})$ or $(\om_{2k},\om_{2k+1})$.]
{\TikzBacktrackEvenFirst{0.65}
\quad
\TikzBacktrackOddFirst{0.65}\label{fig:TikzBacktracks}}
\quad
\subfigure[Quasi-backtrack $(\om_{2k},\om_{2k+1})$ or $(\om_{2k+1},\om_{2k})$.]
{
\TikzQuasiBacktrackEvenFirst{0.65}
\quad
\TikzQuasiBacktrackOddFirst{0.65}
\label{fig:TikzQuasiBacktracks}
}
\vspace{-1mm}

\caption{A pair $(\om_j,\om_{j+1})$ of consecutive steps along a radius $\tau$ is either a quasi-backtrack or a~backtrack cycle if every ideal triangle of $T^o$ is an ordinary triangle.}

\end{figure}

\begin{defin}[Laurent monomial from a $T^o$-path]\label{def:x_tpath}\label{def:laurent_monomial_tpath}
We identify each step with the label of the quasi-arc/boundary edge which it traverses and def\/ine the Laurent monomial
$x(\om)$ corresponding to a complete $(T^o,\ga)$-path $\om$ by
\begin{gather*}
x(\om) = \prod_{i \textnormal{ odd}} x_{\om_i}
\prod_{i \textnormal{ even}} x^{-1}_{\om_i}.
\end{gather*}
\end{defin}

\begin{rem}
Two or more $(T^o,\ga)$-paths may correspond to the same Laurent monomial, e.g., see Example~\ref{example:nine_paths}.
For each $(T^o,\ga)$-path $\om$, the denominator of $x(\om)$, before reducing, is equal to $x_{i_1}, x_{i_2}, \ldots, x_{i_d}$ which corresponds to the
arcs $\tau_{i_1}, \ldots, \tau_{i_d}$ of $T^o$ which cross~$\ga$.
\end{rem}

\begin{theorem}[$T^o$-path formula for $\Cn$, an extension of {\cite[Theorem~3.1]{Sch10}}, {\cite[Theorem~3.2]{ST09}}]
\label{thm:tpath_expansion_formula}
Let $T^o$ be an ideal triangulation of $\Cn$, let $\ga\notin T^o$ be an ordinary arc of $\Cn$, and let $x_\ga$ denote the corresponding element in the cluster algebra which arises from $\Cn$ $($see Theorem {\rm \ref{thm:fst_thm})}. Then
\begin{gather*}
x_\ga = \sum_{\om} x(\om),
\end{gather*}
where the sum is taken over all $(T^o,\ga)$-paths. The formula does not depend on the choice of orientation on~$\ga$.
\end{theorem}

The proof of Theorem~\ref{thm:tpath_expansion_formula} is given in Section~\ref{subsection:proof_of_tpath_expansion_formula}. Note that, since $x_\ell := x_r x_{r\notch}$, Theo\-rem~\ref{thm:tpath_expansion_formula} also provides a formula for the cluster variable associated to every tagged arc of $\Cn$.

\begin{exam}\label{example:T1_Tpath_computation}
The following are the f\/ive $(T^o,\ga)$-paths for the
situation of Fig.~\ref{fig:T1_triangulation}. Note that $\ga$ crosses $1$, $2$, $3$ in order so
$\taui{1}=1, \taui{2}=2$, and $\taui{3}=3$ in this example.
\begin{alignat*}{4}
& 1)~~({b_1}, 1, \underline{2, 2}, 4, 3, {b_3}), \qquad && 3)~~({b_4}, 1, {b_2}, \underline{2, 2}, 3, 4), \qquad && 5)~~({b_4}, \underline{1, 1}, 2, \underline{3, 3}, {b_3}). &\\
& 2)~~({b_1}, 1, \underline{2, 2}, 2, 3, 4), \qquad && 4)~~({b_4}, 1 , {b_2}, 2, 4, 3, {b_3}), \qquad && &
\end{alignat*}
Theorem \ref{thm:tpath_expansion_formula} thus implies that
\begin{gather*}
x_\ga = \frac{x_{b_1} x_2   x_4  x_{b_3}
+ x_{b_1}  x_2  x_2  x_4 +x_{b_4}   x_{b_2}  x_2  x_4
+ x_{b_4}  x_{b_2}  x_4  x_{b_3}
+ x_{b_4}  x_1   x_3  x_{b_3}}{x_1 x_2 x_3},
\end{gather*}
where $x_{b_j}=1$ for each boundary edge $b_j$. 
\end{exam}

\begin{exam}\label{example:nine_paths}
The nine $(T^o,\ga)$-paths from Fig.~\ref{fig:T2_triangulated_polygon_original}
are as follows (with the backtracks underlined):
\begin{alignat*}{3}
& (1)~~(b_4, 	1, 	b_2, 	2, 	b_3, 		\ell, 	\underline{r, 		r}, 	\underline{\ell, 		\ell}, 	b_3), \qquad &&
(6)~~(b_4, 	1, 	b_2, \underline{2, 	2}, 		\underline{\ell, 	\ell},	\rrNon, 	\ell, 	b_3), &\\
& (2)~~(b_4, 	\underline{1,	1}, 	2, 	\underline{\ell, 		\ell}, 	\underline{r, 		r}, 	\underline{\ell, 		\ell}, 	b_3), \qquad && (7)~~(b_4, 	1, 	b_2, \underline{2, 	2}, 		\ell, 	\rrNon,	\underline{\ell, 		\ell}, 	b_3), &\\
& (3)~~(b_1, 	1, 	\underline{2, 	2}, 	2, 		\underline{\ell, 	\ell}, 		\underline{r,	r}, 		\ell, 	2), \qquad &&
(8)~~(b_1,1,	\underline{2, 2}, 	2, 		\underline{\ell, 	\ell}, 		\rrNon, \ell, 	b_3), &\\
& (4)~~(b_1, 	1, 	\underline{2, 	2}, 	b_3, 		\ell, 	\underline{r, 		r},	\underline{\ell, 		\ell}, 	b_3), \qquad &&
(9)~~(b_1, 	1, 	\underline{2, 	2}, 	2, 		\ell, 	{\bf {r}, 	r},	\underline{\ell, 		\ell}, 	b_3). &\\
& (5)~~(b_4, 	1, 	b_2, 	\underline{2, 	2}, 		\underline{\ell, 	\ell}, 		\underline{r,	r}, 		\ell, 	2),\qquad &&&
\end{alignat*}
The last four of these $(T^o,\ga)$-paths  
contain a (counterclockwise) non-backtrack $(\rrNon)$ 
and are illustrated in Fig.~\ref{fig:TtwoTpaths}.
They are drawn so that the backtrack cycles $(\underline{2,2})$ and $(\underline{\ell,\ell})$ are ignored.

We apply Theorem \ref{thm:tpath_expansion_formula} and replace each ${(x_r x_\ell)}/{(x_\ell x_r x_\ell)}$ with $1/x_\ell$ and each  $x_\ell$ with $x_r x_{r\notch}$ to get{\samepage
\begin{gather*}
\begin{split}
& x_\ga=
\big(x_{b_4} x_{ b_2}  x_{b_3}  x_{b_3} +
  x_{b_4} x_1x_r x_{r\notch} x_{b_3} +
 x_{b_1} x_2^3 +
 x_{b_1} x_{b_3}^2 x_2 +
x_{b_4} x_{b_2} x_2^2  \\
& \hphantom{x_\ga=}{}+ 2   x_{b_4} x_{b_2} x_{b_3} x_2 +
 2  x_{b_1} x_{b_3} x_2^2\big)/( x_1 x_2  x_r x_{r\notch} ),
 \end{split}
\end{gather*}
where $x_{b_j}=1$ for each boundary edge~$b_j$.}
\end{exam}

\begin{exam} \label{example:eleven_paths}
Fig.~\ref{fig:TikzNoose} illustrates the four $(T^o,\ga)$-paths for the case where the ideal triangulation $T^o$ and the $\ell$-loop $\ga\notin T^o$ are illustrated in Fig.~\ref{fig:TikzNooseGammaIdealTriangulation}. In this example, no $(T^o,\ga)$-path contains any cycle $(\om_j,\om_{j+1})$.

Fig.~\ref{fig:TikzNotNoose} illustrates three of the f\/ive $(T^o,\ga)$-paths for the case where the ideal triangulation $T^o$ and the arc $\ga\notin T^o$ are illustrated in Fig.~\ref{fig:TikzNotNooseGammaIdealTriangulation}. The backtrack cycle $(\underline{3,3})$ in Fig.~\ref{fig:TikzNotNooseTpathThree} is not drawn.
The other two $(T^o,\ga)$-paths are $(0,1,b_2,2,0,3,b_3)$ and $(b_1,1,\overline{2,2},0,3,b_3)$.
\end{exam}

\newcommand\TikzFourTpathsSize{0.35}
\newcommand\TikzNooseTpathsSize{0.41}
\begin{figure}[t!]
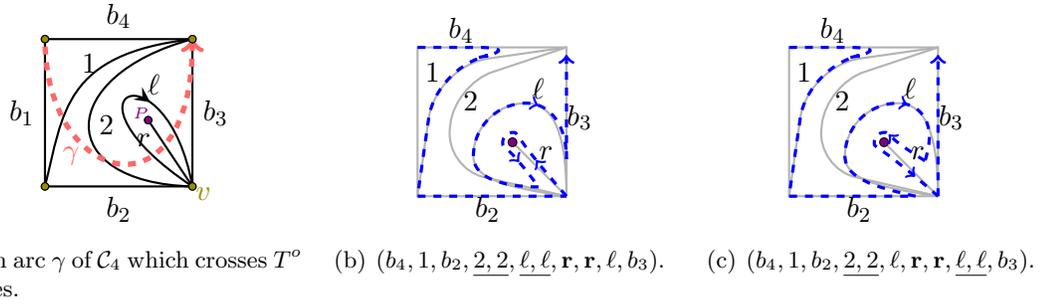
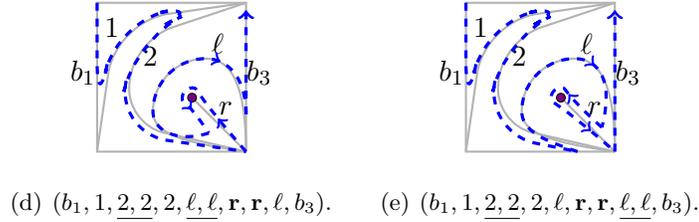

\centering

\subfigure[An arc $\ga$ of $\Cnn{4}$ which crosses~$T^o$ $5$ times.]{
\hspace{6mm}\TikzTtwoIdealTriangulationWithNumberLabels{1.4*\TikzFourTpathsSize}\hspace{6mm}
\label{fig:T2_triangulated_polygon_original}
}\quad
\subfigure[$(b_4,1,b_2,\underline{2,2},\underline{\ell,\ell},\rrNon,\ell,b_3)$.]
{\hspace{7.5mm}\TikzTpathSix{\TikzFourTpathsSize}\hspace{7.5mm}\label{fig:path_six}}
\quad
\subfigure[$(b_4,1,b_2,\underline{2,2},\ell,\rrNon,\underline{\ell,\ell},b_3)$.]
{\hspace{7mm}\TikzTpathSeven{\TikzFourTpathsSize}\hspace{7mm}}\\
\subfigure[$(b_1,1,\underline{2,2},2,\underline{\ell,\ell},\rrNon,\ell,b_3)$.]
{\hspace{6mm}\TikzTpathEight{\TikzFourTpathsSize}\hspace{6mm}}
\quad
\subfigure[$(b_1,1,\underline{2,2},2,\ell,\rrNon,\underline{\ell,\ell},b_3)$.]
{\hspace{6mm}\TikzTpathNine{\TikzFourTpathsSize}\hspace{6mm}}
\caption{The four $(T^o,\ga)$-paths (of $T^o$ and $\ga$ from Fig.~\ref{fig:T2_triangulated_polygon_original}) which contain a (counterclockwise) non-backtrack cycle $(\rrNon)$.
All backtracks $(\underline{2,2})$ and $(\underline{\ell,\ell})$ have been omitted, and steps are not drawn exactly along the arcs/boundary edges for illustration purposes.
}\label{fig:TtwoTpaths}
\end{figure}

\begin{figure}[t!]
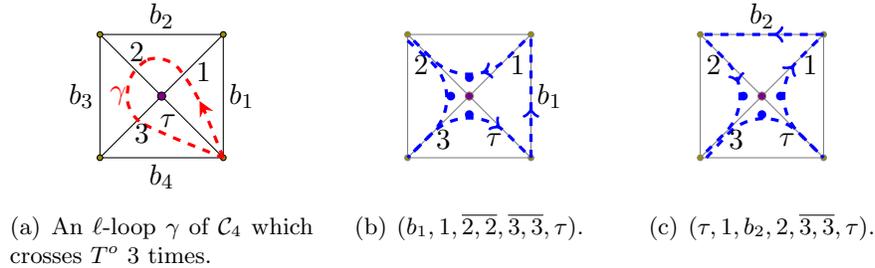
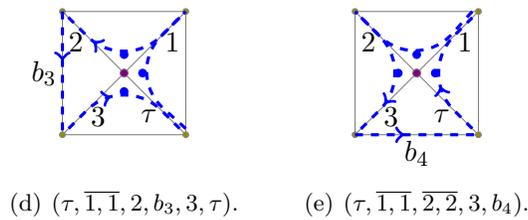
\centering
\subfigure[An $\ell$-loop $\ga$ of $\Cnn{4}$ which crosses $T^o$ $3$ times.]{
\hspace{4mm}\TikzNooseGammaIdealTriangulation{\TikzNooseTpathsSize}\hspace{4mm}
\label{fig:TikzNooseGammaIdealTriangulation}
}\quad
\subfigure[$(b_1,1,\overline{2,2},\overline{3,3},\tau)$.]
{\hspace{2mm}\TikzNooseTpathOne{\TikzNooseTpathsSize}\hspace{2mm}\label{fig:TikzNooseTpathOne}}
\quad \subfigure[$(\tau,1,b_2,2,\overline{3,3},\tau)$.]
{\hspace{2mm}\TikzNooseTpathTwo{\TikzNooseTpathsSize}\hspace{2mm}}  \\
\subfigure[$(\tau,\overline{1,1},2,b_3,3,\tau)$.]
{\hspace{2mm}\TikzNooseTpathThree{\TikzNooseTpathsSize}\hspace{2mm}}
\quad \subfigure[$(\tau,\overline{1,1},\overline{2,2},3,b_4)$.]
{\hspace{2mm}\TikzNooseTpathFour{\TikzNooseTpathsSize}\hspace{2mm}}

\caption{The four $(T^o,\ga)$-paths of the ideal triangulation $T^o$ and the $\ell$-loop $\ga$ of Fig.~\ref{fig:TikzNooseGammaIdealTriangulation}.
Each path contains a quasi-backtrack $(\overline{1,1})$, $(\overline{2,2})$, or $(\overline{3,3})$.
}\label{fig:TikzNoose}
\end{figure}

\begin{figure}[t!]
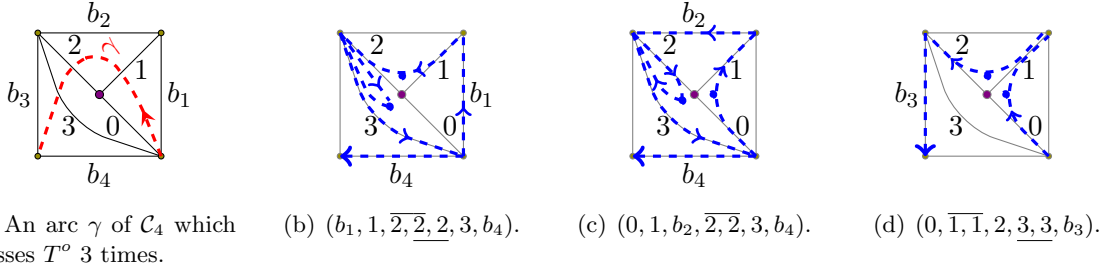
\centering
\subfigure[An arc $\ga$ of $\Cnn{4}$ which crosses $T^o$ $3$ times.]{
\hspace{2mm}\TikzNotNooseGammaIdealTriangulation{\TikzNooseTpathsSize}\hspace{2mm}
\label{fig:TikzNotNooseGammaIdealTriangulation}}
\quad \subfigure[$(b_1,1,\overline{2,\,}\underline{\overline{2},2},3,b_4)$.]
{\hspace{2mm}\TikzNotNooseTpathOne{\TikzNooseTpathsSize}\hspace{2mm}}
\quad \subfigure[$(0,1,b_2,\overline{2,2},3,b_4)$.]
{\hspace{2mm}\TikzNotNooseTpathTwo{\TikzNooseTpathsSize}\hspace{2mm}}
\quad \subfigure[$(0,\overline{1,1},2,\underline{3,3},b_3)$.]
{\hspace{2mm}\TikzNotNooseTpathThree{\TikzNooseTpathsSize}\hspace{2mm}\label{fig:TikzNotNooseTpathThree}}
\caption{Three of the f\/ive $(T^o,\ga)$-paths of the ideal triangulation $T^o$ and the arc $\ga$ of Fig.~\ref{fig:TikzNotNooseGammaIdealTriangulation}.
}\label{fig:TikzNotNoose}
\end{figure}

\begin{figure}[t!]
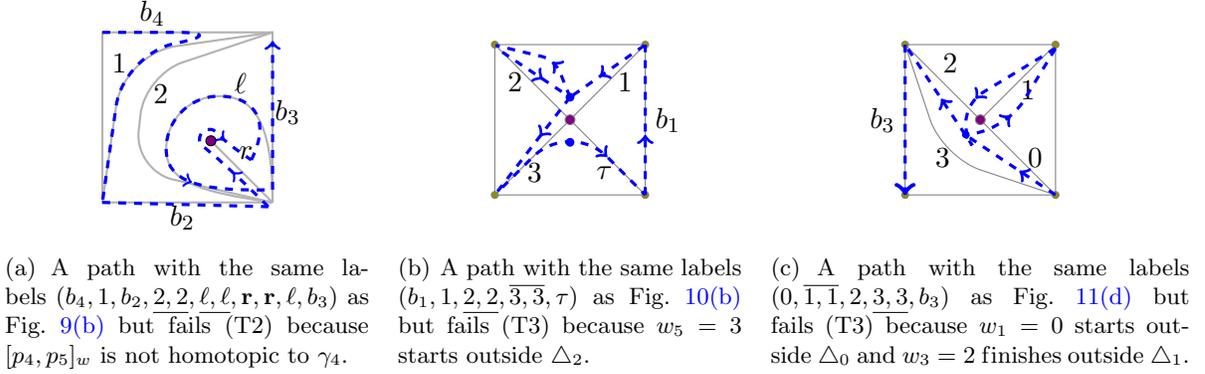
\centering
\renewcommand\TikzNooseTpathsSize{0.5}
\subfigure[A path with the same labels $(b_4,1,b_2,\underline{2,2},\underline{\ell,\ell},\rrNon,\ell,b_3)$ as Fig.~\ref{fig:path_six}
 but fails (T2) because ${[p_4,p_5]_{\om}}$ is not homotopic to~$\ga_4$.]
{\hspace{1mm}\TikzTpathSixBad{0.4}\hspace{1mm}
\label{fig:path_six_bad}
}\quad
\subfigure[A path with the same labels $(b_1,1,\underline{2,2},\overline{3,3},\tau)$ as Fig.~\ref{fig:TikzNooseTpathOne} but fails~(T3) because $\om_5=3$ starts outside~$\blacktri{2}$.]{
\hspace{1mm}\TikzNooseTpathOneBad{\TikzNooseTpathsSize}\hspace{1mm}
\label{fig:TikzNooseTpathOneBad}
}\quad
\subfigure[A path with the same labels $(0,\overline{1,1},2,\underline{3,3},b_3)$ as Fig.~\ref{fig:TikzNotNooseTpathThree} but fails~(T3) because $\om_{1}=0$ starts outside~$\blacktri{0}$ and~$\om_3=2$ f\/inishes outside~$\blacktri{1}$.]{
\hspace{1mm}\TikzNotNooseTpathThreeBad{\TikzNooseTpathsSize}\hspace{1mm}
\label{fig:TikzNotNooseTpathThreeBad}
}

\caption{Examples of \emph{non}-$(T^o,\ga)$-paths for the situations in Figs.~\ref{fig:T2_triangulated_polygon_original}, \ref{fig:TikzNooseGammaIdealTriangulation}, and~\ref{fig:TikzNotNooseGammaIdealTriangulation}.
Each path is homotopic to $\ga$ and satisf\/ies~(T1) but fails~(T2) or~(T3).}
\label{fig:non_tpaths}
\end{figure}

\subsection{A non-backtrack cycle can only go along a self-folded triangle's radius}
We prove the assertions from Remark \ref{rem:pair}.

\begin{defin}[Crossing a self-folded triangle]
Suppose that $T^o$ contains a self-folded triangle with radius $r$.
We say that $\ga$ crosses $r$ in the \emph{counterclockwise direction} (respectively, \emph{clockwise direction}) if
it matches (respectively, if $\ga$ has the opposite orientation of) Fig.~\ref{fig:TikzSelfFoldedTriangleCrossesRadiusNoose}.
\end{defin}

\begin{prop}
\label{prop:pair}
Let $T^o$ be an ideal triangulation of $\Cn$.
Suppose $\ga$ is an ordinary arc which crosses $\tau\in T^o$, and let $\om=(\om_1,\dots,\om_{2d+1})$ be a $(T^o,\ga)$-path. Suppose $(\om_j,\om_{j+1})$ is a pair of steps both going along quasi-arcs associated to~$\tau$.
\begin{enumerate}[$1)$]\itemsep=0pt
\item \label{prop:pair:peripheral}
Suppose $\tau\in T^o$ is an $\ell$-loop $\ell$. Then $(\om_j,\om_{j+1})$ are two opposite orientations of $\ell$, so that $(\om_j,\om_{j+1})$ is a backtrack cycle. See Remark~{\rm \ref{rem:pair}(\ref{rem:pair:peripheral})}.
\item Suppose $T^o$ does not contain any self-folded triangle,
and let $\tau\in T^o$ be a radius between the puncture~$\puncture$ and a marked point $v$ on the boundary.
\begin{enumerate}[$a)$]\itemsep=0pt
\item \label{prop:pair:v_p_v}
Suppose that $\om_j$ begins at $v$, so that $(\om_j,\om_{j+1})$ is a concatenation of two quasi-arcs
$v \leadsto \puncture' \leadsto v$
where $\puncture'$ is in the vicinity of the puncture.
Then $(\om_j,\om_{j+1})$ is a backtrack cycle.
See Remark~{\rm \ref{rem:pair}(\ref{rem:pair:v_p_v})}.

\item \label{prop:pair:p_v_p}
Suppose that $\om_j$ ends at $v$, so that $(\om_j,\om_{j+1})$ is a concatenation of two quasi-arcs
$\puncture' \leadsto v \leadsto \puncture''$
where $\puncture'$ and $\puncture''$ are in the vicinity of the puncture.
Then $\puncture' \neq \puncture''$, and hence $(\om_j,\om_{j+1})$ is a quasi-backtrack.
See Remark~{\rm \ref{rem:pair}(\ref{rem:pair:p_v_p})}.
 \end{enumerate}
 \end{enumerate}
\end{prop}
\begin{rem}
In Examples \ref{example:T1_Tpath_computation} and \ref{example:eleven_paths},
 none of the $(T^o,\ga)$-paths include a non-backtrack cycle along a radius
since $T^o$ contains no self-folded triangle (see Proposition \ref{prop:pair}(\ref{prop:pair:v_p_v})).
In contrast, four of the $(T^o,\ga)$-paths in Example \ref{example:nine_paths} for $(T^o,\ga)$-paths contains a non-backtrack cycle.
Furthermore, no cycle $\puncture' \leadsto v \leadsto \puncture'$ can appear in general (see Proposition \ref{prop:pair}(\ref{prop:pair:p_v_p})).  See Fig.~\ref{fig:non_tpaths} and Example \ref{example17} for examples of non-$(T^o,\ga)$-paths.
\end{rem}

\begin{figure}[t!]
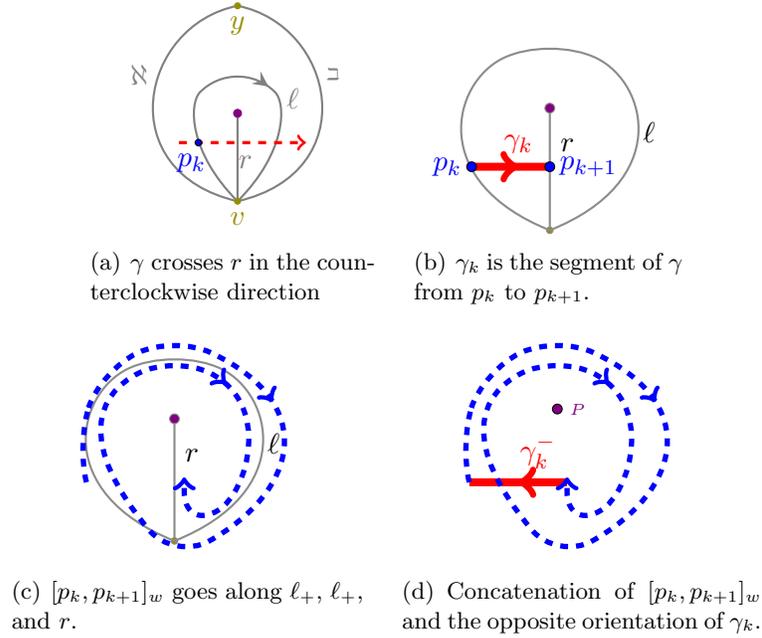
\centering
\newcommand\TikzProveLPlusLPlusIsImpossibleSize{0.65}
\subfigure[$\ga$ crosses $r$ in the counterclockwise direction]{
\TikzCrossingCounterclockwiseLRAB{\TikzProveLPlusLPlusIsImpossibleSize}\label{fig:cross_counterclock_to_prove_Lplus_Lplus_is_impossible}} \quad
\subfigure[$\ga_k$ is the segment of $\ga$ from $p_k$ to $p_{k+1}$.]
{\TikzGammaK{\TikzProveLPlusLPlusIsImpossibleSize}\label{fig:TikzGammaK}}\\
\subfigure[$\pkpkone{\om}$ goes along~$\ell_+$,~$\ell_+$, and~$r$.]
{\hspace{4mm}\TikzEllPlusEllPlus{\TikzProveLPlusLPlusIsImpossibleSize}\hspace{4mm}\label{fig:TikzEllPlusEllPlus}} \quad
\subfigure[Concatenation of $\pkpkone{\om}$ and the op\-po\-site orientation of~$\ga_k$.]
{\hspace{6mm}\TikzGammaKConcatenate{\TikzProveLPlusLPlusIsImpossibleSize}\hspace{6mm}\label{fig:TikzGammaKConcatenate}}

\caption{$(\ell_+, \ell_+)$ is not a valid $(T^o,\ga)$-subpath.}
\end{figure}

\begin{proof}[Proof of Proposition \ref{prop:pair}]
First, we point out that the pair of steps $(\om_j, \om_{j+1})$ must be in
the middle of the $(T^o,\ga)$-path, i.e.,
\begin{gather}
\label{eq:omj_omjplus1_not_the_ends}
 \text{the f\/irst step of $\om$ cannot be $\om_j$ and the last step of $\om$ cannot be $\om_{j+1}$.}
 \end{gather}
 To see this, note that, by (T2), $\om_1$ goes along $\tg{0}$ or $\tg{-1}$, and $\om_{2d+1}$ goes along $\tg{d}$ or $\tg{d+1}$.
But, as illustrated in Fig.~\ref{fig:ways_to_cross_if_k_0}, $\tg{j}\neq \taui{1}$ if $j\in\{0,-1\}$, and $\tg{j}\neq \taui{d}$ if $j\in\{d,d+1\}$.
Since, by~(T1), $\om_2$ goes along $\taui{1}$ and $\om_{2d}$ goes along $\taui{d}$, we must have $1<j$ and $j+1<d$.

We prove case (\ref{prop:pair:peripheral}):
Suppose $T^o$ has a self-folded triangle with radius $r$ and $\ell$-loop $\ell$, and $\tau=\ell$.
For the sake of argument, suppose that $\om_{j}$ and $\om_{j+1}$ are the same orientation of $\ell$.
Note that every arc of $T^o$ other than $\ell$ has two distinct endpoints.

-- First, suppose ${j}=2k$ is even.
Per (\ref{eq:omj_omjplus1_not_the_ends}), another even-indexed step $\om_{2(k+1)}$ follows $(\om_{2k},\om_{2k+1})$. Since $\om_{2k+1}$ ends at $v$, the step $\om_{2k+2}$ (crossing $\ga$) must go from $v$ along a~dif\/ferent arc $\lambda$ to a dif\/ferent marked point.
Using the labels of Figs.~\ref{fig:cross_counterclock_to_prove_Lplus_Lplus_is_impossible}, $\lambda$ can be either~$r$,~$\a$, or~$\b$.
We check all possible cases against~(T1) and~(T2) but only show one of the arguments here.
Suppose $\lambda=r$ and~$\ga$ is a peripheral arc crossing $\ell$, $r$, and $\ell$ again in the counterclockwise direction as in
Fig.~\ref{fig:cross_counterclock_to_prove_Lplus_Lplus_is_impossible} with $j=2k$ and $j+1=2k+1$ so that $\om_{2k}=\ell$ and $\om_{2k+2}=r$.
For contradiction, suppose $(\om_{2k},\om_{2k+1})$ go along the clockwise orientation of $\ell$ twice.
Recall that~$\ga_{k}$ is the segment of~$\ga$ from $p_{k}$ (a point in the interior of $\ell$) and $p_{k+1}$ (a point in the interior of~$r$) as in Fig.~\ref{fig:TikzGammaK}.
By~(T2), the segment $\ga_{k}$ must be homotopic to $\pkpkone{\om}$ (see Def\/inition \ref{def:Tpath} and Fig.~\ref{fig:TikzEllPlusEllPlus}).
However, as illustrated in Fig.~\ref{fig:TikzGammaKConcatenate}, the concatenation of $\pkpkone{\om}$ and the opposite orientation $\ga_k^-$ of $\ga_{k}$ is not contractible.
Hence~$\ga_{k}$ is not homotopic to~$\pkpkone{\om}$.

-- Second, suppose $j=2k+1$ is odd. Again, per (\ref{eq:omj_omjplus1_not_the_ends}), there is another even-indexed step $\om_{2k}$ (crossing $\ga$, and going from a dif\/ferent marked point to $v$ along a peripheral arc) which precedes $(\om_{2k+1},\om_{2k+2})$.
Using the same logic as in the previous paragraph, we show that $(\om_{2k+1},\om_{2k+2})$ must go along two opposite orientations of $\ell$ in order for $\ga_{k}$ to be homotopic to $\pkpkone{\om}$.

We prove case (\ref{prop:pair:v_p_v}):
Suppose $\tau$ is a radius and a side of a regular triangle (hence $T^o$ has no self-folded triangle),
and every arc of $T^o$ has two distinct endpoints.

-- First, suppose ${j}=2k$ is even.
Per (\ref{eq:omj_omjplus1_not_the_ends}), another even-indexed step $\om_{2(k+1)}$ follows $(\om_{2k},\om_{2k+1})$.
Since $\om_{2k+1}$ ends at~$v$, the step $\om_{2k+2}$ (crossing~$\ga$) must be a peripheral step~$\lambda$ to a~dif\/ferent point. But $(\om_{2k},\om_{2k+1})$ would need to be contractible in order for~$\ga_k$ to be homotopic to~$\pkpkone{\om}$.

-- Second, suppose ${j}=2k+1$ is odd.
Again, per (\ref{eq:omj_omjplus1_not_the_ends}), there is another even-indexed step~$\om_{2k}$ (crossing $\ga$, and going from a dif\/ferent marked point to $v$ along a peripheral arc $\lambda$) which precedes $(\om_{2k+1},\om_{2k+2})$.
But $(\om_{2k+1},\om_{2k+2})$ would need to be contractible in order for~$\ga_k$ to be homotopic to $\pkpkone{\om}$.

\looseness=-1
We prove case (\ref{prop:pair:p_v_p}):
First suppose $j=2k+1$ is odd. Then $\puncture'$ is in the interior of~$\blacktri{k}$ since~$\om_{2k+1}$ must start in the interior of~$\blacktri{k}$ by~(T3).
The next step~$\om_{2k+2}$ then goes along~$\tau$ from~$v$ to~$\puncture''$.
Since the last step of~$\om$ is an odd-indexed step, after~$\om_{2k+2}$ there must be another step~$\om_{2k+3}$ which starts at~$\puncture''$.
Again by~(T3), $\om_{2k+3}$ starts in the interior of~$\blacktri{k+1}$.
But $\blacktri{k}$ and $\blacktri{k+1}$ are distinct triangles since $T^o$ has no self-folded triangle.
Hence $\puncture' \neq \puncture''$, and so $(\om_{2k+1},\om_{2k+2})$ is a quasi-backtrack which starts in the interior of~$\blacktri{k}$ and ends in the interior of~$\blacktri{k+1}$.

If $j=2k$ is even, then by a similar argument $(\om_{2k},\om_{2k+1})$ is a quasi-backtrack which starts in the interior of $\blacktri{k-1}$ and ends in the interior of $\blacktri{k}$.
See Fig.~\ref{fig:TikzQuasiBacktracks}.
\end{proof}

\subsection[A $(T^o,\ga)$-path $\om$ on $\Cn$ is uniquely determined by its sequence of labels]{A $\boldsymbol{(T^o,\ga)}$-path $\boldsymbol{\om}$ on $\boldsymbol{\Cn}$ is uniquely determined by its sequence of labels}
\label{subsec:uniquely_determined}
\begin{notation}
Let $r$ and $\ell$ be the radius and $\ell$-loop of a self-folded triangle, and let $v$ be the boundary vertex on the boundary that is adjacent to $\ell$.
Let $\ell_+$ (respectively, $\ell_-$) denote the clockwise (respectively, counterclockwise) orientation along $\ell$.
(We mark each $\ell$ curve in our illustrations with an arrow pointing clockwise to remind the reader that $\ell_+$ denotes the clockwise direction.)

Consider the cycle $(r,r)$
\begin{gather*}
\text{$v \leadsto \puncture' \leadsto v$ along $r$}
\end{gather*}
from Remark \ref{rem:pair}(\ref{rem:pair:v_p_v}).
Let $(\rrBak)$ denote the backtrack cycle and
let $(\rrNon)$ denote a (counterclockwise or clockwise) non-backtrack cycle.
\end{notation}

The following proposition is an analogue of a remark from \cite[Section 3.1]{Sch10}.
\begin{prop}
\label{prop:uniquely_determined}
A $(T^o,\ga)$-path $\om$ on $\Cn$ is uniquely determined by its sequence of labels $(\om_1, \ldots, \om_{2d+1})$, forgetting the orientations of the steps and whether a consecutive pair is a~non-backtrack or a backtrack.
\end{prop}

\begin{exam} \label{example17}
To illustrate Proposition \ref{prop:uniquely_determined},
consider the $(T^o,\ga)$-path of Fig.~\ref{fig:path_six}
\begin{gather*}
\om=(b_4, 1, b_2, 2, 2, \ell_+, \ell_-, \rrNon, \ell_+, b_3).
\end{gather*}
The pair $(\om_8,\om_9)=(\rrNon)$ is a \emph{counterclockwise} non-backtrack along the radius $r$, and $\om_{10}=\ell_+$ goes clockwise around $\ell$.
Consider a dif\/ferent path $\om'$ (see Fig.~\ref{fig:path_six_bad}) which goes along the same sequence of arcs
 such that $(\om'_8,\om'_9)=(\rrNon)$ is a \emph{clockwise} non-backtrack along the radius, and $\om'_{10}=\ell_-$ goes counterclockwise around $\ell$.
Even though $\om'$ is homotopic to $\ga$, the segment $\ga_4$ is not homotopic to $[p_4,p_{5}]_{\om}$, violating Def\/inition~\ref{def:Tpath}(T2).

The paths illustrated in Figs.~\ref{fig:TikzNooseTpathOneBad} and~\ref{fig:TikzNotNooseTpathThreeBad}
go along associated quasi-arcs of the same arcs/edges as the $(T^o,\ga)$-paths of Figs.~\ref{fig:TikzNooseTpathOne} and~\ref{fig:TikzNotNooseTpathThree}, but they fail axiom (T3), which requires each odd-indexed step $\om_{2k+1}$ to start and f\/inish in the interior of $\blacktri{k}$ or at a boundary marked point.
\end{exam}

\begin{proof}[Proof of Proposition \ref{prop:uniquely_determined}]
For short, let the arcs $\taui{1},\dots,\taui{d}$ be denoted by arcs $1,\dots,d$.

First, consider the subsequence $(\om_1,\om_2,\om_3)$.
There are four possibilities:
\begin{enumerate}[i)]\itemsep=0pt
\item $\blacktri{0}$ is an ordinary triangle
(Fig.~\ref{fig:TikzTriangleCrossesOneEdge}).

\item $\blacktri{0}$ is an ideal triangle with two vertices where arc $1$ is the loop
(Fig.~\ref{fig:TikzTwoVertexTriangleCrossesOneEdgeLoop}).

\item $\blacktri{0}$ is an ideal triangle with two vertices where arc $1$ is not the loop
(Fig.~\ref{fig:TikzTwoVertexTriangleCrossesOneEdgeNotLoop}).

\item $\blacktri{0}$ is a self-folded triangle where arc $1$ is the $\ell$-loop
(Fig.~\ref{fig:TikzSelfFoldedTriangleCrossesNoose}).
\end{enumerate}
For each of the f\/irst three cases,
$\blacktri{0}$ has three distinct edges $k$, $\tg{0}$, and $\tg{-1}$.
There are exactly two possible subpaths for $(\om_1,\om_2)$, which are represented by two distinct sequences~$(\tg{0},1)$ and~$(\tg{-1},1)$.

In the fourth case, $\blacktri{0}$ is a self-folded triangle with radius $r$ and $\ell$-loop $\ell$ (see Fig.~\ref{fig:TikzSelfFoldedTriangleCrossesNoose}, top), so that $\tg{0}=\tg{-1}=r$.  However, there are exactly two possible subpaths in this case as well.

\begin{rem}\label{remark:om1_om2}
Suppose $\ga$ is a radius starting at the puncture, as in Fig.~\ref{fig:TikzSelfFoldedTriangleCrossesNoose}.
Let $v$ denote the vertex that $\ell$ is based at.
Due to (T2), the f\/irst step $\om_1$ must go $\puncture' \leadsto v$ along $\tg{0}=\tg{-1}=r$, but there are two valid options for $\om_2$, to go counterclockwise or clockwise along $\ell$, so that there are exactly two possibilities for $(\om_1,\om_2)$:
\begin{enumerate}[{A}1)]\itemsep=0pt
\item \label{remark:om1_om2:rl_min} $(r,\ell_-)$,
\item \label{remark:om1_om2:rl_plus} $(r,\ell_+)$.
\end{enumerate}
\end{rem}

\begin{figure}[t!]
\centering
\subfigure[When $\ga$ starts at the puncture.]{\hspace{5mm}\TikzCrossingNooseOneAndA{0.7}\hspace{5mm}\label{fig:TikzCrossingNooseOneAndA}}
\quad
\subfigure[$(r$,$\ell_-$,$\a$,$\a)$.]{
\TikzRLMinusA{0.7}\label{fig:TikzRLMinusA}}
\\
\subfigure[$(r,\ell_+,\b,\a)$.]{
\TikzRLPlusB{0.7}\label{fig:TikzRLPlusB}}
\quad
\subfigure[$(r,\ell_+,\ell_-,\a)$.]{
\TikzRLPlusLMinusA{0.7}
\label{fig:TikzRLPlusLMinusA}
}

\caption{The subpaths $(\om_1,\dots,\om_4)$
of Remark \ref{remark:om1_om2}(A\ref{remark:om1_om2:rl_min}), (A\ref{remark:om1_om2:rl_plus})
when $\ga$ starts at~$\puncture$ and crosses~$\ell$ and~$\a$. See Fig.~\ref{fig:TikzRL_lift}.}
\label{fig:TikzRL_actual}
\label{fig:inducedmapRL_actual}
\end{figure}

We continue with the proof of Proposition \ref{prop:uniquely_determined} in this fourth case.  In particular, we show that~(A\ref{remark:om1_om2:rl_min}) and (A\ref{remark:om1_om2:rl_plus}) lead to distinct ways for the $(T^o,\ga)$-path to be f\/inished.
These two possible choices  for $(\om_1,\om_2)$ are represented by the same sequence $(\tg{0},1)=(\tg{-1},1)$, but
we claim that the orientation of $\om_2$ determines the possible choices for the next term~$\om_3$.
Suppose that~$\a$ (respectively, $\b$) is the side of $\blacktri{1}$ which lies clockwise (respectively, counterclockwise) of $\ell$, as illustrated in Fig.~\ref{fig:TikzCrossingNooseOneAndA}.

Suppose $\ga$ ends at the vertex $y$. If $\om_2=\ell_-$ is counterlockwise, then $\om_3=\a$
(Fig.~\ref{fig:TikzRLMinusA}), otherwise $\om_3=\b$ (Fig.~\ref{fig:TikzRLPlusB}).

If $\ga$ does not end at $y$, then the second arc that $\ga$ crosses is either $\a$ or $\b$ (say, $\a$).
If $\om_2=\ell_-$ is counterclockwise, then $(\om_3,\om_4)=(\a,\a)$ (Fig.~\ref{fig:TikzRLMinusA}).
If $\om_2=\ell_+$ is clockwise, then $(\om_3,\om_4)=(\b,\a)$ (Fig.~\ref{fig:TikzRLPlusB})
or $(\ell_-,\a)$ (Fig.~\ref{fig:TikzRLPlusLMinusA}).
Hence, the subpath $(\om_1,\om_2,\om_3)$ is either $(r,\ell_-,\a)$, $(r,\ell_+,b)$, or $(r,\ell_+,\ell_-)$.
As these are represented by three distinct subsequenes, each subsequence uniquely determines the f\/irst three steps of $\om$.

By the same logic, the subsequence $(\om_{2d-1},\om_{2d},\om_{d+1})$ uniquely determines the last three steps of $\om$.

Next, consider any triple  $(\om_{2k},\om_{2k+1},\om_{2k+2})$.
There are four possibilities:
\begin{enumerate}[i)]\itemsep=0pt
\item $\blacktri{k}$ is an ordinary triangle (Fig.~\ref{fig:TikzTriangleCrossesTwoEdges}).

\item $\blacktri{k}$ is an ideal triangle with two vertices where neither arc $k$ nor arc $k+1$ are loops (Fig.~\ref{fig:TikzTwoVertexTriangleCrossesTwoEdgesNoLoop}).

\item $\blacktri{k}$ is an ideal triangle with two vertices where one of arcs $k$ and $k+1$ (say, $k+1$) is a loop (Fig.~\ref{fig:TikzTwoVertexTriangleCrossesTwoEdgesKOneisLoop}).

\item $\blacktri{k}$ is a self-folded triangle with radius $r$ and $\ell$-loop $\ell$ where one of arcs $k$ and $k+1$ (say, $k+1$) is the radius $r$ (Fig.~\ref{fig:TikzSelfFoldedTriangleCrossesRadiusNoose}).
\end{enumerate}
For each of the f\/irst three cases, $\blacktri{k}$ has three distinct edges $k$, $k+1$, and $\tg{k}$. There are exactly three legal subpaths for $(\om_{2k},\om_{2k+1},\om_{2k+2})$, and they are represented by three distinct sequences $(k,k+1,k+1)$, $(k,k,k+1)$, and  $(k,\tg{k},k+1)$.  In the case that one of these steps is an $\ell$-loop, it is also clear that (T2) forces a specif\/ic orientation along $\ell$ as part of the $(T^o,\ga)$-path.

In the fourth case, $\ga$ crosses $\ell$, $r$, then $\ell$ in either the counterclockwise or clockwise (say, the former) direction.

\begin{lem}\label{lemma:four_paths}
Suppose $T^o$ contains a self-folded triangle with radius $r$ and $\ell$-loop $\ell$, and suppose that $\ga$ crosses $r$, say, in the counterclockwise direction. Let $v$ be the boundary vertex on the boundary that is adjacent to $\ell$.
Suppose $\ell$, $r$, $\ell$ are the $k$-th, $(k+1)$-th, and $(k+2)$-th arcs crossed by $\ga$ $($see Fig.~{\rm \ref{fig:TikzSelfFoldedTriangleCrossesRadiusNoose})}.
Then there are exactly three possible subpaths for $(\om_{2k},\om_{2k+1},\om_{2(k+1)})$:
\begin{enumerate}[${a}1)$]\itemsep=0pt
\item \label{itm:lrr_bak} $(\ell_-,\rrBak)$: follow $\ell_-$ then backtrack cycle $(\rrBak)$ $($Fig.~{\rm \ref{fig:TikzEllMinusRR})},
\item \label{itm:llr} $(\ell_+,\ell_-,r)$: follow $(\ell_+,\ell_-)$ then follow $r$ from $v$ to $\puncture'$ $($a point in the vicinity of the puncture$)$ $($Fig.~{\rm \ref{fig:TikzEllPlusEllMinusR})},
\item \label{itm:lrr_non} $(\ell_+,\rrNon)$: follow $\ell_+$ then counterclockwise non-backtrack cycle $(\rrNon)$
$($Fig.~{\rm \ref{fig:TikzEllPlusRRnonBT})}.
\end{enumerate}
There are also exactly three possible subpaths for $(\om_{2(k+1)},\om_{2k+3},\om_{2(k+2)})$:
\begin{enumerate}[${b}1)$]\itemsep=0pt
\item \label{itm:rll}  $(r,\ell_-,\ell_+)$: follow $r$ from $\puncture'$ (a point in the vicinity of the puncture) to $v$ then $(\ell_-,\ell_+)$ $($see Fig.~{\rm \ref{fig:TikzRLMinLPlus})},
\item \label{itm:rrl_bak} $(\rrBak,\ell_-)$: backtrack cycle $(\rrBak)$ then $\ell_-$ $($Fig.~{\rm \ref{fig:TikzRREllMinus})},
\item \label{itm:rrl_non} $(\rrNon,\ell_+)$: counterclockwise non-backtrack cycle $(\rrNon)$ then $\ell_+$
$($Fig.~{\rm \ref{fig:TikzRREllPlusnonBT})}.
\end{enumerate}

The $(\lrrll)$-subsequence corresponds to exactly $2$ valid subpaths for $(\om_{2k},\dots,\om_{2k+2})$,
one where  $(\rrBak)$  is a backtrack and the other where $(\rrNon)$ is a non-backtrack.
We get them by combining $(a\ref{itm:lrr_bak})$ with $(b\ref{itm:rll})$ and combining $(a\ref{itm:lrr_non})$ with $(b\ref{itm:rll})$:
\begin{enumerate}[$I)$]\itemsep=0pt
\item \label{item:four_paths_lrrll_backtrack}
$(\ell_-,\rrBak,\ell_-,\ell_+)$:
follow $\ell_-$ then backtrack cycle $(\rrBak)$ then $(\ell_-,\ell_+)$.

\item \label{item:four_paths_lrrll_nonbacktrack}
$(\ell_+,\rrNon,\ell_-,\ell_+)$:
follow $\ell_+$ then counterclockwise non-backtrack cycle $(\rrNon)$
then $(\ell_-,\ell_+)$.

\suspend{enumerate}
Similarly, the $(\llrrl)$-subsequence would correspond to exactly two valid subpaths,
one where $(\rrBak)$ is a backtrack and the other where $(\rrNon)$ is a non-backtrack.
We get them by combining $(a\ref{itm:llr})$ with $(b\ref{itm:rrl_bak})$,
and combining $(a\ref{itm:llr})$ with $(b\ref{itm:rrl_non})$.
\resume{enumerate}[{[$I)$]}]

\item \label{item:four_paths_llrrl_backtrack}
$(\ell_+,\ell_-,\rrBak,\ell_-)$:
follow $\ell_+$, then $\ell_-$ then
backtrack cycle $(\rrBak)$ then $\ell_-$.

\item \label{item:four_paths_llrrl_nonbacktrack}\sloppy
$(\ell_+,\ell_-,\rrNon,\ell_+)$:
follow $\ell_+$, then $\ell_-$ then
counterclockwise non-backtrack cycle~$(\rrNon)$ then~$\ell_+$.
\end{enumerate}
\end{lem}

\newcommand\ThreeSubpathsOnActualScale{0.7}

\begin{figure}[t!]
\centering
\renewcommand*\thesubfigure{(\arabic{subfigure})}
\subfigure[$(\ell_-,\rrBak)$.]{
\TikzEllMinusRR{\ThreeSubpathsOnActualScale}
\label{fig:TikzEllMinusRR}
}\quad
\subfigure[$(\ell_+,\ell_-,r)$.]{
\TikzEllPlusEllMinusR{\ThreeSubpathsOnActualScale}
\label{fig:TikzEllPlusEllMinusR}
}\\
\subfigure[$(\ell_+,\rrNon)$.]{
\TikzEllPlusRRnonBT{\ThreeSubpathsOnActualScale}
\label{fig:TikzEllPlusRRnonBT}
}\quad
\subfigure[$\ga$ crosses $r$ in the counterclockwise direction.]{
\TikzCrossingCounterclockwiseLRAB{\ThreeSubpathsOnActualScale}
\label{fig:cross_counterclock}
}

\caption{The subpaths $(\om_{2k},\om_{2k+1},\om_{2k+2})$
of Lemma \ref{lemma:four_paths}(a\ref{itm:lrr_bak}), (a\ref{itm:llr}), (a\ref{itm:lrr_non})
if $\ga$ crosses arcs $\taui{k}=\ell$, $\taui{k+1}=r$, and $\taui{k+2}=\ell$ in the counterclockwise order. See Fig.~\ref{fig:tilom_2k_2k_plus2}.}
\label{fig:om_2k_2k_plus2}
\end{figure}

\begin{figure}[t!]\centering\renewcommand*\thesubfigure{(\arabic{subfigure})}
\subfigure[$(r,\ell_-,\ell_+)$.]{
\TikzRLMinLPlus{\ThreeSubpathsOnActualScale}
\label{fig:TikzRLMinLPlus}
}\quad
\subfigure[$(\rrBak,\ell_-)$.]{
\TikzRREllMinus{\ThreeSubpathsOnActualScale}
\label{fig:TikzRREllMinus}
}\quad
\subfigure[$(\rrNon,\ell_+)$.]{
\TikzRREllPlusnonBT{\ThreeSubpathsOnActualScale}
\label{fig:TikzRREllPlusnonBT}
}

\caption{The subpaths $(\om_{2k+2},\om_{2k+3},\om_{2k+4})$
of Lemma \ref{lemma:four_paths}(b\ref{itm:rll}), (b\ref{itm:rrl_bak}), (b\ref{itm:rrl_non})
if $\ga$ crosses arcs $\taui{k}=\ell$, $\taui{k+1}=r$, and $\taui{k+2}=\ell$ in the counterclockwise order.
See Fig.~\ref{fig:tilom_2kplus2_2k_plus4}.}
\label{fig:om_2kplus2_2k_plus4}
\end{figure}

\begin{figure}[t!]\centering
\subfigure[Concatenation of $\pkpkone{\om}$ ($\ell_-$,$\rrBak$(backtrack)) and $\ga_k^-$ \emph{is} contr\-actible, so $\pkpkone{\om}$ is homotopic to $\ga_k$.]
{
\TikzEllMinusRRConcatenate{0.8}
\label{fig:TikzEllMinusRRnonBTConcatenate}}\quad
\subfigure[Similarly, $\pkpkone{\om}$ ($\ell_+$,$\rrNon$ (counterclockwise non-backtrack)) is homotopic to~$\ga_k$.]
{\hspace{2mm}
\TikzEllPlusRRnonBTConcatenate{0.8}\hspace{2mm}
\label{fig:TikzEllPlusRRnonBTConcatenate}}\quad
\subfigure[$\pkpkone{\om}$ $(\ell_+$,$\rrBak$(back\-track)) is \emph{not} homotopic to~$\ga_k$.]
{\hspace{2mm}
\TikzEllPlusRRConcatenate{0.8}\hspace{2mm}
\label{fig:TikzEllPlusRRConcatenate}}
\caption{Concatenation of $\pkpkone{\om}$ $(\ell,r,r)$ and the opposite orientation of $\ga_k^-$ of $\ga_k$.}
\end{figure}

\begin{rem}
In Lemma \ref{lemma:four_paths}, note that the subpaths (\ref{item:four_paths_lrrll_backtrack}) and (\ref{item:four_paths_llrrl_backtrack}) are homotopic to a counterclockwise orientation of $\ell$.
The subpaths (\ref{item:four_paths_lrrll_nonbacktrack}) and (\ref{item:four_paths_llrrl_nonbacktrack}) are contractible to $v$.
\end{rem}
\begin{proof}[Proof of Lemma \ref{lemma:four_paths}]
By (T1), since the $k$-th, $(k+1)$-th, $(k+2)$-th arcs crossed by~$\ga$ are~$\ell$,~$r$, and~$\ell$, the sequence
\begin{gather*}
(\om_{2k},\om_{2k+1},\om_{2(k+1)},\om_{2k+3},\om_{2(k+2)})
\end{gather*}
must be
\begin{gather*}(\ell, ~~, r, ~~, \ell).\end{gather*} For this to be a connected path, we must f\/ill in
the odd steps with $\ell$ then $r$ (or $r$ then $\ell$) by~(T2), so that
the 5-term subsequence is either~$(\lrrll)$ or~$(\llrrl)$.

The two subpaths ({a}\ref{itm:lrr_bak}) and ({a}\ref{itm:lrr_non}) (both represented by $(\ell,r,r)$) satisfy (T2):
\begin{enumerate}[{Case a}(\ref{itm:lrr_bak})]\itemsep=0pt
\item[Case ({a}\ref{itm:lrr_bak})]
If $(\om_{2k},\om_{2k+1},\om_{2(k+1)})=(\ell_-,\rrBak)$ (where $(\rrBak)$ is a backtrack cycle),
the concatenation of $\pkpkone{\om}$ and the opposite orientation ${\ga_k}^-$ of $\ga_k$ is contractible, so $\pkpkone{\om}$ and $\ga_k$ are homotopic.
See Fig.~\ref{fig:TikzEllMinusRRnonBTConcatenate}.

\item[Case ({a}\ref{itm:lrr_non})]
If $(\om_{2k},\om_{2k+1},\om_{2(k+1)})=(\ell_+,\rrNon)$ (where $(\rrNon)$ is a counterclockwise non-backtrack cycle),
the concatenation of $\pkpkone{\om}$ and the opposite orientation ${\ga_k}^-$ of $\ga_k$ is contractible, so $\pkpkone{\om}$ and $\ga_k$ are homotopic.
See Fig.~\ref{fig:TikzEllPlusRRnonBTConcatenate}.
\end{enumerate}
The sequence $(\ell,r,r)$ may also represent a subpath $(\ell_+,\rrBak)$ where $(\rrBak)$ is a backtrack or $(\ell_-,\rrNon)$ where $(\rrNon)$ is a non-backtrack cycle, but we claim that these are not valid $T^o$-subpath:
For contradiction, suppose $(\om_{2k}, \om_{2k+1}, \om_{2k+2})$ $=$ $(\ell_+,\rrBak)$ where $(\rrBak)$ is a backtrack.
By~(T2), the segment $\ga_k$ (Fig.~\ref{fig:TikzGammaK}) must be homotopic to~$\pkpkone{\om}$.
However, as illustrated in Fig.~\ref{fig:TikzEllPlusRRConcatenate}, the concatenation of~$\pkpkone{\om}$ and the opposite orientation ${\ga_k}^-$ of $\ga_k$ is not contractible.
Hence~$\ga_k$ is not homotopic to~$\pkpkone{\om}$, and so this subpath does not appear.
Similarly, no subpath $(\ell_-,\rrNon)$ with a non-backtrack cycle~$(\rrNon)$ can appear.

By similar logic,
\begin{enumerate}[--]\itemsep=0pt
\item ({a}\ref{itm:llr}) is the only valid $T^o$-subpath which can be represented by  $(\ell,\ell ,r)$,
\item ({b}\ref{itm:rll}) is the only valid $T^o$-subpath which can be represented by $(r,\ell,\ell)$, and
\item ({b}\ref{itm:rrl_bak}) and ({b}\ref{itm:rrl_non}) are the only two valid $T^o$-subpaths which can be represented by  $(r,r,\ell)$.
\end{enumerate}
This concludes our proof of Lemma \ref{lemma:four_paths}.
\end{proof}

We continue with the proof of Proposition \ref{prop:uniquely_determined}.
Note that the subpaths (\ref{item:four_paths_lrrll_backtrack}) and (\ref{item:four_paths_lrrll_nonbacktrack}) are represented by the same $5$-term subsequence $(\lrrll)$.
We claim that the steps $\om_{2k-1}$ which precede (\ref{item:four_paths_lrrll_backtrack}) and (\ref{item:four_paths_lrrll_nonbacktrack}) go along distinct arcs/edges.
Note that the ideal triangle $\blacktri{k-1}$ is a two-vertex, three-edge triangle (Fig.~\ref{fig:TikzTwoVertexTriangle}).
\begin{enumerate}[\text{Case }1:]\itemsep=0pt
\item
Suppose $k=1$, so that $\ga$ starts at the marked point $y$.
Recall that the side of $\blacktri{0}$ that lies clockwise of $\ell$ is $\tg{0}$, and the side of $\blacktri{0}$ that lies counterclockwise of $\ell$ is $\tg{-1}$ (see Fig.~\ref{fig:TikzTwoVertexTriangleCrossesOneEdgeLoop}). In Fig.~\ref{fig:cross_counterclock}, $\tg{0}=\a$ and $\tg{-1}=\b$.
	\begin{enumerate}[--]\itemsep=0pt
	\item Since (\ref{item:four_paths_lrrll_backtrack}) starts with a counterclockwise $\ell_-$ step,
the f\/irst step $\om_1$ must go along~$\tg{-1}$.
So the f\/irst six-term subsequence is $(\tg{-1},\lrrll)$.
	\item Since (\ref{item:four_paths_lrrll_nonbacktrack}) starts with a clockwise $\ell_+$ step,
the f\/irst step $\om_1$ must go along~$\tg{0}$.
So the f\/irst six-term subsequence is~$(\tg{0},\lrrll)$.
	\end{enumerate}

\item
Suppose $1<k$,
so that $\ga$ crosses arc $k-1$ (in Fig.~\ref{fig:cross_counterclock}, this arc is either~$\a$ or~$\b$.)
Suppose (without loss of generality) that arc $k-1$ lies clockwise of $\ell$ (i.e., arc $k-1$ is~$\a$ in Fig.~\ref{fig:cross_counterclock}).
Then the third arc $\tg{k-1}$ of $\blacktri{k-1}$ lies counterclockwise of $\ell$ (i.e., side~$\tg{k-1}$ is~$\b$ in Fig.~\ref{fig:cross_counterclock}).
Recall that $\om_{2k-2}$ must traverse arc $k-1$ by~(T1).
	\begin{enumerate}[--]\itemsep=0pt
	\item Since (\ref{item:four_paths_lrrll_backtrack}) starts with a counterclockwise $\ell_-$ step,
	the previous step $\om_{2k-1}$ must go along either $\ell$ or $\tg{k-1}$ (i.e., side $\b$ in Fig.~\ref{fig:cross_counterclock}), so that $(\om_{2k-1},\dots,\om_{2(k+2)})$ form a~six-term subsequence
	$(\ell, \lrrll)$ or $(\tg{k-1}, \lrrll)$.
	\item Since (\ref{item:four_paths_lrrll_nonbacktrack}) starts with a clockwise $\ell_+$ step,
	the previous step $\om_{2k-1}$ must go along arc $k-1$ (i.e., arc $\a$ in Fig.~\ref{fig:cross_counterclock}), so that $(\om_{2k-1},\dots,\om_{2(k+2)})$ $=$ $(k-1, \lrrll)$.
	\end{enumerate}
\end{enumerate}
Similarly, the subpaths (\ref{item:four_paths_llrrl_backtrack}) and (\ref{item:four_paths_llrrl_nonbacktrack}) are represented by the same $5$-term subsequence $(\ell, \ell, r$, $r, \ell)$.
We claim that the steps $\om_{2(k+2)+1}$ following (\ref{item:four_paths_llrrl_backtrack}) and (\ref{item:four_paths_llrrl_nonbacktrack}) go along distinct arc/edges.
Note that the ideal triangle $\blacktri{k+2}$ is a two-vertex, three-edge triangle (Fig.~\ref{fig:TikzTwoVertexTriangle}).
\begin{enumerate}[\text{Case }1:]\itemsep=0pt
\item Suppose $d=k+2$, so that $\ga$ ends at the marked point $y$.
Recall that the side of $\blacktri{d}$ that lies clockwise of $\ell$ is $\tg{d}$, and the side of $\blacktri{0}$ that lies counterclockwise of $\ell$ is $\tg{d+1}$ (see Fig.~\ref{fig:TikzTwoVertexTriangleCrossesOneEdgeLoop}).
In Fig.~\ref{fig:cross_counterclock}, $\tg{d}=\a$ and $\tg{d+1}=\b$.
\begin{enumerate}[--]\itemsep=0pt
\item Since (\ref{item:four_paths_llrrl_backtrack}) ends with a counterclockwise $\ell_-$, the last step must go along $\tg{d}$, so that the last six-term subsequence is $(\llrrl, \tg{d})$
\item Since (\ref{item:four_paths_llrrl_nonbacktrack}) ends with a clockwise $\ell_+$, the last step must go along $\tg{d+1}$, so that the last six-term subsequence is $(\llrrl, \tg{d+1})$.
\end{enumerate}

\item Suppose $d>k+2$, so that $\ga$ crosses arc $k+3$ (in Fig.~\ref{fig:cross_counterclock}, this arc is either $\a$ or $\b$.)
Suppose (without loss of generality) that arc $k+3$ lies clockwise of $\ell$ (i.e., arc $k+3$ is~$\a$ in Fig.~\ref{fig:cross_counterclock}).
Then the third arc $\tg{k+2}$ of $\blacktri{k+2}$ lies counterclockwise of $\ell$ (i.e., side~$\tg{k+2}$ is $\b$ in Fig.~\ref{fig:cross_counterclock}).
\begin{enumerate}[--]\itemsep=0pt
\item Since (\ref{item:four_paths_llrrl_backtrack}) ends with a counterclockwise $\ell_-$, the next step $\om_{2k+5}$ must go along arc $k+3$ (i.e., arc $\a$  in Fig.~\ref{fig:cross_counterclock}), so that $(\om_{2k},\dots,\om_{2k+5})=(\llrrl, k+3)$.
\item Since (\ref{item:four_paths_llrrl_nonbacktrack}) ends with a clockwise $\ell_+$, the next step $\om_{2k+5}$ must go along either counterclockwise $\ell_-$ or $\tg{k+2}$ (i.e., arc/edge labeled $\b$ in Fig.~\ref{fig:cross_counterclock}), so that $(\om_{2k},\dots,\om_{2k+5})$ form a six-term subsequence $(\llrrl, \ell)$ or $(\llrrl, \tg{k+2})$.\hfill \qed
\end{enumerate}
\end{enumerate}\renewcommand{\qed}{}
\end{proof}

\section[Proof of Theorem \ref{thm:tpath_expansion_formula}, the $(T^o,\ga)$-path expansion formula for a once-punctured disk]{Proof of Theorem \ref{thm:tpath_expansion_formula}, the $\boldsymbol{(T^o,\ga)}$-path expansion formula \\ for a once-punctured disk}
\label{sec:proof_of_tpath_expansion_formula}

We begin with an outline of our proof. Let $T^o$ be an ideal triangulation of a once-punctured $n$-gon $\Cn$ and let $\ga\notin T^o$ be an oriented ordinary arc (possibly an $\ell$-loop).
\begin{enumerate}[Step 1:]\itemsep=0pt
\item
Following \cite[Section 7]{MSW11}, we construct a triangulated $(d+3)$-gon $\tilS$
(modeled after $\blacktri{0},\dots,\blacktri{d}$). See Section \ref{subsection:triangulated_polygon_and_lifted_arc}.

\item In Section \ref{subsection:define_induced_covering_map}, we give a bijection
$
\overpi\colon  \{(\tilS,\tilga) \textnormal{-paths}  \}
\to \{(T^o,\ga) \textnormal{-paths} \}.
$

\item
Section \ref{subsection:formula_from_perfect_matchings} recaps \cite[Theorem 4.10]{MSW11},
a $T^o$-expansion of $x_\ga$ in terms of perfect matchings of a snake graph $\GT$.

\item
\cite[Lemma 4.5]{MS10}  gives a bijection $F\colon \{$perfect matchings of $\GT\} \to \{ (\tilS,\tilga)$-paths$\}$.
See Section \ref{subsection:proof_of_tpath_expansion_formula}.

\item The composition $\overpi \circ F$ gives a bijection between perfect matchings of $\GT$ and
 $(T^o,\ga)$-paths. Applying $\overpi \circ F$ to \cite[Theorem 4.10]{MSW11} yields Theorem \ref{thm:tpath_expansion_formula}.
\end{enumerate}

\subsection[A triangulated polygon $\tilS$ and a lifted arc $\tilga$]{A triangulated polygon $\boldsymbol{\tilS}$ and a lifted arc $\boldsymbol{\tilga}$}
\label{subsection:triangulated_polygon_and_lifted_arc}
Let $T^o$ be an ideal triangulation of a once-punctured $n$-gon $\Cn$ and let $\ga\notin T^o$ be an oriented ordinary arc (possibly an $\ell$-loop) from the point $s$ to the point $t$. Let $\Sga$ be the union of the ideal triangles $\blacktri{k}$ ($k=0,\dots,d$) crossed by $\ga$. First, for each $k=0,\dots,d$, we build a triangle~$\tilblacktri{k}$ with three distinct labels.

\begin{defin}[triangles $\tilblacktri{k}$]
\label{def:tilde_triangles}
If $\blacktri{k}$ has three distinct sides, its lift $\tilblacktri{k}$ is an ordinary triangle with the edge labels $\tiltaui{k}=\taui{k}$, $\tiltaui{k+1}=\taui{k+1}$, $\tiltg{k}=\tg{k}$, and the same orientation as $\blacktri{k}$.

When $\blacktri{k}$ is a self-folded triangle with radius $r$ and $\ell$-loop $\ell$, its lift $\tilblacktri{k}$ is formed by two lifts of $r$ and a lift of $\ell$, as follows:
\begin{enumerate}[{Case} $k=1$ or $d$:]\itemsep=0pt
\item[Case $k=1$ or $d$:] If $k=0$, then $\tilblacktri{0}$ has edge labels $\tiltaui{1}=\ell$, $\tiltg{0}=r$, $\tiltg{-1}=\ccr$ (arranged in clockwise order). Similarly, if $k=d$, then $\tilblacktri{d}$ has edge labels $\tiltaui{d}=\ell$, $\tiltg{d}=r$, $\tiltg{d+1}=\ccr$ (arranged in clockwise order). See Fig.~\ref{fig:tiltriangle_0}.

\item[Case $0<k<d$:]
Suppose $k=\ell$, $k+1=r$, and $k+2=\ell$. Then $\tilblacktri{k}$ has three distinct edge labels $\tiltaui{k}$, $\tiltaui{k+1}$, and $\tiltg{k}$: If $\ga$ crosses $r$ in counterclockwise (respectively, clockwise) direction,
let $\tiltaui{k+1}=r$ be the label of the edge of $\tilblacktri{k}$ which lies counterclockwise (respectively, clockwise) of $\tiltaui{k}=\ell$. Let $\tiltg{k}=\dot{r}$ be the label of the third edge of $\tilblacktri{k}$.

Similarly, $\tilblacktri{k+1}$ has three distinct edge labels $\tiltaui{k+1}$, $\tiltaui{k+2}$, and $\tiltg{k+1}$: If $\ga$ crosses $r$ in counterclockwise (respectively, clockwise) direction, let $\tiltg{k+1}=\ddot{r}$ be the label of the edge of $\tilblacktri{k+1}$ which lies counterclockwise (respectively, clockwise) of $\tiltaui{k+2}=\ell$. Let $r$ be the label of the third edge of $\tilblacktri{k+1}$. See Figs.~\ref{fig:TikzDoubleTriangleCounterclockwise} and~\ref{fig:TikzDoubleTriangleClockwise}.
\end{enumerate}
\end{defin}

\begin{figure}[t!]
\centering
\subfigure[$\ga$ is a radius starting or ending at the puncture.]
{\TikzCrossingNooseOne{0.7}
\TikzCrossingNooseD{0.7}
\label{fig:TikzCrossingNooseOne_TikzCrossingNooseD}
}\\
\newcommand\localscalesize{0.84}
\subfigure[The triangle $\tilblacktri{0}$ or $\tilblacktri{d}$ when $\ga$ is a radius starting or ending at the puncture.]{
\TikzLiftedTriangleOneTwo{\localscalesize}
\TikzLiftedTriangleD{\localscalesize}
\label{fig:tiltriangle_0}
}\\
\subfigure[$\ga$ crosses $r$ in the counterclockwise direction.]{
\TikzCrossingCounterclockwiseLRAB{0.7}\label{fig:cross_counterclock_later}
}\quad
\subfigure[$\ga$ crosses $r$ in the clockwise direction.]{
\TikzCrossingClockwiseLRAB{0.7}\label{fig:cross_clockwise}
}
\\
\subfigure[The triangles $\tilblacktri{k}$, $\tilblacktri{k+1}$ when~$\ga$ crosses a self-folded triangle's radius~$r$ in the counterclockwise direction.]{\hspace{3.5mm}\TikzDoubleTriangleCounterclockwise{\localscalesize}\hspace{3.5mm}
\label{fig:TikzDoubleTriangleCounterclockwise}
}
\quad
\subfigure[The triangles $\tilblacktri{k}$, $\tilblacktri{k+1}$ when~$\ga$ crosses a self-folded triangle's radius~$r$ in the clockwise direction.]{\hspace{3.5mm}\TikzDoubleTriangleClockwise{\localscalesize}\hspace{3.5mm}
\label{fig:TikzDoubleTriangleClockwise}
}
\caption{The triangles $\tilblacktri{k}$ when $\ga$ crosses an $\ell$-loop.}
\end{figure}

\begin{defin}[triangulated $(d+3)$-gon]
Glue $\tilblacktri{1}$ to $\tilblacktri{0}$ along $\taui{1}$, glue $\tilblacktri{2}$ to $\tilblacktri{1}$ along $\taui{2}$, and so forth to form a triangulated $(d+3)$-gon $\tilS$ with internal edges
\begin{gather*}
\tiltaui{1},\dots,\tiltaui{d}
\end{gather*}
and boundary edges \begin{gather*}
\tiltg{-1},\tiltg{0},\dots,\tiltg{d+1}.
\end{gather*}

Let the lift  $\widetilde{s}=\tilp_0$ of $s=p_0$ be the vertex of $\tilblacktri{0}$ that is adjacent to the lifted edges~$\tiltg{0}$ and~$\tiltg{-1}$, and the lift $\widetilde{t}=\tilp_{d+1}$ of $t=p_{d+1}$ be the vertex of $\tilblacktri{d}$ that is adjacent to the lifted edges~$\tiltg{d}$ and~$\tiltg{d+1}$. If $k=1,\dots,d-1$, let $\tilp_{k}$, $\tilp_{k+1}$ denote the lifts of $p_k$, $p_{k+1}$ on~$\tilblacktri{k}$
which lies on the interior of $\tiltaui{k}$, $\tiltaui{k+1}$ (respectively).

Let $\tilga$ be the arc in $\tilS$ from $\tils$ to $\tilt$.
We call $\tilga$ the \emph{lift} of $\ga$.
Let $\pi\colon \tilS \to \Sga$ denote the covering map of $\tilS \to \Sga$.
See Fig.~\ref{fig:T2_triangulated_polygon_tilS} for the lifts $\tilga$ and $\tilS$ of $\ga$ and $\Sga$ of Fig.~\ref{fig:T2_triangulated_polygon_tilS_original_T}.
\end{defin}

\begin{rem}\label{rem:same_labels}
By construction, $\tiltaui{k}=\taui{k}$ for each $k$.
Furthermore, $\tiltg{k}=\tgk$ unless $\tiltg{k}$ has label $\dot{r}$, $\ddot{r}$, or $\ccr$.
In addition, $\tiltg{0},\tiltg{d} \notin \{ \dot{r}, \ddot{r}, \ccr \}$, so $\tiltg{0}=\tg{0}$, $\tiltg{d}=\tg{d}$.
\end{rem}

\begin{figure}[t!]
\centering
\subfigure[Arc $\ga$ of $\Cnn{4}$ where $d=5$.]
{\TikzTtwoIdealTriangulationWithNumberLabels{0.8}
\label{fig:T2_triangulated_polygon_tilS_original_T}
}
\subfigure[Lifted arc $\tilga$ on the constructed triangulated \mbox{$(5+3)$-gon} $\tilS$.]
{\TtwoTriangulatedPolygon{0.45}\label{fig:T2_triangulated_polygon_tilS}}
\caption{The construction of the triangulated polygon $\tilS$ and lifted arc $\tilga$ for $T^o$ and $\ga$ from Fig.~\ref{fig:T2_triangulated_polygon_tilS_original_T}.}
\end{figure}

\begin{figure}[t!]
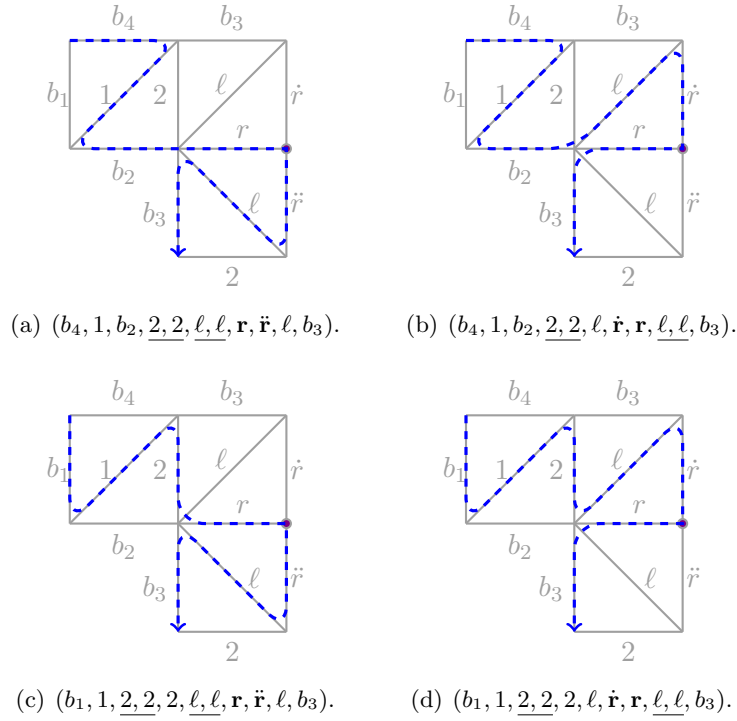

\centering
\newcommand\TikzPolygonTpathSize{0.36}
\subfigure[$(b_4,1,b_2,\underline{2,2},\underline{\ell,\ell},\RddotR,\ell,b_3)$.]
{\hspace{6mm}\TikzPolygonTpathSix{\TikzPolygonTpathSize}\hspace{6mm}}
\subfigure[$(b_4,1,b_2,\underline{2,2},\ell,\dotRR,\underline{\ell,\ell},b_3)$.]
{\hspace{6mm}\TikzPolygonTpathSeven{\TikzPolygonTpathSize}\hspace{6mm}}\\
\subfigure[$(b_1,1,\underline{2,2},2,\underline{\ell,\ell},\RddotR,\ell,b_3)$.]
{\hspace{6mm}\TikzPolygonTpathEight{\TikzPolygonTpathSize}\hspace{6mm}\label{fig:TikzPolygonTpathEight}}
\subfigure[$(b_1,1,\underline{2,2},2,\ell,\dotRR,\underline{\ell,\ell},b_3)$.]
{\hspace{6mm}\TikzPolygonTpathNine{\TikzPolygonTpathSize}\hspace{6mm}}
\caption{The $(\tilTo,\tilga)$-paths corresponding to the four $(T^o,\ga)$-paths from Fig.~\ref{fig:TtwoTpaths}. All backtracks have been omitted.}
\label{fig:TtwoLiftedTpaths}
\end{figure}

\subsection[A bijection $\overpi$ between $(\tilS,\tilga)$-paths and $(T^o,\ga)$-paths]{A bijection $\boldsymbol{\overpi}$ between $\boldsymbol{(\tilS,\tilga)}$-paths and $\boldsymbol{(T^o,\ga)}$-paths}
\label{subsection:define_induced_covering_map}

Keep the same setup as in the previous section, where $T^o$ is an ideal triangulation of $\Cn$ and $\ga\notin T^o$ is an ordinary arc of $\Cn$. As $\tilS$ is a triangulation of a polygon, we can consider a~$(\tilS,\tilga)$-path $\tilom=(\tilom_1,\dots,\tilom_{2d+1})$ as def\/ined in Def\/inition~\ref{def:Tpath}. For example, the $(\tilS,\tilga)$-paths of Fig.~\ref{fig:TtwoLiftedTpaths} correspond to the four $(T^o,\ga)$-paths
of Fig.~\ref{fig:TtwoTpaths}.

\begin{lem}
\label{lemma:induced_covering_map_bijection}
The covering map $\pi\colon \tilS \to \Sga$ which gives
\begin{gather*}
\pi (\tilom_j) =
\begin{cases}
r &\textnormal{if $\tilom_j$ has label $\dot{r}$, $\ddot{r}$, or $\ccr$,}\\
\textnormal{the same label as }\tilom_j & \textnormal{otherwise.}
\end{cases}
\end{gather*}
induces a bijection
\begin{gather*}
\overpi\colon \  \{(\tilS,\tilga) \text{-paths in }  \tilS \}
 \to  \{(T^o,\ga) \text{-paths in } \Sga \},
\\
\hphantom{\overpi\colon}{} \ \overpi (\tilom_1, \tilom_2, \dots, \tilom_{2d+1})  =
(\pi(\tilom_1), \pi( \tilom_2), \dots, \pi( \tilom_{2d+1})).
\end{gather*}
\end{lem}

\begin{figure}[t!]
\centering
\subfigure[When $\ga$ starts at the puncture.]{\hspace{7.5mm}\TikzLiftedTriangleOneTwoCrossesA{1}\hspace{7.5mm}\label{fig:TikzLiftedTriangleOneTwoCrossesA_next_to_lifted_paths}}
\quad
\subfigure[$(r,\ell,\a,\a)$.]
{\TikzNooseOneRLMinusA{1}}
\quad
\subfigure[$(\ccr,\ell,\b,\a)$.]{
\TikzNooseOneRLPlusB{1}}
\quad
\subfigure[$(r,\ell,\ell,\a)$.]{
\TikzNooseOneRLPlusLMinus{1}
}
\caption{The lifts $(\tilom_1,\dots,\tilom_4)$
of Remark \ref{remark:om1_om2}(A\ref{remark:om1_om2:rl_min}), (A\ref{remark:om1_om2:rl_plus})
when $\ga$ starts at the puncture, see Fig.~\ref{fig:TikzRL_actual}.}
\label{fig:inducedmapRL_tilS}
\label{fig:TikzRL_lift}
\end{figure}

\newcommand\ThreeSubpathsOnPolygonScale{0.9}
\begin{figure}[t!]
\renewcommand*\thesubfigure{(\arabic{subfigure})}
\centering
\mbox{
\subfigure[$(\ell,\rrBak)$.]{
\TikzDoubleTriangleKPlusPlus{\ThreeSubpathsOnPolygonScale}}\quad
\subfigure[$(\ell,\ell,r)$.]{
\TikzDoubleTriangleKKPlus{\ThreeSubpathsOnPolygonScale}}\quad
\subfigure[$(\ell,\dotRR)$.]{
\TikzDoubleTriangleKThirdEdgePlus{\ThreeSubpathsOnPolygonScale}}
}

\caption{The lifts $(\tilom_{2k},\tilom_{2k+1},\tilom_{2k+2})$ of the subpaths of Lemma \ref{lemma:four_paths}(a\ref{itm:lrr_bak}), (a\ref{itm:llr}), (a\ref{itm:lrr_non}), see Fig.~\ref{fig:om_2k_2k_plus2}.}
\label{fig:tilom_2k_2k_plus2}

\subfigure[$(r,\ell,\ell)$.]{
\TikzDoubleTriangleRLL{\ThreeSubpathsOnPolygonScale}}\quad
\subfigure[$(\rrBak,\ell)$.]{
\TikzDoubleTriangleRRL{\ThreeSubpathsOnPolygonScale}}\quad
\subfigure[$(r,\ddot{r},\ell)$.]{
\TikzDoubleTriangleRddotRL{\ThreeSubpathsOnPolygonScale}}

\caption{The lifts $(\tilom_{2k+2},\tilom_{2k+3},\tilom_{2k+4})$ of the subpaths of Lemma \ref{lemma:four_paths}(b\ref{itm:rll}), (b\ref{itm:rrl_bak}), (b\ref{itm:rrl_non}), see Fig.~\ref{fig:om_2kplus2_2k_plus4}.}
\label{fig:tilom_2kplus2_2k_plus4}
\end{figure}

\begin{proof}
Let $\tilom$ be a $(\tilS,\tilga)$-path.
Note that, since $\tilS$ is a polygon, there is only one way to concatenate a pair $(\tilom_{j},\tilom_{j+1})$ of steps.
By (T1), every $\tilom_{2k}$ has label $\tiltaui{k}$, and every $\om_{2k}$ has label $\taui{k}$. Per Remark \ref{rem:same_labels}, $\tiltaui{k}$ has the same label as $\taui{k}$, and $\tiltg{j}$ has the same label as $\tg{j}$ if $j\in\{0,d\}$.

By construction, $\tilblacktri{0}$ has three distinct labels, $\taui{1}$, $\tg{0}$, and $\tiltg{-1}$, and $\tilblacktri{d}$ has three distinct labels, $\taui{d}$, $\tg{d}$, and $\tiltg{d+1}$.
Hence, as discussed in Section \ref{subsec:uniquely_determined},
the subpath $(\tilom_1,\tilom_2)$ is either $(\tg{0},\taui{1})$ or $(\tiltg{-1},\taui{1})$, and
the subpath $(\tilom_{2d},\tilom_{2d+1})$ is either $(\taui{d},\tg{d})$ or $(\taui{d},\tiltaui{d+1})$.
\begin{itemize}\itemsep=0pt
\item
If $\blacktri{0}$ has three distinct sides, $\taui{1}$, $\tg{0}$, and $\tg{-1}$, then $\tiltg{-1}=\tg{-1}$, and $\overpi$ maps $(\tilom_1,\tilom_2)$ to the $(T^o,\ga)$-subpath $(\om_1,\om_2)$ with the same labels.

 If $\blacktri{0}$ is a self-folded triangle with radius $r$ and $\ell$-loop $\ell$, then recall that $\tilblacktri{0}$ is an ordinary triangle with edge labels $\ell,\ccr,r$ (in counterclockwise order) as in Fig.~\ref{fig:TikzLiftedTriangleOneTwoCrossesA_next_to_lifted_paths}.
By~(T1) and~(T2), $\tilom_1$ goes along $\tiltg{0}$ or $\tiltg{-1}$ and $\tilom_2$ goes along $\ell$.
Hence, either $(\tilom_1, \tilom_2)=(\ccr, \ell)$ or $(\tilom_1, \tilom_2)=(r, \ell)$ (see Fig.~\ref{fig:inducedmapRL_tilS}).
By Remark \ref{remark:om1_om2}, there are two possible $(T^o,\ga)$-subpaths $(\om_1,\om_2)$, either $(r,\ell_-)$ or $(r,\ell_+)$.
We see that $\overpi$ maps
\begin{enumerate}[]\itemsep=0pt
\item $(r, \ell) \bijection$ Remark \ref{remark:om1_om2}(A\ref{remark:om1_om2:rl_min}): $(r, \ell_-)$,
\item $(\ccr,  \ell) \bijection$ Remark \ref{remark:om1_om2}(A\ref{remark:om1_om2:rl_plus}): $(r, \ell_+)$.
\end{enumerate}

\item Similarly, if $\blacktri{d}$ has three distinct sides, then $(\tilom_{2d},\tilom_{2d+1})$ is mapped to the subpath with the same labels. Otherwise, if $\blacktri{d}$ is a self-folded triangle with radius $r$ and $\ell$-loop $\ell$, then recall that $\tilblacktri{d}$ is an ordinary triangle with edges~$\ell$,~$\ccr$,~$r$ (in counterclockwise order), see Fig.~\ref{fig:tiltriangle_0}.
Then either $(\tilom_{2d}, \tilom_{2d+1})=(\ell,\ccr)$ or $(\tilom_{2d}, \tilom_{2d+1})=(\ell,r)$, and $\overpi$ maps
\begin{enumerate}[]\itemsep=0pt
\item $(\ell,\ccr) \bijection (\ell_-,r)$,
\item $(\ell,r) \bijection (\ell_+,r)$.
\end{enumerate}
\end{itemize}

Similarly, each $\tilblacktri{k}$ has three distinct labels $\taui{k}$, $\taui{k+1}$, and $\tiltg{k}$. As discussed in Section \ref{subsec:uniquely_determined},
the subpath $(\tilom_{2k},\tilom_{2k+1},\tilom_{2k+2})$ is one of $(\taui{k},\taui{k+1},\taui{k+1})$, $(\taui{k},\taui{k},\taui{k+1})$, or $(\taui{k},\tiltg{k},\taui{k+1})$.

\begin{itemize}\itemsep=0pt
\item If $\blacktri{k}$ ($k=1,\dots,d$) has three distinct sides, $\taui{k}$, $\taui{k+1}$, and $\tg{k}$,
then $\tiltg{k}=\tg{k}$. Hence $(\tilom_{2k},\tilom_{2k+1},\tilom_{2k+2})$ is mapped to the $(T^o,\ga)$-subpath $(\om_{2k},\om_{2k+1},\om_{2k+2})$ with the same labels.

Otherwise, suppose $\blacktri{k}$ is a self-folded triangle with radius $r$ and $\ell$-loop $\ell$. Assume that~$\ell$,~$r$,~$\ell$ are the $k$-th, $(k+1)$-th, and $(k+2)$-th arcs crossed by $\ga$ and that $\ga$ crosses them in the counterclockwise direction (see Fig.~\ref{fig:cross_counterclock_later}).
By construction (see Def\/inition~\ref{def:tilde_triangles} and Fig.~\ref{fig:TikzDoubleTriangleCounterclockwise}), $\tilblacktri{k}$ has sides labeled $\tiltaui{k}=\ell$, $\tiltaui{k+1}=r$, and $\tiltg{k}=\dot{r}$ (in counterclockwise order).
The subpath $(\tilom_{2k},\tilom_{2k+1},\tilom_{2k+2})$ is one of $(\tiltaui{k},\tiltaui{k+1},\tiltaui{k+1})=(\ell,r,r)$, $(\tiltaui{k},\tiltaui{k},\tiltaui{k+1})=(\ell,\ell,r)$, or $(\tiltaui{k},\tg{k},\tiltaui{k+1})=(\ell,\dot{r},r)$. See Fig.~\ref{fig:tilom_2k_2k_plus2}.
Per Lemma \ref{lemma:four_paths},
there are three possible $(T^o,\ga)$-subpaths for $(\om_{2k},\om_{2k+1},\om_{2k+2})$.
We see that $\overpi$ maps
\begin{enumerate}[]\itemsep=0pt
\item $(\ell,r,r) \bijection$ Lemma \ref{lemma:four_paths}(a\ref{itm:lrr_bak}): $(\ell_-,\rrBak)$,
\item $(\ell,\ell,r) \bijection$ Lemma \ref{lemma:four_paths}(a\ref{itm:llr}): $(\ell_+,\ell_-,r)$,
\item $(\ell,\dot{r},r) \bijection$ Lemma \ref{lemma:four_paths}(a\ref{itm:lrr_non}): $(\ell_+,\rrNon$ (counterclockwise non-backtrack)$)$.
\end{enumerate}

By construction (see Def\/inition \ref{def:tilde_triangles} and Fig.~\ref{fig:TikzDoubleTriangleCounterclockwise}), $\tilblacktri{k+1}$ has sides labeled $\tiltaui{k+1}=r$,  $\tiltaui{k+2}=\ell$, and $\tiltg{k+1}=\ddot{r}$ (in counterclockwise order).
The subpath $(\tilom_{2k+2},\tilom_{2k+3},\tilom_{2k+4})$ is one of $(\tiltaui{k+1},\tiltaui{k+2},\tiltaui{k+2})=(r,\ell,\ell)$, $(\tiltaui{k+1},\tiltaui{k+1},\tiltaui{k+2})=(r,r,\ell)$, or $(\tiltaui{k+1},\tiltg{k+1},\tiltaui{k+2})=(r,\ddot{r},\ell)$. See Fig.~\ref{fig:tilom_2kplus2_2k_plus4}.
Per Lemma \ref{lemma:four_paths}, there are three possible $(T^o,\ga)$-subpaths for $(\om_{2k+2},\om_{2k+3},\om_{2k+4})$.
We see that $\overpi$ maps
\begin{enumerate}[]\itemsep=0pt
\item $(r,\ell,\ell) \bijection $ Lemma \ref{lemma:four_paths}(b\ref{itm:rll}):  $(r,\ell_-,\ell_+)$,
\item $(r,r,\ell) \bijection $ Lemma \ref{lemma:four_paths}(b\ref{itm:rrl_bak}): $(\rrBak,\ell_-)$ then $\ell_-$,
\item $(r,\ddot{r},\ell) \bijection $ Lemma \ref{lemma:four_paths}(b\ref{itm:rrl_non}): $(\rrNon$ (counterclockwise non-backtrack) $,\ell_+)$.
\hfill {\qed}
\end{enumerate}
\end{itemize}\renewcommand{\qed}{}
\end{proof}

We expect that Lemma \ref{lemma:induced_covering_map_bijection} can be generalized to other punctured surfaces. However, other punctured surfaces can have more complication ideal triangulations (like those containing an ideal triangle of Fig.~\ref{fig:TikzOneVertexTriangle}) that have not been considered in our argument.

\subsection{Perfect matching expansion formula}
\label{subsection:formula_from_perfect_matchings}
We recap the snake graph expansion formula of \cite[Theorem 4.10]{MSW11} for ordinary arcs
of any bordered surface (including once-punctured polygons).
We continue to restrict our attention to the case of the once-punctured $n$-gon $\Cn$.

\begin{figure}[t!]
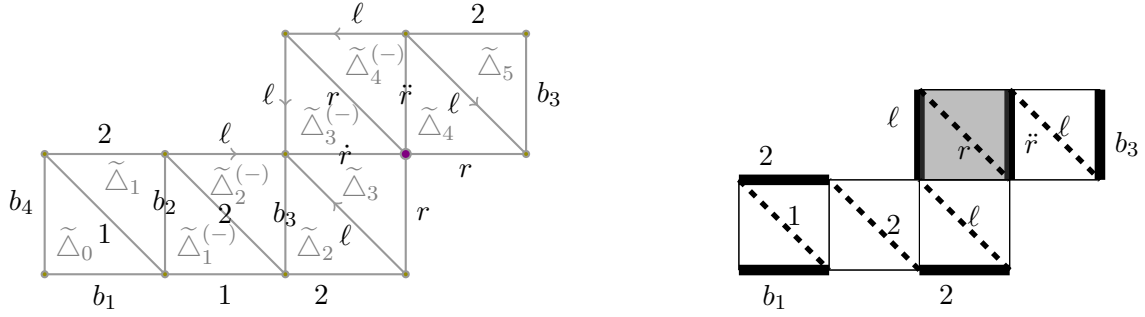
\centering\vspace{-2mm}
\subfigure[The snake graph $\overGT$ constructed from $\tilS$ in Fig.~\ref{fig:T2_triangulated_polygon_tilS}.]{
\TtwoSnakeGraphPlainer{0.4}
\label{fig:T2 MSW snake graph}
}
\quad
\subfigure[A path corresponding to the $(\tilS,\tilga)$-path of Fig.~\ref{fig:TikzPolygonTpathEight}, with odd steps $\PMEight$, see Fig.~\ref{fig:TtwoPerfectMatchingsEight}, and even steps the diagonals $1$, $2$, $\ell$, $r$, and~$\ell$.]
{\TtwoPerfectMatchingsWithDiagonalsEight{0.3}
\label{fig:ToIllustrateBijectionsTtwoPerfectMatchingsWithDiagonalsEight}
}

\caption{The snake graph $\overGT$ of Def\/inition \ref{def:snake_graph} and an example of a path on $\overGT$, see Remark~\ref{remark:bijection_F}.}
\end{figure}

\begin{defin}[snake graph $\GT$]\label{def:snake_graph}
We \emph{unfold} $\tilS$ into a graph $\overGT$, called a \emph{snake graph},
by inserting negative-oriented copies $\tilblacktri{k}^{(-)}$ of $\tilblacktri{k}$ (for $k=1,\dots,d-1$) into $\tilS$.
Fig.~\ref{fig:T2 MSW snake graph} is the graph $\overGT$ that is
built from $\tilS$ of Fig.~\ref{fig:T2_triangulated_polygon_tilS}.
See \cite[Section 4]{MSW11} for details.
A \emph{tile} consists of $\tilblacktri{k}$ and $\tilblacktri{k+1}$ (or $\tilblacktri{k}^{(-)}$ and $\tilblacktri{k+1}^{(-)}$) glued together along a \emph{diagonal} labeled $\tiltaui{k}$.
Let $\GT$ denote the graph obtained from $\overGT$ by removing the diagonal from each tile.
\end{defin}

\begin{defin}[crossing monomial]
Def\/ine the \emph{crossing monomial} of $\ga$ with respect to $T^o$ to~be
\begin{gather*}
\text{cross}(T^o, \ga) = \prod_{k=1}^d x(\tau_{i_k}).
\end{gather*}
\end{defin}

\begin{figure}[t!]
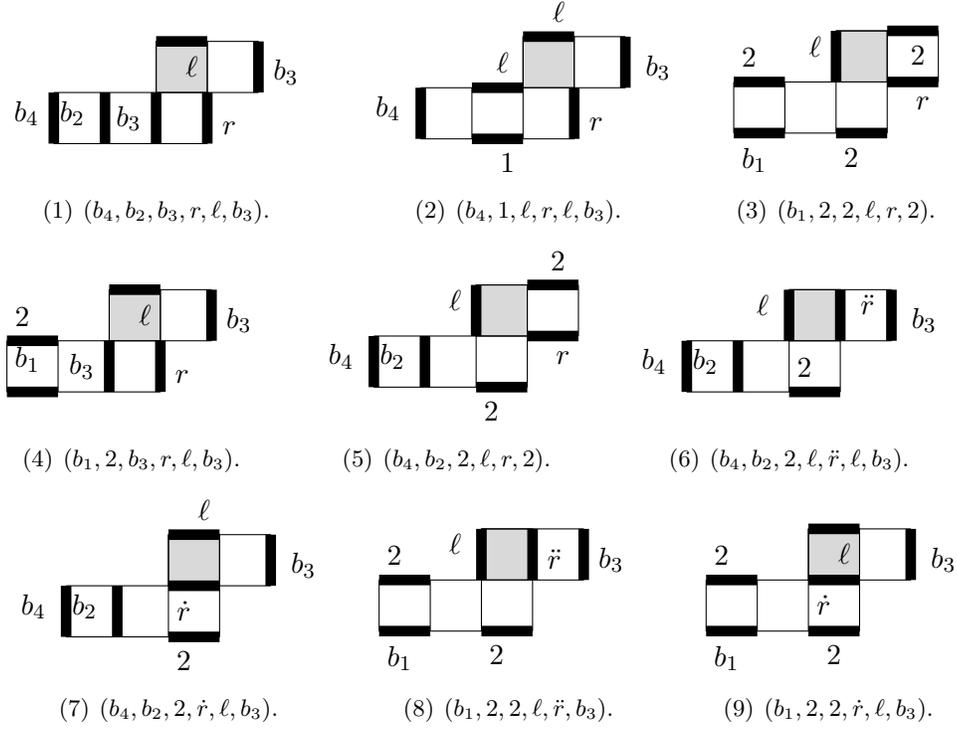
\centering\vspace*{-2mm}
\renewcommand*\thesubfigure{(\arabic{subfigure})}
\newcommand\PMOneThroughNineSize{0.17}

\subfigure[$\PMOne$.]{\TtwoPerfectMatchingsOne{\PMOneThroughNineSize}}\quad
\subfigure[$\PMTwo$.]{ \TtwoPerfectMatchingsTwo{\PMOneThroughNineSize}}\quad
\subfigure[$\PMThree$.]{\TtwoPerfectMatchingsThree{\PMOneThroughNineSize}}\quad
\vspace{-1.5mm}\\
\subfigure[$\PMFour$.]{\TtwoPerfectMatchingsFour{\PMOneThroughNineSize}}\quad
\subfigure[$\PMFive$.]{ \TtwoPerfectMatchingsFive{\PMOneThroughNineSize}}\quad
\subfigure[$\PMSix$.]{\TtwoPerfectMatchingsSix{\PMOneThroughNineSize}}\quad
\vspace{-1.5mm}\\
\subfigure[$\PMSeven$.]{\TtwoPerfectMatchingsSeven{\PMOneThroughNineSize}}\quad
\subfigure[$\PMEight$.]{\TtwoPerfectMatchingsEight{\PMOneThroughNineSize} \label{fig:TtwoPerfectMatchingsEight}}\quad
\subfigure[$\PMNine$.]{ \TtwoPerfectMatchingsNine{\PMOneThroughNineSize}}

\caption{The nine perfect matchings of the snake graph $\GT$ of Fig.~\ref{fig:T2 MSW snake graph}.
Matchings (6)--(9) correspond to the $(T^o,\ga)$-paths of Fig.~\ref{fig:TtwoTpaths}.
Note that these four are the only matchings where their restrictions to the gray-shaded tile are perfect matchings.}
\label{fig:TtwoPerfectMatchings}
\end{figure}

\begin{defin}[perfect matchings and weights]
\label{def:perfect_matchings_and_weights}
A \emph{perfect matching} of a graph $G$ is a subset~$E$ of the edges of~$G$ such that each vertex of~$G$ is incident to exactly one edge of $E$.
See Fig.~\ref{fig:TtwoPerfectMatchings}.
If the edges of a perfect matching $E$ are labeled $\be_{j_1}, \ldots, \be_{j_r}$,
then the \emph{weight}~$x(E)$ of~$E$ is the product $x_{\be_{j_1}} \cdots x_{\be_{j_r}}$.
\end{defin}

\begin{theorem} [{\cite[Theorem 4.10]{MSW11}}, perfect matching expansion formula]
\label{thm:410}\label{thm:MSW_perfect_matching_formula}
Let $T^o$ be an ideal triangulation of any surface~$S$
and let $\ga\notin T^o$ be an ordinary arc $($note: $\ga$ may be an $\ell$-loop$)$.
Let $x_\ga$ denote the element corresponding to $\ga$ in the cluster algebra which arises from~$S$ $($see Theorem~{\rm \ref{thm:fst_thm})}.
Then{\samepage
\begin{gather*}
x_\ga = \sum_E \frac{x(E)}
{\text{\rm cross}(T^o,\ga)},
\end{gather*}
where the sum is over all perfect matchings $E$ of $\GT$.}
\end{theorem}

\subsection[Proof of the $T^o$-path formula for $\Cn$ (Theorem \ref{thm:tpath_expansion_formula})]{Proof of the $\boldsymbol{T^o}$-path formula for $\boldsymbol{\Cn}$ (Theorem \ref{thm:tpath_expansion_formula})}
\label{subsection:proof_of_tpath_expansion_formula}

\begin{rem}[{\cite[Lemma 4.5]{MS10}}]
\label{remark:bijection_F}
There is a bijection
\begin{align*}
F\colon \  \{ \text{perfect matchings of the graph }\GT \}  &\to  \big\{ (\tilS,\tilga)\text{-paths on $\tilS$}\big\}, \\
F\colon \  E  &\mapsto  \tilom_E,
\end{align*}
where the diagonals of $\overGT$ correspond to the even-indexed steps of each $(\tilS,\tilga)$-path (see Fig.~\ref{fig:ToIllustrateBijectionsTtwoPerfectMatchingsWithDiagonalsEight})
and each matching $E$ correspond to the odd-indexed steps of $\tilom_E$.
See Figs.~\ref{fig:PMom_2k_2k_plus2} and~\ref{fig:PMom_2kplus2_2k_plus4}.
\end{rem}

\newcommand\ThreeSubpathsOnPMScale{0.2}
\begin{figure}[t!]
\renewcommand*\thesubfigure{(\arabic{subfigure})}
\centering
\mbox{
\subfigure[edge $r \bijection \tilom_{2k+1}$.]{
\TikzPMLRR{\ThreeSubpathsOnPMScale}}\quad
\subfigure[edge $\ell \bijection \tilom_{2k+1}$.]{
\TikzPMLLR{\ThreeSubpathsOnPMScale}}\quad
\subfigure[edge $\dot{r} \bijection \tilom_{2k+1}$.]{
\TikzPMLRdotR{\ThreeSubpathsOnPMScale}}
}
\caption{The subpaths on $\overGT$ corresponding to the subpaths $(\tilom_{2k},\tilom_{2k+1},\tilom_{2k+2})=(\ell, ~~, r)$ from Fig.~\ref{fig:tilom_2k_2k_plus2}.}
\label{fig:PMom_2k_2k_plus2}

\mbox{
\subfigure[edge $r \bijection \tilom_{2k+3}$.]{
\TikzPMRLL{\ThreeSubpathsOnPMScale}}\quad
\subfigure[edge $\ell \bijection \tilom_{2k+3}$.]{
\TikzPMRRL{\ThreeSubpathsOnPMScale}}\quad
\subfigure[edge $\ddot{r} \bijection \tilom_{2k+3}$.]{
\TikzPMRRddotL{\ThreeSubpathsOnPMScale}}
}
\caption{The subpaths on $\overGT$ corresponding to the subpaths $(\tilom_{2k+2},\tilom_{2k+3},\tilom_{2k+4})=(r, ~~,\ell)$ from Fig.~\ref{fig:tilom_2kplus2_2k_plus4}.}
\label{fig:PMom_2kplus2_2k_plus4}
\end{figure}

The bijections $\overpi$ (of Lemma~\ref{lemma:induced_covering_map_bijection})
and $F$ compose to form a bijection between the perfect matchings of the snake graph~$\GT$ and the $(T^o,\ga)$-paths on~$\Cn$.

\begin{theorem}[snake graph matchings to $(T^o,\ga)$-paths]
\label{thm:overpi_F}
Consider an arc $\ga$ of $\Cn$ and the $(T^o,\ga)$-paths and snake graph $\GT$ corresponding to~$\ga$.
The map
\begin{align*}
\overpi \circ F \colon \  \{ \text{perfect matchings in } \GT \}   & \to \{(T^o,\ga) \text{-paths on }\Cn \},  \\
E  & \mapsto \om_E:= \overpi(\tilom_E)
\end{align*}
is a bijection
where the edges of each matching $E$ correspond to the odd-indexed steps of the $(T^o,\ga)$-path $\om_E$, see Fig.~{\rm \ref{fig:ToIllustrateBijectionsTtwoPerfectMatchingsWithDiagonalsEight}}.
In particular, the weight $x(E)$ of the perfect matching $E$ is equal to the numerator of the Laurent monomial $x(\om_E)$ $($see Definitions {\rm \ref{def:laurent_monomial_tpath}} and {\rm \ref{def:perfect_matchings_and_weights})}.
\end{theorem}

\begin{rem}[matching restrictions and $(T^o,\ga)$-subpaths]\label{remark:correspondence_matchings_subpaths}
The four $(T^o,\ga)$-subpaths of Lem\-ma~\ref{lemma:four_paths}(\ref{item:four_paths_lrrll_backtrack})--(\ref{item:four_paths_llrrl_nonbacktrack}) correspond to the four possible matchings of the triple-tiles of $\GT$.
Compare Figs.~\ref{fig:om_2k_2k_plus2}--\ref{fig:om_2kplus2_2k_plus4}, Figs.~\ref{fig:tilom_2k_2k_plus2}--\ref{fig:tilom_2kplus2_2k_plus4}, and Figs.~\ref{fig:PMom_2k_2k_plus2}--\ref{fig:PMom_2kplus2_2k_plus4}.
\end{rem}

\begin{proof}[Proof of Theorem \ref{thm:tpath_expansion_formula}]
Applying the bijection $\overpi \circ F$ to the formula of \cite[Theorem 4.10]{MSW11} yields{\samepage
\begin{gather*}
x_\ga  =
\frac{1}{\text{cross}(T^o,\ga)} \sum_{\substack{E \text{ perfect matchings}\\ \text{of $\GT$}}} x(E)
 =
\frac{1}{\text{cross}(T^o,\ga)} \sum_{(T^o,\ga)\text{-paths }\om}  \left( \prod_{k=0}^d x(\om_{2k+1}) \right).
\end{gather*}
Since $\prod\limits_{k=1}^d x(\om_{2k}) = \text{cross}(T^o,\ga)$,
this concludes the proof.}
\end{proof}

\begin{exam}\label{example:nine_pm}
The nine perfect matchings of the graph $\GT$ (from Fig.~\ref{fig:T2 MSW snake graph})
are listed in Fig.~\ref{fig:TtwoPerfectMatchings}.
Per Theorem \ref{thm:overpi_F}, they correspond to the nine $(T^o,\ga)$-paths of Example \ref{example:nine_paths}.
\end{exam}

\section[Combinatorial proof of atomic bases for type $D$]{Combinatorial proof of atomic bases for type $\boldsymbol{D}$}
\label{sec:proof_atomic_basis}
\subsection[Atomic bases for the once-punctured $n$-gons (type $D_n$)]{Atomic bases for the once-punctured $\boldsymbol{n}$-gons (type $\boldsymbol{D_n}$)}\label{subsection:atomic_bases_D}

\begin{defin}[cluster monomials]
\label{def:cluster_monomial}
Let $\AAA$ be a coef\/f\/icient-free cluster algebra.
A \emph{cluster monomial} is a monomial in cluster variables all belonging to a single cluster.
\end{defin}

Recall that, in the case that $\AAA$ arises from a surface, a cluster corresponds to a tagged triangulation and a cluster monomial corresponds to a multi-tagged triangulation (see Theorem~\ref{thm:fst_thm}).
See Def\/inition~\ref{def:compatibility_of_tagged_arcs} for a list of compatible pairs of tagged arcs of $\Cn$.

The concepts of positive elements and atomic bases were f\/irst introduced in~\cite{SZ04} for the case of an annulus with one marked point on each boundary.
\begin{defin}[positive elements and atomic bases]
An element  $y\in \AAA$ is called \emph{positive} if
the Laurent expansion of $y$ in the variables of every cluster of $\AAA$ has non-negative coef\/f\/icients.

A $\Z$-linear basis $\B$ of $\AAA$ is called an \emph{atomic basis} of $\AAA$ if
any positive element of $\AAA$ is a~non-negative $\Z$-linear
combination of~$\B$.
Note that, if such a $\B$ exists, $\B$ is the collection of all indecomposable positive elements
(i.e., elements which cannot be written as a sum of positive elements) of $\AAA$, hence it is unique.
\end{defin}
\begin{rem}
To prove that a collection $\B\subset \AAA$ is an atomic basis of $\AAA$, it suf\/f\/ices to verify the following:
\begin{gather}
 \label{eqn:atomic_basis_B_is_pos_elt}\textnormal{every element of $\B$ is a positive element of $\AAA$,}\\
 \label{eqn:atomic_basis_Z_basis}\textnormal{$\B$ is a $\Z$-linear basis of $\AAA$,}\\
 \label{eqn:y_pos_linear_combi}\textnormal{every positive element $y\in\AAA$ can be written as a  $\Z_{\geq 0}$-linear combination of $\B$.}
\end{gather}
\end{rem}

\begin{theorem}[\cite{Cer11}]\label{thm:DT}\label{thm:Cer11_ADE}
If $\AAA$ is a cluster algebra of type $A$, $D$, or $E$, the cluster monomials of~$\AAA$
form the atomic basis of~$\AAA$.
\end{theorem}
\begin{proof}
\cite{HL10,Nak11}, \cite{CK08}, and \cite{Cer11,CL12} give representation theoretic proofs for~(\ref{eqn:atomic_basis_B_is_pos_elt}), (\ref{eqn:atomic_basis_Z_basis}), and~(\ref{eqn:y_pos_linear_combi}), respectively.
\cite{MSW11} and~\cite{MSW13} provide combinatorial arguments for~(\ref{eqn:atomic_basis_B_is_pos_elt}) and~(\ref{eqn:atomic_basis_Z_basis}), respectively, and~\cite{DT13} gives a combinatorial proof for type~$A$ for~(\ref{eqn:y_pos_linear_combi}).
We present in the remainder of this paper a combinatorial proof for type $D$ for (\ref{eqn:y_pos_linear_combi}) which relies on the $T^o$-path formula (Theorem~\ref{thm:tpath_expansion_formula}) and is inspired by~\cite{DT13}.
\end{proof}

Recall that a multi-tagged triangulation $\Si$ is not compatible with a tagged triangulation~$T$ if~$\Si$ contains an arc~$\si$ which is not in~$T$.
Recall also that a \emph{proper Laurent monomial in variables~$u_i$} is a product of the form $u_1^{c_1}\cdots u_r^{c_r}$ where at least one of the~$c_i$ is negative.
\begin{lem}
\label{lemma:proper_laurent}
Consider a once-punctured polygon~$\Cn$.
For every tagged triangulation~$T$ and every multi-tagged triangulation~$\Si$ which is not compatible with~$T$,
the $T$-expansion of~$x_\Si$ is a~sum of proper Laurent monomials.
\end{lem}
\begin{rem}
A general version of Lemma \ref{lemma:proper_laurent} is known as the \emph{proper Laurent monomial property} (following \cite{CL12}) which was proven for cluster algebras from surfaces and any choice of coef\/f\/icients in \cite{CL12}
and then for any skew-symmetric cluster algebras in \cite{CKLP13},
using representation theoretic arguments in both cases.
\end{rem}

\begin{defin}
Let $N$ be a subset of a tagged triangulation $T$. We write the \emph{degree with respect to $N$} to mean the degree with respect to the cluster variables corresponding to $N$.
\end{defin}

\begin{proof}[Outline of proof of Lemma~\ref{lemma:proper_laurent}]
Suppose~$x_\Si$ is a cluster monomial not compatible with a~tagged triangulation $T$, i.e., $\Si$ contains an arc $\si$ which is not in~$T$.
For brevity, suppose that all tagged arcs of $\Si\backslash T$ are peripheral and the corresponding ideal triangulation $T^o$ has no self-folded triangle.
For full details and the rest of the cases, see the remainder of this paper.
\begin{enumerate}[\text{Step }1:]\itemsep=0pt
\item
Choose a tagged arc $\si\in \Si \backslash T$ such that $\si$ is as close as possible to the puncture and, if possible, $\si$ only crosses every arc of $T$ at most once.
 Let $(T,\si)$-$cross$ (respectively, $(T,\si)$-$doublecross$) be the set of arcs of $T$ which $\si$ crosses (respectively, the set of arcs of $T$ which $\si$ crosses twice).

 \item
 For each $(T,\si)$-path $\om$, we compare the number of odd-indexed steps (contributing to the numerator of $x(\om)$) and the number of even-indexed steps (contributing to the denominator of $x(\om)$) that belong to $(T,\si)$-$cross$ and $(T,\si)$-$doublecross$. See Def\/inition~\ref{def:x_tpath}. This allows us to show that each term in the $T$-expansion of~$x_\si$ is of \emph{negative} degree either with respect to the cluster variables corresponding to $(T,\si)$-$cross$ or with respect to the cluster variables corresponding to $(T,\si)$-$doublecross$.

\item
Similarly, for each factor $x_\be$ in the product of $x_\Si$, we consider a $(T,\be)$-path and compare the number of odd steps versus the number of even steps to show that every term in the $T$-expansion of $x_\be$ has \emph{non-positive} degree with respect to both the cluster variables corresponding to $(T,\si)$-$cross$ and the cluster variables corresponding to $(T,\si)$-$doublecross$.

\item It follows, since $\si\in\Si$ and every term in the $T$-expansion of $x_\Si$ has non-positive degree with respect to both $(T,\si)$-$cross$ and $(T,\si)$-$doublecross$ that every term in the $T$-expansion of $x_\Si$ has \emph{negative} degree with respect to either $(T,\si)$-$cross$ or $(T,\si)$-$doublecross$.
  \end{enumerate}
 Since $(T,\ga)$-$cross$ and $(T,\ga)$-$doublecross$ are subsets of $T$, the $T$-expansion of $x_\Si$ is a sum of proper Laurent monomials.
\end{proof}

\begin{theorem}\label{thm:main_thm}
If $\AAA$ is a coefficient-free cluster algebra of type $D$,
every positive element $y\in\AAA$ is equal
to a linear combination $y=\sum_{\Ga} m_\Ga x_\Ga$ of cluster monomials where each $m_\Ga$ is non-negative.
\end{theorem}
\begin{proof}[Proof of Theorem \ref{thm:main_thm}]
We assume (\ref{eqn:atomic_basis_B_is_pos_elt}) and (\ref{eqn:atomic_basis_Z_basis}) for $\B=\{$cluster monomials$\}$.
The following argument appears in \cite{CL12,DT13, SZ04}, and we include it here for completeness.
Let $y$ be a positive element of a cluster algebra of type $D$. Write $y=\sum_{\Ga} m_\Ga x_\Ga$ as a linear combination of cluster monomials.
Then we prove that every $m_\Ga$ is non-negative as follows.

Let $x_\Ga$ be a cluster monomial from the sum, and let $\Ga$ denote the multi-tagged triangulation corresponding to it.
Choose a tagged triangulation $T:=T_\Ga$ that is compatible with the multi-tagged triangulation $\Ga$, i.e., so that $\ga\in T$ for every arc $\ga \in \Ga$.
Consider a dif\/ferent cluster monomial $x_\Si$ from the sum.
We argue that
\begin{gather}\label{eq:main_thm}
\text{$x_\Ga$ does not appear in the $T$-expansion of $x_\Si$.}
\end{gather}
If all the arcs of $\Si$ were in $T$, then $x_\Si$ is its own $T$-expansion, and we are done. Suppose that~$\Si$ has an arc that is not in~$T$.  By Lemma~\ref{lemma:proper_laurent},
every term in the $T$-expansion of $x_\Si$ is a proper Laurent monomial in the cluster corresponding to $T$.
Since $x_\Ga$ is a monomial (i.e., not a proper Laurent monomial) in this cluster, we have proven~(\ref{eq:main_thm}).

Hence, $m_\Ga$ is equal to the coef\/f\/icient of $x_\Ga$ in the $T$-expansion of~$y$. Since $y$ is a positive element, the coef\/f\/icient of $x_\Ga$ in the $T$-expansion of $y$ is non-negative, as required.
\end{proof}

In Section \ref{sec:atomicbases_notations}, we provide notations and lemmas which are helpful toward proving the lemmas in Sections~\ref{subsection:technical_lemmas_allperipheral} and~\ref{subsection:technical_lemmas_radius}.
Section \ref{subsection:technical_lemmas_allperipheral} (respectively, Section~\ref{subsection:technical_lemmas_radius})
discusses lemmas which are needed to prove Lemma~\ref{lemma:proper_laurent}
for the case where all arcs of $\Si\backslash T$ are peripheral (respectively, for the case where $\Si\backslash T$ contains a radius).
Even though a few of the lemmas can also be proven using the existing snake graph formula~\cite[Theorem~4.3]{MSW11}
(recapped in Section~\ref{subsection:formula_from_perfect_matchings}),
many arguments (such as the proofs of Lemmas~\ref{lemma:Lemma2}, \ref{lemma:si_doublecross_nonpos}, and~\ref{lemma:be}) can be done more easily using the $T^o$-path formula.
Finally, we prove Lemma~\ref{lemma:proper_laurent}
in Section~\ref{sec:proper_laurent}:
if~$\Si$ is not compatible with the tagged triangulation $T$, we can choose a tagged arc $\si\in \Si\backslash T$ such that
every term in the $T$-expansion of~$x_\Si$ has negative degree with respect to some subset $N_\si^T$ of $T\backslash \Si$
by Lemmas~\ref{lemma:Lemma2}, \ref{lemma:si_doublecross_nonpos},
\ref{lemma:si}, \ref{lemma:si_notched}, \ref{lemma:be}, and~\ref{lemma:be_notched}.

\subsection{Notations and technical lemmas to prove 
Lemma \ref{lemma:proper_laurent}}\label{sec:atomicbases_notations}

Recall from Remark \ref{remark:every_tagged_arc_belongs_to_1of3_classes}
that every tagged arc of $\Cn$ belongs to one of three classes: plain radius, notched radius, or peripheral.

\begin{rem}
Per Remark \ref{remark:only_need_to_prove_two_cases},
we only need to prove
Lemma~\ref{lemma:proper_laurent}
and the related statements
for two cases: one where $T$ has all plain radii and
one where $T$ has $2$ parallel radii $r$, $r\notch$.
\end{rem}

\begin{defin}\label{def:cross}
Let $T$ be a tagged triangulation, $T^o$ the corresponding ideal triangulation, and~$\lambda$ a tagged arc.
Let $(T^o,\lambda)$-$cross$ denote the arcs of~$T^o$ that cross~$\lambda$, and
let $(T^o,\lambda)$-$doublecross$ denote be the arcs of~$T^o$ that cross $\lambda$ twice.
(Note that, If $T$ contains two parallel radii~$r$,~$r\notch$, then $(T^o,\lambda)$-$cross$ and $(T^o,\lambda)$-$doublecross$ may contain the associated $\ell$-loop.)

If $T$ contains two parallel radii $r$, $r\notch$ (so that $T^o$ contains a self-folded triangle with radius~$r$ and~$\ell$-loop $\ell$), and if $\lambda$ crosses the $\ell$-loop $\ell$, let
\begin{gather*}
(T,\lambda)\text{-}cross  :=
\text{$\{ r\notch$ and the peripheral arcs of $T$ that cross $\lambda\}$, and}\\
(T,\lambda)\text{-}doublecross  :=
 \text{$\{ r\notch$ and the (peripheral) arcs of $T$ that cross $\lambda$ twice$\}$}.
\end{gather*}

Otherwise, if $\lambda$ does not cross an $\ell$-loop that belongs to $T^o$, let
\begin{gather*}
\text{$(T,\lambda)$-$cross:=(T^o,\lambda)$-$cross$ \ and \ $(T,\lambda)$-$doublecross:=(T^o,\lambda)$-$doublecross$.}
\end{gather*}
\end{defin}

\begin{lem}\label{lemma:replace_ell}
Let $T$ be a tagged triangulation with two parallel radii $r,r\notch$ $($so that $T^o$ has a~self-folded triangle with radius $r$ and $\ell$-loop~$\ell)$.
Suppose $\be$ and $\si$ $($possibly $\be=\si)$ are compatible tagged arcs not in $T$, and either $\si$ is a plain radius or $\be$ is peripheral.
Let $w=(\om_1,\dots,\om_{2d+1})$ be a $(T^o,\be)$-path.
Recall that $x_\ell:=x_r\,x_{r\notch}$.

Then $x(\om)$ having negative $($respectively, non-positive$)$ degree with respect to $(T^o,\si)$-$cross$
implies that $x(\om)$ has negative $($respectively, non-positive$)$ degree with respect to $(T,\si)$-$cross$.

Similarly, if $(T^o,\si)$-$doublecross$ is nonempty and if $x(\om)$ has negative $($respectively, non-positive$)$ degree with respect to $(T^o,\si)$-$doublecross$,
then $x(\om)$ has negative $($respectively, non-positive$)$ degree with respect to $(T,\si)$-$doublecross$.
\end{lem}

Note that this statement is not necessarily true if $\si$ is peripheral and $\be$ is a radius.

\begin{proof}
If $\si$ is a radius (Fig.~\ref{fig:TikzCrossingNooseOne_TikzCrossingNooseD}), then it cannot cross $r$, so $r\notin (T^o,\si)$-$cross$.
Since $r\notin (T,\si)$-$cross$,
 replacing the two gradings $(T^o,\si)$-$cross$ and $(T,\si)$-$cross$ is equivalent to
simply replacing~$\ell$ in~$(T^o,\si)$-$cross$ with $r\notch$ in $(T,\si)$-$cross$, and the degree would stay the same.

If $\be$ is a peripheral arc crossing $r$ (Figs.~\ref{fig:cross_clockwise},~\ref{fig:cross_counterclock_later}), then $(T^o,\si)$-$cross$ contains $r$, but $\om$ goes along $r$ exactly twice in a row, per Lemma~\ref{lemma:four_paths}(\ref{item:four_paths_lrrll_backtrack})--(\ref{item:four_paths_llrrl_nonbacktrack}). So $x(\om)$ (as a $T^o$-monomial) has zero degree with respect to $r$, and we can again replace $\ell$ in $(T^o,\si)$-$cross$ with $r\notch$ in $(T,\si)$-$cross$ to arrive at the same conclusion.

Since $r$ is a radius, $\si$ can cross it at most once, so $r\notin (T^o,\si)$-$doublecross$.
Since $r\notin (T,\si)$-$doublecross$, replacing the two gradings $(T^o,\si)$-$doublecross$ and $(T,\si)$-$doublecross$ is equivalent to
simply replacing $\ell$ in $(T^o,\si)$-$doublecross$ with $r\notch$ in $(T,\si)$-$doublecross$, and the degree would stay the same.
\end{proof}

\begin{notation}\label{notation:lies_inside}
Every ordinary peripheral arc $\lambda$ snips the surface into a smaller once-punctured disk (denoted $\outsidedisk{\lambda}$) and
a region not containing the puncture (denoted $\disk{\lambda}$).
See Fig.~\ref{fig:disk_lambda}.

Let $\si$ be a peripheral arc with endpoints $s$ and $t$ (possibly $s=t$), and let $\be$ be another ordinary arc.
We say that $\be$ \emph{lies in} $\disk{\si}$ (respectively, $\outsidedisk{\si}$) if the interior of $\be$ lies entirely in the interior of $\disk{\si}$ (respectively, $\outsidedisk{\si}$).

Let $v$ be an endpoint of $\be$.
We say that $v$ \emph{lies in} $\disk{\si}$ if $v$ is a marked point on the boundary of $\disk{\si}$, possibly $v=s$ or $t$.
We say that $v$ \emph{lies in} $\outsidedisk{\si}$ if $v$ is a marked point in the interior or on the boundary of $\outsidedisk{\si}$, possibly $v=s$ or $t$.
We say that $v$ \emph{lies strictly in} $\disk{\si}$ (respectively,~$\outsidedisk{\si}$) if $v$ lies in $\disk{\si}$ (respectively,~$\outsidedisk{\si}$) and~\emph{$v$ is not equal to~$s$ nor~$t$}.
We say $v$ \emph{lies outside of} $\disk{\si}$ if~$v$ is not a marked point of $\disk{\si}$ (hence $v$ is not equal to~$s$ nor~$t$).
We use the same language when the region is a disk cut out by two radii.
\end{notation}

\begin{figure}[t!]\centering
\mbox{
\subfigure[Once-punctured disk $\outsidedisk{\lambda}$.]
{\OncePuncturedDiskLambda{0.6}
}\qquad
\subfigure[Disk $\disk{\lambda}$.]
{\DiskLambda{0.6}
}
}
\caption{$\lambda$ cuts out a smaller once-punctured disk $\outsidedisk{\lambda}$ and a disk $\disk{\lambda}$, shaded in gray.}\label{fig:disk_lambda}
\end{figure}

We use Lemmas \ref{lemma:tgk_in_Tsi_cross} and \ref{lemma:tgk_in_Tsi_doublecross} to prove Lemma \ref{lemma:pre_be}(\ref{lemma:pre_be:cross_si}),~(\ref{lemma:pre_be:cross_si_doublecross}).
\begin{lem}
\label{lemma:tgk_in_Tsi_cross}
Suppose $\ga$ and $\si$ are distinct compatible ordinary arcs not in~$T^o$.
If $\si$ crosses~$\tg{k}$, then~$\si$ crosses  arc~$k^\ga$ or~$k+1^\ga$,
the $k$-th and $(k+1)$-th arcs crossed by~$\ga$.
\end{lem}
\begin{proof}
To simplify the proof, assume that $0 < k < d$. There are four possibilities for ways of~$\ga$ to cross the ideal triangle~$\blacktri{k}^\ga$, see Fig.~\ref{fig:ways_to_cross_if_k_not_0}. If $\blacktri{k}^\ga$ is self-folded, then $\tg{k}$ is equal to either arc~$k^\ga$ or~$k+1^\ga$, and we are done. If arc $k^\ga$, arc $k+1^\ga$, and $\tg{k}$ are all distinct, then $\ga_k$ (the segment between~$p_k^\ga$ and~$p_{k+1}^\ga$ contained in $\blacktri{k}$) cuts $\blacktri{k}^\ga$ into two regions, and one of the regions is an ideal quadrilateral $Q=\Quad{\blacktri{k}^\ga}$ with sides $\tg{k}$, $\ga_k$, part of arc $k^\ga$, and part of arc $k+1^\ga$ (see Figs.~\ref{fig:TikzTriangleCrossesTwoEdges}, \ref{fig:TikzTwoVertexTriangleCrossesTwoEdgesNoLoop}, and~\ref{fig:TikzTwoVertexTriangleCrossesTwoEdgesKOneisLoop}).
By assumption, $\si$ crosses~$\tg{k}$,
so~$\si$ cuts through~$Q$.
Since $\si$ does not cross $\ga_k$ (because $\si$ and $\ga$ are compatible by assumption), $\si$ has to cross either arc~$k^\ga$ or~$k+1^\ga$.
\end{proof}

\begin{figure}[t!]\vspace*{-2mm}
\centering
\newcommand\TikzPunctureClosestToSize{0.8}

\subfigure[$\puncture$ closest to arc~$k^\ga$.]{
\hspace{2mm}\TikzPunctureClosestToK{0.7}\hspace{2mm}}
\quad
\subfigure[$\puncture$ closest to arc~${k+1^\ga}$.]{
\hspace{5mm}\TikzPunctureClosestToKOne{0.7}\hspace{5mm}}
\quad
\subfigure[$\puncture$ closest to arc~$\tg{k}$.]{
\hspace{4mm}\TikzPunctureClosestToTGK{0.7}\hspace{4mm}}

\vspace{-1.5mm}

\subfigure[$\puncture$ closest to arc~$k^\ga$.]{
\TikzPunctureClosestToKTwoVertexTriangle{0.7}}
\quad
\subfigure[$\puncture$ closest to arc~${k+1^\ga}$.]{
\hspace{2mm}\TikzPunctureClosestToKOneTwoVertexTriangle{0.7}\hspace{2mm}}
\quad
\subfigure[$\puncture$ closest to arc~$\tg{k}$.]{
\hspace{4mm}\TikzPunctureClosestToTGKTwoVertexTriangle{0.7}\hspace{4mm}}

\vspace{-1.5mm}

\subfigure[$\puncture$ closest to arc~$k^\ga$.]{
\hspace{4mm}\TikzPunctureClosestToKWithSigma{\TikzPunctureClosestToSize}\hspace{4mm}
\label{fig:si_crosses_k_twice_ordinary}
}
\quad
\subfigure[$\puncture$ closest to arc~${k+1^\ga}$.]{
\hspace{4mm}\TikzPunctureClosestToKOneWithSigma{\TikzPunctureClosestToSize}\hspace{4mm}
\label{fig:si_crosses_kone_twice_ordinary}
}
\quad
\subfigure[$\puncture$ closest to arc~$\tg{k}$.]{
\TikzPunctureClosestToTGKWithSigma{\TikzPunctureClosestToSize}
\TikzPunctureClosestToTGKWithSigmaOpeningLeft{\TikzPunctureClosestToSize}
\label{fig:si_crosses_k_or_kone_twice_ordinary}
}

\vspace{-1.5mm}

\subfigure[$\puncture$ closest to arc $k^\ga$.]{
\TikzPunctureClosestToKTwoVertexTriangleWithSigma{\TikzPunctureClosestToSize}
\label{fig:si_crosses_k_twice_twovertextriangle}
}
\quad
\subfigure[$\puncture$ closest to arc ${k+1^\ga}$.]{
\TikzPunctureClosestToKOneTwoVertexTriangleWithSigma{\TikzPunctureClosestToSize}
\label{fig:si_crosses_kone_twice_twovertextriangle}
}
\quad
\subfigure[$\puncture$ closest to arc $\tg{k}$.]{
\TikzPunctureClosestToTGKTwoVertexTriangleWithSigma{\TikzPunctureClosestToSize}
\TikzPunctureClosestToTGKTwoVertexTriangleWithSigmaOpeningRight{\TikzPunctureClosestToSize}
\label{fig:si_crosses_k_or_kone_twice_twovertextriangle}
}
\vspace{-1mm}

\caption{Lemma \ref{lemma:tgk_in_Tsi_doublecross}, when $\si$ crosses $\tg{k}$ twice.}\label{fig:when_si_crosses_tgk_twice}

\vspace{-2mm}

\end{figure}

\begin{lem}
\label{lemma:tgk_in_Tsi_doublecross}
Suppose $\ga$ and $\si$ are $($distinct, compatible$)$ peripheral arcs such that~$\ga$ that is contained in the interior of~{\rm $\disk{\si}$} $($i.e., $\si$ is closer than $\ga$ to the puncture, see Fig.~{\rm \ref{fig:when_si_crosses_tgk_twice})}.
If~$\si$ crosses~$\tg{k}$ twice, then~$\si$ crosses arc  $k^\ga$ twice or~$\si$ crosses~$k+1^\ga$ twice.
\end{lem}
\begin{proof}
To simplify the proof, assume $0<k<d$. Let $v$ be the marked point adjacent to both arcs $k^\ga$ and $k+1^\ga$.
Suppose $\si$ crosses $\tg{k}$ twice.
There are four possibilities for ways of $\ga$ to cross the ideal triangle $\blacktri{k}$, see Fig.~\ref{fig:ways_to_cross_if_k_not_0}.
If $\blacktri{k}^\ga$ is self-folded, then $\tg{k}$ is equal to either arc~$k^\ga$ or~$k+1^\ga$, and we are done.

Suppose arc $k^\ga$, arc $k+1^\ga$, and $\tg{k}$ are all distinct. Since $\ga$ is peripheral, $\ga$ cuts the surface~$\Cn$ into two regions, the smaller once-punctured disk~$\outsidedisk{\ga}$ containing~$\si$, and the disk~$\disk{\ga}$ containing~$v$.

We claim that $\tg{k}$, arc $k^\ga$, and arc $k+1^\ga$ are all peripheral (one of them possibly an $\ell$-loop).
First, $\tg{k}$ must be peripheral because a radius can be crossed at most once.
Because~$v$ is in the disk $\disk{\ga}$, $v$ is not the puncture, and so the arcs $k$ and $k+1$ are also peripheral.

The puncture is either closest to arc $k^\ga$, arc $k+1^\ga$, or $\tg{k}$ (see Fig.~\ref{fig:when_si_crosses_tgk_twice}).
Since $\si$ cuts out a~disk containing~$\ga$, if the puncture is closest to arc $k^\ga$ (respectively, to arc $k+1^\ga$), then $\si$ has to cross arc $k^\ga$ (respectively, arc $k+1^\ga$) twice. See Figs.~\ref{fig:si_crosses_k_twice_ordinary} and~\ref{fig:si_crosses_k_twice_twovertextriangle} (respectively, Figs.~\ref{fig:si_crosses_kone_twice_ordinary} and~\ref{fig:si_crosses_kone_twice_twovertextriangle}).
If the puncture is closest to $\tg{k}$, then $\si$ has to cross either arc $k^\ga$ or arc~$k+1^\ga$ twice
(Figs.~\ref{fig:si_crosses_k_or_kone_twice_ordinary} and~\ref{fig:si_crosses_k_or_kone_twice_twovertextriangle}).
\end{proof}

\begin{lem}\label{lemma:pre_be}
Let $T^o$ be an ideal triangulation.
Suppose $\be$ and $\si$ are compatible ordinary arcs, with $\be\notin T^o$, and possibly $\be=\si$.
Let $\om=(\om_1,\dots,\om_{2d+1})$ be a $(T^o,\be)$-path. Suppose we have one of these four scenarios:

\begin{enumerate}[$i)$]\itemsep=0pt
\item Let $N_\si^{T^o}:=\{$all radii of $T^o\}$ and suppose $\be$ is peripheral.
$($Note that $N_\si^{T^o}$ does not depend on $\si$ here.$)$
\label{lemma:pre_be:radii}
\item\label{lemma:pre_be:radii_minus_si}
Let $\si$ be a radius of $T^o$, and
let $N_\si^{T^o}:=\{$all radii of $T^o\}\backslash\{\si\}$,
and suppose $\be$ is peripheral.
\item Let $N_\si^{T^o}$ denote $(T^o,\si)$-$cross$, if $\si\notin T$. \label{lemma:pre_be:cross_si}
\item Let $N_\si^{T^o}$ denote $(T^o,\si)$-$doublecross$ if $(T^o,\si)$-$doublecross$ is non-empty, and suppose that~$\be$ is a peripheral arc such that either $\be=\si$ or $\be$ is contained in the interior of {\rm $\disk{\si}$}.
\label{lemma:pre_be:cross_si_doublecross}
\end{enumerate}
In all four of these scenarios, if $\om_{2k+1} \in N_\si^{T^o}$,
then there is a step $\om_{2k}$ or $\om_{2k+2}$ which goes along an arc of $N_\si^{T^o}$.
\end{lem}
\begin{proof}[Proof of Lemma \ref{lemma:pre_be}]
Suppose $\om_{2k+1}$ belongs to $N_\si^{T^o}$.
To simplify the proof, we only discuss cases where $k \neq 0,d$.
Recall that (T1) requires $\be$ to cross each even step $\om_{2i}$ of $\om$.

First, we prove the assertion in scenarios (\ref{lemma:pre_be:radii}) and  (\ref{lemma:pre_be:radii_minus_si}).
If $\om_{2k+1}$ goes along a radius from the boundary to the puncture, we see that $\om_{2k+2}$ must be a radius, as needed for scenario (\ref{lemma:pre_be:radii}).
Furthermore, since $\be$ crosses $\om_{2k+2}$ but not $\si$,
we see that $\om_{2k+2}$ is a radius that is not equal to $\si$,
proving scenario (\ref{lemma:pre_be:radii_minus_si}).
If $\om_{2k+1}$ goes from the puncture to the boundary, we can repeat a~similar argument to show that $\om_{2k}\in N_\si^{T^o}$.

In scenarios (\ref{lemma:pre_be:cross_si}) and (\ref{lemma:pre_be:cross_si_doublecross}),
recall that the steps $\om_{2k}$, $\om_{2k+1}$, and $\om_{2k+2}$ must go along edges in the $(k+1)$-th ideal triangle $\triangle_k^\be$
crossed by $\be$, and $\be$ crosses $\om_{2k}$ and $\om_{2k+2}$.
If $\om_{2k+1}$ goes along the same arc as $\om_{2k}$ or $\om_{2k+2}$, then we are done.
Hence assume that the ideal triangle $\triangle_k^\be$ has three distinct edges $\om_{2k}=\taui{k}^\be$, $\om_{2k+1}=\tbk$, $\om_{2k+2}=\taui{k+1}^\be$.

We consider scenario (\ref{lemma:pre_be:cross_si}):
If $\be=\si$, then $\om_{2k}$ and $\om_{2k+2}$ cross $\be=\si$, and we are done. Assume that $\be\neq \si$.
By Lemma \ref{lemma:tgk_in_Tsi_cross}, since $\si$ crosses $\tbk$, then $\si$ crosses $\taui{k}^\be$ or $\taui{k+1}^\be$.

Finally, we consider scenario (\ref{lemma:pre_be:cross_si_doublecross}):
By assumption, $\be=\si$ or $\be$ is peripheral arc contained in the interior of $\disk{\si}$.

First, suppose $\be=\si$. Having distinct $\om_{2k}$, $\om_{2k+1}$, and~$\om_{2k+2}$ means that the segment~$\si_k$ cuts~$\triangle_k$ into two regions, and one of the regions contains~$\om_{2k+1}$.
Since $\si$~cuts~$\om_{2k+1}$ twice, $\si$~cuts~$\triangle_k^\si$ three times. Since the surface~$\Cn$ only contains one puncture, this is impossible,
hence either $\om_{2k}=\om_{2k+1}$ or  $\om_{2k+2}=\om_{2k+1}$.

Next, suppose $\be$ is a peripheral arc ($\be\neq \si$) contained in the interior of $\disk{\si}$.
By Lem\-ma~\ref{lemma:tgk_in_Tsi_doublecross} (see Fig.~\ref{fig:when_si_crosses_tgk_twice}), since $\si$ crosses $\om_{2k+1}=\tbk$ twice,
$\si$ crosses $\taui{k}^\be$ or $\taui{k+1}^\be$ twice.
Hence $\si$ crosses $\om_{2k}$ or $\om_{2(k+1)}$, as needed.
\end{proof}

In particular, Lemma \ref{lemma:pre_be} implies Corollary \ref{cor:pre_be}.

\begin{defin}
Let $\om$ be a $T^o$-path.
Let $N$ be a subset of $T^o$.
A consecutive subpath $\om'=(\om_a,\dots,\om_b)$ of $\om$ is called an \emph{$N$-subpath} if each step $\om_j$ of $\om'$ comes from $N$.
We say that~$\om'$ is a \emph{longest $N$-subpath} of $\om$ if the step $\om_{a-1}$ (if any) and the step $\om_{b+1}$ (if any) do not come from $N$.
\end{defin}

\begin{cor}
\label{cor:pre_be}
Suppose we have the same setup as Lemma~{\rm \ref{lemma:pre_be}}.
\begin{enumerate}[$a)$]\itemsep=0pt
\item \label{itm:cor:pre_be:be_does_not_cross} If $\be$ does not cross any arc of $N_\si^{T^o}$, then $\om$ does not contain any arc of $N_\si^{T^o}$,
and the degree of $x(\om)$ with respect to $N_\si^{T^o}$ is zero.

\item \label{itm:cor:pre_be:x_om_has_pos_deg} If the corresponding term $x(\om)$ has positive degree with respect to $N_\si^{T^o}$,
there must be an odd-length $($of length three or greater$)$
longest $N_\si^{T^o}$-subpath $\overom=\om_{2i-1},\dots,\om_{2j+1}$ of $\om$.
\end{enumerate}
\end{cor}

\begin{proof}
a) If $\be$ does not cross any arc of $N_\si^{T^o}$, then by (T1) there is no even-indexed step $\om_{2k}$ from $N_\si^{T^o}$. By Lemma \ref{lemma:pre_be}, there is no odd-indexed step from $N_\si^{T^o}$.

b) Suppose  $x(\om)$ has a positive degree with respect to $N_\si^{T^o}$. Then $|$odd-indexed steps from $N_\si^{T^o}|>|$even-indexed steps from $N_\si^{T^o}|$ by the $T^o$-path formula, and so there is at least one odd-indexed step from $N_\si^{T^o}$. By Lemma \ref{lemma:pre_be}, each odd-indexed step $\om_{2k+1}\in N_\si^{T^o}$ must either follow $\om_{2k}\in N_\si^{T^o}$ or preceed $\om_{2k+2}\in N_\si^{T^o}$, and so the assertion follows.
\end{proof}

\begin{figure}[t!]
\centering
\DiskWithPunctureTauPeripheral{0.85}
\quad
\TikzTauPeripheralSelfFolded{0.85}
\caption{Lemma \ref{lemma:adjacent}, when $\si$ crosses two peripheral arcs $\al$, $\taui{1}$ of $\blacktri{0}$.}
\label{fig:tau_peripheral}
\end{figure}

\begin{lem}
\label{lemma:adjacent}
Suppose $T^o$ is an ideal triangulation and $\si\notin T^o$ is a peripheral arc from~$v_1$ to~$v_1'$ $($with $v_1\neq v_1')$.
Let $\tau:=\taui{1}$, the edge opposite $v_1$ in the first ideal triangle~$\triangle_0$ through which~$\si$ passes.
Let~$\al$ and $\lambda$ denote the two sides $($the arcs $\ts{0}$, $\ts{-1}$, not necessarily in this order$)$ of~$\triangle_0$ adjacent to~$\si$ at~$v_1$.
Suppose $\si$ crosses $\al$ $($see Fig.~{\rm \ref{fig:tau_peripheral})}.
Then:
\begin{enumerate}[$a)$]\itemsep=0pt
\item \label{lemma:adjacent:al_peripheral} $\al$ is a peripheral arc with distinct endpoints.
\item \label{lemma:adjacent:tau_peripheral} $\lambda$ is a peripheral arc $($or boundary edge$)$ with distinct endpoints.
\item \label{lemma:adjacent:crosses_twice} $\taui{1}$ is peripheral $($possibly an $\ell$-loop$)$
and $\taui{1}\in (T^o,\si)$-$doublecross$.
\item \label{lemma:adjacent:only_one} There is no arc from $(T^o,\si)$-$cross$ that is adjacent to $\si$ at $v_1'$,
and $\al$ is the only arc of $(T^o,\si)$-$cross$ that is adjacent to $\si$.
Consequently, if $w=(\om_1,\dots, \om_{2d+1})$ is a  $(T^o,\si)$-path, then either $\om_1\notin (T^o,\si)$-$cross$ or $\om_{2d+1}\notin (T^o,\si)$-$cross$.
\item \label{lemma:adjacent:crosses_all_radii} $\si$ crosses all radii of $T^o$.
\end{enumerate}
\end{lem}

\begin{proof}
Observe that, since $v_1$ lies on the boundary, the f\/irst ideal triangle $\si$ crosses, $\triangle_0$, is not self-folded, and so all its edges $\taui{1}$, $\lambda$, $\al$ are pairwise distinct.
\begin{enumerate}[a)]\itemsep=0pt
\item
First, we prove part (\ref{lemma:adjacent:al_peripheral}):
We know that $\al$ cannot be a radius, because a radius from $v_1$ cannot cut an arc adjacent to $v_1$.
Furthermore, since an $\ell$-loop based at $v_1$ would have to lie entirely in $\outsidedisk{\si}$, no $\ell$-loop based on $v_1$ can cross $\si$,
hence $\al$ is a non-$\ell$-loop peripheral arc.

\item
Second, we prove part (\ref{lemma:adjacent:tau_peripheral}) of this lemma.
\begin{itemize}\itemsep=0pt
\item Suppose for contradiction that $\lambda$ is an $\ell$-loop. Again, since no $\ell$-loop based on $v_1$ can cross $\si$, we see that $\lambda$ and $\si$ do not cross.
Hence $\lambda$ would cut out a region $\disk{\lambda}$ not containing the puncture but containing $\si$ as well as $\al$.
It is impossible for $\al$ and $\si$ to both cross and be adjacent in this region, hence $\lambda$ is not an $\ell$-loop.

\item Suppose for contradiction that $\lambda$ is a radius.
Then $\blacktri{0}$ is an ordinary triangle with the puncture as one of the vertices.
Hence the region outside of $\blacktri{0}$ contains no puncture, so it is impossible for $\si$ to cross $\blacktri{0}$ a second time. Therefore $\si$ cannot cross $\al$, contradicting our assumption.
\end{itemize}
Hence $\lambda$ is a peripheral arc with distinct endpoints or a boundary edge.

\item
We now prove part (\ref{lemma:adjacent:crosses_twice}) of this lemma.
By (\ref{lemma:adjacent:al_peripheral}) and (\ref{lemma:adjacent:tau_peripheral}),
$\tau=\taui{1}$ is also a peripheral arc (possibly an $\ell$-loop).
Note that the puncture must be closer to $\taui{1}$ than to $\al$.
Hence, the second arc crossed by $\si$ is contained outside of $\disk{\tau}$,
while $\al$ lies in $\disk{\tau}$.
Since $\si$ crosses $\al$, $\si$ must cross $\tau$ the second time, as needed for part (\ref{lemma:adjacent:crosses_twice}).

\item
Next, we prove part (\ref{lemma:adjacent:only_one}) of this lemma.
Suppose that $\al'\in T^o$ is adjacent to the other end\-point~$v_1'$ of~$\si$.
Note that $\al$ cuts out a disk $\disk{\al}$ containing~$v_1'$,
and no arc of~$T^o$ in the interior of $\disk{\al}$ can both be adjacent to~$\si$ and belong to~$(T^o,\si)$-$cross$.
But both $\al$ and $\al'$ belong to~$T^o$ so they do not cross, hence $\al'$ is contained in $\disk{\al}$, and so $\al'$ cannot both be adjacent to~$\si$ and belong to $(T^o,\si)$-$cross$.

Furthermore, by def\/inition
the only arcs that are adjacent to $\si$ at $v_1$ are $\ts{0}$ and $\ts{-1}$ (i.e., $\al$ and $\lambda$).
But $\lambda$ is contained in $\disk{\si}$, so $\lambda$ cannot cross $\si$.
\item
Finally, we prove part (\ref{lemma:adjacent:crosses_all_radii}) of this lemma.
If $\taui{1}$ is an $\ell$-loop, then $\si$ must cross the radius (which is the only radius of $T^o$) containing it.
Otherwise, since $\taui{1}\in T^o$ is peripheral, all radii of $T^o$ are contained in $\outsidedisk{\tau}$.
Since the segment of $\si$ that is contained in $\outsidedisk{\tau}$ is homotopic to the segment of the original surface's boundary along $\outsidedisk{\tau}$
and since the two endpoints of~$\si$ lie \emph{strictly} in $\disk{\tau}$,
it follows that $\si$ must cross all the radii.\hfill {\qed}
\end{enumerate}\renewcommand{\qed}{}
\end{proof}

The following lemma is also used many times throughout the rest of this paper.

\begin{lem}
\label{lemma:xw_degree}
Suppose $T^o$ is an ideal triangulation, $\be\notin T^o$ is an ordinary arc,
and $\om=(\om_1,\dots,$ $\om_{2d+1})$ is a $(T^o,\be)$-path.
Let $N$ be a subset of $T^o$.
\begin{enumerate}[$a)$]\itemsep=0pt
\item \label{lemma:xw_degree:be_nonpos_neg}  \label{lemma:xw_degree:be_neg}
If $|$even steps from $N| \geq |$odd steps from $N|$,
then $x(\om)$ has non-positive degree with respect to~$N$.
Furthermore, inequality implies negative degree.
\item  \label{lemma:xw_degree:si_nonpos} $x(\om)$ has non-positive degree with respect to $(T^o,\be)$-$cross$.
\item  \label{lemma:xw_degree:si_neg}
If $\om_1, \om_{2d+1}\notin (T^o,\be)$-$cross$, then $x(\om)$ has negative degree with respect to $(T^o,\be)$-$cross$.
Hence, if $\be$ is not adjacent to any arc that crosses it, $x(\om)$ has negative degree with respect to $(T^o,\be)$-$cross$.
\end{enumerate}
\end{lem}
\begin{proof}
Consider a $(T^o,\be)$-path $\om=(\om_1,\dots,\om_{2d+1})$.
Part (\ref{lemma:xw_degree:be_nonpos_neg}) follows from the $T^o$-path formula.

The corresponding term $x(\om)$ in the $T^o$-expansion of
$x_\be$ has $d+1$ factors in the numerator
and $d$ factors in the denominator.
All the factors in the denominator correspond to the arcs of~$T^o$ that cross $\be$.
By Lemma~\ref{lemma:adjacent}(\ref{lemma:adjacent:only_one}),
at least one of the endpoints of $\be$ is not adjacent to any arc from $(T^o,\be)$-$cross$, so
either the f\/irst step $\om_1$ or the last step $\om_{2d+1}$ of $\om$ does not come from $(T^o,\be)$-$cross$ by~(T2).
Hence, there are at most $d$ odd steps contributing to the degree. Since there are exactly $d$ even steps from $(T^o,\be)$-$cross$,
$x(\om)$ has non-positive degree with respect to the arcs of $T^o$ which cross $\si$, satisfying part~(\ref{lemma:xw_degree:si_nonpos}).

If the f\/irst and last steps $\om_1$, $\om_{2d+1}$ of $\om$ do not come from $(T^o,\si)$-$cross$, there are at most $d-1$ odd steps which contribute to the degree,
and~(\ref{lemma:xw_degree:si_neg}) follows.
Furthermore, if $\be$ is not adjacent to any arc from $(T^o,\si)$-$cross$,
then $\om_1, \om_{2d+1}\notin(T^o,\si)$-$cross$ by~(T2).
\end{proof}

\subsection[Technical lemmas to prove Lemma \ref{lemma:proper_laurent}
for cases where all arcs of $\Si \backslash T$ are peripheral]{Technical lemmas to prove Lemma \ref{lemma:proper_laurent}
for cases\\ where all arcs of $\boldsymbol{\Si \backslash T}$ are peripheral}
\label{subsection:technical_lemmas_allperipheral}

\subsubsection[Proving that $x_\Si$ has non-positive degree with respect to $(T,\si)$-$cross$]{Proving that $\boldsymbol{x_\Si}$ has non-positive degree with respect to $\boldsymbol{(T,\si)}$-$\boldsymbol{cross}$}

\begin{lem}\label{lemma:si_nopos_wrt_To}
Let $T$ be a tagged triangulation and
suppose $\si$ is a tagged peripheral arc not in~$T$.
\begin{enumerate}[$a)$]\itemsep=0pt
\item Then each termin the $T$-expansion of $x_\si$ has non-positive degree with respect to $(T,\si)$-$cross$.
\item If $(T^o,\si)$-$doublecross$ is empty, then each term in the $T$-expansion of $x_\si$ has negative degree with respect to $(T,\si)$-$cross$.
\end{enumerate}
\end{lem}

\begin{proof}
By Lemma \ref{lemma:xw_degree}(\ref{lemma:xw_degree:si_nonpos}), each term in the $T^o$-expansion of $x_\si$ has non-positive degree
with respect to $(T^o,\si)$-$cross$.

By Lemma \ref{lemma:adjacent}(\ref{lemma:adjacent:crosses_twice}), if $(T^o,\si)$-$doublecross$ is empty, then $\si$ is not adjacent to any arc from $(T^o,\si)$-$cross$. Therefore, by Lemma \ref{lemma:xw_degree}(\ref{lemma:xw_degree:si_neg}), $x(\om)$ has negative degree with respect to $(T^o,\si)$-$cross$.

By Lemma \ref{lemma:replace_ell}, both assertions follow.
\end{proof}

\begin{defin}
Let $\Lambda$ be a multi-tagged triangulation containing only peripheral arcs.
We say that $\lambda\in\Lambda$ is \emph{central in $\Lambda$} if $\lambda$ is as close as possibly to the puncture, i.e.,
\begin{gather*}
\text{if $\be\in\Lambda\backslash \{ \lambda \}$, then the disk $\disk{\be}$ does not contain $\lambda$.}
\end{gather*}
\end{defin}
\begin{rem}\label{rem:central_arcs}
Let $\{\lambda_1,\dots,\lambda_r \}$ be the set of central arcs of $\Lambda$.
Then, by def\/inition, every $\be\in\Lambda$ is contained in $\disk{\lambda_j}$ for some $j$.
See Fig.~\ref{fig:choice_of_sigma_lambdas}.
\end{rem}

\begin{figure}[t!]\centering
\subfigure[The peripheral arcs $\lambda_1,\dots,\lambda_r$ as close as possible to the puncture.]
{\ChoiceOfSigmaLambdaKs{0.75}
\label{fig:choice_of_sigma_lambdas}}\qquad
\subfigure[If $(T^o,\si)$-$doublecross$ is non\-empty, then $r=1$.]
{\ChoiceOfSigmaLambdaR{0.75}
\label{fig:choice_of_sigma_lambdas_doublecross}}

\vspace{-1.5mm}

\caption{Proof of Lemma \ref{lemma:proper_laurent}: the case where $\Si\backslash T$ are all peripheral.}
\label{fig:choice_of_sigma_proof}

\vspace{-2mm}

\end{figure}

\begin{figure}[t!]
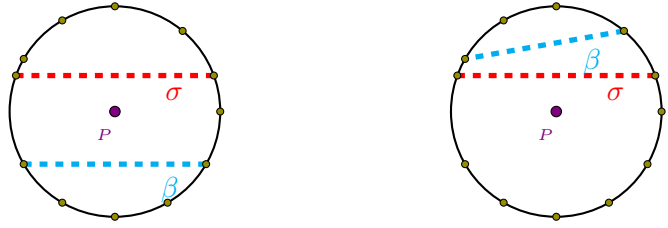
\centering
\subfigure[Correct setup for Lemma \ref{lemma:Lemma2}
but wrong setup for Lemma \ref{lemma:Lemma2_doublecross}: $\si$ is not the only central arc.]
{\ChoiceOfSigmaGoodB{0.7}
\label{fig:choice_of_sigma_good_be_outside_disksi}}
\qquad
\subfigure[Correct setup for both Lemmas \ref{lemma:Lemma2} and \ref{lemma:Lemma2_doublecross}: $\si$ is the only central arc.]
{
\ChoiceOfSigmaGoodA{0.7}
\label{fig:choice_of_sigma_good_be_inside_disksi}}

\vspace{-1.5mm}

\caption{Setups for Lemmas \ref{lemma:Lemma2} and \ref{lemma:Lemma2_doublecross}.}
\label{fig:choice_of_sigma}
\vspace{-2mm}

\end{figure}

\begin{lem}\label{lemma:Lemma2}
Suppose $\Si$ is a multi-tagged triangulation and $T$ is a tagged triangulation.
Suppose $\si\in\Si\backslash T$ is a tagged peripheral arc that is central in $\Si\backslash T$,  i.e.,
\begin{gather}\label{eq:choice_of_sigma}
\text{if $\be\in\Si\backslash T$ is not $\si$, then the disk {\rm $\disk{\be}$} does not contain $\si$ {\rm (}see Fig.~{\rm \ref{fig:choice_of_sigma})}.}
\end{gather}
Let $\be$ be any tagged arc in $\Si$ such that, if $\be\in\Si\backslash T$, then $\be$ is peripheral.
See Figs.~{\rm \ref{fig:choice_of_sigma_good_be_outside_disksi}} and~{\rm \ref{fig:choice_of_sigma_good_be_inside_disksi}}.
Then each term in the $T$-expansion of $x_\be$ has non-positive degree with respect to $(T,\si)$-$cross$.
\end{lem}

\begin{proof}
Let $s^\si$ be the starting point and let $t^\si$ be the f\/inishing point of $\si$. Since $\si$ is not an $\ell$-loop by assumption, $s^\si\neq t^\si$.

If $\be\in T$, then its degree with respect to any of the above gradings is zero since it cannot cross $\si$. If $\be=\si$, we are done by Lemma~\ref{lemma:si_nopos_wrt_To}.

Otherwise, suppose $\be\in\Si\backslash T$ and $\be\neq \si$, and so $\be$ is peripheral by assumption. Let $\om=(\om_1,\dots,\om_{2d+1})$ be a $(T^o,\be)$-path.
First, we prove that
\begin{gather}\label{eq:Lemma2}
\text{the $T^o$-expansion of $x(\om)$ has non-positive degree with respect to $(T^o,\si)$-$cross$.}
\end{gather}
For contradiction, suppose otherwise.
By Corollary \ref{cor:pre_be}(\ref{itm:cor:pre_be:x_om_has_pos_deg}), there must be a longest $(T^o,\si)$-$cross$-subsequence $\overom = (\om_{2i-1}, \dots, \om_{2j+1})$ of $\om$ (of length three or greater).
 Since $\si$ and $\be$ are compatible, there are two possibilities for $\be$:
 either $\be$ lies in the once-punctured disk $\outsidedisk{\si}$ or $\be$ lies in the disk $\disk{\si}$.
 \begin{enumerate}[\text{Case }1:]\itemsep=0pt
 \item
 \emph{First, suppose $\be$ lies in $\outsidedisk{\si}$.} \looseness=-1 See Fig.~\ref{fig:choice_of_sigma_good_be_outside_disksi}.
 Per (\ref{eq:choice_of_sigma}),
 $\be$ has the property that it cuts out the disk $\disk{\be}$ not containing $\si$,
 so no radius can cross both $\be$ and $\si$.
Hence each of the even-indexed arcs $\om_{2i}, \dots, \om_{2j}$ is a peripheral arc
because each of them crosses both~$\be$ and~$\si$.
It then follows that each odd-indexed arc $\om_{2i+1}, \dots, \om_{2j-1}$ is also a peripheral arc.

We claim that $\overom$ begins and f\/inishes in $\outsidedisk{\si}$ (recall that this means possibly at $s^\si$ or $t^\si$):
\begin{itemize}\itemsep=0pt
\item If $i=1$, then $\om_1$ starts
in $\outsidedisk{\si}$ since $\be$ lies in $\outsidedisk{\si}$.
Otherwise, there is a pre\-vious step~$\om_{2i-2}$ which crosses $\be$.
Since $\om_{2i-2}$ does not belong to $(T^o,\si)$-$cross$ by assumption,
$\om_{2i-2}$ ends in $\outsidedisk{\si}$.

\item If $j=d$, then $\om_{2d+1}$ f\/inishes
in $\outsidedisk{\si}$ since $\be$ lies in $\outsidedisk{\si}$.
Otherwise, there is a~next step $\om_{2j+2}$ which crosses $\be$.
Since $\om_{2j+2}$ does not belong to $(T^o,\si)$-$cross$ by assumption,
$\om_{2j+2}$ ends in~$\outsidedisk{\si}$.
\end{itemize}
In fact, by induction starting from $\om_{2i-1}$, every odd-indexed arc of $\overom$ must start from~$\outsidedisk{\si}$.
But this means that $\om_{2j+1}$ is a peripheral arc which begins and ends in $\outsidedisk{\si}$ but also crosses~$\si$, which is impossible.

\item
\emph{Next, suppose that $\be$ lies in {\rm $\disk{\si}$}.} See Figs.~\ref{fig:choice_of_sigma_good_be_inside_disksi} and~\ref{fig:TikzDoubleCrossingActual}.
To simplify our argument,
consider
\begin{gather}
\nonumber\text{a $2$-fold cover $\tilTo$ of $T^o$, a triangulated once-punctured disk containing two lifts} \\
\nonumber\text{of every arc $\tau\in T^o$ and of every boundary marked point of $\Cn$ (see Fig.~\ref{fig:TikzTwoFoldCover}),}\\
\text{every ideal triangle in $\tilTo$ is an ordinary triangle,}  \label{eq:tilTo}\\
\nonumber\text{a lift of $\ga$ is a radius if and only if $\ga$ is a radius of $\Cn$.}
\end{gather}
Since every pair of arcs of $\Cn$ cross each other at most twice,
\begin{gather}\label{eq:lifted_arcs_cross_at_most_once_on_tilTo}
\text{every pair of lifted arcs cross each other at most once on the two-fold cover.}\!\!\!
\end{gather}
Hence, by Lemma \ref{lemma:adjacent}(\ref{lemma:adjacent:crosses_twice}),
\begin{gather}\label{eq:adjacent_lifted_arcs_on_tilTo_do_not_cross}
\text{a pair of adjacent lifted arcs do not cross on the two-fold cover.}
\end{gather}

Choose a lift $\tilsi$ of $\si$ going from $\tils^\si$ to $\tilt^\si$.
Here $\tilsi$ cuts out a disk $\tildisk{\si}$ from $\tilTo$.
Let $\tilbe$ denote the lift of $\be$ which lies in $\tildisk{\si}$,
and let $\tilom$ be the lift of $\om$ which is a $(\tilTo,\tilbe)$-path.
We abuse notation by writing $\overom$ to denote the subpath of $\tilom$ corresponding to $\overom$.

\begin{enumerate}[i)]\itemsep=0pt
\item First, assume that $\overom$ has an earlier step $\om_{2i-2}$ and a later step $\om_{2j+2}$.
We see that
$\overom$ starts from inside $\tildisk{\si}$ (otherwise, $\tilom_{2i-2}$ would have to cross $\tilsi$) and
$\overom$ ends inside $\tildisk{\si}$ (otherwise, $\tilom_{2j+2}$ would have to cross $\tilsi$).
Recall that by inside $\tildisk{\si}$ we mean possibly at $\tils^\si$ or $\tilt^\si$.
In fact, by induction starting from $\om_{2i-1}$, every odd-indexed arc of $\overom$ must start from inside $\tildisk{\si}$.
But this means that $\tilom_{2j+1}$ begins and ends inside $\tildisk{\si}$ and also crosses $\tilsi$.
This requires $\tilsi$ and $\tilom_{2j+1}$ to cross twice, contradicting (\ref{eq:lifted_arcs_cross_at_most_once_on_tilTo}).

\item Second, assume that $i=1$ or $j=d$.
Without loss of generality, assume $i=1$.
There are two cases: either $\be$ and $\si$ are adjacent at $\si$ or not.
\begin{itemize}\itemsep=0pt
\item If $\be$ and $\si$ are adjacent at $s$, then $\tilom_1$ and $\tilsi$ are adjacent and so $\tilom_1$ cannot cross~$\tilsi$ by~(\ref{eq:adjacent_lifted_arcs_on_tilTo_do_not_cross}).
But $\om_1$ crosses~$\si$ by assumption, and
so $\tilom_1$ must cross the other lift~$\tilsi^2$ of~$\si$ which is outside of $\tildisk{\si}$.
Therefore~$\tilom_1$ ends outside of~$\tildisk{\si}$.
\item If $\be$ and $\si$ are not adjacent at $\si$, then $\om_1$ starts inside of $\tildisk{\si}$ because $\tilbe$ lies in $\tildisk{\si}$. Since $\om_1$ crosses $\si$ by assumption, $\tilom_1$ has to cross $\tilsi$ and end outside of $\tildisk{\si}$.
\end{itemize}
So, either way, $\tilom_1$ ends outside of $\tildisk{\si}$.
Hence $\tilom_2, \dots, \tilom_{2j+1}$ forms a zig-zag pattern
crossing $\tilsi$ an even number ($2j$) of times, ending outside of $\tildisk{\si}$.
Since~$\tilbe$ lies in $\tildisk{\si}$, there must be a step $\tilom_{2j+2}$ right after $\overom$ which crosses $\tilom$, contra\-dic\-ting the assumption that $\overom$ is a longest $(\tilTo,\tilsi)$-$cross$-subpath.
\end{enumerate}

This ends the proof for (\ref{eq:Lemma2}).
The conclusion follows from Lemma \ref{lemma:replace_ell}.\hfill{\qed}
\end{enumerate}\renewcommand{\qed}{}
\end{proof}

\begin{figure}[t!]\centering
\subfigure[$T^o$ contains a self-folded triangle with radius $r$ and $\ell$-loop $\ell$.]
{\TikzDoubleCrossingActualCrossingLR{0.65}
}\qquad
\subfigure[$T^o$ has radii $\rho_1,\dots,\rho_f$ ($f\geq 2$).]
{\TikzDoubleCrossingActualRhoOneRhoF{0.65}
}

\caption{$T^o$ when $(T^o,\si)$-$doublecross = \{\theta_1,\dots,\theta_h\}$ is nonempty.}
\label{fig:TikzDoubleCrossingActual}

\subfigure[$T^o$ contains a self-folded triangle with radius $r$ and $\ell$-loop $\ell$.]
{
\TikzTwoFoldCoverCrossingLR{0.65}
}\quad
\subfigure[$T^o$ has radii $\rho_1,\dots,\rho_f$ ($f\geq 2$).]{
\TikzTwoFoldCoverCrossingRhoOneRhoF{0.65}
}
\caption{Two-fold cover $\tildisk{\si}$ of $T^o$ when $(T^o,\si)$-$doublecross = \{\theta_1,\dots,\theta_h\}$ is nonempty.}
\label{fig:TikzTwoFoldCover}
\end{figure}

\subsubsection[Proving that $x_\Si$ has non-positive degree with respect to $(T,\si)$-$doublecross$,
and $x_\si$ has negative degree with respect to $(T,\si)$-$cross$ and $(T,\si)$-$doublecross$]{Proving that $\boldsymbol{x_\Si}$ has non-positive degree with respect\\ to $\boldsymbol{(T,\si)}$-$\boldsymbol{doublecross}$,
and $\boldsymbol{x_\si}$ has negative degree\\ with respect to $\boldsymbol{(T,\si)}$-$\boldsymbol{cross}$ and $\boldsymbol{(T,\si)}$-$\boldsymbol{doublecross}$}

\begin{lem}\label{lemma:Lemma2_doublecross}\label{lemma:si_doublecross_nonpos}
Suppose $\Si$ is a multi-tagged triangulation and $T$ is a tagged triangulation.
Let $\si\in\Si\backslash T$ be a tagged peripheral arc.
\begin{enumerate}[$1)$]\itemsep=0pt
\item \label{lemma:Lemma2_doublecross:nonpos}
Assume $(T,\si)$-$doublecross$ is non-empty, and
$\si$ is central in $\Si\backslash T$ and no other arc is central in $\Si\backslash T$.
Let $\be\in \Si$ $($possibly $\be = \si)$.
\label{lemma:si_doublecross_nonpos:nonpos} Then each term in the $T$-expansion of $x_\be$ has non-positive degree with respect to $(T,\si)$-$doublecross$.
\item
\label{lemma:si_doublecross_nonpos:neg}
Each term in the $T$-expansion of $x_\si$ has negative degree with respect to $(T,\si)$-$cross$ or $(T,\si)$-$doublecross$.
\end{enumerate}
\end{lem}

\begin{proof}
 If $\be\in T$, then its degree with respect to $(T^o,\si)$-$doublecross$ is zero since it cannot cross~$\si$.
 Suppose $\be\in\Si\backslash T$. By Remark \ref{rem:central_arcs}, $\be=\si$ or $\be$ lies in $\disk{\si}$, i.e., if $\be \neq \si$, then $\si$ is closer than $\be$ to the puncture. See Notation \ref{notation:lies_inside}.

 Let $\om^\be=(\om_1,\dots,\om_{2e+1})$ be a $(T^o,\be)$-path.
 If $\be$ does not cross any arc from $(T^o,\si)$-doublecross, then by Corollary \ref{cor:pre_be}(\ref{itm:cor:pre_be:be_does_not_cross}) the degree of $\om$ with respect to $(T^o,\si)$-doublecross is zero. Suppose $\be$ crosses at least an arc from $(T^o,\si)$-doublecross. See Figs.~\ref{fig:TikzDoubleCrossingActual}.

Suppose $x(\om^\be)$ has positive degree with respect to the $\theta_k$'s.
Then by Corollary \ref{cor:pre_be}(\ref{itm:cor:pre_be:x_om_has_pos_deg}) there is a longest $\{ \theta_k$'s$\}$-subpath $\overom^\be=(\om_{2i-1},\dots,\om_{2j+1})$ of $\om$.

Let arc $\taui{1}, \dots, \taui{d}$ be the arcs crossed by $\si$, in order.

Let $\theta_1,\dots,\theta_h$ ($h \geq 1$) denote the arcs in that are crossed by $\si$ twice
(ordered so that $\theta_1$ is the f\/irst and last arc from $(T^o,\si)$-$doublecross$ crossed by~$\si$).
Note that $\theta_1,\dots,\theta_{h-1}$ are contained in the disk $\disk{\theta_h}$ cut by~$\theta_h$,
see Fig.~\ref{fig:TikzDoubleCrossingActual}.
Let $\rho_1,\dots,\rho_f$ ($f\geq 1$) be the radii of~$T^o$, in order of their intersections with $\si$.
Recall that, by Lemma~\ref{lemma:adjacent}(\ref{lemma:adjacent:crosses_all_radii}), $\si$ crosses all radii of~$T^o$.
Let $\taui{1},\dots,\taui{m}$ (if any) be the peripheral arcs which~$\si$ exactly once, right before~$\si$ crosses~$\theta_1$ for the f\/irst time.
Let $\taui{m+h+f+h+1}, \dots, \taui{d}$ (if any) be the peripheral arcs which~$\si$ exactly once, right after~$\si$ crosses~$\theta_1$ for the second time.

Consider a $2$-fold cover $\tilTo$ of $T^o$, a triangulated once-punctured disk containing two lifts of every arc $\tau\in T^o$ and of every boundary marked point of~$\Cn$, see (\ref{eq:tilTo}) and Fig.~\ref{fig:TikzTwoFoldCover}.

Recall that (\ref{eq:adjacent_lifted_arcs_on_tilTo_do_not_cross}) gives that
\begin{gather*}
\text{a pair of adjacent lifted arcs on $\tilTo$ do not cross.}
\end{gather*}

Let $s^\si$ be the starting point of $\si$ and let $t^\si$ be the f\/inishing point of $\si$.
By assumption, $\si$ is not an $\ell$-loop, so $s\neq t$.
Choose a lift $\tilsi$ of $\si,$ running from a lift $\tils$ of $s$ to a lift $\tilt$ of $t$.
Let $\tildisk{\si}$ denote the disk that is cut out by $\tilsi$.
Let $\tilsi^2$ be the other lift of $\si$, running from $\tils^2$ to $\tilt^2$.
Note that $\tilsi$ and $\tilsi^2$ never coincide, and, in particular,
\begin{gather}\label{eq:lifts_of_si_share_no_common_endpoint}
\tilsi \text{ and } \tilsi^2 \text{ share no common endpoint.}
\end{gather}
Let $\tildiskplus{\si}{2}$ be the disk that is cut out by $\tilsi$.

Let $\rho_1^1, \dots, \rho_f^1$ denote the lifts of $\rho_1, \dots, \rho_f$ which are crossed by $\tilsi$.
Let $\rho_1^2, \dots, \rho_f^2$ denote the other lifts of $\rho_1, \dots, \rho_f$, lying entirely outside of $\tildisk{\si}$.
Let $\theta_1^1, \dots, \theta_h^1$ (respectively, $\theta_1^2, \dots, \theta_h^2$)
denote the lifts of $\theta_1, \dots, \theta_h$ which $\tilsi$ crosses f\/irst (respectively, second).
In general, let $\tau^1$ denote the lift of $\tau\in T^o$ in
the region cut out by $\rho_f^1$ and $\rho_f^2$ which contains $\tils$, and let $\tau^2$ denote the lift of $\tau\in T^o$ in the other region.

Since $\si$ cuts $\taui{1}, \dots, \taui{m}$, $\theta_1, \dots, \theta_h$, $\rho_1, \dots, \rho_f$, $\theta_h, \dots, \theta_1$, $\taui{m+h+f+h+1}, \dots, \taui{d}$, in this order, we see that $\tilsi$ cuts
\begin{gather*}
\text{$\taui{1}^1, \dots, \taui{m}^1$,
$\theta_1^1, \dots, \theta_h^1$,
$\rho_1^1, \dots, \rho_f^1$,
$\theta_1^2, \dots, \theta_h^2$,
$\taui{m+h+f+h+1}^2, \dots, \taui{d}^2$},
\end{gather*}
 in this order.

If $\tau\neq\theta_k$, then
by assumption $\tau\in T^o$ cuts $\si$ at most once, so
\begin{gather}
\text{if $\tau\neq\theta_k$, then $\tilsi$ cuts at most one of the two lifts $\tau^1$, $\tau^2$ of $\tau$,}\nonumber\\
\text{and a lift of $\tau$ cuts at most one of the two lifts $\tilsi$, $\tilsi^2$ of $\si$.}\label{eq:tils_cuts_once}
\end{gather}

Since $\si$ cuts each arc $\theta_k$ twice,
$\tilsi$ cuts both lifts $\theta_k^1$ and $\theta_k^2$ of $\theta_k$.
Due to (\ref{eq:adjacent_lifted_arcs_on_tilTo_do_not_cross}),
\begin{gather}\label{eq:tils_not_adjacent_to}
\text{$\tilsi$ cannot be adjacent to $\theta_k^1$ or $\theta_k^2$.}
\end{gather}

Let $\tilbe$ be the lift of $\be$ that is contained in $\tildisk{\si}$.
Let $\tilw^\be$ denote the $(\tilTo,\tilbe)$-path corresponding to $\om^\be$, and
we abuse notation by writing $\overom^\be$ to refer to the subpath of $\tilw^\be$ corresponding to $\overom^\be$.

If $\overom^\be$ contains both $\theta_a^1$ and $\theta_b^2$ for some $1\leq a,b \leq h$,
then $\overom^\be$ must contain the radius/radii $\rho_1^1, \dots, \rho_f^1$, contradicting that the fact that $\overom^\be$  is a $\{\theta_k$'s$\}$-subpath. Without loss of generality, assume each step of~$\overom^\be$ goes along~$\theta_k^1$ for some~$k$.

Note that, because $\theta_k$ crosses $\si$ twice,
\begin{gather*}
\text{every $\theta_k^1$ crosses both $\tilsi$ and $\tilsi^2$.}
\end{gather*}
Hence, since every step of $\overom^\be$ is $\theta_k^1$ for some $k$,
\begin{gather}\label{eq:barw_for_beta_cannot_start_at_tils}
\text{every step of $\overom^\be$ begins \emph{strictly} inside $\tildisk{\si}$ or $\tildiskplus{\si}{2}$.}
\end{gather}
Recall that starting \emph{strictly} inside $\tildisk{\si}$ (respectively, $\tildiskplus{\si}{2}$) means at a marked point of $\tildisk{\si}$
(respectively, $\tildiskplus{\si}{2}$) that is not $\tils$ or $\tilt$ (respectively, $\tils^2$ or $\tilt^2$).

If $\be=\si$, then $\overom_1$ starts at $\tils$ and $\overom_{2d+1}$ f\/inishes at $\tilt$, so
by (\ref{eq:barw_for_beta_cannot_start_at_tils})
there must be an earlier step $\tilom_{2i-2}$ and a later step $\om_{2j+2}$ of $\tilom$.
If $\be\neq \si$,
then it is possible to have $i=1$ (if $\be$ is not adjacent to $\si$ at its starting point)
or $j=e$ (if $\be$ is not adjacent to $\si$ at its f\/inishing point).

We claim that
$\overom^\be$ starts strictly inside $\tildisk{\si}$.
For the sake of argument, suppose that $\overom^\be=(\tilom_{2i-1}, \dots, \tilom_{2j+1})$ starts strictly inside of $\tildiskplus{\si}{2}$.
 \begin{enumerate}[\text{Case }1:]\itemsep=0pt
\item First, suppose that $\tilom^\be$ has an earlier step $\tilom_{2i-2}$.
Since $\tilom_{2i-2}$ must cross $\tilbe$, which lies in $\tildisk{\si}$, we see that $\tilom_{2i-2}$ must cross both $\tilsi$ and $\tilsi^2$ to get from strictly inside $\tildisk{\si}$ to strictly inside $\tildiskplus{\si}{2}$. Since $\tilom_{2i-2}$ is not a lift of one of the $\theta_k$'s, this contradicts~(\ref{eq:tils_cuts_once}).

\item
Second, assume that $\be\neq \si$, and we have $i=1$.
Since $\tilbe$ lies in $\tildisk{\si}$ by assumption, either $\tilom_1$ starts at $\tils^\si$ or $\tilt^\si$, or $\tilom_1$ starts \emph{strictly} inside of $\tildisk{\si}$ (see Notation~\ref{notation:lies_inside}).
By~(\ref{eq:barw_for_beta_cannot_start_at_tils}), $\tilom_1$ starts strictly in $\tildisk{\si}$.

\end{enumerate}
Either way, $\overom^\be$ starts strictly in $\tildisk{\si}$.
Per (\ref{eq:barw_for_beta_cannot_start_at_tils}),
we see that
 $(\tilom_{2i-1}, \dots, \tilom_{2j+1})$ forms a~zip-zag pattern bouncing back and forth between (strictly) $\tildisk{\si}$ and (strictly) $\tildiskplus{\si}{2}$.
So, by induction starting from $\tilom_{2i-1}$, every odd-indexed arc of $\overom$ must start strictly in $\tildisk{\si}$ and ends strictly in $\tildiskplus{\si}{2}$.
Hence the last step $\tilom_{2j+1}$ ends strictly in~$\tildiskplus{\si}{2}$.

The step $\tilom_{2j+1}$ cannot be the last step of $\tilom^\be$ because the f\/inishing point of $\tilbe$ lies (strictly) outside of $\tildiskplus{\si}{2}$ per (\ref{eq:lifts_of_si_share_no_common_endpoint}).
Hence $\tilom^\be$ has a next step $\tilom_{2j+2}$.
Since $\tilom_{2j+2}$
must cross $\tilbe$, which lies in $\tildisk{\si}$, we see that $\tilom_{2j+2}$
 must cross both $\tilsi$ and $\tilsi^2$ to go strictly from $\tildiskplus{\si}{2}$ to strictly in $\tildisk{\si}$.
 Since $\tilom_{2j+2}$ is not a lift of one of the $\theta_k$'s,
this contradicts (\ref{eq:tils_cuts_once}).

\medskip

We prove part (\ref{lemma:si_doublecross_nonpos:neg}):

Suppose $\om=(\om_1,\dots,\om_{2d+1})$ is a $(T^o,\si)$-path.
First, we claim that the $T^o$-expansion of~$x(\om)$ has negative degree with respect to either
$(T^o,\si)$-$cross$ or $(T^o,\si)$-$doublecross$:
By Lemma~\ref{lemma:adjacent}(\ref{lemma:adjacent:only_one}),
either $\om_1\notin (T^o,\si)$-$cross$ or $\om_{2d+1}\notin (T^o,\si)$-$cross$.
If necessary, reverse the orientation of $\om$ so that $\om_{2d+1}\notin (T^o,\si)$-$cross$.
If not all the odd steps $\om_1, \om_3, \dots, \om_{2d-1}$ cross~$\si$, then~$x(\om)$ has negative degree with respect to $(T^o,\si)$-$cross$
by Lemma~\ref{lemma:xw_degree}(\ref{lemma:xw_degree:si_neg}).
Hence suppose that all steps of~$\om$ except~$\om_{2d+1}$ come from $(T^o,\si)$-$cross$.

We claim that $\om_{2h+1}\notin (T^o,\si)$-$doublecross$:

Since $\om_1\in (T^o,\si)$-$cross$,
Lemma \ref{lemma:adjacent}(\ref{lemma:adjacent:crosses_twice}) gives us $\om_2=\theta_1$.
Since $\om_1$ crosses $\si$, its lift~$\tilw_1$ (which lies in $\tildiskplus{\theta_1}{1}$) must cross a lift of $\si$.
Because $\tilw_1$ is adjacent to $\tilsi$, they cannot cross by~(\ref{eq:adjacent_lifted_arcs_on_tilTo_do_not_cross}),
so $\tilw_1$ must cross the other lift $\tilsi^2$ of $\si$ which is located (strictly) outside of $\tildisk{\si}$.
Hence~$\tilw_1$ goes from $\tils$ to outside of~$\tildisk{\si}$,
and $\tilw_2$ goes along $\theta_1^1$ from outside to strictly inside of~$\tildisk{\si}$.
By (T1) and (T2), the steps $\tilw_2, \tilw_3, \dots, \tilw_{2h-1}, \tilw_{2h}$ make a zig-zag pattern along the arcs $\theta_1^1, \dots, \theta_h^1$,
so that $\tilw_{2h}$ goes along $\theta_h^1$ and ends strictly in~$\tildisk{\si}$.
Hence, by (T2), $\tilw_{2h+1}$ either goes along $\theta_h^1$ or $\rho_1^1$.
If  $\tilw_{2h+1}$ goes along $\theta_h^1$, it goes from strictly inside to outside of~$\tildisk{\si}$.
Since the point of $\rho_1^1$ that is adjacent to $\theta_1^1$ is strictly in $\tildisk{\si}$,
it is impossible for~$\tilw_{2(h+1)}$ to go along~$\rho_1^1$.
This contradicts the fact that~$\tilw_{2(h+1)}$ must go along~$\rho_1^1$.
Hence $\tilw_{2h+1}$ goes along~$\rho_1^1$,
as needed to show that $w_{2h+1}\notin (T^o,\si)$-$doublecross$.

By part (\ref{lemma:si_doublecross_nonpos:nonpos}), the odd-indexed steps $\om_{2(h+f+1)-1}, \dots,\om_{2(h+f+h)+1}$ cannot all belong to  $(T^o,\si)$-$doublecross$, so there are at most $h$ odd-indexed steps coming from  $\om_{2(h+f+1)-1}, {\dots},\om_{2(h+f+h)+1}$.
As there are at most  $(h-1)+h$ odd-indexed steps from $\om_3, \om_5, \dots, \om_{2h-1}$, there are at most $(h-1)+h$ odd-indexed steps of $\om$ coming from $(T^o,\si)$-$doublecross$ total.
There are exactly $2h$ even steps of $\om$ from $(T^o,\si)$-$doublecross$ by (T1).
Hence, by Lemma \ref{lemma:xw_degree}(\ref{lemma:xw_degree:be_nonpos_neg}), the term $x(\om)$ has negative degree with respect to $(T^o,\si)$-$doublecross$, as required.

Per Lemma \ref{lemma:replace_ell}, the conclusion follows.
\end{proof}

\begin{figure}[t!]\vspace{-2mm}
\centering
\subfigure[$4$-fold cover of $T$, which has at least two radii $1$, $2$, from Fig.~\ref{fig:wheel}, and $\si$ is a radius.]
{\hspace{6mm}
\FourCoverWheel{1}\hspace{6mm}
\label{fig:T_4fold_peripheral_beta}}
\qquad
\subfigure[$4$-fold cover of $T^o$ from Fig.~\ref{fig:SelfFoldedLocalTriangulation}, where
$T^o$ has a self-folded triangle $r$,~$\ell$. Here $\si$ is a radius.]
{\hspace{6mm}
\FourCoverSelfFolded{1}\hspace{6mm}
\label{fig:selffoldedT_4fold_peripheral_beta}}

\caption{$4$-fold cover $\tilTo$ of $T^o$ where $\be$ is a peripheral arc, and $\si$ is a peripheral arc (respectively, a radius), and $\tildisk{\si}$ is a disk cut out by the peripheral
$\tilsi$ (respectively, the fundamental domain area bounded by the two lifts $\tilsi$ and $\tilsi'$ of $\si$.)}
\end{figure}

\subsection[Technical lemmas to prove Lemma \ref{lemma:proper_laurent} for cases where $\Si$ has a radius not in $T$]{Technical lemmas to prove Lemma \ref{lemma:proper_laurent} for cases\\ where $\boldsymbol{\Si}$ has a radius not in $\boldsymbol{T}$}
\label{subsection:technical_lemmas_radius}

\subsubsection[$\Si\backslash T$ contains a plain radius $\si$]{$\boldsymbol{\Si\backslash T}$ contains a plain radius $\boldsymbol{\si}$}

\begin{lem}[Corollary of Lemma \ref{lemma:xw_degree}(\ref{lemma:xw_degree:si_neg})]\label{lemma:si}
Let $T$ be a tagged triangulation and let
$\si\notin T$ be a~plain radius,
see Figs.~{\rm \ref{fig:TechnicalLemmaUnnotchedRadius}}
and~{\rm \ref{fig:TechnicalLemmaUnnotchedRadiusSelfFolded}}.
Then every term in the $T$-expansion of $x_\si$ has negative degree with respect to
\begin{gather*}
N_\si^T:=(T,\si)\text{-}cross.
\end{gather*}
As defined in Definition {\rm \ref{def:cross}},
\begin{gather*}
\begin{split}
&  (T,\si)\text{-}cross\\
&  =
\begin{cases}
(T^o,\si)\text{-}cross &\text{if $T$ has no parallel radii}\\
\{ r\notch,\text{ and the $($peripheral$)$ arcs of $T$ that cross }\si \} &\text{if $T$ has parallel radii $r$, $r\notch$.}
\end{cases}
\end{split}
\end{gather*}
\end{lem}

\begin{proof}
Since $\si$ is a radius, it is not adjacent to any arc from $(T^o,\si)$-$cross$.
Hence, by Lem\-ma~\ref{lemma:xw_degree}(\ref{lemma:xw_degree:si_neg}),
\begin{gather}
\text{each term in the $T^o$-expansion of $x_\si$ has negative degree}\nonumber\\
\qquad\text{with respect to $(T^o,\si)$-$cross$,}\label{eq:radius_si_has_neg_deg_wrt_To_cross}
\end{gather}
and therefore, by Lemma \ref{lemma:replace_ell}, with respect to $N_\si^{T}:=(T,\si)$-$cross$.
\end{proof}

\begin{figure}[t!]\vspace*{-15mm}
\centering
\subfigure[$T$ has no parallel radii.]
{
\TikzTechnicalLemmaUnnotchedRadius{0.43}
\label{fig:TechnicalLemmaUnnotchedRadius}
}\qquad
\subfigure[$T^o$ has a self-folded triangle $r$,~$\ell$.]
{\hspace{3mm}
\TikzTechnicalLemmaUnnotchedRadiusSelfFolded{0.43}\hspace{3mm}
\label{fig:TechnicalLemmaUnnotchedRadiusSelfFolded}
}\qquad
\subfigure[$\ell_{1}$ in dotted line is the loop surrounding the radius $\rho_1$.]
{
\TikzMinilemmaRadiusGammaPolygonRhoOneLoop{0.43}
\label{fig:TikzMinilemmaRadiusGammaPolygonRhoOneLoop}
}

\caption{Lemmas \ref{lemma:si} and \ref{lemma:be} assume that $\Si\backslash T$ contains a plain radius $\si$.}
\end{figure}

\begin{lem}\label{lemma:be}
Assume the same setup of $T$, $\si$, and $N_\si^T$ as in Lemma~{\rm \ref{lemma:si}} above,
and let $\be$ be a~tagged arc that is compatible with $\si$.
Then each term in the $T$-expansion of~$x_\be$ has non-positive degree with respect to $N_\si^T = (T,\si)$-$cross$.
\end{lem}

\begin{proof}[Proof of Lemma \ref{lemma:be}]
Let $\rho_1, \dots, \rho_f$ denote the radii of~$T^o$
and let~$\si$ run from the puncture to the boundary such that the radii $\rho_1$, $\rho_f$ and peripheral arc $\taui{1}$ (the f\/irst arc crossed by~$\si$)
form the f\/irst triangle that $\si$ crosses.
If $T^o$ has a self-folded triangle $r$, $\ell$, then $f=1$ and $\rho_1=r=\rho_f$ and $\taui{1}=\ell$.
See Figs.~\ref{fig:TechnicalLemmaUnnotchedRadius} and~\ref{fig:TechnicalLemmaUnnotchedRadiusSelfFolded}.

If $\be = \si$, we are done by the previous Lemma~\ref{lemma:si}.
Otherwise, since $\be$ must be compatible with~$\si$, there are only four possibilities:
$\be$ is an arc of $T$,
another plain radius,
the notched radius $\si\notch$,
or a peripheral arc.
\begin{enumerate}[\text{Case }1:]\itemsep=0pt
\item
\emph{First, suppose $\be\in T$ is compatible with $\si$.} Then $\be$ cannot cross $\si$ and $\be$ cannot be a~notched radius, so $\be\notin (T,\si)$-$cross$.

\item
\emph{Second, suppose $\be$ is another plain radius $($not in $T^o)$.}
Let $\om =(\om_1,\dots,\om_{2d+1})$ be a~$(T^o,\be)$-path which runs from the puncture to the boundary.
Then $\om_1$ is a radius, and hence $\om_1$ does not cross $\si$.
Note that $\rho_1$ and $\rho_f$ (possibly $\rho_1=\rho_f$) bound a region (say, $R$) containing $\si$ and all the (peripheral) arcs in $(T^o,\si)$-$cross$.
If $\be$ is outside of~$R$, then $\om$ does not contain any arc of $(T^o,\si)$-$cross$.
If $\be$ is contained in~$R$, then the even-indexed arcs of $\om$ that are crossed by $\si$ form a consecutive string
\begin{gather*}
\om_2=\taui{1}, \dots, \om_{2j}=\taui{j}.
\end{gather*}
By Lemma \ref{lemma:pre_be}(\ref{lemma:pre_be:cross_si}),
the only odd-indexed step/s that may belong to $(T^o,\si)$-$cross$ are $\om_3, \dots, \om_{2j+1}$ since $\om_1$ is a radius (and hence does not cross $\si$).
Hence $\om$ has at $j$ even-indexed steps and at most $j$ odd-indexed steps from $(T^o,\si)$-$cross$,
so, by Lemma~\ref{lemma:xw_degree}(\ref{lemma:xw_degree:be_nonpos_neg}),
$x(\om)$ has non-positive degree with respect to $(T^o,\si)$-$cross$, and hence, by Lemma \ref{lemma:replace_ell},
with respect to $(T,\si)$-$cross$.

\item
\emph{Third, suppose $\be = \si\notch$.} First, suppose $T$ has parallel radii $r$, $r\notch$. By~(\ref{eq:radius_si_has_neg_deg_wrt_To_cross}), each term of the $T^o$-expansion of $x_\si$ has negative degree with respect to $(T^o,\si)$-$cross = \{ \ell,$ (peripheral) arcs of $T$ which cross $\si\}$. Each term $x(\om)$ (as a~$T^o$-monomial) corresponding to a $(T^o,\si)$-path $\om$ has degree $+1$ with respect to $r$, so $x(\om)|_{r\leftrightarrow r\notch}$ (as a~${T\notch}^o$-monomial) has degree~$+1$ with respect to $r\notch$ and negative degree with respect to $(T^o,\si)$-$cross = \{ \ell,$ (peripheral) arcs of $T$ which cross $\si\}$, and, therefore, non-positive degree (as a $T$-monomial) with respect to $(T,\si)$-$cross$.
Since $x_{\si\notch} = x_{\si} |_{ {r}\leftrightarrow {r\notch}}$ by \cite[Proposition~3.15]{MSW11}, every term in the $T$-expansion of $x_{\si\notch}$ has non-positive degree with respect to $(T,\si)$-$cross$, as needed.

Next, suppose $f \geq 2$, and we have, by \cite[Proposition~3.15]{MSW11},
\begin{gather}\label{eqn:to_get_ga_notched_replace_radius}
x_{\si\notch} = x_{\si} |_{ {\rho_i}\leftrightarrow {\rho_i}\notch}
= x_{\si} |_{ {\rho_1}\leftrightarrow {\rho_1}\notch,  {\rho_f}\leftrightarrow {\rho_f}\notch }.
\end{gather}
The second equality is due to the fact that no $(T,\si)$-path would contain the radii $\rho_2, \dots, \rho_{f-1}$, as these radii are not adjacent to~$\si$.
Let $\ell_{1}$ (respectively, $\ell_{f}$) be the loop surrounding the radius~$\rho_1$ (respectively,~$\rho_f$),
see Fig.~\ref{fig:TikzMinilemmaRadiusGammaPolygonRhoOneLoop}.

Since only the starting endpoint (and not the f\/inish endpoint) of $\si$ is adjacent to $\rho_1$ and $\rho_f$,
we see that every term in the $T$-expansion of $x_\si$ has degree $+1$ with respect to $\{ \rho_1, \rho_f \}$,
so
\begin{gather}
\text{each term in the $T$-expansion of $\si\notch$ has degree $+1$}\nonumber\\
\qquad\text{with respect to }
\{ x_{\rho_1\notch}=x_{\ell_1} / x_{\rho_1}, \, x_{\rho_f\notch}=x_{\ell_f} / x_{\rho_f}\}.\label{eqn:si_notch_deg_1_wrt_peripheral_arcs}
\end{gather}
Furthermore, we observe that each term in the $T$-expansions of $x_{\ell_1}$ and $x_{\ell_f}$ has at most degree $+1$ with respect to $(T,\si)$-$cross$.
Since $x_{\rho_i\notch} = x_{\ell_i}/ x_{\rho_i}$ for $i=1,f$,
\begin{gather}
\text{each term in the $T$-expansions of $x_{\rho_1\notch}$ and $x_{\rho_f\notch}$
has at most degree $+1$}\nonumber \\
\qquad\text{with respect to $(T,\si)$-$cross$.}
\label{eqn:rho_1_rho_f_deg_1_wrt_peripheral_arcs}
\end{gather}
By Lemma \ref{lemma:si},
every term in the $T$-expansion of $x_\si$ has negative degree with respect to $(T,\si)$-$cross$.
Combining this fact with (\ref{eqn:to_get_ga_notched_replace_radius}),
(\ref{eqn:si_notch_deg_1_wrt_peripheral_arcs}), and
(\ref{eqn:rho_1_rho_f_deg_1_wrt_peripheral_arcs}),
we see that each term in the $T$-expansion of $x_{\si\notch}$ has non-positive degree with respect to $(T,\si)$-$cross$.

\item
\emph{Fourth, suppose $\be$ is a peripheral arc $($not in $T)$ that does not cross $\si$.}
Consider a~$4$-fold cover $\tilTo$ of $T^o$, a disk with one puncture ($\tilP$) and $4$ lifts for every arc and marked point of~$T^o$.
Let $\tilbe$ be a lift of $\be$.
Since~$\be$ and~$\si$ do not cross,
there are two lifts of~$\si$, say,~$\tilsi$, $\tilsi'$, which bound exactly one fundamental domain containing $\tilbe$.
We denote this fundamental domain by $\tildisk{\si}$.
Let $\tilt$ and $\tilt'$ be the boundary endpoints of~$\tilsi$ and~$\tilsi'$, and let~$\tilP$ be the puncture.
See Fig.~\ref{fig:T_4fold_peripheral_beta} (if~$T$ has no parallel radii) or~\ref{fig:selffoldedT_4fold_peripheral_beta} (if~$T^o$ contains a self-folded triangle~$r$,~$\ell$).

Let $\om = (\om_1,\dots, \om_{2d+1})$ \looseness=-1 denote a $(T^o,\be)$-path and also (by abuse of notation) its corresponding lifted $(\tilTo,\tilbe)$-path.
Per Corollary~\ref{cor:pre_be}(\ref{itm:cor:pre_be:x_om_has_pos_deg}), in order for $x(\om)$ to have a~positive degree with respect $(T^o,\si)$-$cross$,
there must be a~subsequence $\overom$ of consecutive edges that cross $\si$, say, $\overom=\om_{2i-1}, \dots, \om_{2j+1}$,
where the edges directly prior to $\overom$ and directly after $\overom$ do not cross $\si$.
Since no radius of $T^o$ can cross $\si$, all edges in $\overom$ are peripheral arcs.

We claim $\overom$ starts and f\/inishes strictly inside $\tildisk{\si}$, i.e.,
at a marked point of $\disk{\si}$ that is none of $\tilP$, $\tilt$, and $\tilt'$.

If $i=1$, then $\om_1$ starts strictly in $\disk{\si}$ because $\om_1$ is not adjacent to $\tilsi$ nor $\tilsi'$.
Similarly, if $j=d$, then $\om_{2d+1}$ ends strictly $\disk{\si}$ because $\om_{2d+1}$ is not adjacent to $\tilsi$ nor $\tilsi'$.

Otherwise,
$\overom$ must start strictly in $\tildisk{\si}$ because $\om_{2i-2}$ must cross $\tilbe$ but not $\tilsi$ nor $\tilsi'$, and
$\overom$ must end strictly in $\tildisk{\si}$ because $\om_{2j+2}$ must cross $\tilbe$ but not $\tilsi$ nor $\tilsi'$.

By induction,   every odd edge of $\overom$ starts strictly in $\tildisk{\si}$.
But this means that $\om_{2j+1}$ starts and ends strictly $\tildisk{\si}$ even though $\om_{2j+1}$ cuts one of $\tilsi$ and $\tilsi'$, which is impossible.

Hence, $x(\om)$ has non-positive degree with respect to $(T^o,\si)$-$cross$, and, by Lemma \ref{lemma:replace_ell},
also with respect to $(T,\si)$-$cross$, as needed.\hfill{\qed}
\end{enumerate}\renewcommand{\qed}{}
\end{proof}

\subsubsection[$\Si\backslash T$ contains no plain radius but there is a notched radius $\si\notch\in \Si\backslash T$]{$\boldsymbol{\Si\backslash T}$ contains no plain radius but there is a notched radius $\boldsymbol{\si\notch\in \Si\backslash T}$}

The following lemma is helpful toward proving Lemmas \ref{lemma:si_notched} and \ref{lemma:be_notched}.

\begin{lem}\label{lemma:pre_si_be_notched}
Suppose $T$ is a tagged triangulation with no parallel radii, so that $T=T^o$, and~$\rho$ is a radius of~$T$.
Then every term in the $T$-expansion of $x_{\rho\notch}$ has negative degree with respect to $\{$all the radii of~$T\}$.
\end{lem}

\begin{proof}[Proof of Lemma \ref{lemma:pre_si_be_notched}]
Let $\rho_1, \dots, \rho_f$ (with $f \geq 2$) denote all the (plain) radii (in consecutive order) of $T$
such that $\rho = \rho_1$.
Consider the loop $\ell_1$ surrounding the radius $\rho_1$ which crosses $\rho_2$ f\/irst and crosses $\rho_f$ last
(see Fig.~\ref{fig:TikzMinilemmaRadiusGammaPolygonRhoOneLoop}).
Consider a $(T,\ell_1)$-path $\om=(\om_1,\dots,\om_{2d+1})$. Here $d=f-1$.

We f\/irst show that $x(\om)$ has non-negative degree with respect to $\{\rho_1, \dots, \rho_f\}$.
Since the arcs that are crossed by $\ell_1$ are precisely $\rho_2, \dots, \rho_f$,
the even steps of $\om$ are
$\rho_2, \rho_3, \dots, \rho_f$ by (T1).
By Lemma \ref{lemma:xw_degree}(\ref{lemma:xw_degree:be_nonpos_neg}),
The only way for $x(\om)$ to have positive degree with respect to $\{\rho_1, \dots, \rho_f\}$
is if all the odd steps of $\om$ come from $\{\rho_1, \dots, \rho_f\}$.
By induction, $\om_1=\rho_1, \om_3=\rho_2, \dots, \om_{2d-1}=\rho_f$ such that each one
goes from the boundary to the puncture.
But this requires $\om_{2d}=\rho_f$ to go from the puncture to the boundary.
Hence $\om_{2d+1}$ goes from boundary to boundary, so $\om_{2d+1}\notin \{\rho_1, \dots, \rho_f\}$.
Hence each term in the $T$-expansion of $x_{\ell_1}$ has non-negative degree with respect to $\{\rho_1, \dots, \rho_f\}$.

Since
$
x_{\rho_1}\notch = x_{\ell_1} / x_{\rho_1}$,
every term in the $T$-expansion of $x_{\rho_1}\notch$ has negative degree with respect to $\{\rho_1, \dots, \rho_f\}$.
\end{proof}

\begin{figure}[t!]
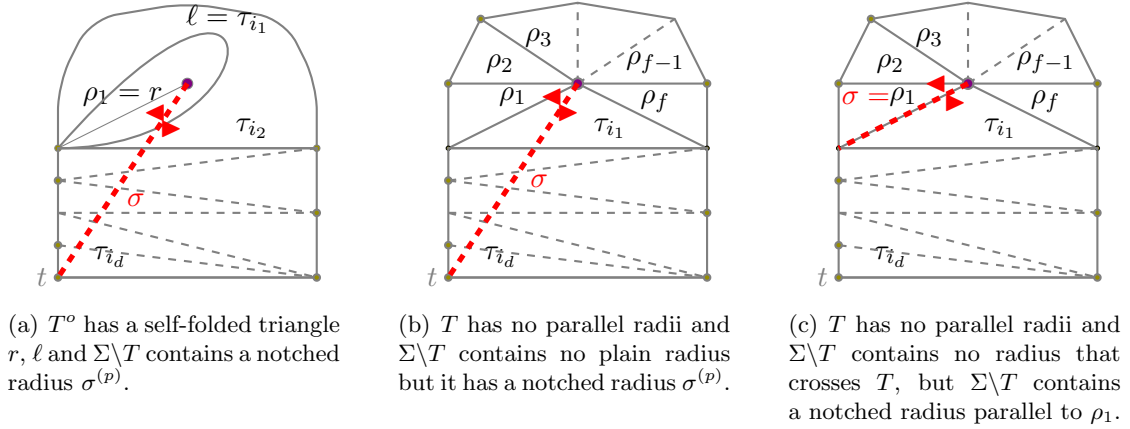
\vspace*{-15mm}
\centering
\subfigure[$T^o$ has a self-folded triangle $r$, $\ell$ and $\Si\backslash T$ contains a notched radius~$\si\notch$.]
{
\TikzTechnicalLemmaNotchedRadiusSelfFolded{0.43}
\label{fig:TechnicalLemmaNotchedRadiusSelfFolded}
}\qquad
\subfigure[$T$ has no parallel radii and $\Si\backslash T$ contains no plain radius but it has a~notched radius~$\si\notch$.]
{
\TikzTechnicalLemmaNotchedRadius{0.43}
\label{fig:TechnicalLemmaNotchedRadius}
}\qquad
\subfigure[$T$ has no parallel radii and $\Si\backslash T$ contains no radius that crosses $T$,
but $\Si\backslash T$ contains a~notched radius parallel to~$\rho_1$.]
{
\TikzTechnicalLemmaRhoOne{0.43}
\label{fig:TechnicalLemmaRhoOne}
}

\vspace{-2mm}

\caption{Three setups for Lemmas \ref{lemma:si_notched} and \ref{lemma:be_notched}.}
\vspace{-1mm}
\end{figure}

\vspace{-1mm}

\begin{lem}\label{lemma:T:A_1_Ga}\label{lemma:si_notched}
Suppose $\Si$ is a multi-tagged triangulation and $T$ is a tagged triangulation
such that $\Si\backslash T$ does not contain any plain radius but it contains a notched radius~$\si\notch$.
There are three cases to consider, Figs.~{\rm \ref{fig:TechnicalLemmaNotchedRadiusSelfFolded}},~{\rm \ref{fig:TechnicalLemmaNotchedRadius}}, and~{\rm \ref{fig:TechnicalLemmaRhoOne}}.

If $T$ has parallel radii $r$, $r\notch$, we denote{\samepage
\begin{gather*}
N_\si^T:=\{r, \text{ and the $($peripheral$)$ arcs of $T$ that cross $\si\}$},
\end{gather*}
where $\si\notch \in \Si\backslash T$
is a notched radius such that $\si\notin T$ $($see Fig.~{\rm \ref{fig:TechnicalLemmaNotchedRadiusSelfFolded})}.}

If $T$ has no parallel radii, we denote
\begin{gather*}
  N_\si^T:=
\!\begin{cases}
\!\{ \text{all radii of }  T \}
&\text{if $\si\notch \in \Si\backslash T$ is a notched radius such that $\si\notin T$}\\[-1pt]
&\text{$($see Fig.~{\rm \ref{fig:TechnicalLemmaNotchedRadius})}},\\
\!\{ \text{all radii of } T \} \backslash \{ \si \}\!\!
&\text{if $\Si \backslash T$ contains no notched radius that crosses $T$, but}\\[-1pt]
&\text{$\si\notch \!\in\! \Si\backslash T$ is a notched radius where $\si\!\in\! T\!$ $($see Fig.\,{\rm \ref{fig:TechnicalLemmaRhoOne}).}}
\end{cases}
\end{gather*}
Then each term in the $T$-expansion of $x_{\si\notch}$
has negative degree with respect to $N_\si^T$.
\end{lem}

\begin{proof}[Proof of Lemma \ref{lemma:si_notched}]
There are three cases to prove,
Figs.~\ref{fig:TechnicalLemmaNotchedRadiusSelfFolded},~\ref{fig:TechnicalLemmaNotchedRadius}, and~\ref{fig:TechnicalLemmaRhoOne}.
\begin{enumerate}[{Fig.~\ref{fig:TechnicalLemmaNotchedRadiusSelfFolded}}:]\itemsep=0pt
\item[Fig.~\ref{fig:TechnicalLemmaNotchedRadiusSelfFolded}:]
In this case $T$ has parallel radii $r$, $r\notch$.
The result follows from the proof of Lem\-ma~\ref{lemma:si} by switching the roles of $r$, $\si$ with $r\notch$, $\si\notch$.
\end{enumerate}
For the remainder of this proof, assume $T$ has no parallel radii. Let $\rho_1, \dots, \rho_f$ denote the radii of $T$ with $f \geq 2$.
\begin{enumerate}[Fig.~\ref{fig:TechnicalLemmaNotchedRadius}]\itemsep=0pt
\item[Fig.~\ref{fig:TechnicalLemmaNotchedRadius}:]
$\Si \backslash T$ contains a notched radius
 $\si\notch$ such that $\si\notin T$. Let $N_\si^T:=\{$all radii of $T\}$, per above assumption.
Assume $\si$ is adjacent to $\rho_1$ and $\rho_f$.
First, recall that (\ref{eqn:to_get_ga_notched_replace_radius}) gives us
\begin{gather*}
x_{\si\notch} = x_{\si} |_{ {\rho_i}\leftrightarrow {\rho_i}\notch}
= x_{\si} |_{ {\rho_1}\leftrightarrow {\rho_1}\notch,  {\rho_f}\leftrightarrow {\rho_f}\notch }.
\end{gather*}
Second, observe that every term in the $T$-expansion of $x_\si$ has degree $+1$ with respect to $\{\rho_1, \rho_f\}$
and does not include any step from $\{\rho_2,\dots,\rho_{f-1}\}$.
Third, by Lemma \ref{lemma:pre_si_be_notched}, every term in the $T$-expansions of $x_{\rho_1}\notch$ and $x_{\rho_f}\notch$
has negative degree with respect to $\{\rho_1, \dots, \rho_f\}$.
Combing these three facts, every term in the $T$-expansion of $x_{\si\notch} = x_{\si} |_{ {\rho_1}\leftrightarrow {\rho_1}\notch,  {\rho_f}\leftrightarrow {\rho_f}\notch}$
has negative degree with respect to $\{\rho_1,\dots,\rho_f\}$.

\item[Fig.~\ref{fig:TechnicalLemmaRhoOne}:]
$\Si \backslash T$ contains no radius that crosses $T$,
but $\si\notch \in \Si\backslash T$ is a notched radius where $\si\in T$.
Let $\rho_1$ denote the radius $\si$ of $T$.
Let $N_\si^T:=\{ \rho_2, \dots,\rho_f \}$, per above assumption.
Consider the loop $\ell_1$ around $\rho_1$ (see Fig.~\ref{fig:TikzMinilemmaRadiusGammaPolygonRhoOneLoop})
and a $(T,\ell_1)$-path $\om$.
Observe that $\om$ has even steps $\rho_2, \rho_3, \dots, \rho_f$, and the f\/irst and last steps of $\om$ are not in $\{ \rho_2, \dots,\rho_f \}$, so $x(\om)$ has negative degree with respect to  $\{ \rho_2, \dots,\rho_f \}$ by Lemma \ref{lemma:xw_degree}(\ref{lemma:xw_degree:be_nonpos_neg}).
Since
\begin{gather*}
x_{{\rho_1}\notch} = \frac{x_{\ell_1}}{x_{\rho_1}},
\end{gather*}
each term in the $T$-expansion of $x_{\rho_1\notch}$ also has negative degree with respect to $\{ \rho_2, \dots,\rho_f \}$.\hfill{\qed}
\end{enumerate}\renewcommand{\qed}{}
\end{proof}

\begin{lem}\label{lemma:be_notched}
Assume the same setup and three different cases as Lemma~{\rm \ref{lemma:si_notched}} above.
If $\be\in\Si$, then each term in the $T$-expansion of $x_\be$ is of non-positive degree with respect to~$N_\si^{T}$.
\end{lem}
\begin{proof}
As in the proof of Lemma \ref{lemma:si_notched}, there are three cases to prove,
Figs.~\ref{fig:TechnicalLemmaNotchedRadiusSelfFolded},~\ref{fig:TechnicalLemmaNotchedRadius}, and~\ref{fig:TechnicalLemmaRhoOne}.

We prove the case of Fig.~\ref{fig:TechnicalLemmaNotchedRadiusSelfFolded} where
$T$ has parallel radii $r$, $r\notch$
 and $N_\si^T$$=\{r$, and the arcs of $T$ that cross $\si \}$
by switching the roles of $r$, $\si$ with $r\notch$, $\si\notch$ in the proof of Lemma \ref{lemma:be}.

For the remainder of this proof, we prove the two setups of
Figs.~\ref{fig:TechnicalLemmaNotchedRadius} and~\ref{fig:TechnicalLemmaRhoOne} for when $T$ has no parallel radii.
Let $\rho_1, \dots, \rho_f$ denote the radii of~$T$ with~$f \geq 2$.

If $\be=\si\notch$, we are done by the previous Lemma \ref{lemma:si_notched}.
Otherwise, since $\be$ is compatible with~$\si\notch$, there are four possibilities:
$\be$ is an arc of $T$, a notched radius that crosses~$T$, a~notched radius $\rho_k\notch$ parallel to a radius of~$T$, or a peripheral arc.
\begin{enumerate}[\text{Case }1:]\itemsep=0pt
\item
\emph{First, suppose $\be\in T$}.
Since $\be$ is compatible with the notched radius $\si\notch$,
we see that $\be\notin N_\si^T$ since either $N_\si^T = \{$all radii of $T\}$ (if $\si\notin T$) or  $N_\si^T = \{$all radii of $T\} \backslash \{ \si \}$ (if~$\si$ is a radius of~$T$).

\item
\emph{Second, suppose $\be$ is a notched radius that crosses $T$}.
The case of Fig.~\ref{fig:TechnicalLemmaNotchedRadius} where $N_\si^T=\{$all radii of $T \}$ is done by Lemma \ref{lemma:si_notched} (for the same case).
We do not need to consider
the case of Fig.~\ref{fig:TechnicalLemmaRhoOne} because it is assumed that $\Si \backslash T$ does not contain any notched radius that crosses $T$.

\item
\emph{Third, suppose that $\be$ is a notched radius $\rho_k\notch$ not equal to $\si\notch$}.
The case of Fig.~\ref{fig:TechnicalLemmaNotchedRadius} where $N_\si^T=\{$all radii of $T \}$ is done by Lemma \ref{lemma:pre_si_be_notched}.

To prove the case of Fig.~\ref{fig:TechnicalLemmaRhoOne}, let $\rho_1$ denote $\si$.
Suppose $\be$ is a notched arc $\rho_k\notch$ with $k \neq 1$, and consider the loop $\ell_k$ around $\rho_k$ and a $(T,\ell_k)$-path $\om$.
We claim that $x(\om)$ has at most degree $+1$ with respect to $\{ \rho_2, \dots,\rho_f \}$.
To prove this,
note that the $f-1$ even arcs of $\om$ are $\rho_{k+1}, \dots, \rho_f, \rho_1,\dots, \rho_{k-1}$, i.e., all the radii of $T$ except for $\rho_k$.
For~$x(\om)$ to have degree $+2$ or more with respect to $\{ \rho_2, \dots,\rho_f \}$,
we need all $f$ odd arcs of~$\om$ to be in $\{ \rho_2, \dots,\rho_f \}$.
But either the step right before or right after the even-step along $\rho_1$ of $\om$ needs to be either $\rho_1$ or another (peripheral) arc.
Hence $x(\om)$ has at most degree $+1$  with respect to $\{ \rho_2, \dots,\rho_f \}$.
Since{\samepage
\begin{gather*}
x_{{\rho_k}\notch} = \frac{x_{\ell_k}}{x_{\rho_k}},
\end{gather*}
each term of $x_{\rho_k}$ has at most degree $0$ with respect to $\{ \rho_2, \dots,\rho_f \}$, as needed.}

\item
\emph{Fourth, suppose $\be$ is a peripheral arc not crossing $\si$.}
Let $\om = (\om_1, \dots, \om_{2d+1})$ be a~$(T,\be)$-path.
Suppose for contradiction that $x(\om)$ is of positive degree with respect to~$N_\si^T$.
By Corollary~\ref{cor:pre_be}(\ref{itm:cor:pre_be:x_om_has_pos_deg}),
\begin{gather}
\text{we must have a subsequence
$\overom=\om_{2i-1}, \dots, \om_{2j+1}$ of $\om$}\nonumber\\
\text{(of length three or greater) where all the steps of $\overom$ belong to $N_\si^T$}\nonumber\\
\text{while the step before and the one after do not belong to $N_\si^T$.}\label{eqn:lemma:be_notch:subsequence}
\end{gather}
We claim that $\overom$ must starts at the boundary (as opposed to the puncture).
Otherwise, $\om$ would contain an earlier step $\om_{2i-2}$ which must be a radius of $T$ going from the boundary to the puncture.
For the case where $\si$ is not a radius of $T$ and $N_\si^T=\{$all radii of $T\}$ (Fig.~\ref{fig:TechnicalLemmaNotchedRadius}),
this shows that $\om_{2i-2}\in N_\si^T$.
For the case where $\si=\rho_1$ and $N_\si^T = \{\rho_2,\dots,\rho_f\}$ (Fig.~\ref{fig:TechnicalLemmaRhoOne}),
since $\be$ cannot cross $\rho_1$, we see that $\om_{2i-2}$ has to come from $\{ \rho_2, \dots,\rho_f \}=N_\si^T$.
For both cases, this contradicts the assumption (\ref{eqn:lemma:be_notch:subsequence}) that no step before $\overom$ comes from $N_\si^T$,
hence $\overom$ must starts at the boundary.
Similarly, $\overom$ must end at the boundary (as opposed to the puncture).

Since every step of $\overom$ is a radius, it follows by induction that every odd step of $\overom$ starts at the boundary.
But this means that $\om_{2j+1}$ is a radius that begins and ends at the boundary, which is impossible.
Hence $x(\om)$ has non-positive degree with respect to~$N_\si^T$.\hfill{\qed}
\end{enumerate}\renewcommand{\qed}{}
\end{proof}

\subsection{Proof of Lemma \ref{lemma:proper_laurent}}\label{sec:proper_laurent}

\begin{proof}[Proof of Lemma \ref{lemma:proper_laurent}]
Suppose $x_\Si$ is a cluster monomial not compatible with a tagged triangulation $T$.

\emph{Suppose that all tagged arcs of $\Si \backslash T$ are peripheral.}
We shall choose a tagged arc $\si\in \Si \backslash T$ that is central and, if possible, crosses every arc of $T$ at most once, as follows. Let $\{\lambda_1,\dots, \lambda_r\}$ be the set of central arcs in $\Si \backslash T$ (see Fig.~\ref{fig:choice_of_sigma_lambdas}). If $r=1$, then choose $\si:=\lambda_r$.

If $r>1$, we choose $\si\in\{\lambda_k$'s$\}$ such that $(T^o,\si)$-$doublecross = \varnothing$. To see that we can do this, suppose that one of these $\lambda_k$ (say, $\lambda_r$) cuts an arc $\tau\in T^o$ twice, see Fig.~\ref{fig:choice_of_sigma_lambdas_doublecross}. Since $\tau$ must cut out a simply-connected region containing all of $\lambda_1,\dots,\lambda_{r-1}$, each of $\lambda_1,\dots,\lambda_{r-1}$ crosses every arc of $T^o$ at most once.

Hence, we satisfy the setup for
\begin{gather*}
\begin{split}
&\text{Lemma \ref{lemma:Lemma2}:} \ \text{$\si$ is central in }\Si\backslash T,\text{ and for}\\
& \text{Lemma \ref{lemma:Lemma2_doublecross}:} \ \text{If $(T^o,\si)$-$doublecross$ is nonempty, then $\si$ is the only central arc in }\Si\backslash T.
\end{split}
\end{gather*}
Lemma \ref{lemma:si_doublecross_nonpos}(\ref{lemma:si_doublecross_nonpos:neg})
tells us that each term in the $T$-expansion of $x_\si$
is of negative degree with respect to either $(T,\si)$-$cross$ or $(T,\si)$-$doublecross$.
At the same time, Lemmas \ref{lemma:Lemma2} and \ref{lemma:Lemma2_doublecross}(\ref{lemma:Lemma2_doublecross:nonpos}) tell us that
every term in the $T$-expansion of $x_\Si$ has non-positive degree with respect to both $(T,\si)$-$cross$ and $(T,\si)$-$doublecross$.
It follows, since $\si\in\Si$, that every term in the $T$-expansion of $x_\Si$ has negative degree with respect to either $(T,\si)$-$cross$ or $(T,\si)$-$doublecross$.
 Since $(T,\si)$-$cross$ and $(T,\si)$-$doublecross$ are subsets of $T$, the $T$-expansion of $x_\Si$ is a sum of proper Laurent monomials.

\looseness=-1
\emph{Suppose now that $\Si \backslash T$ contains a radius $\si$.}
A very similar argument appears in the proof of \cite[Proposition~2.3]{DT13} but we repeat it here for completeness.
Lemma \ref{lemma:si} (if $\si$ is a plain radius)
or Lemma \ref{lemma:si_notched} (if all radii of~$\Si \backslash T$ are notched)
tells us that each term in the $T$-expansion of~$x_\si$
is of negative degree with respect to a subset~$N_\si^T$.
At the same time,
Lemma~\ref{lemma:be} (if~$\si$ is a plain radius)
or Lemma~\ref{lemma:be_notched} (if all radii of $\Si \backslash T$ are notched)
tells us that each term in the $T$-expansion of the other factors of~ $x_\Si$
are of non-positive degree with respect to the same grading~$N_\si^T$.
It follows that each term in the $T$-expansion of $x_\Si$ has negative degree with respect to this grading.
 Since~$N_\si^T$ is a subset of~$T$, the $T$-expansion of~$x_\Si$ is a sum of proper Laurent monomials.
\end{proof}

\subsection{Further directions}
 It is not known for a general cluster algebra whether its set of indecomposable positive elements forms an atomic basis.
Examples of cluster algebras where this fails to happen are those of rank~$2$ where the exchange matrix looks like
\begin{math}
\begin{pmatrix}0  & b \cr -c & 0 \cr
\end{pmatrix}
\end{math} with $bc \geq 5$ \cite{LLZ13}.

Per \cite[Conjecture~1.5]{MSW13}, we are exploring the existence of atomic bases for type $\tilD_{n-1}$ cluster algebras
(which arise from $(n-3)$-gons with $2$ punctures).
Toward this goal, we would like to extend
the $T^o$-path expansion formula to general bordered surfaces with more than 1 puncture.
Additionally, the $T^o$-path formula presented in this paper was only done for ordinary arcs in the setting of ideal triangulations~$T^o$. We are working on proving a similar formula for tagged arcs and for tagged triangulations~$T$.

\subsection*{Acknowledgements}

The authors would like to thank Pasha Pylyavskyy, Vic Reiner, and Peter Webb for looking over an early version of this paper, Hugh Thomas for helpful discussions, and the referees for many useful comments.
The authors were supported by NSF Grants DMS-1067183 and DMS-1148634.

\pdfbookmark[1]{References}{ref}
\LastPageEnding

\end{document}

%% file: figures.tex
\usepackage{tikz}
\usetikzlibrary{arrows,shapes,trees}
\usetikzlibrary{backgrounds}
\usetikzlibrary{calc,decorations.pathmorphing,patterns}
\usetikzlibrary{decorations.markings}

\tikzstyle{NooseArrowAtPosStyleTwo}=[decoration={ markings, mark=at position 0.2 with {\arrow[line width=1.5pt]{stealth};}}, postaction={decorate}]
\tikzstyle{NooseArrowAtPosStyleThree}=[decoration={ markings, mark=at position 0.3 with {\arrow[line width=1.5pt]{stealth};}}, postaction={decorate}]
\tikzstyle{NooseArrowAtPosStyleFour}=[decoration={ markings, mark=at position 0.4 with {\arrow[line width=1.5pt]{stealth};}}, postaction={decorate}]
\tikzstyle{NooseArrowAtPosStyleFive}=[decoration={ markings, mark=at position 0.5 with {\arrow[line width=1.5pt]{stealth};}}, postaction={decorate}]
\tikzstyle{NooseArrowAtPosStyleSix}=[decoration={ markings, mark=at position 0.6 with {\arrow[line width=1.5pt]{stealth};}}, postaction={decorate}]
\tikzstyle{NooseArrowAtPosStyleSeven}=[decoration={ markings, mark=at position 0.7 with {\arrow[line width=1.5pt]{stealth};}}, postaction={decorate}]

\newcommand\PathArrowLineWidth{1.8pt}
\newcommand\PathArrowStyle{to}
\tikzstyle{PathArrowAtPosStyleOne}=[decoration={ markings, mark=at position 0.1 with {\arrow[line width=\PathArrowLineWidth]{\PathArrowStyle};}}, postaction={decorate}]

\tikzstyle{PathArrowAtPosStyleTwo}=[decoration={ markings, mark=at position 0.2 with {\arrow[line width=\PathArrowLineWidth]{\PathArrowStyle};}}, postaction={decorate}]

\tikzstyle{PathArrowAtPosStyleThree}=[decoration={ markings, mark=at position 0.3 with {\arrow[line width=\PathArrowLineWidth]{\PathArrowStyle};}}, postaction={decorate}]

\tikzstyle{PathArrowAtPosStyleFour}=[decoration={ markings, mark=at position 0.4 with {\arrow[line width=\PathArrowLineWidth]{\PathArrowStyle};}}, postaction={decorate}]
\tikzstyle{PathArrowAtPosStyleFourSmall}=[decoration={ markings, mark=at position 0.4 with {\arrow[line width=0.7*\PathArrowLineWidth]{\PathArrowStyle};}}, postaction={decorate}]

\tikzstyle{PathArrowAtPosStyleFive}=[decoration={ markings, mark=at position 0.5 with {\arrow[line width=\PathArrowLineWidth]{\PathArrowStyle};}}, postaction={decorate}]
\tikzstyle{PathArrowAtPosStyleSix}=[decoration={ markings, mark=at position 0.6 with {\arrow[line width=\PathArrowLineWidth]{\PathArrowStyle};}}, postaction={decorate}]
\tikzstyle{PathArrowAtPosStyleSixSmall}=[decoration={ markings, mark=at position 0.6 with {\arrow[line width=0.7*\PathArrowLineWidth]{\PathArrowStyle};}}, postaction={decorate}]

\tikzstyle{PathArrowAtPosStyleSeven}=[decoration={ markings, mark=at position 0.7 with {\arrow[line width=\PathArrowLineWidth]{\PathArrowStyle};}}, postaction={decorate}]
\tikzstyle{PathArrowAtPosStyleSevenSmall}=[decoration={ markings, mark=at position 0.7 with {\arrow[line width=0.7*\PathArrowLineWidth]{\PathArrowStyle};}}, postaction={decorate}]

\tikzstyle{PathArrowAtPosStyleEight}=[decoration={ markings, mark=at position 0.8 with {\arrow[line width=\PathArrowLineWidth]{\PathArrowStyle};}}, postaction={decorate}]
\tikzstyle{PathArrowAtPosStyleEightSmall}=[decoration={ markings, mark=at position 0.8 with {\arrow[line width=0.7*\PathArrowLineWidth]{\PathArrowStyle};}}, postaction={decorate}]

\tikzstyle{PathArrowAtPosStyleNine}=[decoration={ markings, mark=at position 0.9 with {\arrow[line width=\PathArrowLineWidth]{\PathArrowStyle};}}, postaction={decorate}]
\tikzstyle{PathArrowAtPosStyleEnd}=[decoration={ markings, mark=at position 1 with {\arrow[line width=\PathArrowLineWidth]{\PathArrowStyle};}}, postaction={decorate}]

\tikzstyle{PathArrowAtPosStyleReversedFive}=[decoration={ markings, mark=at position 0.5 with {\arrow{to reversed};}}, postaction={decorate}]

\tikzstyle{NonOrientedTriangulationStyle}=[style=thick]

\tikzstyle{TPathStyle}=[blue, very thick,
	style=dashed]

\tikzstyle{orientedTriangulationStyle}=[style=thick]

\tikzstyle{inducedmapTPathStyle}=[
	blue, very thick,
	->, style=dashed, rounded corners = 20pt,
	decoration={ markings, mark=at position .6 with {\arrow[blue, line width=2.5pt]{to};}}]
	
\tikzstyle{inducedmapTPathonCnStyle}=[
	blue, very thick, 
	style=dashed, rounded corners = 30pt,
	decoration={ markings, mark=at position .6 with {\arrow[blue, line width=2.5pt]{to};}}]

\tikzstyle{RedGammaStyle}=[very thick, dashed, red]
\tikzstyle{CyanGammaStyle}=[very thick, dashed, cyan]
\tikzstyle{CyanBetaStyle}=[very thick, dashed, cyan]
\tikzstyle{CyanTilBetaStyle}=[line width=2pt, dashed, cyan]
\tikzstyle{RedTilSigmaStyle}=[line width=2pt, dashed, red]
\tikzstyle{RedSigmaStyle}=[very thick, dashed, red]
\newcommand\TikzPunctureSize{0.04}
\newcommand\TikzPointSize{0.03}

\newcommand \MacroLiftedTriangleOneTwoNoArrow {
\coordinate (P) at (0,0); \coordinate (v1) at (\Tikzn,0); \coordinate (v2) at (0,\Tikzn); \coordinate (y) at (\Tikzn,\Tikzn);

\draw (P) --(v1);
\draw (P) -- (v2);

\draw (v2) -- (y);
\draw (v1) -- (y);

\draw (v2) -- (v1);

\foreach \pos in {(v1),(v2)}
	{\draw [fill=olive] \pos circle (\TikzPointSize*\Tikzn);}

\foreach \pos in {(P),(y)}
	{\draw [fill=violet] \pos circle (\TikzPunctureSize*\Tikzn);}
}

\newcommand \MacroLiftedTriangleOneTwo {
\coordinate (P) at (0,0); \coordinate (v1) at (\Tikzn,0); \coordinate (v2) at (0,\Tikzn); \coordinate (y) at (\Tikzn,\Tikzn);

\draw (P) --(v1) node[pos=0.5, below] {$r$};
\draw (P) -- (v2) node[pos=0.5,left] {$\ccr$};

\draw (v2) -- (y) node[pos=0.5,above]{$\a$};
\draw (v1) -- (y) node[pos=0.5,right]{$\b$};

\draw[NooseArrowAtPosStyleSix] (v2) -- (v1) node[pos=0.7,right] {$\ell$};

\foreach \pos in {(v1),(v2),(y)}
	{\draw [fill=olive] \pos circle (\TikzPointSize*\Tikzn);}

\foreach \pos in {(P)}
	{\draw [fill=violet] \pos circle (\TikzPunctureSize*\Tikzn);}
}

\newcommand\TikzLiftedTriangleOneTwoCrossesA[1]{
\begin{tikzpicture}[scale=#1, DoubleTriangleStyle]
\newcommand\Tikzn{2}
\MacroLiftedTriangleOneTwo

\draw [RedGammaStyle, ->] (P) --(0.8*\Tikzn,1.2*\Tikzn) node[pos=0.4,left] {$\ga$};

\begin{scope}[gray]
\draw (\Tikzn * 0.4,\Tikzn * 0.22) node {$\tilblacktri{0}$};
\draw (\Tikzn * 0.75,\Tikzn *0.72) node {$\tilblacktri{1}$};
\end{scope}
\end{tikzpicture}
}

\newcommand\TikzLiftedTriangleOneTwo[1]{
\begin{tikzpicture}[scale=#1, DoubleTriangleStyle]
\newcommand\Tikzn{2}
\MacroLiftedTriangleOneTwo

\draw [RedGammaStyle, <-] (0.5*\Tikzn,0.8*\Tikzn) -- (0.1,0.1);

\begin{scope}[gray]
\draw (\Tikzn * 0.4,\Tikzn * 0.22) node {$\tilblacktri{0}$};
\draw (\Tikzn * 0.75,\Tikzn *0.72) node {$\tilblacktri{1}$};
\end{scope}
\end{tikzpicture}
}
\newcommand\TikzLiftedTriangleD[1]{
\begin{tikzpicture}[scale=#1, DoubleTriangleStyle]
\newcommand\Tikzn{2}
\MacroLiftedTriangleOneTwo

\draw [RedGammaStyle, ->] (0.5*\Tikzn,0.7*\Tikzn) -- (0.1,0.1);

\begin{scope}[gray]
\draw (\Tikzn * 0.4,\Tikzn * 0.22) node {$\tilblacktri{d}$};
\draw (\Tikzn * 0.75,\Tikzn *0.72) node {$\tilblacktri{d-1}$};
\end{scope}
\end{tikzpicture}
}
\newcommand\TikzCrossingNooseOne[1]{
\begin{tikzpicture}[scale={#1}]
\newcommand\Tikzn{2}
\MacroCrossingLRAB
\draw[RedGammaStyle, <- ] (-0.2*\Tikzn,1.6*\Tikzn) -- (P) node[blue,left, pos=0.4] {$p_1$};
\newcommand\TikzHeight{1.3}
\draw[fill=blue] (-0.1*\Tikzn,\Tikzn*\TikzHeight) circle (0.2*\Tikzn ex);
\end{tikzpicture}
}
\newcommand\TikzCrossingNooseOneAndA[1]{
\begin{tikzpicture}[scale={#1}]
\newcommand\Tikzn{2}
\MacroCrossingLRAB
\draw[RedGammaStyle, -> ]
(P) -- (-0.6*\Tikzn,2.1*\Tikzn) node[pos=0.45,left]{$\ga$};
\newcommand\TikzHeight{1.3}
\end{tikzpicture}
}

\newcommand\TikzCrossingNooseD[1]{
\begin{tikzpicture}[scale={#1}]
\newcommand\Tikzn{2}
\MacroCrossingLRAB
\draw[RedGammaStyle, -> ] (-0.2*\Tikzn,1.6*\Tikzn) -- (P) node[blue,left, pos=0.4] {$p_d$};
\newcommand\TikzHeight{1.3}
\draw[fill=blue] (-0.1*\Tikzn,\Tikzn*\TikzHeight) circle (0.2*\Tikzn ex);
\end{tikzpicture}
}

\newcommand\MacroCrossingLRAB{
\begin{scope}[thick, gray]
	\coordinate (v) at (0,0);
	\coordinate (P) at (0,0.9*\Tikzn);
	\coordinate (y) at (0,2*\Tikzn);
	\draw[thick] (0,0) -- (P);
	\draw (0.08*\Tikzn,0.4*\Tikzn) node {$r$};
	\draw (v) node[below,olive] {$v$}; \draw (y) node[below,olive] {$y$};
	
\begin{scope}[gray]
	\newcommand\TikzLowerPointerLeft{(-\Tikzn*0.5,\Tikzn*0.5)}
	\newcommand	\TikzLowerPointerRight{(\Tikzn*0.5,\Tikzn*0.5)}
	\newcommand\TikzUpperPointerLeft{(-\Tikzn*0.5,\Tikzn)} \newcommand\TikzUpperPointerRight{(\Tikzn*0.5,\Tikzn)}

	\draw[style = thick, rounded corners = 20pt,NooseArrowAtPosStyleSix]	
	(v) .. controls \TikzLowerPointerLeft and \TikzUpperPointerLeft  .. (0,\Tikzn*1.5)
	.. controls \TikzUpperPointerRight and \TikzLowerPointerRight .. (v) node[pos=0.3,right=3pt] {$\ell$};
\end{scope}
	
	\newcommand\TikzLowerPointA{(-1.3*\Tikzn,0.3*\Tikzn)} \newcommand\TikzUpperPointA{(-\Tikzn,1.8*\Tikzn)}
	\newcommand\TikzLowerPointB{(1.3*\Tikzn,0.3*\Tikzn)} \newcommand\TikzUpperPointB{(\Tikzn,1.8*\Tikzn)}
	\draw (v) .. controls \TikzLowerPointA and \TikzUpperPointA .. (y) node[pos=0.6,left] {$\a$};
	\draw (v) .. controls \TikzLowerPointB and \TikzUpperPointB .. (y) node[pos=0.6,right] {$\b$};
\end{scope}
	
	\draw [fill=violet, violet] (P) circle (\TikzPunctureSize*\Tikzn);
	\foreach \pos in {(v),(y)}
	{\draw [fill=olive, olive] \pos circle (\TikzPointSize*\Tikzn );}
}

\newcommand\TikzCrossingCounterclockwiseLRAB[1]{
\begin{tikzpicture}[scale={#1}]
\newcommand\Tikzn{2}
\MacroCrossingLRAB

\newcommand\TikzHeight{0.6}
\draw [RedGammaStyle, dashed, ->] (-\Tikzn *0.6, \Tikzn * \TikzHeight) -- (\Tikzn *0.7, \Tikzn * \TikzHeight) node[blue,below, pos=0.1] {$p_k$};
\draw[fill=blue] (-0.4*\Tikzn,\Tikzn*\TikzHeight) circle (0.2*\Tikzn ex);
\end{tikzpicture}
}

\newcommand\TikzCrossingClockwiseLRAB[1]{
\begin{tikzpicture}[scale={#1}]
\newcommand\Tikzn{2}
\MacroCrossingLRAB
\newcommand\TikzHeight{0.6}
\draw [RedGammaStyle, dashed, <-] (-\Tikzn *0.7, \Tikzn * \TikzHeight) -- (\Tikzn *0.7, \Tikzn * \TikzHeight)
node[blue,below, pos=0.9] {$p_k$};
\draw[fill=blue] (0.4*\Tikzn,\Tikzn*\TikzHeight) circle (0.2*\Tikzn ex);
\end{tikzpicture}
}

\newcommand\MacroDoubleCrossingActualThetaOneThetaH{
\draw[gray, thick, rounded corners=5pt] (v) -- (-2.3*\Tikzn,0) -- (-2.3*\Tikzn,2*\Tikzn) -- (y)
(-2*\Tikzn,0) -- (-1*\Tikzn,2*\Tikzn) node [black,pos=0.5, left] {$\theta_1$};
}
\newcommand\TikzDoubleCrossingActualCrossingLR[1]{
	\begin{tikzpicture}[scale={#1}]
	\newcommand\Tikzn{2}
	\MacroCrossingLRAB
	\MacroDoubleCrossingActualThetaOneThetaH
	
	\newcommand\TikzHeightBelow{0.6}
	\newcommand\TikzHeightExtraBelow{0.3}
	\newcommand\TikzHeightAbove{1.5}
	\draw [RedGammaStyle, dashed, ->, rounded corners=15pt] (-\Tikzn *2, \Tikzn * \TikzHeightBelow) -- (\Tikzn *0.2, \Tikzn * \TikzHeightBelow) -- (\Tikzn *0.2, \Tikzn * \TikzHeightAbove) -- (-2*\Tikzn, \Tikzn*\TikzHeightAbove)
	node[pos=0.8, above] {$\si$};
	
	\draw [CyanBetaStyle, dashed, ->] (-\Tikzn *2, \Tikzn * \TikzHeightExtraBelow) -- (\Tikzn *0.45, \Tikzn * \TikzHeightExtraBelow)
	node[pos=0.2, below] {$\ga$};
	
	\end{tikzpicture}
}

\newcommand\TikzDoubleCrossingActualRhoOneRhoF[1]{
\begin{tikzpicture}[scale={#1}]
\newcommand\Tikzn{2}
\begin{scope}[thick, gray]
	\coordinate (v) at (0,0);
	\coordinate (P) at (0,0.9*\Tikzn);
	\coordinate (y) at (0,2*\Tikzn);
	\draw[thick] (0,0) -- (y) ;
	\draw[black] (-0.1*\Tikzn,0.4*\Tikzn) node {$\rho_1$};
	\draw[black] (-0.2*\Tikzn,1.3*\Tikzn) node {$\rho_f$};
	\draw (v) node[below,white] {$v$};
	
	\newcommand\TikzLowerPointerLeft{(-\Tikzn*0.5,\Tikzn*0.5)}
	\newcommand\TikzLowerPointerRight{(\Tikzn*0.5,\Tikzn*0.5)}
	\newcommand\TikzUpperPointerLeft{(-\Tikzn*0.5,\Tikzn)} \newcommand\TikzUpperPointerRight{(\Tikzn*0.5,\Tikzn)}
	
	\newcommand\TikzLowerPointA{(-1.3*\Tikzn,0.3*\Tikzn)} \newcommand\TikzUpperPointA{(-\Tikzn,1.8*\Tikzn)}
	\newcommand\TikzLowerPointB{(1.3*\Tikzn,0.3*\Tikzn)} \newcommand\TikzUpperPointB{(\Tikzn,1.8*\Tikzn)}
	\draw (v) .. controls \TikzLowerPointA and \TikzUpperPointA .. (y) node[pos=0.5,black] {$\theta_h$};
	\draw (v) .. controls \TikzLowerPointB and \TikzUpperPointB .. (y);
\end{scope}
	
	\draw [fill=violet, violet] (P) circle (\TikzPunctureSize*\Tikzn);
	\foreach \pos in {(v),(y)}
	{\draw [fill=olive, olive] \pos circle (\TikzPointSize*\Tikzn );}

\MacroDoubleCrossingActualThetaOneThetaH

\newcommand\TikzHeightBelow{0.6}
\newcommand\TikzHeightExtraBelow{0.3}
\newcommand\TikzHeightAbove{1.5}
\draw [RedGammaStyle, dashed, ->, rounded corners=15pt] (-\Tikzn *2, \Tikzn * \TikzHeightBelow) -- (\Tikzn *0.5, \Tikzn * \TikzHeightBelow) -- (\Tikzn *0.5, \Tikzn * \TikzHeightAbove) -- (-2*\Tikzn, \Tikzn*\TikzHeightAbove)
node[pos=0.8, above] {$\si$};
\draw [CyanBetaStyle, dashed, ->] (-\Tikzn *2, \Tikzn * \TikzHeightExtraBelow) -- (\Tikzn *0.3, \Tikzn * \TikzHeightExtraBelow)
node[pos=0.2, below] {$\ga$};
\end{tikzpicture}
}

\newcommand\MacroTwoFoldCoverCrossingRL{
\begin{scope}[thick, gray]
	\coordinate (v) at (0,0);
	\coordinate (P) at (0,1*\Tikzn);
	\coordinate (y) at (0,2*\Tikzn);
	\draw[thick] (0,0) -- (y);
	\draw[black] (0.0*\Tikzn,0.6*\Tikzn) node {$r^1$};
	\draw[black] (-0.0*\Tikzn,1.6*\Tikzn) node {$r^2$};
	\draw (v) node[below,olive] {$v^1$}; \draw (y) node[above,olive] {$v^2$};
	
	\newcommand\TikzLowerPointerLeft{(-\Tikzn*0.6,\Tikzn*0.5)}
	\newcommand\TikzLowerPointerRight{(\Tikzn*0.6,\Tikzn*0.5)}
	\newcommand\TikzUpperPointerLeft{(-\Tikzn*0.6,1.5*\Tikzn)}
	\newcommand\TikzUpperPointerRight{(\Tikzn*0.6,1.5*\Tikzn)}
	
	\draw[style = thick, NooseArrowAtPosStyleThree,NooseArrowAtPosStyleSeven]	
	(v) .. controls \TikzLowerPointerLeft and \TikzUpperPointerLeft  .. (y) node[black,pos=0.7] {$\ell^1$}
	(y) .. controls \TikzUpperPointerRight and \TikzLowerPointerRight .. (v) node[black,pos=0.5] {$\ell^2$};
	
	\newcommand\TikzLowerPointA{(-1.3*\Tikzn,0.3*\Tikzn)} \newcommand\TikzUpperPointA{(-0.7*\Tikzn,1.8*\Tikzn)}
	\newcommand\TikzLowerPointB{(0.7*\Tikzn,0.2*\Tikzn)} \newcommand\TikzUpperPointB{(1.3*\Tikzn,1.7*\Tikzn)}
	\draw (v) .. controls \TikzLowerPointA and \TikzUpperPointA .. (-1.4,2*\Tikzn) node[black,pos=0.6] {$\a^1$};
	\draw (1.4,0*\Tikzn)  .. controls \TikzLowerPointB and \TikzUpperPointB .. (y) node[black,pos=0.6] {$\a^2$};
	\draw[fill=olive] (-1.4,2*\Tikzn) circle (\TikzPointActualSize) node [olive, above] {$y^2$};
	\draw (-0.8,2.1*\Tikzn) node {$\b$};
	\draw[fill=olive] (1.4,0*\Tikzn) circle (\TikzPointActualSize) node [olive, below] {$y^1$};
	\draw (0.8,-0.1*\Tikzn) node {$\b$};
	
\end{scope}
	
	\draw [fill=violet, violet] (P) circle (\TikzPunctureActualSize);
	\foreach \pos in {(v),(y)}
	{\draw [fill=olive, olive] \pos circle (\TikzPointSize*\Tikzn );}
}

\newcommand\TikzTwoFoldCoverCrossingLR[1]{
\begin{tikzpicture}[scale={#1}]
\newcommand\Tikzn{2}
\MacroTwoFoldCoverCrossingRL
\newcommand\TikzWidth{2}
\draw[gray, thick,rounded corners=5pt] (v) -- (-\TikzWidth*\Tikzn,0) -- (-\TikzWidth*\Tikzn,2*\Tikzn) -- (y)
(-1.8*\Tikzn,0) -- (-1*\Tikzn,2*\Tikzn) node [black,pos=0.7, left] {$\theta_1^1$};
\draw[gray, thick,rounded corners=5pt] (v) -- (\TikzWidth*\Tikzn,0) -- (\TikzWidth*\Tikzn,2*\Tikzn) -- (y)
(1.8*\Tikzn,2*\Tikzn) -- (1*\Tikzn,0) node [black,pos=0.3, left] {$\theta_1^2$};

\newcommand\TikzHeightBelow{0.4}
\newcommand\TikzHeightAbove{0.8}
\draw [RedTilSigmaStyle, dashed, ->] (-\Tikzn *2, \Tikzn * \TikzHeightAbove) -- (\Tikzn *0.7, \Tikzn * \TikzHeightAbove) -- (2*\Tikzn, \Tikzn*\TikzHeightAbove) node[pos=0.7,above] {$\tilsi$};
\draw[fill=olive] (-\Tikzn *2, \Tikzn * \TikzHeightAbove) circle (\TikzPunctureActualSize) node[olive,left] {$\tils$};
\draw[fill=olive] (2*\Tikzn, \Tikzn*\TikzHeightAbove) circle (\TikzPunctureActualSize) node[olive,right] {$\tilt$};

\draw [CyanTilBetaStyle, dashed, ->] (-\Tikzn *1.5, \Tikzn * \TikzHeightBelow) -- (0.5*\Tikzn, \Tikzn*\TikzHeightBelow) node[pos=0.2, below] {$\tilbe$};

\end{tikzpicture}
}

\newcommand\MacroTwoFoldCoverCrossingRhoOneRhoF{
\begin{scope}[thick, gray]
	\coordinate (v) at (0,0);
	\coordinate (P) at (0,1*\Tikzn);
	\coordinate (y) at (0,2*\Tikzn);
	\draw[thick] (0,0) -- (y);
	\draw[black] (0.0*\Tikzn,0.6*\Tikzn) node {$\rho_1^1$};
	\draw[black] (-0.0*\Tikzn,1.6*\Tikzn) node {$\rho_1^2$};
	
	\draw  (-0.7*\Tikzn,2*\Tikzn)  --
	(0.7*\Tikzn,0)
	node[black,pos=0.25] {$\rho_f^2$} node[black,pos=0.8] {$\rho_f^1$};
	\newcommand\TikzLowerPointerLeft{(-\Tikzn*0.6,\Tikzn*0.5)}
	\newcommand\TikzLowerPointerRight{(\Tikzn*0.6,\Tikzn*0.5)}
	\newcommand\TikzUpperPointerLeft{(-\Tikzn*0.6,1.5*\Tikzn)}
	\newcommand\TikzUpperPointerRight{(\Tikzn*0.6,1.5*\Tikzn)}
	
	\newcommand\TikzLowerPointA{(-0.8*\Tikzn,0.3*\Tikzn)} \newcommand\TikzUpperPointA{(-0.8*\Tikzn,1.8*\Tikzn)}
	\newcommand\TikzLowerPointB{(1*\Tikzn,0.3*\Tikzn)} \newcommand\TikzUpperPointB{(0.5*\Tikzn,1.8*\Tikzn)}
	\draw (v) .. controls \TikzLowerPointA and \TikzUpperPointA .. (-0.7*\Tikzn,2*\Tikzn) node[black,pos=0.6] {$\theta_{h}^1$};
	\draw (0.7*\Tikzn,0) .. controls \TikzLowerPointB and \TikzUpperPointB .. (y) node[black,pos=0.6] {$\theta_{h}^2$};
\end{scope}
	
	\draw [fill=violet, violet] (P) circle (\TikzPunctureActualSize);
	\foreach \pos in {(v),(y),(0.7*\Tikzn,0), (-0.7*\Tikzn,2*\Tikzn)}
	{\draw [fill=olive, olive] \pos circle (\TikzPointSize*\Tikzn );}
}

\newcommand\TikzTwoFoldCoverCrossingRhoOneRhoF[1]{
\begin{tikzpicture}[scale={#1}]
\newcommand\Tikzn{2}
\MacroTwoFoldCoverCrossingRhoOneRhoF
\newcommand\TikzWidth{2}
\draw[gray, thick,rounded corners=5pt] (v) -- (-\TikzWidth*\Tikzn,0) -- (-\TikzWidth*\Tikzn,2*\Tikzn) -- (y)
(-1.8*\Tikzn,0) -- (-1*\Tikzn,2*\Tikzn) node [black,pos=0.7, left] {$\theta_1^1$};
\draw[gray,thick, rounded corners=5pt] (v) -- (\TikzWidth*\Tikzn,0) -- (\TikzWidth*\Tikzn,2*\Tikzn) -- (y)
(1.8*\Tikzn,2*\Tikzn) -- (1*\Tikzn,0) node [black,pos=0.25, left] {$\theta_1^2$};

\newcommand\TikzHeightBelow{0.4}
\newcommand\TikzHeightAbove{0.8}
\draw [RedTilSigmaStyle, dashed, ->] (-\Tikzn *2, \Tikzn * \TikzHeightAbove) -- (\Tikzn *0.7, \Tikzn * \TikzHeightAbove) -- (2*\Tikzn, \Tikzn*\TikzHeightAbove) node[pos=0.7,above] {$\tilsi$};
\draw[fill=olive] (-\Tikzn *2, \Tikzn * \TikzHeightAbove) circle (\TikzPunctureActualSize) node[olive,left] {$\tils$};
\draw[fill=olive] (2*\Tikzn, \Tikzn*\TikzHeightAbove) circle (\TikzPunctureActualSize) node[olive,right] {$\tilt$};
\draw [CyanTilBetaStyle, dashed, ->] (-\Tikzn *1.5, \Tikzn * \TikzHeightBelow) -- (0.3*\Tikzn, \Tikzn*\TikzHeightBelow) node[pos=0.2, below] {$\tilbe$};

\draw (-1.4,2*\Tikzn) node [white, above] {$y^2$};
\draw (1.4,0*\Tikzn) node [white, below] {$y^1$};
\end{tikzpicture}
}

\tikzstyle{LiftedPathStyle} = [blue, line width=2pt, dashed]

\newcommand\TikzNooseOneRLMinusA[1]{
\begin{tikzpicture}[scale=#1, DoubleTriangleLightStyle]
\newcommand\Tikzn{2}
\MacroLiftedTriangleOneTwoNoArrow
\draw [LiftedPathStyle, PathArrowAtPosStyleTwo, PathArrowAtPosStyleFour, PathArrowAtPosStyleEnd, rounded corners=5pt]
(P) --(v1) node[pos=0.4, below] {$r$}
-- (0.1,0.9*\Tikzn) node[pos=0.6,right] {$\ell$}
(0.1,0.9*\Tikzn)
-- (\Tikzn,0.9*\Tikzn) node[pos=0.3,above=3pt]{$\a$}
-- (\Tikzn,1.1*\Tikzn)
-- (0.6*\Tikzn,1.1*\Tikzn);
\end{tikzpicture}
}

\newcommand\TikzNooseOneRLPlusB[1]{
\begin{tikzpicture}[scale=#1, DoubleTriangleLightStyle]
\newcommand\Tikzn{2}
\MacroLiftedTriangleOneTwoNoArrow
\draw [LiftedPathStyle, PathArrowAtPosStyleTwo, PathArrowAtPosStyleFive, PathArrowAtPosStyleEnd, rounded corners=5pt]
 (P) -- (v2) node[pos=0.3,left] {$\ccr$}
-- (v1) node[pos=0.6,right] {$\ell$}
-- (y) node[pos=0.5,right]{$\b$}
-- (0.3,\Tikzn) node[pos=0.5,above]{$\a$};
\end{tikzpicture}
}

\newcommand\TikzNooseOneRLPlusLMinus[1]{
\begin{tikzpicture}[scale=#1, DoubleTriangleLightStyle]
\newcommand\Tikzn{2}
\MacroLiftedTriangleOneTwoNoArrow
\draw [LiftedPathStyle, PathArrowAtPosStyleOne, PathArrowAtPosStyleSeven, PathArrowAtPosStyleEnd, rounded corners=5pt]
(P) -- (0,0.8*\Tikzn) node[pos=0.5,left] {$\ccr$}
-- (v1)  -- (0.2,\Tikzn)  node[pos=0.6,right] {$\ell$}
 -- (y) node[pos=0.5,above]{$\a$};
\end{tikzpicture}
}

\newcommand\MacroCrossingLRABPlain{
\begin{scope}[thick, gray]
	\coordinate (v) at (0,0);
	\coordinate (P) at (0,0.9*\Tikzn);
	\coordinate (y) at (0,2*\Tikzn);
	\draw[thick] (0,0) -- (P);
	\draw (0.08*\Tikzn,0.4*\Tikzn) node {$r$};
	\draw (v) node[below,olive] {$v$}; \draw (y) node[below,olive] {$y$};
	
	\newcommand\TikzLowerPointerLeft{(-\Tikzn*0.5,\Tikzn*0.5)}
	\newcommand	\TikzLowerPointerRight{(\Tikzn*0.5,\Tikzn*0.5)}
	\newcommand\TikzUpperPointerLeft{(-\Tikzn*0.5,\Tikzn)} \newcommand\TikzUpperPointerRight{(\Tikzn*0.5,\Tikzn)}

	\draw[style = thick, rounded corners = 20pt]	
	(v) .. controls \TikzLowerPointerLeft and \TikzUpperPointerLeft  .. (0,\Tikzn*1.5)
	.. controls \TikzUpperPointerRight and \TikzLowerPointerRight .. (v);
	
	\newcommand\TikzLowerPointA{(-1.3*\Tikzn,0.3*\Tikzn)} \newcommand\TikzUpperPointA{(-\Tikzn,1.8*\Tikzn)}
	\newcommand\TikzLowerPointB{(1.3*\Tikzn,0.3*\Tikzn)} \newcommand\TikzUpperPointB{(\Tikzn,1.8*\Tikzn)}
	\draw (v) .. controls \TikzLowerPointA and \TikzUpperPointA .. (y) ;
	\draw (v) .. controls \TikzLowerPointB and \TikzUpperPointB .. (y) ;
\end{scope}
	
	\draw [fill=violet, violet] (P) circle (\TikzPunctureSize*\Tikzn);
	\foreach \pos in {(v),(y)}
	{\draw [fill=olive, olive] \pos circle (\TikzPointSize*\Tikzn );}
}

\newcommand\TikzRLMinusA[1]{
\begin{tikzpicture}[scale={#1}] 	
\newcommand\Tikzn{2}
	\MacroCrossingLRABPlain

	\newcommand\TikzLowerPointerLeft{(-\Tikzn*0.5,\Tikzn*0.5)}
	\newcommand	\TikzLowerPointerRight{(\Tikzn*0.5,\Tikzn*0.5)}
	\newcommand\TikzUpperPointerLeft{(-\Tikzn*0.5,\Tikzn)} \newcommand\TikzUpperPointerRight{(\Tikzn*0.5,\Tikzn)}
	\newcommand\TikzLowerPointA{(-1.3*\Tikzn,0.3*\Tikzn)} \newcommand\TikzUpperPointA{(-\Tikzn,1.8*\Tikzn)}

\draw [LiftedPathStyle, PathArrowAtPosStyleOne,PathArrowAtPosStyleFour,rounded corners=15pt]
	(P) -- (0.1,0.1)
	(0.1,0.1) .. controls \TikzLowerPointerRight and \TikzUpperPointerRight .. (0,\Tikzn*1.5)
	node[pos=0.7,right=3pt] {$\ell_-$}
	.. controls \TikzUpperPointerLeft and \TikzLowerPointerLeft .. (-0.3,0.2);
	
\draw [LiftedPathStyle, PathArrowAtPosStyleNine,rounded corners=15pt]
	 (-0.3,0.2) .. controls \TikzLowerPointA and \TikzUpperPointA .. (y) node[pos=0.4,left] {$\a$};
\draw [LiftedPathStyle, PathArrowAtPosStyleEnd,rounded corners=15pt]
	(y) --  (-0.4*\Tikzn,2.2*\Tikzn) -- (-1*\Tikzn,1.5*\Tikzn);
	
\draw[LiftedPathStyle] (0.08*\Tikzn,0.4*\Tikzn) node {$r$};	

\end{tikzpicture}
}

\newcommand\TikzRLPlusB[1]{
\begin{tikzpicture}[scale={#1}] 	
\newcommand\Tikzn{2}
	\MacroCrossingLRABPlain

	\newcommand\TikzLowerPointerLeft{(-\Tikzn*0.5,\Tikzn*0.5)}
	\newcommand	\TikzLowerPointerRight{(\Tikzn*0.5,\Tikzn*0.5)}
	\newcommand\TikzUpperPointerLeft{(-\Tikzn*0.5,\Tikzn)} \newcommand\TikzUpperPointerRight{(\Tikzn*0.5,\Tikzn)}
	\newcommand\TikzLowerPointA{(-1.3*\Tikzn,0.3*\Tikzn)} \newcommand\TikzUpperPointA{(-\Tikzn,1.8*\Tikzn)}
	\newcommand\TikzLowerPointB{(1.3*\Tikzn,0.3*\Tikzn)} \newcommand\TikzUpperPointB{(\Tikzn,1.8*\Tikzn)}
	
\draw [LiftedPathStyle, PathArrowAtPosStyleOne, PathArrowAtPosStyleSeven, rounded corners=15pt]
	(P) -- (0,0.3)
	(0,0.3) .. controls \TikzLowerPointerLeft and \TikzUpperPointerLeft  .. (0,\Tikzn*1.5)
	.. controls \TikzUpperPointerRight and \TikzLowerPointerRight .. (0.2,0) node[pos=0.3,right=3pt] {$\ell_+$};

\draw [LiftedPathStyle, PathArrowAtPosStyleNine, rounded corners=15pt]
	 (0.2,0) .. controls \TikzLowerPointB and \TikzUpperPointB .. (y) node[pos=0.6,right] {$\b$};
	
\draw [LiftedPathStyle, PathArrowAtPosStyleNine]
	 (y) .. controls \TikzUpperPointA and \TikzLowerPointA .. (-0.8,0.3) node[pos=0.6,left] {$\a$};	

\draw[LiftedPathStyle] (0.08*\Tikzn,0.4*\Tikzn) node {$r$};
\end{tikzpicture}
}

\newcommand\TikzRLPlusLMinusA[1]{
\begin{tikzpicture}[scale={#1}] 	
\newcommand\Tikzn{2}
	\MacroCrossingLRABPlain

	\newcommand\TikzLowerPointerLeft{(-\Tikzn*0.5,\Tikzn*0.5)}
	\newcommand	\TikzLowerPointerRight{(\Tikzn*0.5,\Tikzn*0.5)}
	\newcommand\TikzUpperPointerLeft{(-\Tikzn*0.5,\Tikzn)} \newcommand\TikzUpperPointerRight{(\Tikzn*0.5,\Tikzn)}
	\newcommand\TikzLowerPointA{(-1.2*\Tikzn,0.3*\Tikzn)} \newcommand\TikzUpperPointA{(-\Tikzn,1.8*\Tikzn)}

\draw [LiftedPathStyle, PathArrowAtPosStyleOne, PathArrowAtPosStyleFour, rounded corners=15pt]
	(P) -- (0,0.3)
	(0,0.3) .. controls \TikzLowerPointerLeft and \TikzUpperPointerLeft  .. (0,\Tikzn*1.4)
	.. controls \TikzUpperPointerRight and \TikzLowerPointerRight .. (0.3,0.2); 

	\renewcommand\TikzLowerPointerLeft{(-\Tikzn*0.7,\Tikzn*0.5)}
	\renewcommand\TikzLowerPointerRight{(\Tikzn*0.7,\Tikzn*0.5)}
	\renewcommand\TikzUpperPointerLeft{(-\Tikzn*0.8,\Tikzn)}
	\renewcommand\TikzUpperPointerRight{(\Tikzn*0.8,\Tikzn)}
	
\draw [LiftedPathStyle,PathArrowAtPosStyleFour,rounded corners=15pt]
	(0.3,0.2) .. controls \TikzLowerPointerRight and \TikzUpperPointerRight .. (0,\Tikzn*1.6)
	node[pos=0.7,right=3pt] {$\ell$}
	.. controls \TikzUpperPointerLeft and \TikzLowerPointerLeft .. (-0.5,0.3);

\draw [LiftedPathStyle, PathArrowAtPosStyleNine,rounded corners=15pt]
	 (-0.5,0.3) .. controls \TikzLowerPointA and \TikzUpperPointA .. (y) node[pos=0.6,left] {$\a$};

\draw[LiftedPathStyle] (0.08*\Tikzn,0.4*\Tikzn) node {$r$};
\end{tikzpicture}
}

\tikzstyle{DoubleTriangleStyle}=[thick]
\tikzstyle{DoubleTriangleLightStyle}=[gray!90]

\newcommand\MacroDoubleTriangle{
\newcommand\Tikzn{2}
\coordinate (P) at (0,\Tikzn);

\draw (0,0) -- (P) node[pos=0.3, left] {$r$};

\draw (-\Tikzn,0) -- (-\Tikzn,\Tikzn) node[pos=0.6,left]{$\b$}; 
\draw (0,0) -- (-\Tikzn,0) node[pos=0.5,below]{$\a$}; 
\draw (\Tikzn,0) -- (\Tikzn,\Tikzn) node[pos=0.6,right]{$\a$}; 
\draw (0,0) -- (2,0) node[pos=0.5,below]{$\b$}; 

\draw [fill=violet] (P) circle (\TikzPunctureSize*\Tikzn);
\foreach \pos in {(0,0),(-\Tikzn,\Tikzn) ,(\Tikzn,\Tikzn), (-\Tikzn,0) ,(\Tikzn,0)}
	{\draw [fill=olive] \pos circle (\TikzPointSize*\Tikzn);}
}

\newcommand\MacroDoubleTriangleWithArrow{
\MacroDoubleTriangle
\draw[NooseArrowAtPosStyleSeven] (0,0) -- (-\Tikzn,\Tikzn) node[pos=0.3,left] {$\ell$}; 
\draw[NooseArrowAtPosStyleFour] (\Tikzn,\Tikzn) -- (0,0) node[pos=0.7,right] {$\ell$}; 
}

\newcommand\MacroDoubleTriangleNoArrow{
\MacroDoubleTriangle
\draw (0,0) -- (-\Tikzn,\Tikzn) node[pos=0.3,left] {$\ell$}; 
\draw (\Tikzn,\Tikzn) -- (0,0) node[pos=0.7,right] {$\ell$}; 
}

\newcommand\TikzDoubleTriangleCounterclockwise[1]{
\begin{tikzpicture}[scale=#1, DoubleTriangleStyle]
\MacroDoubleTriangleWithArrow
\draw (-\Tikzn,\Tikzn) -- (P) node[pos=0.5,above]{$\dot{r}$};
\draw (\Tikzn,\Tikzn) -- (P) node[pos=0.5,above]{$\ddot{r}$};
\draw [RedGammaStyle, ->] (-\Tikzn *0.65, \Tikzn *0.5) -- (\Tikzn *0.8, \Tikzn *0.5);

\begin{scope}[gray]
\draw (-\Tikzn * 0.35,\Tikzn * 0.72) node {$\tilblacktri{k}$};
\draw (\Tikzn * 0.35,\Tikzn *0.72) node {$\tilblacktri{k+1}$};
\end{scope}

\end{tikzpicture}
}

\newcommand\TikzDoubleTriangleClockwise[1]{
\begin{tikzpicture}[scale=#1, DoubleTriangleStyle]
\MacroDoubleTriangleWithArrow

\draw (-\Tikzn,\Tikzn) -- (P) node[pos=0.5,above]{$\ddot{r}$};
\draw (\Tikzn,\Tikzn) -- (P) node[pos=0.5,above]{$\dot{r}$};

\draw [RedGammaStyle, <-] (-\Tikzn *0.8, \Tikzn *0.5) -- (\Tikzn *0.7, \Tikzn *0.5);

\begin{scope}[gray]
\draw (-\Tikzn * 0.35,\Tikzn * 0.72) node {$\tilblacktri{k+1}$};
\draw (\Tikzn * 0.35,\Tikzn *0.72) node {$\tilblacktri{k}$};
\end{scope}

\end{tikzpicture}
}

\newcommand\TikzDoubleTriangleKPlusPlus[1]{
\begin{tikzpicture}[scale=#1, DoubleTriangleLightStyle]
\MacroDoubleTriangleNoArrow
\draw (-\Tikzn,\Tikzn) -- (P) node[pos=0.5,above]{$\ddot{r}$};
\draw (\Tikzn,\Tikzn) -- (P) node[pos=0.5,above]{$\dot{r}$};

\coordinate (pk) at (-0.6*\Tikzn,0.6*\Tikzn);
\coordinate (pk1) at (0.2,0.3*\Tikzn);
\draw [LiftedPathStyle, PathArrowAtPosStyleTwo, PathArrowAtPosStyleEnd, rounded corners=5pt] (pk) -- (-0.1,0) -- ($(P)+(-0.1,0)$) -- (pk1);

\end{tikzpicture}
}
\newcommand\TikzEllMinusRR[1]{ %
\begin{tikzpicture}[scale={#1}] 	
	\MacroChubbySelfFoldedTriangle
	\MacroEllMinusRRLong
\end{tikzpicture}
}

\newcommand\TikzDoubleTriangleKKPlus[1]{
\begin{tikzpicture}[scale=#1, DoubleTriangleLightStyle]
\MacroDoubleTriangleNoArrow
\draw (-\Tikzn,\Tikzn) -- (P) node[pos=0.5,above]{$\ddot{r}$};
\draw (\Tikzn,\Tikzn) -- (P) node[pos=0.5,above]{$\dot{r}$};

\coordinate (pk) at (-0.2+-0.5*\Tikzn,0.4*\Tikzn);
\coordinate (pk1) at (0,0.9*\Tikzn);
\draw [LiftedPathStyle, PathArrowAtPosStyleSix, PathArrowAtPosStyleEnd, rounded corners=5pt] (pk) -- (-\Tikzn+0.2,\Tikzn) -- (0,0) -- (pk1);

\end{tikzpicture}
}

\newcommand\TikzEllPlusEllMinusR[1]{ %
\begin{tikzpicture}[scale={#1}] 	
	\MacroChubbySelfFoldedTriangle
	\MacroEllPlusEllMinusR
\end{tikzpicture}
}

\newcommand\TikzDoubleTriangleKThirdEdgePlus[1]{
\begin{tikzpicture}[scale=#1, DoubleTriangleLightStyle]
\MacroDoubleTriangleNoArrow
\draw (-\Tikzn,\Tikzn) -- (P) node[pos=0.5,above]{$\ddot{r}$};
\draw (\Tikzn,\Tikzn) -- (P) node[pos=0.5,above]{$\dot{r}$};

\coordinate (pk) at (-0.3*\Tikzn,0.3*\Tikzn);
\coordinate (pk1) at (0,0.1*\Tikzn);
\draw [LiftedPathStyle, PathArrowAtPosStyleTwo, PathArrowAtPosStyleSix, PathArrowAtPosStyleEnd, rounded corners=5pt] (pk) -- (-\Tikzn,\Tikzn) -- (P) -- (pk1);

\end{tikzpicture}
}

\newcommand\TikzDoubleTriangleRRL[1]{
\begin{tikzpicture}[scale=#1, DoubleTriangleLightStyle]
\MacroDoubleTriangleNoArrow
\draw (-\Tikzn,\Tikzn) -- (P) node[pos=0.5,above]{$\ddot{r}$};
\draw (\Tikzn,\Tikzn) -- (P) node[pos=0.5,above]{$\dot{r}$};

\coordinate (pk2) at (0.6*\Tikzn,0.6*\Tikzn);
\draw [LiftedPathStyle, PathArrowAtPosStyleSix, PathArrowAtPosStyleEnd, rounded corners=5pt]  (-0.2,0.3) -- ($(P)+(-0.1,-0.2)$) -- (0.2,0.1) -- (pk2);
\end{tikzpicture}
}

\newcommand\TikzDoubleTriangleRLL[1]{
\begin{tikzpicture}[scale=#1, DoubleTriangleLightStyle]
\MacroDoubleTriangleNoArrow
\draw (-\Tikzn,\Tikzn) -- (P) node[pos=0.5,above]{$\ddot{r}$};
\draw (\Tikzn,\Tikzn) -- (P) node[pos=0.5,above]{$\dot{r}$};

\draw [LiftedPathStyle, PathArrowAtPosStyleTwo, PathArrowAtPosStyleEnd, rounded corners=15pt]
(0,0.9*\Tikzn) --  (0,0.2) -- (\Tikzn-0.2,\Tikzn) -- (0.7,0.5) ;
\end{tikzpicture}
}

\newcommand\TikzRREllMinus[1]{ %
\begin{tikzpicture}[scale={#1}] 	
	\MacroChubbySelfFoldedTriangle
	\MacroRREllMinus
\end{tikzpicture}
}

\newcommand\TikzDoubleTriangleRddotRL[1]{
\begin{tikzpicture}[scale=#1, DoubleTriangleLightStyle]
\MacroDoubleTriangleNoArrow
\draw (-\Tikzn,\Tikzn) -- (P) node[pos=0.5,above]{$\ddot{r}$};
\draw (\Tikzn,\Tikzn) -- (P) node[pos=0.5,above]{$\dot{r}$};

\draw [LiftedPathStyle, PathArrowAtPosStyleTwo, PathArrowAtPosStyleEnd, rounded corners=5pt] (0,0) -- (0,\Tikzn)  -- (\Tikzn,\Tikzn) -- (0.3,0.3) ;

\end{tikzpicture}
}

\newcommand\TikzRREllPlusnonBT[1]{ %
\begin{tikzpicture}[scale={#1}] 	
	\MacroChubbySelfFoldedTriangle
	\MacroRREllPlusnonBT
\end{tikzpicture}
}
\newcommand\TikzRLMinLPlus[1]{ %
\begin{tikzpicture}[scale={#1}] 	
	\MacroChubbySelfFoldedTriangle
	\MacroRLMinLPlus
\end{tikzpicture}
}

\newcommand\WheelLocalTriangulation[1]{ 
\begin{tikzpicture}[scale=#1]
	\path (0:0) coordinate (origin);
	\path (1*50:2) coordinate (1);
	\path (2*50:2) coordinate (2);
	\path (3*50:2) coordinate (3);
	\path (4*50:2) coordinate (4);
	\path (5*50:2) coordinate (5);
	\path (6*50:2) coordinate (6);
	\path (7*50:2) coordinate (7);
	
	\draw [fill=gray!20, style=dashed] (origin) circle (2);
	
	\begin{scope}[NonOrientedTriangulationStyle]
	\draw  [fill=white] (1) -- (2) -- (3) -- (4) -- (5) -- (6) -- (7) -- cycle;
	\draw [postaction={decorate}] (1) -- (origin);
	\draw [postaction={decorate}] (2) -- (origin);
	\draw [postaction={decorate}] (3) -- (origin);
	\draw [postaction={decorate}] (4) -- (origin);
	\draw [postaction={decorate}] (5) -- (origin);
	\draw [postaction={decorate}] (6) -- (origin);
	\draw [postaction={decorate}] (7) -- (origin);
	\draw [fill=violet] (origin) circle (1.2*\TikzPunctureActualSize);
	\draw (origin) node[above=4pt, right=2pt, violet] {{$\puncture$}};
	\node at (0:3) {}; \node at (180:3) {};
	\end{scope}
\end{tikzpicture}
  }

\newcommand \SelfFoldedLocalTriangulation[1] { 
	\begin{tikzpicture}[scale={#1}]
	\node at (-4,1) {};  \node at (4,1) {};
	\coordinate (v1L) at (-2.5,1);
	\coordinate (v1R) at (2.5,1);
	\coordinate (v2L) at (-2.5,3.5);
	\coordinate (v2R) at (2.5,3);

	\begin{scope}[NonOrientedTriangulationStyle]
		\draw [style = dashed, fill=gray!20]  
		(0,0) .. controls (-3,2) and (-2,3.5) .. (0,4)  -- (v2L) -- (v1L) -- (0,0);
		
		\draw [style = dashed, fill=gray!20]  
		(0,0) .. controls (3,2) and (2,3.5) .. (0,4)  -- (v2R) -- (v1R) -- (0,0);
	
		\draw[style = thick, postaction={decorate}, fill=white] (0,0) -- (0,2.2);
		\draw (0.2,1) node {$r$};
		\draw[style = thick, postaction={decorate}] (0,0) .. controls (-3,4) and (3,4) .. (0,0)
		node[pos = 0.7, above]{$\ell$};
		\draw[style = thick, postaction={decorate}] (0,0) .. controls (-3,2) and (-2,3.5) .. (0,4);
		\draw (-2.1,2.7) node {$\a$};
		\draw[style = thick, postaction={decorate}] (0,0) .. controls (3,2) and (2,3.5) .. (0,4);
		\draw (2.05,2) node {$\b$};
	\end{scope}

	\coordinate (Ps) at (0,2.2);  \coordinate (xs) at (0,0); \coordinate (ys) at (0,4);
	\foreach \pos in {(0,0),(0,4)}
		{\draw[fill=olive] \pos circle (1.2*\TikzPointActualSize);}
	\draw [fill=violet] (Ps) circle (1.2*\TikzPunctureActualSize);  \draw[violet] (0,2.2) node[above] {{$\puncture$}}; 
		
	\end{tikzpicture}
}

 \newcommand\MacroTtwoTpathChubby{
	 \begin{scope}[thick, gray!60]
	 \draw[rounded corners = 20pt] (-4,4)  -- (0,7) node[pos=0.8,above,black]{$1$} -- (4,4);
	
	 \draw[rounded corners = 15pt] (0,0) --
	 (-2.5,4) 
	 --  (0,6) node[pos=0.8,right,black]{$2$}  -- (4,4);
	
	 \draw[thick] (0,0) -- (0,2.9) node[pos=0.8, right=2pt,black]{$r$};
	
	\draw[style = thick, rounded corners = 20pt]
	(0,0) .. controls (-1.5,2) and (-2,3) .. (0,5) ..
	controls (2,3)  and (2,2) .. (0,0) node[pos=0.4, above=2pt,black] {$\ell$};
	
	\draw
	(0,0) -- 
	(4,4) --
	(0,8) --
	node[pos=0.5, left] {}
	 (-4,4) --
	(0,0);	
	\end{scope}
	\draw [fill=violet] (0,2.9) circle (1ex);
 }

\newcommand\MacroCounterClockRadius[4]{
\draw[TPathStyle, 
rounded corners = 5pt,
	decoration={ markings, mark=at position .3 with {\arrow[line width=1pt]{to};},
	mark=at position .8 with {\arrow[line width=1pt]{to};}}, postaction={decorate}]
	({#1},{#2}) .. controls (0.2,3) and (0.4,3) .. (0,3.5) 
	 .. controls (-0.4,3) and (-0.2,3) .. ({#3},{#4});
}

\newcommand\MacroClockLoop[4]{
	\draw[TPathStyle, 
	rounded corners = 20pt
	,decoration={ markings, mark=at position .73 with {\arrow[blue, line width=1pt]{to};}},
	postaction={decorate}]
	({#1},{#2}) .. controls (-1.5,2) and (-2,3) .. (0,5) ..
	controls (2,3)  and (2,2) .. ({#3},{#4}); 
}

\newcommand\MacroClockRadius[4]{
\draw[TPathStyle, 
rounded corners = 5pt,
	decoration={ markings, mark=at position .3 with {\arrow[line width=1pt]{to};},
	mark=at position .8 with {\arrow[line width=1pt]{to};}}, postaction={decorate}]
	({#1},{#2}) .. controls  (-0.2,3) and (-0.4,3) .. (0,3.5)
	.. controls (0.4,3) and (0.2,3)  .. ({#3},{#4});
}

\newcommand\MacroCounterClockLoop[4]{
	\draw[TPathStyle, 
	rounded corners = 20pt
	,decoration={ markings, mark=at position .73 with {\arrow[blue, line width=1pt]{to};}},
	postaction={decorate}]
	({#1},{#2}) .. controls (2,2) and (2,3) .. (0,5) ..
	controls  (-2,3) and (-1.5,2)  ..   ({#3},{#4}); 
}

\newcommand\TikzTpathSix[1]{
	\begin{tikzpicture} [scale = {#1}]
	\begin{scope}[rotate=45]
	\MacroTtwoTpathChubby
	\begin{scope}[TPathStyle]
	\draw[rounded corners = 20pt]  (0,8) -- (3,5) -- (0,7) -- (-4,4);
	\draw (-4,4) -- (0.2*0,-0.1*0);
	
	\MacroCounterClockRadius{0.2*0}{-0.1*0}{-0.5}{1}	
	\MacroClockLoop{-0.5}{1}{1}{1}

	\draw[->] (1,1) -- (3.8,3.8);
	\end{scope}
	\end{scope}
	\node at (-4,6.3) {$b_4$};
	\node at (-3,-0.5) {$b_2$};
	\node at (0.5,3) {$b_3$};
	\end{tikzpicture}
}

\newcommand\TikzTpathSixBad[1]{
	\begin{tikzpicture} [scale = {#1}]
	\begin{scope}[rotate=45]
	\MacroTtwoTpathChubby
	\begin{scope}[TPathStyle]
	\draw[rounded corners = 20pt]  (0,8) -- (3,5) -- (0,7) -- (-4,4) ;
	\draw (-4,4) -- (-0.2,0);
	
	\MacroClockRadius{-0.2}{0}{0.5}{1.5}	
	\MacroCounterClockLoop{0.5}{1.5}{0.3}{0.3}

	\draw[->] (0.3,0.3) -- (3.8,3.8);
	\end{scope}
	\end{scope}
	\node at (-4,6.3) {$b_4$};
	\node at (-3,-0.5) {$b_2$};
	\node at (0.5,3) {$b_3$};
	\node[white] at (-8,3){}; \node[white] at (2,3){};
	\end{tikzpicture}
}

\newcommand\TikzTpathSeven[1]{
	\begin{tikzpicture} [scale = {#1}]
	\begin{scope}[rotate=45]
	\MacroTtwoTpathChubby
	\begin{scope}[TPathStyle]
	\draw[rounded corners = 20pt]  (0,8) -- (3,5) -- (0,7) -- (-4,4);
	\draw  (-4,4) -- (-0.6,0.6);

	\MacroClockLoop{-0.6}{0.6}{0.6}{1.3}
	\MacroCounterClockRadius{0.6}{1.3}{-0.2*0}{0}
	\draw[->] (0,0) -- (3.8,3.8);
	\end{scope}
	\end{scope}
	\node at (-4,6.3) {$b_4$};
	\node at (-3,-0.5) {$b_2$};
	\node at (0.5,3) {$b_3$};
	\end{tikzpicture}
}

\newcommand\MacroTikzTpathEight[1]{
	\begin{scope} [scale = {#1}]
	\begin{scope}[rotate=45]
	\MacroTtwoTpathChubby
	\begin{scope}[TPathStyle]
	\draw[rounded corners = 20pt]  (0,8) -- (-3,5) -- (0,7) -- (3,5) --
	(0,6) -- (-3,4) -- (0.1*0,-0.1*0);
	
	\MacroCounterClockRadius{0.1*0}{-0.1*0}{-0.5}{1.5}
	\MacroClockLoop{-0.5}{1.5}{1}{1}
	\draw[->] (1,1) -- (3.8,3.8);
	\end{scope}
	\end{scope}
	\end{scope}
}

\newcommand\TikzTpathEight[1]{
	\begin{tikzpicture} [scale = {#1}]
	\MacroTikzTpathEight{1}
	\node at (-6.2,3) {$b_1$};
	\node at (0.5,3) {$b_3$};
	\node at (-4,6.3) {};
	\node at (-3,-0.5) {};
	\end{tikzpicture}
}

\newcommand\TikzTpathNine[1]{
	\begin{tikzpicture} [scale = {#1}]
	\begin{scope}[rotate=45]
	\MacroTtwoTpathChubby
	\begin{scope}[TPathStyle]
	\draw[rounded corners = 20pt]  (0,8) -- (-3,5) -- (0,7) -- (3,5) --
	(0,6) -- (-3,4) -- (-1,1);
	\MacroClockLoop{-1}{1}{0.3}{1}
	\MacroCounterClockRadius{0.3}{1}{0}{0}
	\draw[->] (0,0) -- (3.8,3.8);
	\end{scope}
	\end{scope}
	\node at (-6.2,3) {$b_1$};
	\node at (0.5,3) {$b_3$};
	\node at (-4,6.3) {};
	\node at (-3,-0.5) {};
	\end{tikzpicture}
}

\newcommand\MacroTtwoLiftedPolygon{ 
	\begin{scope}[gray!75,style=thick]
	\tilTtwoIdealTriangulationCoordinates
	\draw 
	 (x1) -- (-4,4) node[pos = 0.6, left=2pt]{$\ell$};
	\draw 
	(-4,4) -- (0,4); 
	\draw (-2,4.4) node {$\dot{r}$};
	\draw 
	(0,0) -- (0,4) node[pos=0.6, above]
	{$r$}; 
	\draw 
	(4,4) -- (0,0) node[pos = 0.5, right=2pt]{$\ell$};
	\draw 
	(4,4) -- (0,4) ;
	\draw (2,4.4) node {$\ddot{r}$};

	\draw (x2) -- (t1) node[pos=0.5, above] {$b_3$};
	\draw (t1) -- (x1) node[pos=0.5, left] {$2$};
	\draw (x1) -- (t2) node[pos=0.6, left] {$b_3$}; 
	\draw (t2) -- (x3) node[pos=0.5, below] {$2$}; 
	
	\draw (t1) -- (s) node[pos=0.5, above] {$b_4$}; 
	\draw (s) -- (z); 
	\draw (-2,-4.4) node {$b_1$};
	\draw (z) -- (x1) node[pos=0.5, below] {$b_2$}; 
	\draw (t1) -- (z) node[pos=0.5, left] {$1$};

	\draw [fill=violet] (P) circle (1ex);
	\end{scope}
}

\newcommand\rdPolygonTpath{10}
\newcommand \TikzPolygonTpathSix[1]{
\begin{tikzpicture}[scale={#1}, rotate=-90]
	\MacroTtwoLiftedPolygon
	\begin{scope}[TPathStyle]
	\draw [->, rounded corners = \rdPolygonTpath pt] (-4,-4) -- (-4,0) -- (0,-4) -- (0,0) [rounded corners = 0] --  (0,4)  [rounded corners= \rdPolygonTpath pt] --  (4,4) -- (0,0) -- (4,0);
	\end{scope}
	
	\end{tikzpicture}
}
\newcommand \TikzPolygonTpathSeven[1]{
\begin{tikzpicture}[scale={#1}, rotate=-90]
	\MacroTtwoLiftedPolygon
	\begin{scope}[TPathStyle]
	\draw [->, rounded corners = \rdPolygonTpath pt] (-4,-4) -- (-4,0) -- (0,-4) -- (0,0) -- (-4,4)  [rounded corners=0] -- (0,4)  [rounded corners=\rdPolygonTpath pt] -- (0,0) --(4,0);
	\end{scope}
	
	\end{tikzpicture}
}

\newcommand \MacroTikzPolygonTpathEight[1]{
\begin{scope}[scale={#1}, rotate=-90]
	\MacroTtwoLiftedPolygon
	\begin{scope}[TPathStyle]
	\draw [TPathStyle,->, rounded corners = \rdPolygonTpath pt] (-4,-4) -- (0,-4) -- (-4,0) -- (0,0)  [rounded corners=0] -- (0,4)  [rounded corners=\rdPolygonTpath pt] -- (4,4) -- (0,0) -- (4,0);
	\end{scope}
	\end{scope}
}
\newcommand \TikzPolygonTpathEight[1]{
\begin{tikzpicture}[scale={#1}]
\MacroTikzPolygonTpathEight{1}
\end{tikzpicture}
}

\newcommand \TikzPolygonTpathNine[1]{
\begin{tikzpicture}[scale={#1}, rotate=-90]
	\MacroTtwoLiftedPolygon
	\begin{scope}[TPathStyle]
	\draw [->, rounded corners = \rdPolygonTpath pt] (-4,-4) -- (0,-4) -- (-4,0) -- (0,0) -- (-4,4)  [rounded corners=0]-- (0,4)  [rounded corners=\rdPolygonTpath pt] -- (0,0) -- (4,0);
	\end{scope}
	\end{tikzpicture}
}

\newcommand\MacroChubbySelfFoldedTriangleTau{
	\draw[thick] (0,0) -- (0,2.9) node[pos=0.4, left=4pt]{$\tau$};
	\draw[style = thick, rounded corners = 20pt]
	(0,0) .. controls (-2,1) and (-1.5,3) .. (0,4.2) ..
	controls (2,3)  and (2.5,1) .. (0,0);
	\draw [fill=violet] (0,2.9) circle (1ex);
	\draw [olive] (0,0) node[below] {$v$};
}

\newcommand \TikzBacktrackNonbacktrackBacktrack[1] {
\begin{tikzpicture}[scale={#1} ]
	\MacroChubbySelfFoldedTriangleTau
	\draw[blue,style = dashed, line width = 1pt, rounded corners = 5pt,
	decoration={ markings, mark=at position .3 with {\arrow[line width=1pt]{to};},
	mark=at position .8 with {\arrow[line width=1pt]{to};}}, postaction={decorate}]
 	(-0.1,0) .. controls (-0.2,2) and (-0.3,2.2) .. (0,2.5) node[pos = 0.4, left]{$\om_j$}
	node[pos=0.8,style=thick,fill=blue,blue,circle,minimum width=4pt,inner sep=0pt,draw] {}
	.. controls (0.2,2) and (0.1,1) .. (0.1,0) node[pos = 0.5, right]{$\om_{j+1}$};
	\draw[blue] (-0.3,2.5) node {$\puncture'$};
	\draw [fill=olive] (0,0) circle (.6ex);
	\draw [violet] (-0.2,3.1) node {$\puncture$};
\end{tikzpicture}
}

\newcommand\MacroQuasiBacktrackBackground{
\newcommand\Tikzn{2}
\newcommand\Tikzwidth{1.2*\Tikzn}
\coordinate (P) at (0,2.9);
\coordinate (v) at (0,0);
\coordinate (left_top) at (-\Tikzwidth,4);
\coordinate (right_top) at (\Tikzwidth,4);
\coordinate (left_bottom) at (-\Tikzwidth,0);
\coordinate (right_bottom) at (\Tikzwidth,0);
\coordinate (left_mid) at (-0.8*\Tikzwidth,0.3*\Tikzn);
\coordinate (right_mid) at (0.8*\Tikzwidth,0.3*\Tikzn);
\draw [gray!40,fill=gray!40, rounded corners=2pt] (P) -- (left_top) -- (right_top) -- cycle;
\draw [gray!40,fill=gray!40, rounded corners=10pt] (left_bottom) -- (left_top) -- (left_mid) -- (v);
\draw [gray!40,fill=gray!40, rounded corners=10pt] (right_bottom) -- (right_top) -- (right_mid) -- (v) ;
\draw[dashed, rounded corners =10pt] (left_top) -- (left_mid) -- (v);
\draw[dashed, rounded corners =10pt] (right_top) -- (right_mid) -- (v);

\draw [rounded corners=2pt] (P) -- (left_top) ;
\draw[dashed, rounded corners =2pt] (P) -- (right_top) -- (left_top);

	\draw[thick] (0,0) -- (0,2.9);
\draw[rounded corners=2pt] (v) -- (-\Tikzwidth,0) -- (-\Tikzwidth,4) -- (\Tikzwidth,4) -- (\Tikzwidth,0) -- (v);

\draw [fill=olive] (v) circle (.6ex); \draw [fill=olive] (left_top) circle (.6ex);
	\draw [fill=violet] (P) circle (1ex);
	\draw [olive] (v) node[below] {$v$}; 
}

\newcommand\TikzBacktrackOddFirst[1]{
\begin{tikzpicture}[scale={#1},NonOrientedTriangulationStyle]	
	\MacroQuasiBacktrackBackground
	
	\begin{scope}
	\tikzset{anchor=west}
  	\node[rotate=90,red] at (0.2,1.5) {$k+1$};
	\end{scope}

	\draw[blue,style = dashed, line width = 1pt, rounded corners = 1pt,
	decoration={ markings, mark=at position .4 with {\arrow[line width=1pt]{to reversed};}},
	decoration={ markings, mark=at position .7 with {\arrow[line width=1pt]{to reversed};}},
	postaction={decorate}]
	(v) .. controls (-0.1,2) and (-0.3,2.5) .. (-0.5,2.8) node[pos = 0.2, left]{$\om_{2k+1}$}
	node[pos=1,style=thick,fill=blue,blue,circle,minimum width=4pt,inner sep=0pt,draw] {}
	.. controls (-0.7,2.3) and (-0.3,2) .. (-0.1,0) node[pos = 0.8, right=2pt]{$\om_{2k+2}$};
	\draw[blue] (-0.8,2.9) node {$\puncture'$};

	\draw[gray] (-0.5*\Tikzn,1*\Tikzn) node {$\blacktri{k}$};
	\draw[gray] (0.5*\Tikzn,1*\Tikzn) node {$\blacktri{k+1}$};
\end{tikzpicture}
}

\newcommand\TikzBacktrackEvenFirst[1]{
\begin{tikzpicture}[scale={#1},NonOrientedTriangulationStyle]
\MacroQuasiBacktrackBackground

	\draw (-0.15,2.3) node[red] {$k$};

	\draw[blue,style = dashed, line width = 1pt, rounded corners = 1pt,
	decoration={ markings, mark=at position .4 with {\arrow[line width=1pt]{to};}},
	decoration={ markings, mark=at position .7 with {\arrow[line width=1pt]{to};}},
	postaction={decorate}]
	(v) .. controls (0.1,2) and (0.3,2.5) .. (0.5,2.7) node[pos = 0.3, left]{$\om_{2k}$}
	node[pos=1,style=thick,fill=blue,blue,circle,minimum width=4pt,inner sep=0pt,draw] {}
	.. controls (0.7,2.3) and (0.5,2) .. (0.1,0) node[pos = 0.8, right]{$\om_{2k+1}$};
	\draw[blue] (0.8,3) node {$\puncture'$};

	\draw[gray] (-0.5*\Tikzn,1*\Tikzn) node {$\blacktri{k-1}$};
	\draw[gray] (0.5*\Tikzn,1*\Tikzn) node {$\blacktri{k}$};
\end{tikzpicture}
}

\newcommand\TikzQuasiBacktrackOddFirst[1]{
\begin{tikzpicture}[scale={#1},NonOrientedTriangulationStyle]
	\begin{scope}
	\tikzset{anchor=west}
  	\node[rotate=90,red] at (0.2,1.5) {$k+1$};
	\end{scope}

	\MacroQuasiBacktrackBackground
	
        \draw[blue,style = dashed, line width = 1pt, rounded corners = 5pt,
	decoration={ markings, mark=at position .6 with {\arrow[line width=1pt]{to};}},
	postaction={decorate}]
	(-0.3,2.7) .. controls (-0.2,2.5) and (-0.2,2) .. (-0.1,0) node[pos = 0.7, left]{$\om_{2k+1}$}
	node[pos=0,style=thick,fill=blue,blue,circle,minimum width=4pt,inner sep=0pt,draw] {};
	\draw[blue] (-0.6,2.8) node {$\puncture'$};

	\draw[blue,style = dashed, line width = 1pt, rounded corners = 5pt,
	decoration={ markings, mark=at position .4 with {\arrow[line width=1pt]{to};}},
	postaction={decorate}]
	(0.1,0) .. controls (0.1,2) and (0.7,2.5) .. (0.8,2.7) node[pos = 0.3, right]{$\om_{2k+2}$}
	node[pos=1,style=thick,fill=blue,blue,circle,minimum width=4pt,inner sep=0pt,draw] {};
	\draw[blue] (1.2,3) node {$\puncture''$};

	\draw [fill=olive] (v) circle (.6ex); \draw [fill=olive] (left_top) circle (.6ex);
	\draw [fill=violet] (P) circle (1ex);
	\draw [olive] (v) node[below] {$v$}; 
	
\end{tikzpicture}
}

\newcommand\TikzQuasiBacktrackEvenFirst[1]{
\begin{tikzpicture}[scale={#1},NonOrientedTriangulationStyle]
	\MacroQuasiBacktrackBackground
	\draw (0.15,2.3) node[red] {$k$};
        \draw[blue,style = dashed, line width = 1pt, rounded corners = 5pt,
	decoration={ markings, mark=at position .6 with {\arrow[line width=1pt]{to};}},
	postaction={decorate}]
	(-0.3,2.7) .. controls (-0.2,2.5) and (-0.1,2.4) .. (-0.1,0) node[pos = 0.7, left]{$\om_{2k}$}
	node[pos=0,style=thick,fill=blue,blue,circle,minimum width=4pt,inner sep=0pt,draw] {};
	\draw[blue] (-0.6,2.8) node {$\puncture'$};

	\draw[blue,style = dashed, line width = 1pt, rounded corners = 5pt,
	decoration={ markings, mark=at position .4 with {\arrow[line width=1pt]{to};}},
	postaction={decorate}]
	(0.1,0) .. controls (0.1,2) and (0.6,2.5) .. (0.6,2.7) node[pos = 0.3, right]{$\om_{2k+1}$}
	node[pos=1,style=thick,fill=blue,blue,circle,minimum width=4pt,inner sep=0pt,draw] {};
	\draw[blue] (1,3.1) node {$\puncture''$};

	\draw[gray] (-0.5*\Tikzn,1*\Tikzn) node {$\blacktri{k-1}$};
	\draw[gray] (0.5*\Tikzn,1*\Tikzn) node {$\blacktri{k}$};
\end{tikzpicture}
}

\newcommand \TikzBacktrackNonbacktrackCounterclock[1] {
\begin{tikzpicture}[scale={#1}]
	\MacroChubbySelfFoldedTriangleTau
	\draw[blue,style = dashed, line width = 1pt, rounded corners = 5pt,
	decoration={ markings, mark=at position .3 with {\arrow[line width=1pt]{to};},
	mark=at position .8 with {\arrow[line width=1pt]{to};}}, postaction={decorate}]
	(0.1,0) .. controls (0.2,3) and (0.4,3) .. (0,3.5) node[pos = 0.4, right]{$\om_{j+1}$}
	node[pos=0.9,style=thick,fill=blue,blue,circle,minimum width=4pt,inner sep=0pt,draw] {}
	 .. controls (-0.4,3) and (-0.2,3) .. (-0.1,0) node[pos = 0.6, left]{$\om_j$};
	 \draw[blue] (0.4,3.6) node {{$\puncture'$}};
	\draw [fill=olive] (0,0) circle (.6ex);
	\draw [violet] (-0.4,3.1);
	\node at (3.5,0) {}; \node at (-3.5,0) {};
\end{tikzpicture}
}

\newcommand \TikzBacktrackNonbacktrackClock[1] {
\begin{tikzpicture}[scale={#1}]
	\MacroChubbySelfFoldedTriangleTau
	\draw[blue,style = dashed, line width = 1pt, rounded corners = 5pt,
	decoration={ markings, mark=at position .3 with {\arrow[line width=1pt]{to};},
	mark=at position .8 with {\arrow[line width=1pt]{to};}}, postaction={decorate}]
	(-0.1,0) .. controls (-0.2,3) and (-0.4,3) .. (0,3.5) node[pos = 0.4, left]{$\om_j$}
	.. controls (0.4,3) and (0.2,3) .. (0.1,0) node[pos = 0.6, right]{$\om_{j+1}$}
	node[pos=0.1,style=thick,fill=blue,blue,circle,minimum width=4pt,inner sep=0pt,draw] {};
	\draw[blue] (0.4,3.6) node {{$\puncture'$}};
	\draw [fill=olive] (0,0) circle (.6ex);
	\draw [violet] (-0.4,3.1);
	\node at (3,0) {}; \node at (-3,0) {};
\end{tikzpicture}
}

\newcommand\TikzGammaK[1]{
\begin{tikzpicture}[scale={#1}] 	
	\MacroChubbySelfFoldedTriangle
			
	\draw [red, line width = 3pt,
	decoration={ markings, mark=at position .6 with {\arrow[line width=3pt]{to};}}, postaction={decorate}] (-1.6,1.3) -- (0,1.3)
	node[pos=0.6, above]{{$\ga_k$}};
	\draw [fill=blue] (-1.6,1.3) circle (.6ex) node[left, blue] {$p_k$};
	\draw [fill=blue] (0,1.3) circle (.6ex) node[right, blue] {$p_{k+1}$}; 	
\end{tikzpicture}
}

\newcommand\MacroEllPlusEllPlus{
	\draw[->,blue,style = dashed, line width = 2pt, rounded corners = 15pt,
	decoration={ markings, mark=at position .3 with {\arrow[line width=2pt]{to};},
	mark=at position .8 with {\arrow[line width=2pt]{to};}}, postaction={decorate}]
	(-1.8,1.2)  -- (-2,2) -- 
	(-1.5,3.5) --
	(0,4.2) .. 	controls (3,2.8)  and (2.5,1) ..
	(0,-0.4)  .. controls (-2,2) and (-1.6,3) ..
	(0,3.8) .. 	controls (2,2.8)  and (1.5,1) ..
	(0.2,0.2) -- (0.2,1.3) ;
}

\newcommand\TikzEllPlusEllPlus[1]{ 
\begin{tikzpicture}[scale={#1}] 	
	\MacroChubbySelfFoldedTriangle
	\MacroEllPlusEllPlus
\end{tikzpicture}
}

\newcommand\TikzGammaKConcatenate[1]{
\begin{tikzpicture}[scale={#1}] 	
	\draw [fill=violet] (0,2.7) circle (0.6ex) node[right=1pt, violet] {$\puncture$};
	 \draw [red, line width = 3pt,
	decoration={ markings, mark=at position .5 with {\arrow[line width=3pt]{to};}}, postaction={decorate}]  (0.2,1.2) -- (-1.8,1.2)
	node[pos=0.3, above]{{$\ga_k^-$}};
	\MacroEllPlusEllPlus
\end{tikzpicture}
}
\newcommand\MacroChubbySelfFoldedTriangle{
    \begin{scope}[gray!90]
	\draw[thick]
	(0,0) -- (0,2.5) node[black,pos=0.7, right]{$r$};
	\draw[fill=violet] (0,2.5) circle (0.6ex);
	\draw[fill=olive] (0,0) circle (0.4ex);
	\draw[style = thick, rounded corners = 20pt]
	(0,0) .. controls (-2.5,1) and (-2,3) .. (0,4) ..
	controls (2,3)  and (2.5,1) .. (0,0) node[black,pos=0.5, right] {$\ell$};
    \end{scope}
}
\newcommand\MacroEllMinusRRLong{
	\draw[<-,blue,style = dashed, line width = 2pt, rounded corners = 15pt,
	decoration={ markings, mark=at position .5 with {\arrow[line width=2pt]{to reversed};}}, postaction={decorate}]
	(0.2,0.8) --  (-0.1,2.4) --
	(-0.4,0.3)  .. controls (-1.8,1.3) and (-2,2.8) .. (0,3.9) ..  
	controls (2,2.8)  and (1.5,1) .. (1,0.8) ;
}
\newcommand\MacroEllPlusEllMinusR{
		\draw[->,blue,style = dashed, line width = 2pt,
		rounded corners = 15pt,
	decoration={ markings, mark=at position .3 with {\arrow[line width=2pt]{to};},
	mark=at position .85 with {\arrow[line width=2pt]{to};}},
	postaction={decorate}]
	(-1.8,1.5)  -- (-2,2) --
	(-1.5,3.5) -- 	(0,4.1) .. 	controls (2.5,2.8)  and (2,1.5) .. (0.6,0.4)
	 ..  controls  (1.5,1) and  (2,2.8)  ..
	(0,3.8) .. controls (-1.6,3) and (-2.3,2)  ..
	(0,0) -- (0,2.2) ;
}
\newcommand\MacroRLMinLPlus{
		\draw[->,blue,style = dashed, line width = 2pt, rounded corners = 15pt,
	decoration={ markings, mark=at position .06 with {\arrow[line width=2pt]{to};},
	mark=at position .15 with {\arrow[line width=2pt]{to};}},
	postaction={decorate}]
	(-0,2) -- (-0,0) .. controls (1.5,1) and (2,2.8)  .. (0,3.8)
	 ..
	 controls (-1.6,3) and (-2.3,2)  ..
	  (-1,0.6)  ..
	controls (-2.2,2)  and  (-1.7,3.3)  .. (0,4.1) .. controls (1.7,3.3) and (2.2,2) .. (1.6,1.3);
}

\newcommand\MacroEllPlusRR{
	\draw[->,blue,style = dashed, line width = 2pt, rounded corners = 15pt,
	decoration={ markings, mark=at position .2 with {\arrow[line width=2pt]{to};}}, postaction={decorate}]
	(-1.4,1.2)  .. controls (-1.7,2) and (-1.6,3) .. (0,3.8) ..  
	controls (2,2.8)  and (1.5,1) .. (0.3,0.2) 
	-- (0.1,2.5) -- (-0.3,1);
}
\newcommand\MacroRREllMinus{
	\draw[<-,blue,style = dashed, line width = 2pt, rounded corners = 15pt,
	decoration={ markings, mark=at position .5 with {\arrow[line width=2pt]{to reversed};}}, postaction={decorate}]
	(-1.4,1.2)  .. controls (-1.7,2) and (-1.6,3) .. (0,3.8) ..  
	controls (2,2.8)  and (1.5,1) .. (0.2,0.2) 
	-- (0.1,2.5) -- (-0.3,1);
}

\newcommand\TikzEllPlusRRConcatenate[1]{
\begin{tikzpicture}[scale={#1}] 	
	\MacroChubbySelfFoldedTriangle
	 \draw [red, line width = 3pt,
	decoration={ markings, mark=at position .5 with {\arrow[line width=3pt]{to};}}, postaction={decorate}]  (-0.3,1.1) -- (-1.4,1.1)
	node[pos=0.5, above]{{$\ga_k^-$}};
	\MacroEllPlusRR
\end{tikzpicture}
}

\newcommand\MacroEllMinusRR{
	\draw[->,blue,style = dashed, line width = 2pt, rounded corners = 10pt,
	decoration={ markings, mark=at position .2 with {\arrow[line width=2pt]{to};}}, postaction={decorate}]
	(-1.6,1)  --  (-1.3,0.4) -- (0.2,0)
	-- (0.2,2.5) -- (-0.3,1);
}
\newcommand\TikzEllMinusRRConcatenate[1]{
\begin{tikzpicture}[scale={#1}] 	
	\MacroChubbySelfFoldedTriangle
	 \draw [red, line width = 3pt,
	decoration={ markings, mark=at position .5 with {\arrow[line width=3pt]{to};}}, postaction={decorate}]  (-0.3,1.1) -- (-1.6,1.1)
	node[pos=0.5, above]{{$\ga_k^-$}};
	\MacroEllMinusRR
\end{tikzpicture}
}

\newcommand\MacroEllPlusRRnonBT{
	\draw[->,blue,style = dashed, line width = 2pt, rounded corners = 15pt,
	decoration={ markings, mark=at position .2 with {\arrow[line width=2pt]{to};}}, postaction={decorate}]	
	(-1.4,1.2)  -- (-1.8,1.8) --  
	(-1.5,3.2) --  
	(0,3.8) .. controls (2.6,2.8)  and (1.5,1) .. 
	(0.4,0.2) .. controls (0.3,1) and (0.5,2.8) .. 
	(-0.4,3.3) .. controls (-0.3, 2) and (-0.2,1.5) .. 
	(-0.3,1); 
}

\newcommand\MacroRREllPlusnonBT{
	\draw[->,blue,style = dashed, line width = 2pt, rounded corners = 15pt,
	decoration={ markings, mark=at position .3 with {\arrow[line width=2pt]{to};}},
	decoration={ markings, mark=at position .6 with {\arrow[line width=2pt]{to};}}, postaction={decorate}]	
	(0.3,0.6) .. controls (0.3,1) and (0.5,2.8) .. 
	(-0.4,3.3) .. controls (-0.3, 2) and (-0.2,1.5) .. 
	(-0.2,0) 
	--
	(-1.4,0.7)  -- (-1.8,1.8) --  
	(-1.5,3.2) --  
	(0,3.8) .. controls (2.6,2.8)  and (1.5,1) .. 
	(0.8,0.4);
}

\newcommand\TikzEllPlusRRnonBTConcatenate[1]{
\begin{tikzpicture}[scale={#1}] 	
	\MacroChubbySelfFoldedTriangle
	 \draw [red, line width = 3pt,
	decoration={ markings, mark=at position .5 with {\arrow[line width=3pt]{to};}}, postaction={decorate}]  (-0.3,1.1) -- (-1.4,1.1)
	node[pos=0.5, above]{{$\ga_k^-$}};
	\MacroEllPlusRRnonBT
\end{tikzpicture}
}
\newcommand\TikzEllPlusRRnonBT[1]{
\begin{tikzpicture}[scale={#1}] 	
	\MacroChubbySelfFoldedTriangle
	\MacroEllPlusRRnonBT
\end{tikzpicture}
}

\newcommand \MacroTtwoTriangulationCoordinates{
 \coordinate (t) at (0,0);
 \coordinate(s) at (4,0);
  \coordinate(vL) at (0,4);  \coordinate(vR) at (4,4);
   \coordinate(x) at (0,4);  \coordinate(z) at (4,4);
 }

\newcommand \MacroTtwoTriangulationVertices{
 \draw [fill=olive] (0,0) circle (.6ex);
\draw [fill=olive] (0,4) circle (.6ex);
 \draw (-0.3,4.2) node[olive] {$v$};
\draw [fill=olive] (4,0) circle (.6ex);
\draw [fill=olive] (4,4) circle (.6ex);
\draw [fill=violet] (1.2,2.2) circle (.6ex); 
 \draw (1.4,2) node[violet] {$\puncture$};
  }

 \newcommand \MacroTtwoIdealTriangulationWithNumberLabels {
 \MacroTtwoTriangulationCoordinates
\begin{scope}[NonOrientedTriangulationStyle]
\draw[postaction={decorate}]  (s) -- (vR) node[pos=0.5,left] {$b_1$};
\draw[postaction={decorate}]  (vR) -- (vL) node[pos=0.5,below]{$b_2$};
\draw[postaction={decorate}]  (vL) -- (t) node[pos=0.5,right] {$b_3$};
\draw[postaction={decorate}]  (t) --(s) node[pos=0.5,above] {$b_4$};
\draw [postaction={decorate}] (t) ..  controls(3.5,0.4) and (3.5,2) .. (vR) node[left, pos=0.3] {$1$};
\draw [postaction={decorate}] (t) ..  controls (3.7,1) and (3.8,3.9)   .. (vL) node[right, pos=0.5] {$2$};
\draw [NooseArrowAtPosStyleSix] (vL) ..  controls (3.9,1.4)  and (0.7,0.1)  .. (vL) node[pos=0.7, right,above] {$\ell$}; 
\draw[postaction={decorate}] (vL) -- (1.2,2.2) node[left, pos=0.7] {$r$};
\end{scope}
 }
  \newcommand\TikzTtwoIdealTriangulationWithNumberLabels[1]{ 
 	\begin{tikzpicture} [scale = {#1}, rotate=180] 
	\begin{scope}
	\MacroTtwoIdealTriangulationWithNumberLabels
	\draw [<-,decorate, line width=2 pt, dashed, red!60] (0,0) ..  controls (0,5) and (4,4)   .. (4,0)
	node[pos=0.6, left=2pt] {{$\ga$}}; 
	\MacroTtwoTriangulationVertices
	\draw [fill=violet] (1.2,2.2) circle (.6ex); 
	\draw [fill=olive] (t) circle (.6ex); 
	\draw [fill=olive] (vL) circle (.6ex); 
	\draw [fill=olive] (s) circle (.6ex); 
	\draw [fill=olive] (vR) circle (.6ex); 
	\end{scope}
\end{tikzpicture}
}
\newcommand  \TikzTtwoTagged[1]{
 \begin{tikzpicture} [scale =#1]
 \begin{scope}[rotate=180]
\MacroTtwoTriangulationCoordinates
\draw[style = thick]  (0,0) rectangle (4,4);
\draw [style = thick] (0,0) ..  controls(3.5,0.4) and (3.5,2) .. (4,4);
\draw [style = thick] (0,0) ..  controls (3.7,1) and (3.8,3.9)   .. (0,4);
\draw [style = very thick, gray] (0,4) ..  controls (0.7,0.1) and (3.9,1.4)   .. (0,4) node[pos=0.75, left] {$\ell$}; 
\draw[style = thick] (0,4) -- (1.2,2.2);
\MacroTtwoTriangulationVertices
\end{scope}
\end{tikzpicture}\quad
\begin{tikzpicture}[scale=#1] 
\begin{scope}[rotate=180]
\MacroTtwoTriangulationCoordinates
\draw[style = thick]  (0,0) rectangle (4,4);
\draw [style = thick] (0,0) ..  controls(3.5,0.4) and (3.5,2) .. (4,4); 
\draw [style = thick] (0,0) ..  controls (3.7,1) and (3.8,3.9)   .. (0,4); 
\draw[style = thick] (0,4) .. controls (0.8,1) and (1,2) .. (1.2,2.2);  
\draw[style = very thick, gray] (0,4) -- (1.2,2.2) node[pos=0.5, left=1pt] {$\iota(\ell)$};
\path[fill, gray] (1,2.5) -- (1.1,2.8) -- (1.3,2.5) -- cycle; 
\path[fill, gray] (1,2.5) -- (0.9,2.2) -- (0.7,2.5) -- cycle; 
\MacroTtwoTriangulationVertices
\end{scope}
\end{tikzpicture}
 }
\newcommand \TikzTone{ 
\draw[style = thick]  (0,0) rectangle (4,4);
\draw [style = thick] (0,0) ..  controls(3.5,0.4) and (3.5,2) .. node[pos =0.4,right=2pt] {$1$}  (4,4);
\draw [style = thick] (0,0) ..  controls (3.7,1) and (3.8,3.9)   .. node[pos =0.4,left=2pt] {$2$}  (0,4);
\draw[style = thick] (0,4) -- node[pos =0.3,right=1pt] {$3$}  (1.2,2.2);
\draw[style = thick] (1.2,2.2) -- node[pos =0.4,right=1pt] {$4$}  (0,0);
\draw [fill=violet] (1.2,2.2) circle (.7ex); 
\draw [fill=olive] (0,0) circle (.6ex);
\draw [fill=olive] (0,4) circle (.6ex);
\draw [fill=olive] (4,0) circle (.6ex);
\draw [fill=olive] (4,4) circle (.6ex);
\draw (4,2) node [right = 5pt] {{$b_1$}};
\draw  (2,4) node [above = 5pt] {{$b_2$}};
\draw  (0,2) node [left = 5pt] {{$b_3$}};
\draw  (2,0) node [below = 5pt] {{$b_4$}};
 }
\newcommand  \TikzToneWithGamma[1]{
 \begin{tikzpicture} [scale=#1]
\TikzTone 
\draw [red] (4,0) node [right = 5pt] {{$s$}};
\draw [red] (0,0) node [left = 5pt] {{$t$}};
\draw [<-,decorate, line width=3 pt, dashed, red] (0,0) ..  controls (0,5) and (4,4)   .. (4,0)
node[pos=0.7, above=2pt] {{$\ga$}}; 
\end{tikzpicture}
  }
\newcommand \TikzToneBlackTriangles[1]{
\begin{tikzpicture} [scale=#1]
\MacroTikzToneTgk 
\draw[gray] (3.3,0.5) node {{$\blacktri{0}$}};
\draw[gray]  (3.2,3.2) node {{$\blacktri{1}$}};
\draw[gray]  (2,2.5) node {{$\blacktri{2}$}};
\draw[gray]  (0.5,2) node {{$\blacktri{3}$}};
\end{tikzpicture}
}
\newcommand \MacroTikzToneTgk{
\begin{scope}[gray]
\draw [style = thick]  (0,0) rectangle (4,4);
\draw [style = thick] (0,0) ..  controls(3.5,0.4) and (3.5,2) .. node[black,pos =0.4,right=2pt] {$\taui{1}$}  (4,4);
\draw [style = thick] (0,0) ..  controls (3.7,1) and (3.8,3.9)   .. node[black,pos =0.4,left=2pt] {$\taui{2}$}  (0,4);
\draw [style = thick] (0,4) -- node[black, pos =0.4,right=2pt] {$\taui{3}$}  (1.2,2.2);
\draw [style = thick] (1.2,2.2) -- node[black, pos =0.4,right=2pt] {$\tg{2}$}
node [pos=0.6, right, black] {$=\tg{d}$}  (0,0);
\draw [fill=violet] (1.2,2.2) circle (.7ex); 
\draw [fill=olive] (0,0) circle (.6ex);
\draw [fill=olive] (0,4) circle (.6ex);
\draw [fill=olive] (4,0) circle (.6ex);
\draw [fill=olive] (4,4) circle (.6ex);
\draw (4,2) node [right = 2pt, black]  {$\tg{0}$};
\draw (1.5,4) node [above = 2pt, black] {$\tg{1}$};
\draw  (0,2) node [left = 2pt,black]  {$\tg{d+1}$};
\draw  (2,0) node [below = 2pt,black] {$\tg{-1}$};
\end{scope}
}

\newcommand\DiskWithPunctureTauPeripheral[1]{
\begin{tikzpicture}[scale=#1]
\begin{scope}[xscale=-1]
	\path (0:0) coordinate (origin);
	\path (180:2) coordinate (4);
	\path (255:2) coordinate (v1);
	\path (300:2) coordinate (v1');
	\path (7*50:2) coordinate (7);
	
	\draw [fill=violet] (origin) circle (.6ex); \draw (origin) node[above=2pt, violet] {$\puncture$};
	\draw[thick] (origin) circle (2);
	\draw[thick] (v1) -- (7) node[pos=0.4,below]{$\al$}; 
	\draw[thick, dotted] (v1) -- (4) node[pos=0.4,right]{$\lambda$};

	\draw [->, line width=3 pt, dashed, red] (v1) ..  controls (-2.5,2.5) and (1.5,2.5)   .. (v1')
	node[pos=0.6, left=2pt] {{$\si$}}; 
	
	\draw [fill=olive] (v1) circle (.4ex); \draw (v1) node[below=3pt] {$v_1$};
	\draw [fill=olive] (v1') circle (.4ex); \draw (v1') node[below=3pt] {$v_1'$};	
	\draw [fill=olive] (7) circle (.4ex) node[left=3pt] {$v_2$};
	\draw[thick, rounded corners=20pt] (4) -- (270:0.8) -- (7) node[pos=0.2, above] {$\taui{1}$}; 
	\draw [fill=olive] (4) circle (.4ex) node[right=3pt] {$v_0$};
	\draw (290:1.1) node[right=3pt] {${\blacktri{0}}$};
\end{scope}
\end{tikzpicture}
  }

\newcommand\TikzTauPeripheralSelfFolded[1]{
\begin{tikzpicture}[scale=#1]
\begin{scope}[xscale=-1]
	\path (0:0) coordinate (origin);
	\path (180:2) coordinate (4);
	\path (255:2) coordinate (v1);
	\path (300:2) coordinate (v1');
	\path (7*50:2) coordinate (7);
	\path (80:2) coordinate (80);
	
	\draw [fill=violet] (origin) circle (.6ex); \draw (origin) node[above=2pt, violet] {$\puncture$};
	\draw[thick] (origin) circle (2);
	
	\draw[thick, rounded corners=18pt] (v1) --  (320:1.6) -- (30:1.7)
	node[pos=0.5,left]{$\al$} -- (80) ; 
	
	\draw[thick, dotted, rounded corners=18pt] (v1) -- (210:1.8) -- (160:1.8) -- (100:1.8)
	node[pos=0.5,right]{$\lambda$} -- (80);

	\draw [->, line width=3 pt, dashed, red] (v1) ..  controls (-2.5,2.2) and (1.5,2.2)   .. (v1')
	node[pos=0.65, above] {{$\si$}}; 
	
	\draw [fill=olive] (v1) circle (.4ex); \draw (v1) node[below=3pt] {$v_1$};
	\draw [fill=olive] (v1') circle (.4ex); \draw (v1') node[below=3pt] {$v_1'$};	
	
	\draw[thick, rounded corners=15pt] (80) --  (175:0.8) -- (270:0.8) node[pos=0.6, left] {$\taui{1}$} -- (0:1) -- (80) ; 
	
	\draw [fill=olive] (80) circle (.4ex) node[above=3pt] {$v_0=v_2$};

	\draw (290:1.1) node[right=3pt] {${\blacktri{0}}$};
\end{scope}
\end{tikzpicture}
  }

\newcommand\MacroPuncturedDiskMarkedPoints[1]{
	\draw[fill=olive] (0:2) circle ({#1});
	\draw[fill=olive] (20:2) circle ({#1});
	\draw[fill=olive] (50:2) circle ({#1});
	\draw[fill=olive] (90:2) circle ({#1});
	\draw[fill=olive] (120:2) circle ({#1});
	\draw[fill=olive] (150:2) circle ({#1});
	\draw[fill=olive] (160:2) circle ({#1});
	\draw[fill=olive] (210:2) circle ({#1});
	\draw[fill=olive] (240:2) circle ({#1});
	\draw[fill=olive] (270:2) circle ({#1});
	\draw[fill=olive] (300:2) circle ({#1});
	\draw[fill=olive] (330:2) circle ({#1});
}
\newcommand\DiskLambda[1]{
\begin{tikzpicture}[scale=#1]
	\draw[fill=violet] (0:0) circle (\TikzPunctureActualSize);
	\draw (0:0) node[left=4pt,below=3pt,violet] {$\puncture$};	
	\draw[fill,gray!20]  (20:2) -- (30:2) -- (40:2) -- (50:2) -- (70:2) -- (90:2) -- (110:2) -- (130:2) -- (150:2) -- (160:2) --  cycle;
	\draw[thick] (0:0) circle (2);
	\draw[line width=2pt,dashed, cyan] (20:2) -- (160:2) node[pos=0.2,below,cyan] {${\lambda}$};
	\draw (90:1.2) node[cyan] {${\disk{\lambda}}$};
	\MacroPuncturedDiskMarkedPoints{.4ex}
\end{tikzpicture}
}

\newcommand\OncePuncturedDiskLambda[1]{
\begin{tikzpicture}[scale=#1]
	\coordinate (vleft) at (150:2); \coordinate (vright) at (30:2);
	\draw[fill=gray!20] (0:0) circle (2);
	\draw[fill,white]  (vright)  -- (40:2) -- (50:2) -- (70:2) -- (90:2) -- (110:2) -- (130:2) -- (150:2) -- (vleft) --  cycle;
	\draw[black,thick] (0:0) circle (2);
	\draw[line width=2pt,dashed, cyan] (vright) -- (vleft) node[pos=0.4,above,cyan] {${\lambda}$};
	\draw (-90:1) node[cyan] {${\outsidedisk{\lambda}}$};
	\MacroPuncturedDiskMarkedPoints{.4ex}
	\draw[fill=violet] (0:0) circle (\TikzPunctureActualSize) (0:0) node[left=4pt,above,violet] {$\puncture$};	
	
	\node at (0:3) {};\node at (180:3) {};
\end{tikzpicture}
}

\newcommand\MacroPuncturedDisk{
	\draw[fill=violet] (0:0) circle (.6ex);
	\draw (0:0) node[left=4pt,below=3pt,violet] {$\puncture$};
	\draw[thick] (0:0) circle (2);	
}
\newcommand\MacroPuncturedDiskNoP{
	\draw[fill=violet] (0:0) circle (.6ex);
	\draw (0:0);
	\draw[thick] (0:0) circle (2);	
}

\newcommand\ChoiceOfSigmaLambdaKs[1]{
\begin{tikzpicture}[scale=#1]
	\MacroPuncturedDiskNoP
	\draw[line width=2pt,dashed, red] (20:2) -- (160:2) node[pos=0.3, below,red] {${\lambda_r}$};
	\draw[line width=2pt,dashed, cyan] (50:2) -- (150:2);
	\draw[line width=2pt,dashed, cyan, rounded corners] (90:2) -- (100:1.6) -- (150:2);
	\draw[line width=2pt,dashed, red, rounded corners=20pt] (20:2) -- (0:1) -- (330:2) node[pos=0.8, left,red] {${\lambda_1}$};
	\draw[line width=2pt,dashed, red, rounded corners=20pt] (210:2) -- (240:0.7) node[pos=0.4, above,red] {${\lambda_2}$} -- (300:2);
	\draw[line width=2pt,dashed, cyan, rounded corners=20pt] (210:2) -- (240:1) -- (270:2);
	\MacroPuncturedDiskMarkedPoints{.4ex}
	\node at (0:3) {};\node at (180:3) {};
\end{tikzpicture}
}

\newcommand\ChoiceOfSigmaLambdaR[1]{
\begin{tikzpicture}[scale=#1]
	\MacroPuncturedDiskNoP
	\draw[line width=2pt,dashed, red] (20:2) -- (160:2) node[pos=0.45, above,red] {${\si:=\lambda_r}$};
	\draw[line width=2pt,dashed, cyan] (50:2) -- (150:2);
	\draw[line width=2pt,dashed, cyan, rounded corners] (90:2) -- (100:1.6) -- (150:2);
	\draw[line width=2pt, rounded corners=15pt] (120:2) -- (220:1) node[pos=0.9, left] {${\tau}$} -- (320:1) -- (50:2);
	\MacroPuncturedDiskMarkedPoints{.4ex}
	\node at (0:3) {};\node at (180:3) {};
\end{tikzpicture}
}

\newcommand\ChoiceOfSigmaGoodA[1]{
\begin{tikzpicture}[scale=#1]
	\MacroPuncturedDisk
	\draw[line width=2pt,dashed, red] (20:2) -- (160:2) node[pos=0.2, below,red] {${\si}$};
	\draw[line width=2pt,dashed, cyan] (50:2) -- (150:2) node[pos=0.2, below,cyan] {${\be}$};
	\MacroPuncturedDiskMarkedPoints{.4ex}
	\node at (0:3) {}; \node at (180:3) {};
\end{tikzpicture}
}
\newcommand\ChoiceOfSigmaGoodB[1]{
\begin{tikzpicture}[scale=#1]
	\MacroPuncturedDisk
	\draw[line width=2pt,dashed, red] (20:2) -- (160:2) node[pos=0.2, below,red] {${\si}$};
	\draw[line width=2pt,dashed, cyan] (330:2) -- (210:2) node[pos=0.2, below,cyan] {${\be}$};
	\MacroPuncturedDiskMarkedPoints{.4ex}
	\node at (0:3) {}; \node at (180:3) {};
\end{tikzpicture}
}
\newcommand\MacroFourCover {
\begin{scope}[gray]
	\path (0:0) coordinate (origin);
	\path (0:2) coordinate (v0);
	\path (40:2) coordinate (v40);
	\path (70:2) coordinate (v70);
	\path (90:2) coordinate (v90);
	\path (130:2) coordinate (v130);
	\path (160:2) coordinate (v160);
	\path (180:2) coordinate (v180);
	\path (270:2) coordinate (v270);
	
	\draw [fill=violet] (origin) circle (0.6ex); \draw (origin) node[left=4pt,below=3pt, violet] {$\puncture$};
	\draw[thick] (origin) circle (2);
	
	\draw [fill=olive] (v0) circle (.4ex);
	\draw [fill=olive] (v40) circle (.4ex);
	\draw [fill=olive] (v90) circle (.4ex);
	\draw [fill=olive] (v180) circle (.4ex);
	\draw [fill=olive] (v130) circle (.4ex);
	\draw [fill=olive] (v270) circle (.4ex);
\end{scope}
}

\newcommand\FourCoverWheel[1]{
\begin{tikzpicture}[scale=#1,gray]
	\MacroFourCover
		
	\draw[thick] (v0) .. controls (1.5,0.3)  and (1.4,0.8) .. (v40);
	\draw[thick] (v40)  .. controls (0.6,1.3)  and (0.4,1.4) .. (v90);
	
	\begin{scope}[rotate=90]
		\path (0:2) coordinate (u0);
		\path (90:2) coordinate (u90);
		\path (40:2) coordinate (u40);
		\draw[thick] (u0) .. controls (1.5,0.3)  and (1.4,0.8) .. (u40);
		\draw[thick] (u40)  .. controls (0.6,1.3)  and (0.4,1.4) .. (u90);
	\end{scope}
	
	\draw[thick] (origin) -- (v0) node[pos=0.6,below]{$1$};
	\draw[thick] (origin) -- (v180) node[pos=0.6,below]{$1$};
	\draw[thick] (origin) -- (v90) node[pos=0.4, left]{$1$};
	\draw[thick] (origin) -- (270:2) node[pos=0.6, left]{$1$};
	
	\draw[thick] (origin) -- (v40) node[pos=0.6,below]{$2$};
	\draw[thick] (origin) -- (v130) node[pos=0.6,below]{$2$};
	
	\draw [fill=olive] (v70) circle (.4ex);
	\draw[line width=2 pt, dashed, red] (origin) -- (v70); \draw (v70) node[right, above=3pt, red] {{$\tilsi$}};
	
	\draw [fill=olive] (v160) circle (.4ex);
	\draw[line width=2 pt, dashed, red] (origin) -- (v160); \draw (v160) node[left=5pt, above=3pt, red] {{$\tilsi'$}};
	
	\draw [line width = 1.5pt, dashed, cyan, rounded corners = 20pt] (v70) -- (90:1) -- (u40);
	\draw (130:2) node[left, above=3pt, cyan] {{$\tilbe$}};
\end{tikzpicture}
  }
\newcommand\FourCoverSelfFolded[1]{
\begin{tikzpicture}[scale=#1,gray]
	\MacroFourCover
		
	\draw[thick] (v0) .. controls (1.5,0.3)  and (1.4,0.8) .. (v40) node[pos=0.5,right]{$\b$};
	\draw[thick] (v40)  .. controls (0.6,1.3)  and (0.4,1.4) .. (v90) node[pos=0.5,right]{$\a$};
	
	\begin{scope}[rotate=90]
		\path (0:2) coordinate (u0);
		\path (90:2) coordinate (u90);
		\path (40:2) coordinate (u40);
		\draw[thick] (u0) .. controls (1.5,0.3)  and (1.4,0.8) .. (u40) node[pos=0.5,left]{$\b$};
		\draw[thick] (u40)  .. controls (0.6,1.3)  and (0.4,1.4) .. (u90) node[pos=0.3]{$\a$};
	\end{scope}
	
	\draw[thick] (origin) -- (v0) node[pos=0.4,below]{$r$};
	\draw[thick] (origin) -- (v180) node[pos=0.4,below]{$r$};
	\draw[thick] (origin) -- (v90) node[pos=0.4, left]{$r$};
	\draw[thick] (origin) -- (270:2) node[pos=0.6, left]{$r$};
	
	\draw[thick] (v0) .. controls (0.5,0.4)  and (0.1,1.2) .. (v90) node[pos=0.3,right]{$\ell$};
	\draw[thick] (v180)  .. controls (-0.5,0.4)  and (-0.1,1.2) .. (v90) node[pos=0.3,right]{$\ell$};
	
	\draw [fill=olive] (v70) circle (.4ex);
	\draw[line width=2 pt, dashed, red] (origin) -- (v70); \draw (v70) node[right, above=3pt, red] {{$\tilsi$}};
	
	\draw [fill=olive] (v160) circle (.4ex);
	\draw[line width=2 pt, dashed, red] (origin) -- (v160); \draw (v160) node[left=5pt, above=3pt, red] {{$\tilsi'$}};
	
	\draw [line width = 1.5pt, dashed, cyan, rounded corners = 20pt] (v70) -- (90:1) -- (u40);
	\draw (130:2) node[left, above=3pt, cyan] {{$\tilbe$}};
	
	\draw (v0) node[right,olive] {$v$}; 
	\draw (v90) node[above=2pt,olive] {$v$}; 
	\draw (v180) node[left,olive] {$v$}; 
	\draw (270:2) node[below=2pt,olive] {$v$}; 
	
	\draw (v40) node[right, olive] {$y$}; \draw (v130) node[left, olive] {$y$};
\end{tikzpicture}
  }

\newcommand\MacroMiniLemmaPeripheralPolygon[2]{
\begin{scope}[scale={#1}, {#2},gray]
\coordinate(theta1A) at (-4,-2);
\coordinate(theta1B) at (4,-2);
\coordinate(theta2A) at (-4,-1); \coordinate(theta3A) at (-4,0);  \coordinate(theta4A) at (-4,1);
\coordinate(theta2B) at (4,-1); \coordinate(theta3B) at (4,0);  \coordinate(theta4B) at (4,1);
\coordinate(thetahA) at (-4,2);
\coordinate(thetahB) at (4,2);

\draw (theta1A) -- (theta1B) node[pos=0.2, above,black] {$\taui{d}$}; 
\draw (thetahA) -- (thetahB); 
\draw (theta1A) -- (thetahA) -- (thetahB) -- (theta1B) -- (theta1A);

\draw[style=dashed] (theta2A) -- (theta1B);
\draw[style=dashed] (theta3A) -- (theta1B);
\draw[style=dashed] (theta3A) -- (theta3B);
\draw[style=dashed] (theta4A) -- (theta3B);
\draw[style=dashed] (theta4A) -- (thetahB);

\draw [fill=olive] (theta1A) circle (.6ex) node[left] {$t$}; 
\draw [fill=olive] (theta1B) circle (.6ex); 

\draw [fill=olive] (theta2A) circle (.6ex); 
\draw [fill=olive] (theta4A) circle (.6ex); 
\draw [fill=olive] (theta3B) circle (.6ex); 

\end{scope}
}

\newcommand \MacroMiniLemmaWheelPolygon[2]{
\begin{scope}[scale={#1}, {#2},gray]
\coordinate(rho1A) at (-4,-2); \coordinate(rho2A) at (-4,0); \coordinate(rho3A) at (-3,2);
\coordinate(rho4A) at (0,2.5); \coordinate(rho5A) at (3,2);
\coordinate(P) at (0,0); \coordinate(rhof1A) at (4,0); \coordinate(rhofA) at (4,-2);
\draw (P) -- (rho1A) node[pos=0.5, above,black] {$\rho_1$};
\draw (P) -- (rho2A) node[pos=0.6, above,black] {$\rho_2$};
\draw (P) -- (rho3A) node[pos=0.4, above,black] {$\rho_3$};
\draw[style=dashed] (P) -- (rho4A); \draw[style=dashed] (P) -- (rho5A);
\draw (P) -- (rhof1A) node[pos=0.6, above,black] {$\rho_{f-1}$};
\draw (P) -- (rhofA) node[pos=0.6, above,black] {$\rho_{f}$};
\draw (rho1A) -- (rho2A) -- (rho3A) -- (rho4A) -- (rho5A) -- (rhof1A) -- (rhofA);

\draw[fill=violet] (P) circle(1 ex);
\draw [fill=olive] (rho2A) circle (.6ex);
\draw [fill=olive] (rho3A) circle (.6ex);
\draw [fill=olive] (rhof1A) circle (.6ex);
\end{scope}
}

\newcommand\MacroTechnicalLemma{
\MacroMiniLemmaWheelPolygon{1}{NonOrientedTriangulationStyle}
	\draw [fill=olive] (rho1A) circle (.4ex);
	\draw [fill=olive] (rhofA) circle (.4ex);
	\draw (1,-2) node[above,black] {$\taui{1}$};
	\begin{scope}[yshift=-4cm]
		\MacroMiniLemmaPeripheralPolygon{1}{NonOrientedTriangulationStyle}
	\end{scope}
}

\newcommand\MacroTechnicalLemmaSelfFolded{
         \begin{scope}[gray]
         \coordinate(rho1A) at (-4,-2); \coordinate(rho2A) at (-4,0); \coordinate(rho3A) at (-3,2);
	\coordinate(rho4A) at (0,2.5); \coordinate(rho5A) at (3,2);
	\coordinate(P) at (0,0); \coordinate(rhof1A) at (4,0); \coordinate(rhofA) at (4,-2);

	\draw (P) -- (rho1A) node[pos=0.5, above,black] {$\rho_1=r$};
	\draw[NonOrientedTriangulationStyle, rounded corners = 10pt] (rho1A) -- (rho2A) -- (rho3A) -- (rho4A) -- (rho5A) -- (rhof1A) -- (rhofA);
	\draw[fill=violet] (P) circle(1 ex); 
	\draw (2,-2) node[above,black] {$\taui{2}$};

	\coordinate (startL1) at (rho1A);
	\draw[NonOrientedTriangulationStyle] (rho1A)
	.. controls ($(P)+(3,6)$) and ($(P)+(3,-2)$) .. (rho1A) node[pos=0.5, above=3pt,black] {$\ell=\taui{1}$};
	\draw [fill=olive] (rho1A) circle (.6ex);
	\draw [fill=olive] (rhofA) circle (.6ex);
	\begin{scope}[yshift=-4cm]
		\MacroMiniLemmaPeripheralPolygon{1}{NonOrientedTriangulationStyle}
	\end{scope}
	\end{scope}
}

\newcommand\TikzMinilemmaRadiusGammaPolygonRhoOneLoop[1]{
\begin{tikzpicture}[scale=#1]
	\MacroTechnicalLemma
	\coordinate (startL1) at (rho1A);
	\draw[ inducedmapTPathonCnStyle, 
	postaction={decorate} ] (rho1A)
	.. controls ($(P)+(3,6)$) and ($(P)+(3,-2)$) .. (rho1A) node[pos=0.5, above=5pt] {$\ell_1$};
\end{tikzpicture}
}

\newcommand\TikzTechnicalLemmaUnnotchedRadius[1]{
\begin{tikzpicture}[scale=#1]
	\MacroTechnicalLemma
	\draw[line width=2 pt, dashed, red,PathArrowAtPosStyleSeven] (P) -- (theta1A) node[pos=0.5,right=2pt, red] {{$\si$}};
\end{tikzpicture}
}

\newcommand\TikzTechnicalLemmaUnnotchedRadiusSelfFolded[1]{
\begin{tikzpicture}[scale=#1]
	\MacroTechnicalLemmaSelfFolded
	\draw[line width=2 pt, dashed, red,PathArrowAtPosStyleSeven] (P) -- (theta1A) node[pos=0.5,right=2pt, red] {{$\si$}};
\end{tikzpicture}
}

\newcommand\TikzTechnicalLemmaNotchedRadiusSelfFolded[1]{
\begin{tikzpicture}[scale=#1]
	\MacroTechnicalLemmaSelfFolded
	\draw[line width=2 pt, dashed, red] (P) -- (theta1A) node[pos=0.6,right=2pt, red] {{$\si$}};
	\begin{scope}[yshift=-0.5cm,xshift=-0.2cm]
	\draw (-0.8,-0.4) node[red] {{$\blacktriangleleft$}};
	\draw (-0.3,-0.9) node[red] {{$\blacktriangleright$}};
	\end{scope}
\end{tikzpicture}
}

\newcommand\TikzTechnicalLemmaNotchedRadius[1]{
\begin{tikzpicture}[scale=#1]
	\MacroTechnicalLemma
	\draw[line width=2 pt, dashed, red] (P) -- (theta1A) node[pos=0.5,right=2pt, red] {{$\si$}};
	\draw (-0.8,-0.4) node[red] {{$\blacktriangleleft$}};
	\draw (-0.3,-0.9) node[red] {{$\blacktriangleright$}};
\end{tikzpicture}
}

\newcommand\TikzTechnicalLemmaRhoOne[1]{
\begin{tikzpicture}[scale=#1]
	\MacroTechnicalLemma
	\draw[line width=2 pt, dashed, red] (P) -- (rho1A) node[pos=0.79,above=7pt, red] {{$\si=$}};
	\draw (-1,0) node[red] {{$\blacktriangleleft$}};
	\draw (-0.4,-0.6) node[red] {{$\blacktriangleright$}};
\end{tikzpicture}
}

\newcommand \tilTtwoIdealTriangulationCoordinates {
	\coordinate (P) at (0,4);
	\coordinate (rightd) at (0.2,0); \coordinate (leftd) at (-0.2,0);
	\coordinate (downd) at (0,-0.2); \coordinate (upd) at (0,0.2);
	\coordinate (x1) at (0,0); \coordinate (x2) at (-4,4); \coordinate (x3) at (4,4);
	\coordinate (t1) at (-4,0); \coordinate (s) at (-4,-4); \coordinate (z) at (0,-4);  \coordinate (t2) at (4,0);
}

\newcommand \tilTtwoIdealTriangulationVertices {
\tilTtwoIdealTriangulationCoordinates
\draw [fill=violet] (P) circle (1ex); 
\draw [fill=olive] (x2) circle (.6ex); 
\draw [fill=olive] (x3) circle (.6ex); 
\draw [fill=olive] (t1) circle (.6ex); 
\draw [fill=olive] (t2) circle (.6ex); 
\draw [fill=olive] (x1) circle (.6ex); 

\draw [fill=olive] (s) circle (.6ex); 
\draw [fill=olive] (z) circle (.6ex); 
}

\newcommand \TtwoTriangulatedPolygon[1]{
\begin{tikzpicture}[scale={#1},NonOrientedTriangulationStyle]
\begin{scope}[rotate=-90]
	\tilTtwoIdealTriangulationCoordinates

	\draw[style = thick, NooseArrowAtPosStyleSix]
	(x1) -- (-4,4) node[pos = 0.6, left=2pt]{$\ell$};
	\draw[style = thick, postaction={decorate}]
	(-4,4) -- (0,4) node[pos=0.5, right]
	{$\tiltg{3}=\dot{r}$};
	\draw[style = thick, postaction={decorate}]
	(0,0) -- (0,4) node[pos=0.4, above]
	{$\tiltaui{4}= r$};
	\draw[style = thick, NooseArrowAtPosStyleFive]
	(4,4) -- (0,0) node[pos = 0.3, right=2pt]{$\ell$};
	\draw[style = thick, postaction={decorate}]
	(4,4) -- (0,4) node[pos=0.5, right]
	{$\tiltg{4}= \ddot{r}$};

	\draw [postaction={decorate}] (x2) -- (t1) node[pos=0.3, above] {$b_3$};
	\draw [postaction={decorate}] (t1) -- (x1) node[pos=0.3, left] {$2$};
	\draw [postaction={decorate}] (x1) -- (t2) node[pos=0.6, left] {$\tiltg{6}=b_3$};
	\draw [postaction={decorate}] (t2) -- (x3) node[pos=0.5, below] {$\tiltg{5} = 2$};
	
	\draw [postaction={decorate}] (t1) -- (s) node[pos=0.5, above] {$\tiltg{-1}=b_4$};
	\draw [postaction={decorate}] (s) -- (z) node[pos=0.3, left] {$\tiltg{0}=b_1$};
	\draw [postaction={decorate}] (z) -- (x1) node[pos=0.5, below] {$\tiltg{1}=b_2$};
	\draw [postaction={decorate}] (t1) -- (z) node[pos=0.3, left] {$1$};
	
	\begin{scope}[gray!65]
	\draw (-3,1) node {$\tilblacktri{2}$};
	\end{scope}
	
	\draw
	[red!60, line width = 2.5pt, ->, style=dashed, rounded corners = 10pt,
	decoration={ markings, mark=at position .6 with {\arrow[red!80, line width=2.5pt]{to};}}]
	(s) -- (-2,-2) 
	 -- (-2,2)
	-- (-0,2) -- (2,2) node[pos=0.3,left] {$\tilga$}  -- ($(t2)+(-0.1,0.1)$);

	\tilTtwoIdealTriangulationVertices
	\end{scope}
	\node at (10,0) {};
	\node at (-10,0) {};
	\node[gray!65] at (-3.2,2.2) {$\tilblacktri{0}$};
	\node[gray!65] at (3.1,1.6) {$\tilblacktri{3}$};
	\node[gray!65] at (-1.4,1) {$\tilblacktri{1}$};
	\node[gray!65] at (3.1,-1.2) {$\tilblacktri{4}$};
	\node[gray!65] at (2.1,-3.2) {$\tilblacktri{5}$};
	\end{tikzpicture}
}
\newcommand \tilTtwoSnakeGraphCoordinates {
	\coordinate (P) at (8,0);
	\coordinate (rightd) at (0.2,0); \coordinate (leftd) at (-0.2,0);
	\coordinate (downd) at (0,-0.2); \coordinate (upd) at (0,0.2);
	\coordinate (x1a) at (0,0); \coordinate (x2) at (4,0);
	\coordinate (x1b) at (4,4); \coordinate (x1c) at (8,-4);
	\coordinate (x3a) at (8,4);  \coordinate (x3b) at (12,0);
	\coordinate (t1a) at (-4,0);  \coordinate (t1b) at (4,-4);
	\coordinate (s) at (-4,-4); \coordinate (z) at (0,-4);
	\coordinate (t2) at (12,4);
}

\newcommand \tilTtwoSnakeGraphVertices {
\tilTtwoSnakeGraphCoordinates
\draw [fill=violet] (P) circle (1ex); 
\draw [fill=olive] (x2) circle (.6ex); 
\draw [fill=olive] (x3a) circle (.6ex); 
\draw [fill=olive] (x3b) circle (.6ex); 
\draw [fill=olive] (t1a) circle (.6ex); 
\draw [fill=olive] (t1b) circle (.6ex); 
\draw [fill=olive] (t2) circle (.6ex); 
\draw [fill=olive] (x1a) circle (.6ex); 
\draw [fill=olive] (x1b) circle (.6ex); 
\draw [fill=olive] (x1c) circle (.6ex); 
\draw [fill=olive] (s) circle (.6ex); 
\draw [fill=olive] (z) circle (.6ex); 
}

\tikzstyle{orientedAtEllTriangulationStyle}=[style=thick,
decoration={ markings, mark=at position .6 with {\arrow{to};}} ]
\newcommand\TtwoSnakeGraphPlainer[1]{
\begin{tikzpicture}[scale={#1},orientedAtEllTriangulationStyle,gray!80]
	\tilTtwoSnakeGraphCoordinates
	\draw [postaction={decorate}] (x3a) -- (x3b) node[black,pos=0.4, below] {$\ell$}; 
	\draw (t2) -- (x3a) node[black,pos=0.4, above] {$2$}; 
	\draw (x3b) -- (t2) node[black,pos=0.5, right] {$b_3$}; 
	\draw (x3b) -- (P) node[black,pos=0.5, below] {$r$}; 

	\draw (x1b) -- (P) node[black,pos=0.4, below] {$r$}; 
	\draw [postaction={decorate}] (x1b) -- (x2) node[black,pos=0.5, left] {$\ell$}; 
	\draw [postaction={decorate}] (x3a) -- (x1b) node[black,pos=0.4, above] {$\ell$}; 
	\draw (x3a) -- (P) node[black,pos=0.5] {$\ddot{r}$}; 
	
	\draw [postaction={decorate}] (x1c) -- (x2) node[black,pos=0.5, below] {$\ell$}; 
	\draw (t1b) -- (x1c) node[black,pos=0.3, below] {$2$}; 
	\draw (x1c) -- (P) node[black,pos=0.5, right] {$r$}; 
	\draw (x2) -- (P) node[black,pos=0.5] {$\dot{r}$}; 

	\draw (t1b) -- (x1a) node[black,pos=0.5] {$2$}; 
	\draw (x2) -- (t1b) node[black,pos=0.5] {$b_3$}; 
	\draw (t1b) -- (z) node[black,pos=0.5, below] {$1$}; 
	\draw [postaction={decorate}] (x1a) -- (x2) node[black,pos=0.5, above] {$\ell$}; 
	
	\draw (t1a) -- (z) node[black,pos=0.5, below] {$1$}; 
	\draw (t1a) -- (s) node[black,pos=0.4, left] {$b_4$}; 
	\draw (s) -- (z) node[black,pos=0.5, below] {$b_1$};
	\draw (z) -- (x1a) node[black,pos=0.6] {$b_2$};
	\draw (t1a) -- (x1a) node[black,pos=0.5, above] {$2$}; 
	
	\begin{scope}[gray]
	\draw (-3,-3) node {$\tilblacktri{0}$};
	\draw (-1.4,-0.8) node {$\tilblacktri{1}$};
	\draw (1.4,-3) node {$\tilblacktri{1}^{(-)}$};
	\draw (2.5,-1) node {$\tilblacktri{2}^{(-)}$};
	\draw (5,-3) node {$\tilblacktri{2}$};
	\draw (6.5,-1) node {$\tilblacktri{3}$};
	\draw (5.5,1) node {$\tilblacktri{3}^{(-)}$};
	\draw (7,3) node {$\tilblacktri{4}^{(-)}$}; 
	\draw (9,1) node {$\tilblacktri{4}$}; 
	\draw (11,3) node {$\tilblacktri{5}$};
	\end{scope}

	\tilTtwoSnakeGraphVertices
	\end{tikzpicture}
}

\newcommand\MacroTtwoPerfectMatchingsWithDiagonals{

	\tilTtwoSnakeGraphCoordinates
	\MacroTtwoPerfectMatchings
	\begin{scope}[line width=2pt, dashed]
	\draw (x3a) -- (x3b) node[pos=0.4, right] {$\ell$}; 
	
	\draw [postaction={decorate}] (x1b) -- (P) node[pos=0.5, below] {$r$};

	\draw [postaction={decorate}] (x1c) -- (x2) node[pos=0.6, right] {$\ell$};

	\draw [postaction={decorate}] (t1b) -- (x1a) node[pos=0.5, right] {$2$};
	
	\draw [postaction={decorate}] (t1a) -- (z) node[pos=0.4, right] {$1$};
	\end{scope}
}

\newcommand\MacroTtwoPerfectMatchings{

	\tilTtwoSnakeGraphCoordinates
	\draw (t2) -- (x3a) ;
	\draw (x3b) -- (t2) ;
	\draw (x3b) -- (P) ;
	
	\draw [fill=gray, fill opacity=0.3] (x3a)  -- (x1b) -- (x2) -- (P) -- cycle;
	\draw (x1b) -- (x2);
	\draw (x3a) -- (x1b);
	\draw  (x3a) -- (P);
	
	\draw (t1b) -- (x1c);
	\draw (x1c) -- (P) ;
	\draw (x2) -- (P);

	\draw(x2) -- (t1b) ;
	\draw (t1b) -- (z) node[white,pos=0.1, below] {$1$};
	\draw (x1a) -- (x2) ;
	
	\draw  (t1a) -- (s);
	\draw (s) -- (z);
	\draw (z) -- (x1a);
	\draw (t1a) -- (x1a); 
	
	\node at (t1b) {};
	\node at (x1b) {};
}

\newcommand\TtwoPerfectMatchingsOne[1]{
\begin{tikzpicture}[scale={#1}]
	\MacroTtwoPerfectMatchings
	\begin{scope}[line width=4pt]
	\draw  (t1a) -- (s) node[pos=0.4, left] {$b_4$};
	\draw (z) -- (x1a) node[pos=0.6, left=2pt] {$b_2$};
	\draw(x2) -- (t1b) node[pos=0.5, left] {$b_3$};
	\draw (x1c) -- (P) node[pos=0.3, right] {$r$};
	\draw (x3a) -- (x1b) node[pos=0.3, below] {$\ell$};
	\draw (x3b) -- (t2) node[pos=0.4, right] {$b_3$};
	\end{scope}
\end{tikzpicture}
}

\newcommand\TtwoPerfectMatchingsThree[1]{
\begin{tikzpicture}[scale={#1}]
	\MacroTtwoPerfectMatchings
	\begin{scope}[line width=4pt]
	 \draw  (t1a) -- (x1a) node[pos=0.3, above] {$2$};
	\draw (s) -- (z) node[pos=0.4, below] {$b_1$};
	\draw (t1b) -- (x1c) node[pos=0.3, below] {$2$};
	\draw (x1b) -- (x2) node[pos=0.3, left] {$\ell$};
	\draw (t2) -- (x3a) node[pos=0.4, below] {$2$};
	\draw (x3b) -- (P) node[pos=0.3, below] {$r$};
	\end{scope}
\end{tikzpicture}
}

\newcommand\TtwoPerfectMatchingsFour[1]{
\begin{tikzpicture}[scale={#1}]
	\MacroTtwoPerfectMatchings
	\begin{scope}[line width=4pt]
	 \draw (t1a) -- (x1a) node[pos=0.3, above] {$2$};
	\draw (s) -- (z) node[pos=0.4, above=2pt] {$b_1$};
	\draw(x2) -- (t1b) node[pos=0.5, left] {$b_3$};
	\draw (x1c) -- (P) node[pos=0.3, right] {$r$};
	\draw (x3a) -- (x1b) node[pos=0.3, below] {$\ell$};
	\draw (x3b) -- (t2) node[pos=0.4, right] {$b_3$};
	\end{scope}
\end{tikzpicture}
}

\newcommand\TtwoPerfectMatchingsFive[1]{
\begin{tikzpicture}[scale={#1}]
	\MacroTtwoPerfectMatchings
	\begin{scope}[line width=4pt]
	\draw  (t1a) -- (s) node[pos=0.4, left=2pt] {$b_4$};
	\draw (z) -- (x1a) node[pos=0.6, left=2pt] {$b_2$};
	\draw (t1b) -- (x1c) node[pos=0.3, below] {$2$};
	\draw (x1b) -- (x2) node[pos=0.3, left] {$\ell$};
	\draw (t2) -- (x3a) node[pos=0.4, above] {$2$};
	\draw (x3b) -- (P) node[pos=0.3, below] {$r$};
	\end{scope}
\end{tikzpicture}
}

\newcommand\TtwoPerfectMatchingsTwo[1]{
\begin{tikzpicture}[scale={#1}]
	\MacroTtwoPerfectMatchings
	\begin{scope}[line width=4pt]
	\draw  (t1a) -- (s) node[pos=0.4, left=2pt] {$b_4$};
	\draw (t1b) -- (z) node[pos=0.3, below] {$1$};
	\draw (x1a) -- (x2) node[pos=0.6, above] {$\ell$};
	\draw (x1c) -- (P) node[pos=0.3, right] {$r$};
	\draw (x3a) -- (x1b) node[pos=0.3, above] {$\ell$};
	\draw (x3b) -- (t2) node[pos=0.4, right=2pt] {$b_3$};
	\end{scope}
\end{tikzpicture}
}

\newcommand\TtwoPerfectMatchingsSix[1]{
\begin{tikzpicture}[scale={#1}]
	\MacroTtwoPerfectMatchings
	\begin{scope}[line width=4pt]
	\draw  (t1a) -- (s) node[pos=0.4, left=2pt] {$b_4$};
	\draw (z) -- (x1a) node[pos=0.6, left=2pt] {$b_2$};
	\draw (t1b) -- (x1c) node[pos=0.3, above] {$2$};
	\draw (x1b) -- (x2) node[pos=0.3, left=2pt] {$\ell$};
	\draw  (x3a) -- (P) node[pos=0.3, right=2pt] {$\ddot{r}$};
	\draw (x3b) -- (t2) node[pos=0.4, right=2pt] {$b_3$};
	\end{scope}
\end{tikzpicture}
}

\newcommand\TtwoPerfectMatchingsSeven[1]{
\begin{tikzpicture}[scale={#1}]
	\MacroTtwoPerfectMatchings
	\begin{scope}[line width=4pt]
	\draw  (t1a) -- (s) node[pos=0.4, left=2pt] {$b_4$};
	\draw (z) -- (x1a) node[pos=0.6, left=2pt] {$b_2$};
	\draw (t1b) -- (x1c) node[pos=0.3, below] {$2$};
	\draw (x2) -- (P) node[pos=0.3, below] {$\dot{r}$};
	\draw (x3a) -- (x1b) node[pos=0.3, above] {$\ell$};
	\draw (x3b) -- (t2) node[pos=0.4, right=2pt] {$b_3$};
	\end{scope}
\end{tikzpicture}
}

\newcommand\MacroTtwoPerfectMatchingsEight[1]{
\begin{scope}[scale={#1}]
	\MacroTtwoPerfectMatchings
	\begin{scope}[line width=4pt]
	\draw (s) -- (z) node[pos=0.4, below] {$b_1$};
	\draw [postaction={decorate}] (t1a) -- (x1a) node[pos=0.3, above] {$2$};
	\draw (t1b) -- (x1c) node[pos=0.3, below] {$2$};
	\draw (x1b) -- (x2) node[pos=0.3, left=2pt] {$\ell$};
	\draw  (x3a) -- (P) node[pos=0.6, right] {$\ddot{r}$};
	\draw (x3b) -- (t2) node[pos=0.4, right] {$b_3$};
	\end{scope}
\end{scope}
}

\newcommand\TtwoPerfectMatchingsWithDiagonalsEight[1]{
\begin{tikzpicture}[scale={#1}]
\MacroTtwoPerfectMatchingsEight{1}
\MacroTtwoPerfectMatchingsWithDiagonals
\node at (-8,0) {};\node at (6,0) {};
\end{tikzpicture}
}

\newcommand\TtwoPerfectMatchingsEight[1]{
\begin{tikzpicture}[scale={#1}]
\MacroTtwoPerfectMatchingsEight{1}
\end{tikzpicture}
}

\newcommand\TtwoPerfectMatchingsNine[1]{
\begin{tikzpicture}[scale={#1}]
	\MacroTtwoPerfectMatchings
	\begin{scope}[line width=4pt]
	\draw (s) -- (z) node[pos=0.4, below] {$b_1$};
	\draw [postaction={decorate}] (t1a) -- (x1a) node[pos=0.3, above] {$2$};
	\draw (t1b) -- (x1c) node[pos=0.5, below] {$2$};
	\draw (x2) -- (P) node[pos=0.3, below]{$\dot{r}$};
	\draw (x3a) -- (x1b) node[pos=0.3, below]{$\ell$};
	\draw (x3b) -- (t2) node[pos=0.4, right=2pt] {$b_3$};
	\end{scope}
\end{tikzpicture}
}

\newcommand\TikzPMLRR[1]{
\begin{tikzpicture}[scale={#1}] 	
\MacroTripleTilesOrientationDiagonalsLR{1}
\draw [tilTwoSnakeGraphMatchingStyle] (x1c) -- (P) node[pos=0.3, right] {$r$};
\end{tikzpicture}
}
\newcommand\TikzPMLLR[1]{
\begin{tikzpicture}[scale={#1}] 	
\MacroTripleTilesOrientationDiagonalsLR{1}
\draw [tilTwoSnakeGraphMatchingStyle] (x1b) -- (x2) node[pos=0.3, left] {$\ell$};
\end{tikzpicture}
}
\newcommand\TikzPMLRdotR[1]{
\begin{tikzpicture}[scale={#1}] 	
\MacroTripleTilesOrientationDiagonalsLR{1}
\draw [tilTwoSnakeGraphMatchingStyle] (x2) -- (P) node[pos=0.7, below] {$\dot{r}$};
\end{tikzpicture}
}

\newcommand\TikzPMRLL[1]{
\begin{tikzpicture}[scale={#1}] 	
\MacroTripleTilesOrientationDiagonalsRL{1}
\draw [tilTwoSnakeGraphMatchingStyle] (x3a) -- (x1b) node[pos=0.3, above] {$\ell$};
\end{tikzpicture}
}
\newcommand\TikzPMRRL[1]{
\begin{tikzpicture}[scale={#1}] 	
\MacroTripleTilesOrientationDiagonalsRL{1}
\draw [tilTwoSnakeGraphMatchingStyle] (x3b) -- (P) node[pos=0.3, below] {$r$};
\end{tikzpicture}
}
\newcommand\TikzPMRRddotL[1]{
\begin{tikzpicture}[scale={#1}] 	
\MacroTripleTilesOrientationDiagonalsRL{1}
\draw [tilTwoSnakeGraphMatchingStyle] (x3a) -- (P) node[pos=0.3, left] {$\ddot{r}$};
\end{tikzpicture}
}

\tikzstyle{tilTwoSnakeGraphMatchingStyle}=[line width=3pt, blue]
\newcommand\MacroTripleTilesOrientationDiagonalsLR[1]{ 

	\tilTtwoSnakeGraphCoordinates
	\MacroTtwoPerfectMatchings
	\begin{scope}[tilTwoSnakeGraphMatchingStyle, dashed]
	
	\draw [postaction={decorate}] (x1b) -- (P) node[pos=0.5, above] {$r$};

	\draw [postaction={decorate}] (x1c) -- (x2) node[pos=0.6, below] {$\ell$};
	
	\end{scope}
}

\newcommand\MacroTripleTilesOrientationDiagonalsRL[1]{

	\tilTtwoSnakeGraphCoordinates
	\MacroTtwoPerfectMatchings
	\begin{scope}[tilTwoSnakeGraphMatchingStyle, dashed]
	\draw (x3a) -- (x3b) node[pos=0.4, right] {$\ell$}; 
	
	\draw [postaction={decorate}] (x1b) -- (P) node[pos=0.5, below] {$r$};

	\end{scope}
}

\newcommand\MacroTriangle{
\draw[thick]
(0,0) -- (4,0) -- (2,4) -- (0,0);
\draw [fill=olive] (0,0) circle (.6ex);
\draw [fill=olive] (4,0) circle (.6ex);
\draw [fill=olive] (2,4) circle (.6ex);
}
\newcommand\TikzTriangle[1]{
\begin{tikzpicture}[scale={#1} ]
\MacroTriangle
\node at (4.5,0) {}; \node at (-0.5,0) {};
\end{tikzpicture}
}

\newcommand\TikzTriangleCrossesTwoEdges[1]{
\begin{tikzpicture}[scale={#1} ]
\MacroTriangle
\draw (2,-0.5) node[] {$\tg{k}$};
\draw (4,2.8) node[] {$k+1$};
\draw (0.6,2.8) node[] {$k$};
\draw[TPathStyle]
(0.5,2) -- (3.4,2) node[pos=0.5, below]{$\ga$};
\end{tikzpicture}
}

\newcommand\TikzTriangleCrossesOneEdge[1]{
\begin{tikzpicture}[scale={#1} ]
\MacroTriangle
\draw[TPathStyle,->]
(0,0) -- (3.7,2.1) node[pos=0.5, above]{$\ga$};
\draw (3.1,3) node[] {$1$};
\draw (2,-0.5) node[] {$\tg{0}$};
\draw (0.2,2.8) node[] {$\tg{-1}$};

\begin{scope}[yshift=-6cm]
\MacroTriangle
\draw[TPathStyle, <-]
(0.2,0.1) -- (3.6,2) node[pos=0.5, above]{$\ga$};
\draw (3.1,3) node[] {$d$};
\draw (2,-0.5) node[] {$\tg{d}$};
\draw (0.2,3) node[] {$\tg{d+1}$};
\end{scope}

\end{tikzpicture}
}

\newcommand\MacroTwoVertexTriangle{ 
\draw[style = thick]
	(0,0) .. controls (-1.7,1) and (-1.6,3) .. (0,3.5) ..
	controls (1.6,3)  and (1.7,1) .. (0,0);
	\draw [fill=olive] (0,0) circle (.6ex);
	\draw [fill=olive] (0,3.5) circle (.6ex);

\draw[style = thick, rounded corners = 20pt] 
	(0,0) .. controls (-0.8,1) and (-0.5,1.7) .. (0,2.5) ..
	controls (0.5,1.7)  and (0.8,1) .. (0,0);
}
\newcommand\TikzTwoVertexTriangle[1]{
\begin{tikzpicture}[scale={#1} ]
\MacroTwoVertexTriangle
\end{tikzpicture}
}

\newcommand\TikzTwoVertexTriangleCrossesTwoEdgesNoLoop[1]{
\begin{tikzpicture}[scale={#1} ]
\MacroTwoVertexTriangle
\draw[TPathStyle]
(-1.4,2.5) -- (1.4,2.5) node[pos=0.5, above]{$\ga$};
\draw (0.1,1.5) node[] {$\tg{k}$};
\draw (-1.4,2.3) node[] {$k$};
\draw (1.4,2.3) node[] {$k+1$};
\end{tikzpicture}
}

\newcommand\TikzTwoVertexTriangleCrossesTwoEdgesKOneisLoop[1]{
\begin{tikzpicture}[scale={#1} ]
\MacroTwoVertexTriangle
\draw[TPathStyle]
(-1.3,1.2) -- (-0.2,1.2) node[pos=0.4, above]{$\ga$}
node[black,pos=0,above=24pt]{$k$}
node[black,pos=1,above=18pt]{$k+1$};
\end{tikzpicture}
}

\newcommand\TikzTwoVertexTriangleCrossesOneEdgeLoop[1]{
\begin{tikzpicture}[scale={#1} ]
\MacroTwoVertexTriangle
\draw[TPathStyle,->] (0,3.5) -- (0.1,1.6) node[pos=0.5, left]{$\ga$};
\draw (-0.3,1.3) node[] {$1$};
\draw (-1.7,2.3) node[] {$\tg{0}$};
\draw (1.8,2.3) node[] {$\tg{-1}$};
\begin{scope}[yshift=-4cm]
\MacroTwoVertexTriangle
\draw[TPathStyle,<-] (0,3.4) -- (0.2,1.7) node[pos=0.5, left]{$\ga$};
\draw (-0.3,1.3) node[] {$d$};
\draw (-1.7,2.3) node[] {$\tg{d}$};
\draw (1.8,2.3) node[] {$\tg{d+1}$};
\end{scope}
\end{tikzpicture}
}

\newcommand\TikzTwoVertexTriangleCrossesOneEdgeNotLoop[1]{
\begin{tikzpicture}[scale={#1} ]
\MacroTwoVertexTriangle
\draw[->,TPathStyle,rounded corners=15pt]
(0,0) -- (-1,1) -- (-0.7,2.5) -- (1.3,3) node[pos=0.4, above]{$\ga$};
\draw (0.1,1.5) node[] {$\tg{0}$};
\draw (-1.8,2.3) node[] {$\tg{-1}$};
\draw (1.4,2.3) node[] {$1$};
\begin{scope}[yshift=-4cm]
\MacroTwoVertexTriangle
\draw[PathArrowAtPosStyleReversedFive,TPathStyle,rounded corners=15pt]
(0,0) -- (-1,1) -- (-0.7,2.5) -- (1.2,3) node[pos=0.4, above]{$\ga$};
\draw (0,1.6) node[] {$\tg{d}$};
\draw (-1.8,2.8) node[] {$\tg{d+1}$};
\draw (1.5,2.3) node[] {$d$};
\end{scope}
\end{tikzpicture}
}

\newcommand\MacroVeryThinNoose{
\draw[style = thick, rounded corners = 10pt] 
	(0,0) .. controls (-0.6,1) and (-0.4,1.7) .. (0,2.2) ..
	controls (0.4,1.7)  and (0.6,1) .. (0,0);
}

\newcommand\MacroOneVertexTriangle{
\draw[style = thick, rounded corners = 30pt]
	(0,0) .. controls (-2.7,1) and (-2.4,2.5) .. (0,3.5) ..
	controls (2.4,2.5)  and (2.7,1) .. (0,0);
	
	\draw [fill=olive] (0,0) circle (.6ex);
	
\begin{scope}[rotate=30]
\MacroVeryThinNoose
\end{scope}

\begin{scope}[rotate=-35]
\MacroVeryThinNoose
\end{scope}
}

\newcommand\TikzOneVertexTriangle[1]{
\begin{tikzpicture}[scale={#1} ]
\MacroOneVertexTriangle
\end{tikzpicture}
}

\newcommand\MacroSelfFoldedTriangleNoLabel{
	\draw[thick]
	(0,0) -- (0,2);
	
	\draw [fill=violet] (0,2) circle (1.2*\TikzPunctureActualSize);
	
	\draw[style = thick, rounded corners = 20pt,NooseArrowAtPosStyleSix]	
	(0,0) .. controls (-1.8,1) and (-1.5,2.5) .. (0,3.5) ..
	controls (1.5,2.5)  and (1.8,1) .. (0,0); 
	\draw [fill=olive] (0,0) circle (.6ex);
}
\newcommand\TikzSelfFoldedTriangle[1]{
\begin{tikzpicture}[scale={#1} ]
\MacroSelfFoldedTriangleNoLabel
\draw (-1.6,2) node {$\ell$};
\draw (0.3,1) node {$r$};
\end{tikzpicture}
}

\newcommand\TikzSelfFoldedTriangleCrossesRadiusNoose[1]{
\begin{tikzpicture}[scale={#1} ]
\MacroSelfFoldedTriangleNoLabel
\draw (-1.3,2.8) node {$k$};
\draw (0.2,1.4) node {$k+1$};
\draw (1.3,3.2) node {$k+2$};
\draw[TPathStyle,->]
(-1.5,1) -- (1.6,1) 
node[pos=0.05, below]{$\ga$};
\end{tikzpicture}
}

\newcommand\TikzSelfFoldedTriangleCrossesNoose[1]{
\begin{tikzpicture}[scale={#1} ]
\MacroSelfFoldedTriangleNoLabel
\draw[->,TPathStyle]
(0,2) -- (0.5,3.6) node[pos=0.5, left]{$\ga$};
\draw (-0.1,1.2) node {$\tg{-1}\,\,\, \tg{0}$};
\draw (1.5,2.3) node {$1$};

\begin{scope}[black,yshift=-4cm]
\MacroSelfFoldedTriangleNoLabel
\draw[<-,TPathStyle]
(0,2.2) -- (0.5,3.5) node[pos=0.5, left]{$\ga$};
\draw (-0.1,1.2) node {$\tg{d+1}\,\,\, \tg{d}$};
\draw (1.5,2.3) node {$d$};
\end{scope}

\end{tikzpicture}
}

\newcommand\TikzPunctureActualSize{0.1}
\newcommand\TikzPointActualSize{0.07}

\newcommand\MacroPunctureClosestToKOrKOneTwoVertexTriangle{
	\begin{scope}[rotate=-90]\MacroTwoVertexTriangle\end{scope}
	\draw (1.8,-1.5) node[black]{$\tg{k}$};
	\draw (-0.3,0) node[olive] {$v$};
	\draw[fill=violet] (1.5,0) circle(\TikzPunctureActualSize) node[violet, above] {$\puncture$};
	\draw[CyanGammaStyle] (0.3,-1.2) -- (1.5,1.7) node[left,pos=0.9]{$\ga$};
}
\newcommand\MacroPunctureClosestToKOrKOneTwoVertexTriangleWithSigma{
	\draw[RedSigmaStyle, rounded corners=10pt] (1,-1.6) --
	(1.1,0) -- (1.5,1) -- (2.3,0.8) -- (2.5,-1.6) node[red,pos=0.6,right]{$\si$};
}

\newcommand\TikzPunctureClosestToKTwoVertexTriangle[1]{
\begin{tikzpicture}[scale={#1}]
\MacroPunctureClosestToKOrKOneTwoVertexTriangle
\draw (2.8,1.3) node[black] {$k+1$};
\end{tikzpicture}
}
\newcommand\TikzPunctureClosestToKTwoVertexTriangleWithSigma[1]{
\begin{tikzpicture}[scale={#1}]
\MacroPunctureClosestToKOrKOneTwoVertexTriangle
\MacroPunctureClosestToKOrKOneTwoVertexTriangleWithSigma
\draw (2.8,1.3) node[black] {$k+1$};
\end{tikzpicture}
}
\newcommand\TikzPunctureClosestToKOneTwoVertexTriangle[1]{
\begin{tikzpicture}[scale={#1}]
\begin{scope}[xscale=-1]\MacroPunctureClosestToKOrKOneTwoVertexTriangle\end{scope}
\draw (-2.7,1.3) node[black] {$k$};
\end{tikzpicture}
}
\newcommand\TikzPunctureClosestToKOneTwoVertexTriangleWithSigma[1]{
\begin{tikzpicture}[scale={#1}]
	\begin{scope}[xscale=-1]
	\MacroPunctureClosestToKOrKOneTwoVertexTriangle
	\MacroPunctureClosestToKOrKOneTwoVertexTriangleWithSigma
	\end{scope}
	\draw (-2.7,1.3) node[black] {$k$};
\end{tikzpicture}
}

\newcommand\MacroPunctureClosestToTGKTwoVertexTriangle{
	\begin{scope}\MacroTwoVertexTriangle\end{scope}
	\draw (-0.6,2) node[black]{$\tg{k}$};
	\draw (0,3.8) node[olive] {$v$};
	\draw[fill=violet] (0,1.5) circle(\TikzPunctureActualSize) node[violet, below=2pt] {$\puncture$};
	\draw[CyanGammaStyle] (-1.4,2.9) -- (1.4,2.9) node[pos=0.9, above]{$\ga$};
}
\newcommand\TikzPunctureClosestToTGKTwoVertexTriangle[1]{
\begin{tikzpicture}[scale={#1}]
	\MacroPunctureClosestToTGKTwoVertexTriangle
\end{tikzpicture}
}
\newcommand\MacroDrawGammaSigmaCrossingTwiceOpeningLeft{
	\draw[RedSigmaStyle, rounded corners=10pt] (-1.5,2.3) -- (0.6,2.5) -- (1,1.8) -- (0.7,1) -- (-1.5,1);
}
\newcommand\TikzPunctureClosestToTGKTwoVertexTriangleWithSigma[1]{
\begin{tikzpicture}[scale={#1}]
	\MacroPunctureClosestToTGKTwoVertexTriangle
	\MacroDrawGammaSigmaCrossingTwiceOpeningLeft
\end{tikzpicture}
}
\newcommand\TikzPunctureClosestToTGKTwoVertexTriangleWithSigmaOpeningRight[1]{
\begin{tikzpicture}[scale={#1}]
	\MacroPunctureClosestToTGKTwoVertexTriangle
	\begin{scope}[xscale=-1]\MacroDrawGammaSigmaCrossingTwiceOpeningLeft\end{scope}
\end{tikzpicture}
}

\newcommand\MacroTriangleLabelsKAndKOne{
\draw (-0.55*\TikzBase,0.7*\TikzHeight) node[black] {$k$};
\draw (0.8*\TikzBase,0.7*\TikzHeight) node[black] {$k+1$};
\draw(0,0.15*\TikzHeight) node[black] {$\tg{k}$};
}

\newcommand\MacroTriangleWithGamma{
\newcommand\TikzHeight{1.8}
\newcommand\TikzBase{1.4}
\draw[thick] (-\TikzBase,0) -- (\TikzBase,0) -- (0,\TikzHeight)  -- (-\TikzBase,0);
\draw [fill=olive] (-\TikzBase,0) circle (.6ex);
\draw [fill=olive] (\TikzBase,0)  circle (.6ex);
\draw [fill=olive] (0,\TikzHeight)  circle (.6ex) node[olive,above=2pt] {$v$};
\draw[CyanGammaStyle] (-0.8*\TikzBase,0.8*\TikzHeight) -- (0.8*\TikzBase,0.8*\TikzHeight) node[cyan, pos=0.2, above] {$\ga$};
}
\newcommand\TikzPunctureClosestToK[1]{
\begin{tikzpicture}[scale={#1}]
\MacroTriangleWithGamma
\draw[fill=violet] (-0.9*\TikzBase,0.4*\TikzHeight) circle (\TikzPunctureActualSize);
\MacroTriangleLabelsKAndKOne
\end{tikzpicture}
}
\newcommand\TikzPunctureClosestToKOne[1]{
\begin{tikzpicture}[scale={#1}]
\MacroTriangleWithGamma
\draw[fill=violet] (0.9*\TikzBase,0.4*\TikzHeight) circle (\TikzPunctureActualSize);
\MacroTriangleLabelsKAndKOne
\end{tikzpicture}
}
\newcommand\TikzPunctureClosestToTGK[1]{
\begin{tikzpicture}[scale={#1}]
\MacroTriangleWithGamma
\draw[fill=violet] (0,-0.2*\TikzBase) circle (\TikzPunctureActualSize);
\MacroTriangleLabelsKAndKOne
\end{tikzpicture}
}

\newcommand\MacroOrdinaryTriangleSigmaClosingLeft{
\coordinate (P) at (-0.9*\TikzBase,0.4*\TikzHeight) ;
\draw[fill=violet] (P) circle (\TikzPunctureActualSize);
\draw[RedSigmaStyle, rounded corners=10pt] (0,-0.5) -- (-1.2*\TikzBase,0.4*\TikzHeight)
-- (-0.8*\TikzBase,0.7*\TikzHeight)  -- (0.5*\TikzBase,-0.5) ;
\draw[red] (0*\TikzBase,0.3*\TikzHeight) node {$\si$};
}

\newcommand\TikzPunctureClosestToKWithSigma[1]{
\begin{tikzpicture}[scale={#1}]
\MacroTriangleWithGamma
\MacroOrdinaryTriangleSigmaClosingLeft
\end{tikzpicture}
}
\newcommand\TikzPunctureClosestToKOneWithSigma[1]{
\begin{tikzpicture}[scale={#1}]
\MacroTriangleWithGamma
\begin{scope}[xscale=-1]\MacroOrdinaryTriangleSigmaClosingLeft\end{scope}
\end{tikzpicture}
}
\newcommand\MacroOrdinaryTriangleSigmaOpeningRight{
\draw[fill=violet] (0,-0.2*\TikzBase) circle (\TikzPunctureActualSize);
\draw[RedSigmaStyle, rounded corners=10pt]  (0.9*\TikzBase,0.7*\TikzHeight) --(0.5*\TikzBase,0.6*\TikzHeight)  --
 (-0.4*\TikzBase,-0.3*\TikzHeight) --  (0.3*\TikzBase,-0.3*\TikzHeight) --
 (0.7*\TikzBase,0.3*\TikzHeight) -- (1*\TikzBase,0.4*\TikzHeight);
 \node[red] at (0,0.45*\TikzBase) {$\si$};
 }

\newcommand\TikzPunctureClosestToTGKWithSigma[1]{
\begin{tikzpicture}[scale={#1}]
\MacroTriangleWithGamma
\MacroOrdinaryTriangleSigmaOpeningRight
\end{tikzpicture}
}
\newcommand\TikzPunctureClosestToTGKWithSigmaOpeningLeft[1]{
\begin{tikzpicture}[scale={#1}]
\MacroTriangleWithGamma
\begin{scope}[xscale=-1]\MacroOrdinaryTriangleSigmaOpeningRight\end{scope}
\end{tikzpicture}
}
\newcommand\MacroNooseTriangulation[1]{
\newcommand\Tikzn{2}
\begin{scope}[#1]
\coordinate (P) at (0,0);
\draw (\Tikzn,\Tikzn) -- (\Tikzn,-\Tikzn) -- (-\Tikzn,-\Tikzn) -- (-\Tikzn,\Tikzn)  -- cycle;
\foreach \pos in {(\Tikzn,\Tikzn),(\Tikzn,-\Tikzn),(-\Tikzn,\Tikzn),(-\Tikzn,-\Tikzn)}
		{\draw (P) -- \pos;
		\draw[fill=olive] \pos circle (1.3*\TikzPointActualSize);}
\draw[fill=violet] (P) circle (1.3*\TikzPunctureActualSize);
\end{scope}
}

\newcommand\TikzNooseGammaIdealTriangulation[1]{
\begin{tikzpicture}[scale={#1}]
\MacroNooseTriangulation{black}
\draw[RedSigmaStyle, NooseArrowAtPosStyleTwo,rounded corners=8pt] (\Tikzn,-\Tikzn) -- (0.6*\Tikzn,0) -- (0.1*\Tikzn,0.7*\Tikzn)
-- (-0.5*\Tikzn,0.5*\Tikzn)
-- (-0.6*\Tikzn,-0.3*\Tikzn) -- (0,-0.6*\Tikzn) -- (\Tikzn,-\Tikzn);
\node at (0.4*\Tikzn,0.4*\Tikzn) [right]{$1$};
\node at (-0.4*\Tikzn,0.4*\Tikzn) [above]{$2$};
\node at (-0.6*\Tikzn,-0.6*\Tikzn) [right]{$3$};
\node at (0.4*\Tikzn,-0.4*\Tikzn) [left]{$\tau$};
\node at (1.3*\Tikzn,0) {$b_1$}; \node at (-1.3*\Tikzn,0) {$b_3$};
\node at (0,1.3*\Tikzn) {$b_2$}; \node at (0,-1.3*\Tikzn) {$b_4$};
\node at (-0.7*\Tikzn,0) [red] {$\ga$};
\end{tikzpicture}
}
\newcommand\MacroDrawBlueDot[2]{
\node[style=thick,fill=blue,blue,circle,minimum width=2.5pt,inner sep=0pt,draw] at (#1,#2) {};
}
\newcommand\MacroNooseTPathOne{
\draw [TPathStyle, PathArrowAtPosStyleFourSmall] (\Tikzn,-\Tikzn) -- (\Tikzn,\Tikzn);
\draw [TPathStyle,rounded corners=20pt,PathArrowAtPosStyleFourSmall] (0.9*\Tikzn,0.9*\Tikzn) -- (P) -- (-0.9*\Tikzn,0.9*\Tikzn);
\MacroDrawBlueDot{0}{0.3*\Tikzn}
\draw [TPathStyle,rounded corners=20pt] (-1*\Tikzn,0.9*\Tikzn) -- (P) -- (-\Tikzn,-\Tikzn);
\MacroDrawBlueDot{-0.3*\Tikzn}{0}
\draw [TPathStyle,rounded corners=20pt,PathArrowAtPosStyleSevenSmall] (-0.9*\Tikzn,-0.9*\Tikzn) -- (P) -- (\Tikzn,-\Tikzn) ;
\MacroDrawBlueDot{0}{-0.3*\Tikzn}
}
\newcommand\MacroNooseTpathDrawLabels{
\node at (0.5*\Tikzn,0.5*\Tikzn) [right]{$1$};
\node at (-0.5*\Tikzn,0.5*\Tikzn) [left]{$2$};
\node at (-0.7*\Tikzn,-0.7*\Tikzn) [right]{$3$};
\node at (0.7*\Tikzn,-0.7*\Tikzn) [left]{$\tau$};
\begin{scope}[white]
\node at (1.3*\Tikzn,0) {$b_1$}; \node at (0,1.3*\Tikzn) {$b_2$};\node at (-1.3*\Tikzn,0) {$b_3$};\node at (0,-1.3*\Tikzn) {$b_4$};
\end{scope}
}

\newcommand\TikzNooseTpathOneBad[1]{
\begin{tikzpicture}[scale={#1}]
\MacroNooseTriangulation{gray}
\MacroNooseTpathDrawLabels
\coordinate (P1) at (0,0.3*\Tikzn);
\draw [TPathStyle, PathArrowAtPosStyleFourSmall] (\Tikzn,-\Tikzn) -- (\Tikzn,\Tikzn);
\draw [TPathStyle,rounded corners=5pt,PathArrowAtPosStyleSixSmall] (0.99*\Tikzn,0.99*\Tikzn) -- (P1);
\draw [TPathStyle,rounded corners=5pt,PathArrowAtPosStyleFourSmall] (P1) --  (-0.2*\Tikzn,0.8*\Tikzn) -- (-0.99*\Tikzn,0.99*\Tikzn);
\MacroDrawBlueDot{0}{0.3*\Tikzn}
\draw [TPathStyle,rounded corners=20pt,PathArrowAtPosStyleSixSmall] (-0.9*\Tikzn,0.9*\Tikzn) -- (P1) ;
\draw [TPathStyle,rounded corners=20pt,PathArrowAtPosStyleFourSmall]
(P1) -- (-\Tikzn,-\Tikzn);
\draw [TPathStyle,rounded corners=20pt,PathArrowAtPosStyleSevenSmall] (-0.9*\Tikzn,-0.9*\Tikzn) -- (P) -- (\Tikzn,-\Tikzn) ;
\MacroDrawBlueDot{0}{-0.3*\Tikzn}

\node at (1.3*\Tikzn,0) {$b_1$};

\node at (1.8*\Tikzn,0) {}; \node at (-1.8*\Tikzn,0) {};
\end{tikzpicture}
}
\newcommand\TikzNooseTpathOne[1]{
\begin{tikzpicture}[scale={#1}]
\MacroNooseTriangulation{gray}

\MacroNooseTpathDrawLabels
\MacroNooseTPathOne
\node at (1.3*\Tikzn,0) {$b_1$};
\end{tikzpicture}
}
\newcommand\TikzNooseTpathTwo[1]{
\begin{tikzpicture}[scale={#1}]
\MacroNooseTriangulation{gray}

\begin{scope}[rotate=90]
\MacroNooseTPathOne
\end{scope}
\MacroNooseTpathDrawLabels
\node at (0,1.3*\Tikzn) {$b_2$};

\end{tikzpicture}
}
\newcommand\TikzNooseTpathThree[1]{
\begin{tikzpicture}[scale={#1}]
\MacroNooseTriangulation{gray}

\begin{scope}[rotate=180]
\MacroNooseTPathOne
\end{scope}
\MacroNooseTpathDrawLabels
\node at (-1.3*\Tikzn,0) {$b_3$};

\end{tikzpicture}
}
\newcommand\TikzNooseTpathFour[1]{
\begin{tikzpicture}[scale={#1}]
\MacroNooseTriangulation{gray}

\begin{scope}[rotate=-90]
\MacroNooseTPathOne
\end{scope}
\MacroNooseTpathDrawLabels
\node at (0,-1.3*\Tikzn) {$b_4$};
\end{tikzpicture}
}

\newcommand\MacroNotNooseDrawArcThree[1]{
\draw [rounded corners=15pt,#1] (-\Tikzn,\Tikzn) -- (-0.5*\Tikzn, -0.5*\Tikzn) -- (\Tikzn,-\Tikzn);
}
\newcommand\MacroNotNooseTriangulation[1]{
\newcommand\Tikzn{2}
\begin{scope}[#1]
\coordinate (P) at (0,0);
\draw (\Tikzn,\Tikzn) -- (\Tikzn,-\Tikzn) -- (-\Tikzn,-\Tikzn) -- (-\Tikzn,\Tikzn)  -- cycle;

\MacroNotNooseDrawArcThree{}

\foreach \pos in {(\Tikzn,\Tikzn),(\Tikzn,-\Tikzn),
(-\Tikzn,\Tikzn)
}
		{\draw (P) -- \pos;
		\draw[fill=olive] \pos circle (1.3*\TikzPointActualSize);}
\draw[fill=violet] (P) circle (1.3*\TikzPunctureActualSize);
\draw[fill=olive] (-\Tikzn,-\Tikzn) circle (1.3*\TikzPointActualSize);
\end{scope}
}

\newcommand\TikzNotNooseGammaIdealTriangulation[1]{
\begin{tikzpicture}[scale={#1}]
\MacroNotNooseTriangulation{black}
\draw[RedSigmaStyle, NooseArrowAtPosStyleTwo,rounded corners=8pt] (\Tikzn,-\Tikzn) -- (0.6*\Tikzn,0) -- (0.1*\Tikzn,0.7*\Tikzn)
-- (-0.5*\Tikzn,0.5*\Tikzn)
-- (-\Tikzn,-\Tikzn);
\node at (0.4*\Tikzn,0.4*\Tikzn) [right]{$1$};
\node at (-0.4*\Tikzn,0.5*\Tikzn) [above]{$2$};
\node at (-0.5*\Tikzn,-0.5*\Tikzn) {$3$};
\node at (0.5*\Tikzn,-0.5*\Tikzn) [left]{$0$};
\node at (1.3*\Tikzn,0) {$b_1$}; \node at (-1.3*\Tikzn,0) {$b_3$};
\node at (0,1.3*\Tikzn) {$b_2$}; \node at (0,-1.3*\Tikzn) {$b_4$};
\node[rotate=50] at (0.2*\Tikzn,0.74*\Tikzn) [red] {$\ga$};
\end{tikzpicture}
}
\newcommand\MacroNotNooseTpathDrawLabels{
\node at (0.4*\Tikzn,0.4*\Tikzn) [right]{$1$};
\node at (-0.4*\Tikzn,0.5*\Tikzn) [above]{$2$};
\node at (-0.5*\Tikzn,-0.5*\Tikzn) {$3$};
\node at (0.5*\Tikzn,-0.5*\Tikzn) [right]{$0$};
\node at (1.3*\Tikzn,0) {}; \node at (-1.3*\Tikzn,0) {};

\begin{scope}[white]
\node at (1.3*\Tikzn,0) {$b_1$}; \node at (0,1.3*\Tikzn) {$b_2$};\node at (-1.3*\Tikzn,0) {$b_3$};\node at (0,-1.3*\Tikzn) {$b_4$};
\end{scope}
}

\newcommand\TikzNotNooseTpathThreeBad[1]{
\begin{tikzpicture}[scale={#1}]
\MacroNotNooseTriangulation{gray}

\MacroNotNooseTpathDrawLabels

\begin{scope}[TPathStyle]
\coordinate (P1) at (-0.2*\Tikzn,-0.2*\Tikzn);

\draw[PathArrowAtPosStyleFourSmall] (\Tikzn,-\Tikzn) -- (P1);
\draw[PathArrowAtPosStyleFourSmall, rounded corners=15pt] (P1)  -- (0.3*\Tikzn,-0.2*\Tikzn) -- (\Tikzn,0.95*\Tikzn);

\draw[rounded corners=15pt,PathArrowAtPosStyleSixSmall] (\Tikzn,0.95*\Tikzn) -- (-0.4*\Tikzn,0.1*\Tikzn) 
-- (P1) ;

\MacroDrawBlueDot{-0.2*\Tikzn}{-0.2*\Tikzn}

\draw[rounded corners=15pt,PathArrowAtPosStyleFourSmall] (P1) 
-- (-\Tikzn,\Tikzn);

\draw [TPathStyle, PathArrowAtPosStyleEnd] (-\Tikzn,\Tikzn) -- (-\Tikzn,-\Tikzn); \node[black] at (-1.3*\Tikzn,0) {$b_3$};

\end{scope}
\node at (-2.3*\Tikzn,0) {}; \node at (2.3*\Tikzn,0) {};
\end{tikzpicture}
}
\newcommand\TikzNotNooseTpathOne[1]{
\begin{tikzpicture}[scale={#1}]
\MacroNotNooseTriangulation{gray}

\MacroNotNooseTpathDrawLabels
\begin{scope}[TPathStyle]
\draw [PathArrowAtPosStyleFourSmall] (\Tikzn,-\Tikzn) -- (\Tikzn,\Tikzn); \node[black] at (1.3*\Tikzn,0) {$b_1$};
\draw [TPathStyle,rounded corners=20pt,PathArrowAtPosStyleFourSmall] (0.9*\Tikzn,0.9*\Tikzn) -- (P) -- (-0.9*\Tikzn,0.9*\Tikzn);
\MacroDrawBlueDot{0}{0.3*\Tikzn}

\draw [PathArrowAtPosStyleEightSmall] (-\Tikzn,\Tikzn) --
(-0.2*\Tikzn,-0.1*\Tikzn) ;
\draw [rounded corners=5pt]
(-0.2*\Tikzn,-0.2*\Tikzn) --
(-0.5*\Tikzn,-0.1*\Tikzn) --
(-\Tikzn,\Tikzn);

\MacroDrawBlueDot{-0.2*\Tikzn}{-0.2*\Tikzn}

\MacroNotNooseDrawArcThree{PathArrowAtPosStyleSevenSmall}
\draw [PathArrowAtPosStyleEnd]  (\Tikzn,-\Tikzn) -- (-\Tikzn,-\Tikzn); \node[black] at (0,-1.3*\Tikzn) {$b_4$};
\end{scope}
\end{tikzpicture}
}
\newcommand\TikzNotNooseTpathTwo[1]{
\begin{tikzpicture}[scale={#1}]
\MacroNotNooseTriangulation{gray}
\MacroNotNooseTpathDrawLabels

\begin{scope}[TPathStyle]
\draw [rounded corners=20pt,PathArrowAtPosStyleSevenSmall] (\Tikzn,-\Tikzn) -- (P) -- (\Tikzn,\Tikzn);

\MacroDrawBlueDot{0.3*\Tikzn}{0}

\draw [PathArrowAtPosStyleFourSmall] (\Tikzn,\Tikzn) -- (-\Tikzn,\Tikzn); \node[black] at (0,1.3*\Tikzn) {$b_2$};

\draw [rounded corners=5pt,PathArrowAtPosStyleSevenSmall] (-\Tikzn,\Tikzn) --
(-0.2*\Tikzn,0.2*\Tikzn) --
(-0.2*\Tikzn,-0.1*\Tikzn) ;
\draw [rounded corners=5pt,PathArrowAtPosStyleFourSmall]
(-0.2*\Tikzn,-0.1*\Tikzn) --
(-0.5*\Tikzn,0) --
(-\Tikzn,\Tikzn);

\MacroDrawBlueDot{-0.2*\Tikzn}{-0.1*\Tikzn}

\MacroNotNooseDrawArcThree{PathArrowAtPosStyleSevenSmall}
\draw [PathArrowAtPosStyleEnd]  (\Tikzn,-\Tikzn) -- (-\Tikzn,-\Tikzn); \node[black] at (0,-1.3*\Tikzn) {$b_4$};
\end{scope}
\end{tikzpicture}
}
\newcommand\TikzNotNooseTpathThree[1]{
\begin{tikzpicture}[scale={#1}]
\MacroNotNooseTriangulation{gray}
\MacroNotNooseTpathDrawLabels

\begin{scope}[TPathStyle]
\draw [rounded corners=20pt,PathArrowAtPosStyleFourSmall] (\Tikzn,-\Tikzn) -- (P) -- (\Tikzn,0.95*\Tikzn);

\MacroDrawBlueDot{0.3*\Tikzn}{0}

\draw [TPathStyle,rounded corners=20pt,PathArrowAtPosStyleEightSmall] (0.9*\Tikzn,0.99*\Tikzn) -- (P) -- (-0.9*\Tikzn,0.9*\Tikzn);
\MacroDrawBlueDot{0}{0.3*\Tikzn}

\draw [TPathStyle, PathArrowAtPosStyleEnd] (-\Tikzn,\Tikzn) -- (-\Tikzn,-\Tikzn); \node[black] at (-1.3*\Tikzn,0) {$b_3$};

\end{scope}
\end{tikzpicture}
}
